\documentclass[a4paper,11pt]{article}
\usepackage[numbers]{natbib}

\usepackage[a4paper,left=2cm,right=2cm,top=2cm,bottom=2cm]{geometry}
\usepackage[colorlinks=true, urlcolor=blue]{hyperref} 

\usepackage{amsmath,amssymb,mathtools}
\usepackage{hyperref}
\usepackage{subcaption}
\usepackage{booktabs}
\usepackage{makecell}
\usepackage{siunitx}
\usepackage{verbatim}
\usepackage{xcolor}

\usepackage{authblk}

\definecolor{kitgreen}{RGB}{0,150,130}

\title{{A Structured Review of Reduced Order Modeling for Domain Decomposition Problems: State of the Art and Perspectives}}

\author[1]{Shenhui Ruan \footnote{shenhui.ruan@kit.edu}}
\author[1]{Andreas G. Class \footnote{andreas.class@kit.edu}}
\author[2]{Gianluigi Rozza \footnote{grozza@sissa.it}}

\affil[1]{\small Institute for Thermal Energy Technology and Safety, Karlsruhe Institute of Technology \\

    Kaiserstra\ss{}e 12, 76131, Karlsruhe, Germany}

\affil[2]{\small Mathematics Area, mathLab, International School for Advanced Studies\\

via Bonomea 265, 34136, Trieste, Italy.}

\date{} 

\usepackage{multirow}
\usepackage{pdflscape}
\usepackage{array}
\usepackage{tabularray}

\renewcommand{\arraystretch}{1.3}
\newcommand{\RBsnap}{\renewcommand{\arraystretch}{1} \begin{tabular}[l]{@{}l@{}}Snapshots \\ \& \\ RB \end{tabular}}

\newcommand{\MP}{\renewcommand{\arraystretch}{1} 
    \begin{tabular}[c]{@{}l@{}}Multi- \\ physics\end{tabular}}

\begin{document}
\let\WriteBookmarks\relax
\def\floatpagepagefraction{1}
\def\textpagefraction{.001}
\maketitle

\begin{abstract}
    Reduced Order Models (ROMs) have been regarded as an efficient alternative to conventional high-fidelity Computational Fluid Dynamics (CFD) for accelerating the design and optimization processes in engineering applications. Many industrial geometries feature repeating subdomains or contain sub-regions governed by distinct physical phenomena, making them well-suited to \emph{Domain Decomposition} (DD) techniques. The integration of ROM and DD is promising to further reduce computational costs by constructing local ROMs and assembling them into global solutions. Due to the complexity and necessity of coupling ROMs, many approaches have been proposed in recent years. This review provides a concise overview of existing methodologies combining ROM and DD. We categorize existing methods into intrusive (projection-based) and non-intrusive (data-driven) frameworks. Various strategies for generating local reduced bases and coupling them across subdomains are illustrated. Particular emphasis is placed on intrusive techniques, including equations, numerical algorithms, and practical implementations. The non-intrusive framework is also discussed, highlighting its general procedures, basic formulations, and underlying principles. Finally, we summarise the state of the literature, identify open challenges, and present perspectives on future implementation from an engineering viewpoint.

\end{abstract}

\tableofcontents


\section{Introduction and background}
\label{sec:intro}

Numerical analysis, and especially high-fidelity three-dimensional simulations, have recently become a practical tool for the design and optimization of various industrial applications, such as power plants \cite{roelofs2019towards, courtessole2024experimental} and energy storage systems \cite{chen2017experimental}. In comparison with experiments, which are expensive and time-consuming for large-scale facilities, simulations are economically efficient. They are practical for improving efficiency and estimating the safety threshold in many applications. Over the past several decades, many numerical schemes have been developed, such as the finite element, finite volume, spectral element, and lattice Boltzmann methods \cite{moukalled2016finite, ern2004theory, canuto2007spectral, zhong2024stochastic}, to fulfill specific demands of analysis in various fields.

Nevertheless, even with the fast development of hardware, multi-query and multi-scale simulation of large-scale systems is still challenging. Therefore, even today, a large range of parametric studies, uncertainty, and sensitivity analyses are limited by computational resources. In spite of this, the demand for fast modeling and real-time control is increasing in almost every industry. To overcome the limitation, \emph{Reduced Order Models} (ROMs) are proposed. They can accelerate simulations by several orders of magnitude. Thus, they have been regarded as efficient techniques to obtain high-fidelity solutions (hereinafter called results of Full Order Models (FOMs) in contrast to ROMs) in complex systems \cite{rozza2022advanced, rozza2024real}.

The fundamental philosophy of ROMs is that high-dimensional solutions over a range of parameters can be well approximated using a few sampled solutions. The implementation of ROMs consists of two steps: (i) a set of expensive \emph{offline} calculations is carried out for optimally selected conditions; (ii) real-time \emph{online} calculations are then performed for numerous cases \cite{quarteroni2015reduced}.

In the offline phase, FOM solutions (also known as \emph{snapshots}) corresponding to the sampled parameters are collected. Since the snapshot matrix can be extremely large, dimensionality reduction techniques are employed to extract a few dominant features while preserving the essential dynamics of the original flow field. The reduction of datasets can be achieved utilizing various techniques, including traditional methods like \emph{Proper Orthogonal Decomposition} (POD) and the iterative \emph{Greedy} algorithm \cite{quarteroni2015reduced, hesthaven2016certified}, as well as more recent \emph{neural network}-based techniques \cite{brunton2024data}.

Once the reduced representations are available, ROMs are later constructed based on them. Over the past several decades, many ROM construction techniques have been developed, generally classified as (i) intrusive and (ii) non-intrusive methods \cite{bai2021non, padula2024brief, rozza2022advanced, benner2020snapshot}. The former manipulates the governing \emph{Partial Differential Equations} (PDEs) using \emph{Galerkin projection} to derive a reduced system. In this case, the dominant physical features are known as basis functions because they can span a \emph{Reduced Basis} (RB) of high-fidelity solutions. The latter is purely data-driven, which adopts a set of sampled data to \emph{train} a surrogate model (i.e., ROM), which is then used to \emph{predict} solutions for unknown parameters. In both frameworks, the resulting ROMs are validated against an unseen dataset to evaluate their predictive capability. 

It is worth introducing several widely used open-source platforms for implementing ROMs: 1) pyMOR (python Model Order Reduction) \cite{milk2016pymor}; 2) ITHACA-FV (In real Time Highly Advanced Computational Applications for Finite Volumes) \cite{Stabile2017CAF}; 3) EZyRB (Easy Reduced Basis method) \cite{demo18ezyrb}; 4) RBniCS (Reduced Basis implementation in FEniCS) \cite{rozza2024rbnics}; 5) ARGOS (Advanced Reduced Groupware Online Simulation) \cite{sissa2021ARGOS}; 6) PINA (Physics Informed Neural network for Advanced modeling) \cite{coscia2023physics}; 7) pyLOW (python Low-Order Modeling) \cite{eiximeno2025pylom}.

Now, we turn our attention to the treatment of the geometries in ROMs, a fundamental topic of this review. Classic ROMs approximate the entire computational domain, which is known as \emph{Global ROMs}. However, this direct method has certain drawbacks, as indicated in \cite{benner2020snapshot, heaney2022ai, wentland2024role}:

\begin{itemize}
    \item The RB is optimized to capture the most necessary flow dynamics across the entire model. However, representing numerous local phenomena requires a large number of modes, which weakens the advantage of ROM \footnote{That issue is called a slow decay of the \emph{Kolmogorov n-width} in the frame of ROM \cite{quarteroni2015reduced}.}.
    \item High computational costs may arise in the offline stage due to the need for generating solutions in huge computational domains and performing dimensional reduction on extensive datasets.
    \item Increasing the dimension of the RB in global field approximation may compromise robustness, stability, and accuracy.
\end{itemize}

These shortcomings have led the research community to integrate \emph{Domain Decomposition} (DD) techniques into ROMs. To reduce computational complexity, this approach divides the physical domain into subdomains. First, reduced spaces are computed locally. Then, the \emph{local fields} are \emph{glued} together to construct an accurate global approximation. Compared to global ROMs, the local approach (referred to as local ROM or DD-ROM hereafter) offers several advantages \cite{lovgren2006reduced, antonietti2016discontinuous, benner2020snapshot, riffaud2021dgdd, sambataro2022component}:

\begin{itemize}
    \item Low-dimensional representations are computed at the local level, leading to a better representation of dominant variations in each subdomain.
    \item Different parameterization can be applied locally, enhancing the flexibility.
    \item Computational costs decrease in both high-fidelity CFD simulations, dimensionality reduction of large snapshot datasets, and the training stage of ROMs.
    \item Different numbers of modes or distinct reduced subspaces can be assigned to each subdomain, reducing overall computational costs in the prediction stage while enhancing robustness and stability.
\end{itemize}

It is worth noting that alternative domain decomposition techniques exist. Instead of geometric partitioning, these approaches employ snapshot clustering using classification methods (e.g., the k-means algorithm), as shown in \cite{amsallem2010towards, amsallem2012nonlinear, washabaugh2012nonlinear, riffaud2020reduced}. In this work, we focus exclusively on methods based on spatial domain decomposition.

We emphasize that the local ROM is well-suited for various realistic applications. Firstly, geometries consisting of repeating patterns are widely employed in many practical applications. Examples include fission \cite{batta2017cfd, batta2024cfd, ruan2024local} and fusion \cite{buhler2024geometric} nuclear reactors, heat exchangers \cite{bergman2011fundamentals}, fuel cells \cite{fan2021recent}, thermal storage systems \cite{roos2021thermocline}, and biological models \cite{pegolotti2021model}, among others. Due to symmetry and periodicity, physical phenomena within these repeating structures often exhibit similar behavior. Thus, domain decomposition can be applicable to the treatment of such geometries.

Moreover, for large-scale simulations, different regions may require varying levels of accuracy. FOMs frequently utilize domain decomposition with varying mesh resolutions across partitions, necessitating corresponding adjustments in ROMs to handle non-conformal mesh interfaces. Consequently, domain decomposition is often applied at both the FOM and ROM levels \cite{zappon2024reduced}. Additionally, coupling FOMs with ROMs also relies on spatial partitioning \cite{cinquegrana2011hybrid, barnett2022schwarz}.

This technique is particularly advantageous in multi-physics simulations, where different regions are governed by distinct PDEs. A key example is Fluid-Structure Interaction (FSI) \cite{discacciati2023localized, prusak2024time, zappon2023efficient}, where localized solution methodologies naturally align with partitioned domains.

Given its relevance to our experience, we emphasize that spatial partitioning is widely adopted in energy-related applications, such as the subchannel approach and the Coarse-Grid CFD (CG-CFD) method. The subchannel analysis \cite{du2019thermal, todreas2021nuclear2}, commonly used in nuclear engineering, divides the computational domain into macro Control Volumes (CVs), where interactions and small-scale effects are incorporated through empirically derived correlations. This partitioning strategy can be combined with local ROM techniques for improved nuclear engineering simulations \cite{ruan2024local}. The Coarse-Grid method \cite{Class2010Coarse}, on the other hand, employs a highly coarse discretization of the spatial domain, capturing sub-grid effects through source terms. These terms are estimated via interpolated data from high-resolution simulations, making the approach effective for modeling flow in tube bundles \cite{viellieber2015coarse} and wind farms \cite{class2014two}. 

Despite the use of subchannel analysis and CG-CFD, these methods often produce low-resolution results, limiting their predictive capabilities. In contrast, ROMs achieve high-fidelity solutions with minimal computational cost, making them a compelling alternative for large-scale simulations.

The development of ROMs for domain decomposition problems began around the early 2000s and gained significant traction in the past decade. They are gaining research interest and show promising potential for scientific and engineering applications. They have been used to model various PDEs, including the Laplace equation \cite{maday2002reduced}, the time-dependent heat equation \cite{ohlberger2017true}, convection-diffusion \cite{mu2019domain}, solid mechanics \cite{eftang2014port}, fluid dynamics, particularly the Navier-Stokes equations \cite{prusak2023optimisation}, and multiphysics problems \cite{corigliano2013domain}.

Two comprehensive reviews for this topic have been published. Andreas Buhr et al. \cite{buhr2020localized} published a thorough review of intrusive localized ROM approaches in 2019. Another review regarding non-intrusive procedures was published in 2021 by Alexander Heinlein et al. \cite{heinlein2021combining}, which discusses a combination of Machine Learning (ML) and domain decomposition. 

However, various studies have appeared in recent years. Advanced techniques based on neural networks have been intensively applied to the topic. We observe that a recent in-depth review of the field is lacking. Therefore, this paper provides a structured overview of existing work on local ROMs from an engineering and practical perspective. 

Our goal is to categorize both classical and advanced techniques developed over the past decades and present fundamental aspects of their implementation. We first identify three common aspects that apply to all methods: domain decomposition, parameterization, and dimensionality reduction. As these three items are not unique, they are classified and discussed separately. Then, we review both intrusive and non-intrusive frameworks, with each category further subdivided according to the mathematical formulation and/or numerical algorithm. This hierarchical classification is well-suited for a structured comparison of various approaches, highlighting their unique characteristics and applications.

The review is organized as follows. Section \ref{sec:classappr} presents the overall outline. That denotes the classification of existing techniques with respect to domain decomposition, parameterization, snapshot reduction, and coupling algorithms. Section \ref{sec:dd_local} presents strategies for domain decomposition. Treatments for parameter space are briefly addressed in Section \ref{sec:parameterization}. Various procedures for generating snapshots and computing local RBs are illustrated in Section \ref{sec:snapshots_local_RB}. Section \ref{sec:projection_based_coupling_algorithms} introduces projection-based coupling methods, including monolithic and iterative algorithms. The explanation and overview considering pure data-driven approaches based on \emph{Artificial Intelligence} (AI) are presented in Section \ref{sec:Data_driven_Techniques}. Finally, Section \ref{sec:summary_conclusions} summarizes the review and presents conclusions on the state of the art and challenges in \emph{Model Order Reduction} (MOR) techniques for domain decomposition problems.

\section{Classification and hierarchy}
\label{sec:classappr}
Shortly summarized, three general aspects arose and need to be considered before the construction of local ROM: (i) partition strategies, (ii) parameterization techniques, and (iii) local dimensionality reduction procedures. Additionally, various algorithms have been proposed to reconstruct global fields by coupling local approximations. Therefore, to better organize the review, the categorization and hierarchy for the observed references are presented as follows. 

\subsection{Domain decomposition strategies}

The domain can be decomposed into \emph{overlapping} or \emph{non-overlapping} parts, which mainly depends on the coupling techniques (see Section \ref{subsec:decomposition_with_without_overlaps}). Also, the high-resolution mesh can be \emph{conforming} or \emph{non-conforming} at partition interfaces (shown in Section \ref{subsec:conforming_non-conforming}). Complex structures can be decomposed into many subdomains of different shapes. However, for geometries with repeating patterns, one can adopt several types of \emph{generic} (or \emph{archetype}) blocks to assemble the whole region. A geometric transformation will map the \emph{reference} subdomains into their final position and shape. More specific explanations are presented in Section \ref{subsec:individual_generic_decomposition}.

\subsection{parameterization}
ROMs are built based on a set of high-fidelity simulations. Each solution is carried out considering some conditions, which we interpret as parameters, sampled from a specific space. More precisely, one should consider how to select values within the specified range. The selection is crucial for the quality and performance of ROM, for both global and local frameworks. Thus, \emph{sampling strategies} are presented and briefly explained in Section \ref{subsec:Sampling_strategies}. Note that, for the domain decomposition configuration, each subdomain can be parametrized separately, which enhances flexibility for complex domains. 

Identifying the parameters is also important. Physical properties and boundary conditions can be interpreted directly, as they are already involved in the governing PDEs. Additionally, the shape of computational domains can also be considered a parameter. That is necessary when ROMs are exploited for tasks like shape optimization. Nevertheless, generating geometrical samples is not as straightforward as the physical ones. Thus, Section \ref{subsec:Geometrical_parameterization} illustrates two widely used techniques for achieving the task, namely, \emph{interpolation-based} methods and a more general approach called \emph{Free Form Deformation} (FFD). Two widely used interpolation algorithms are discussed in particular: \emph{Radial Basis Function} (RBF) and \emph{Inverse Distance Weighting} (IDW). We remind that the geometric parameterization can be utilized globally and locally. If a model is decomposed into subdomains, it is possible to deform each one individually.

The dimension of the parameter space is another aspect in this context. For simple conditions, sampling can be performed along every coordinate. However, for a complex scenario with numerous parameters, too many samples result in prohibitively high costs. Hence, Section \ref{subsec:Parameter_space_reduction} provides two approaches for \emph{dimensionality reduction of the parameter space}: \emph{active subspaces} and \emph{AI-based generative model}.

\subsection{Dimensionality reduction}

After decomposition and sampling, we can begin simulations, collect datasets, and extract dominant features. We remind that the snapshots can be collected either globally or locally, which in turn affects the dimensionality reduction. Be aware that the RB is adopted for the following description. However, the illustration also applies to data-driven frameworks.

In the case of the \emph{global snapshots}, the RBs can be computed in three ways: (i) the global RB is computed and then divided into many local parts following the domain decomposition (as presented in Section \ref{subsubSec:localized_global_RB}); (ii) the global snapshot matrix is split into several submatrices corresponding to the partition. The local modes are extracted from each sub-snapshot. Discussion and notation are shown in Section \ref{sec:snapshots_local_RB}; (iii) the entire model is assembled by a few (or a single) generic parts. The global solutions of each archetype are stacked into separate datasets, and the RBs are constructed with regard to each reference shape.

The so-called \emph{localized training} strategy is utilized for generating subdomain-level FOM solutions, in which independent simulations are carried out in a small system. The computational domain consists of several generic blocks, which are much smaller than the original model. Thus, the cost of the offline stage is significantly reduced. We highlight that the technique assumes local features still accurately represent the total dynamics. \emph{Subdomain-specific solutions} are collected as snapshots. A RB is computed for each archetype block and then is transformed into the instantiations. That is addressed in Section \ref{sec:snapshots_local_RB}.

\begin{figure}[h]
    \centering
    \includegraphics[width=\linewidth]{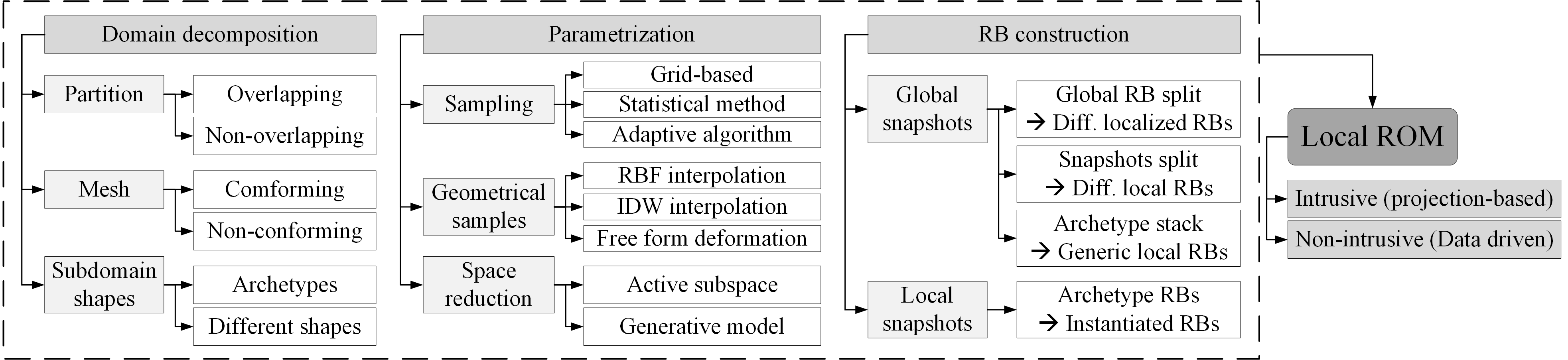}
    \caption{The procedures and classification for preliminaries of constructing ROMs. Abbreviation: RBF (Radial basis function) and IDW (Inverse Distance Weighting).}
    \label{fig:classification_preliminary}
\end{figure}

In short, the three pre-stages for constructing local ROMs are illustrated in Fig. \ref{fig:classification_preliminary}.

\subsection{Coupling techniques}
Once the local RBs are available, the next step is to \emph{glue} them to approximate global solutions. This coupling has gained attention from the research community in recent years. Therefore, various methods have been developed to achieve it. We highlight that the fundamental philosophy of the algorithms is to minimize discontinuities at the adjacent interfaces while satisfying the governing equations and boundary (and initial) conditions. 

Similar to standard ROMs, local ROM approaches can be constructed employing projection-based (Section \ref{sec:projection_based_coupling_algorithms}) or data-driven (Section \ref{sec:Data_driven_Techniques}) methodologies. 

Regarding how reduced systems are constructed and solved, projection-based approaches can be classified into two groups: (i) monolithic schemes (Section \ref{subsec:monolithic}) and (ii) iterative algorithms (Section \ref{subsec:iterative}). The two groups are descried in Section \ref{sec:projection_based_coupling_algorithms}. 

Data-driven techniques are becoming popular due to the rapid development of \emph{Artificial Intelligence} (AI). Moreover, it is a non-intrusive methodology that reduces the complexity of handling PDEs, enhancing the framework's adaptability to different governing equations. Recently, various interpolation/surrogate algorithms, as well as \emph{Neural Networks}, have been intensively employed for constructing local ROMs (see Section \ref{sec:Data_driven_Techniques}). 

Given the diversity of projection-based and data-driven coupling methods, two flowcharts regarding their classifications are presented in Sections \ref{sec:projection_based_coupling_algorithms} and \ref{sec:Data_driven_Techniques}, respectively. We also remind that separate descriptions for multiphysics problems are involved, due to their specific characteristics in coupling.


\section{Preliminaries}
Now, we will explain the three preliminaries: (i) decomposition, (ii) parameterization, and (iii) snapshot collection and dimensionality reduction at the local level.

\subsection{Domain decomposition}
\label{sec:dd_local}
Domain decomposition has been widely used in high-fidelity simulations for parallel computing and multi-physics problems. In a well-known book, Alfio Quarteroni et al. \cite{quarteroni1999domain} discussed the fundamental topics of the theory, numerical analysis, and implementation of domain decomposition approaches.

In the following paragraphs, we focus on the implementation of domain decomposition procedures corresponding to local ROM, and the three classification aspects are clarified with descriptions and sketches.

\subsubsection{Decomposition with/without overlaps}
\label{subsec:decomposition_with_without_overlaps}

There are two main types of partitions: with or without overlaps. The strategy followed to achieve the decomposition into subdomains mainly depends on the coupling algorithms. The latter will be discussed in later sections. Here, we will briefly illustrate the basic characteristics of the two decomposition ideologies. 

\begin{figure}[h]
    \centering
    \begin{subfigure}[b]{0.45\textwidth}
        \centering
        \includegraphics[width=0.5\linewidth]{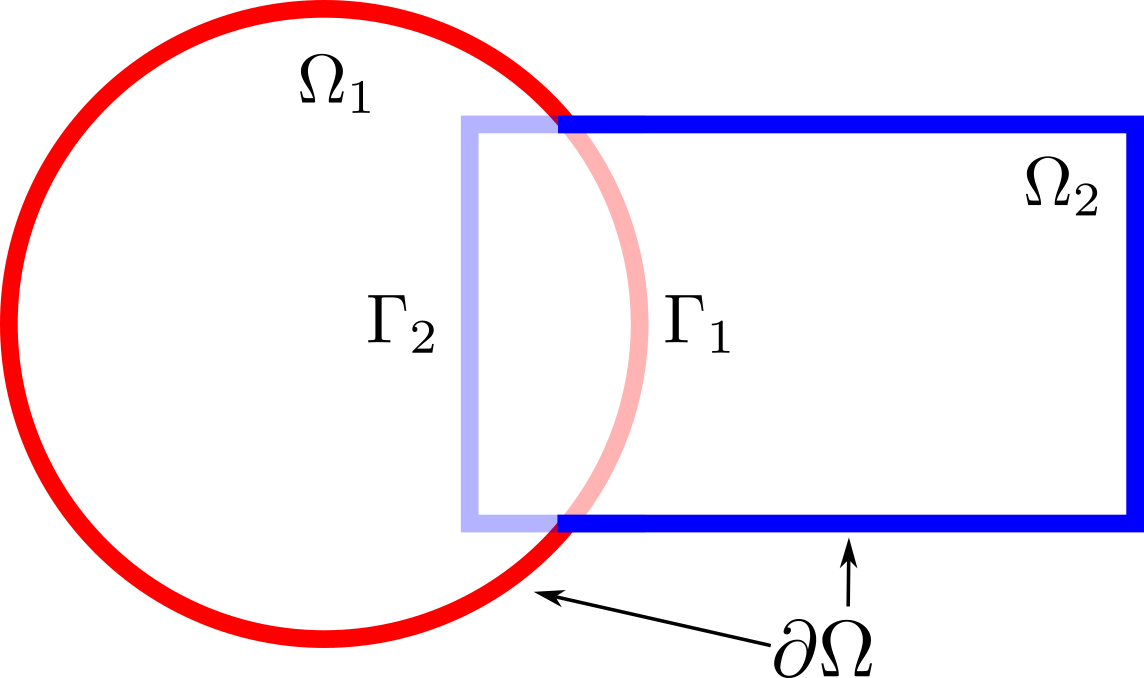}
        \caption{Overlapping decomposition}
        \label{fig:overlapping_decomposition}
    \end{subfigure}
    \begin{subfigure}[b]{0.45\textwidth}
        \centering
        \includegraphics[width=0.5\linewidth]{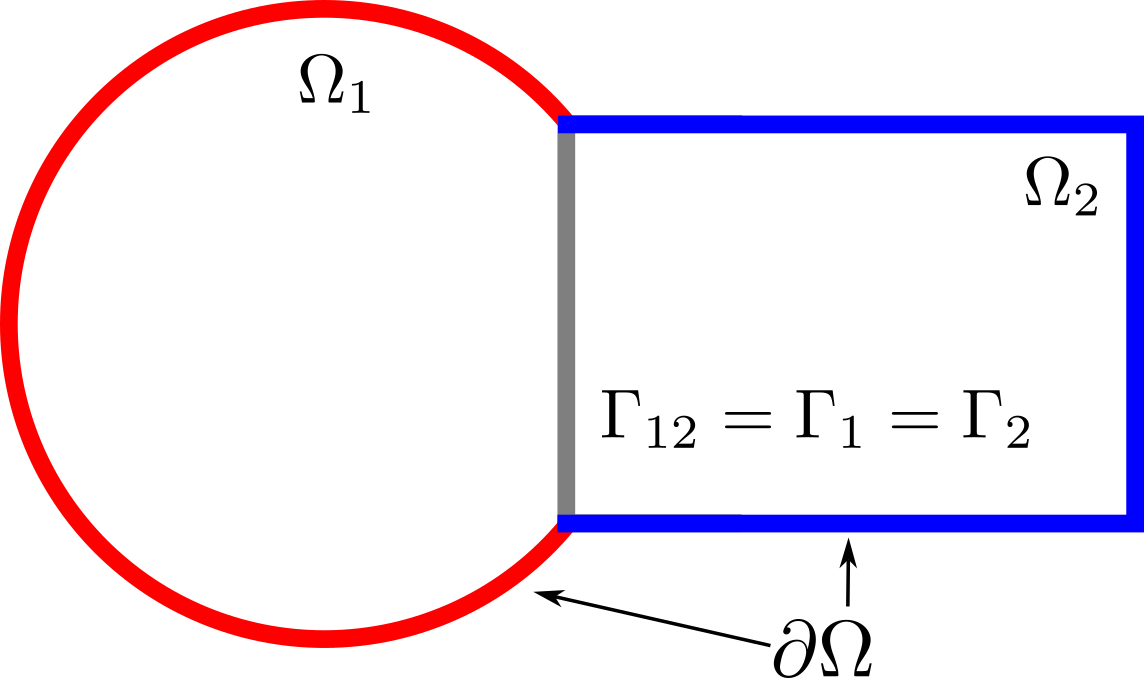}
        \caption{Non-overlapping decomposition}
        \label{fig:nonoverlapping_decomposition}
    \end{subfigure}
    \caption{Schematic of domain decomposition: (a) overlapping; (b) non-overlapping. The global domain $\Omega$ is divided into $\Omega_1$ and $\Omega_2$ parts  with boundary $\partial \Omega$. For overlapping, two interfaces exist, $\Gamma_1$ and $\Gamma_2$; without it a common interaction $\Gamma_{[12]}$.}
    \label{fig:overlap_nonoverlap_scheme}
\end{figure}

The overlapping and non-overlapping decomposition is schematized in Fig. \ref{fig:overlap_nonoverlap_scheme}, where we assume a domain $\Omega$ and its boundary $\partial \Omega$, which is divided into pieces $\Omega_1$ and $\Omega_2$. For the overlapping condition, this division results in two interfaces $\Gamma_1$ and $\Gamma_2$. The procedure without overlapping creates a shared face $\Gamma_{[12]} = \partial \Omega_1 \cap \partial \Omega_2$ \footnote{The subscript $_{[12]}$ is a combined index to designate a single face.}.

The geometry can be split into multiple partitions. For two adjacent subdomains $\Omega_m$ and $\Omega_n$ with overlaps, interfaces are denoted as $\Gamma_{[mn]} = \partial \Omega_m \cap \Omega_n$ and $\Gamma_{[nm]} = \partial \Omega_n \cap \Omega_m$. Nevertheless, a single common internal face is named $\Gamma_{[mn]} = \partial \Omega_m \cap \partial \Omega_n$ for non-overlapping scenarios.

\begin{figure}[h]
    \centering
    \includegraphics[width=0.8\linewidth]{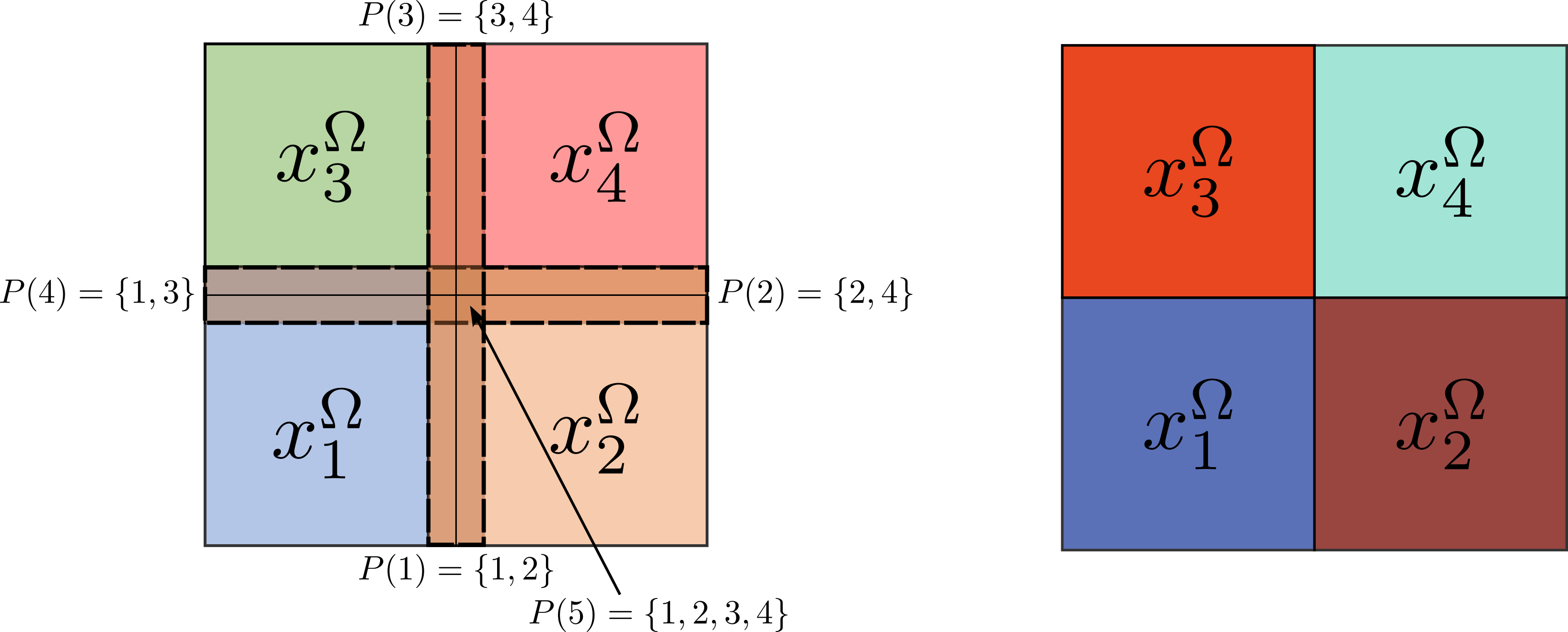}
    \caption{Domain decomposition into $2 \times 2$ subdomains. Left: overlapping partition, redrawn based on \cite{diaz2024fast}. Right: non-overlapping partition, redrawn based on \cite{discacciati2023localized}.}
    \label{fig:2_2_subdomain_decomposition}
\end{figure}

We claim that both methodologies have distinct advantages and limitations. Increasing the shared regions enhances the stability and convergence of coupling techniques \cite{quarteroni2009numerical}. However, overlaps introduce challenges in partitioning. Fig. \ref{fig:2_2_subdomain_decomposition} (left) illustrates a $2 \times 2$ domain decomposition, where the overlapping relationships become increasingly complex \cite{diaz2024fast}. In particular, the central sub-region represents the interaction of all four neighboring subdomains, significantly complicating the construction of Reduced Order Models (ROMs).

In contrast, the non-overlapping approach offers greater flexibility and adaptability, simplifying the subdomain generation process. As depicted in Fig. \ref{fig:2_2_subdomain_decomposition} (right), interfaces exist only between adjacent subdomains, reducing geometric complexity. This approach also allows more flexibility for mesh discretization, where adjacent subdomains can employ either conforming or non-conforming meshes \cite{antonietti2016discontinuous}.

Furthermore, a shared region implies that the same physical phenomena probably occur in the adjoining subdomains \cite{sambataro2022component}. However, non-overlapping techniques enable the coupling of subdomains with different governing PDEs. Despite this advantage, the simplicity of the non-overlapping approach is not free. The approaches are currently limited by the extensive complexities of developing robust algorithms for coupling subdomains at the connected interfaces \cite{discacciati2024overlapping}.

\subsubsection{Conforming and non-conforming meshes}
\label{subsec:conforming_non-conforming}

\begin{figure}[h]
    \centering
    \begin{subfigure}[b]{0.45\textwidth}
        \centering
        \includegraphics[width=0.5\linewidth]{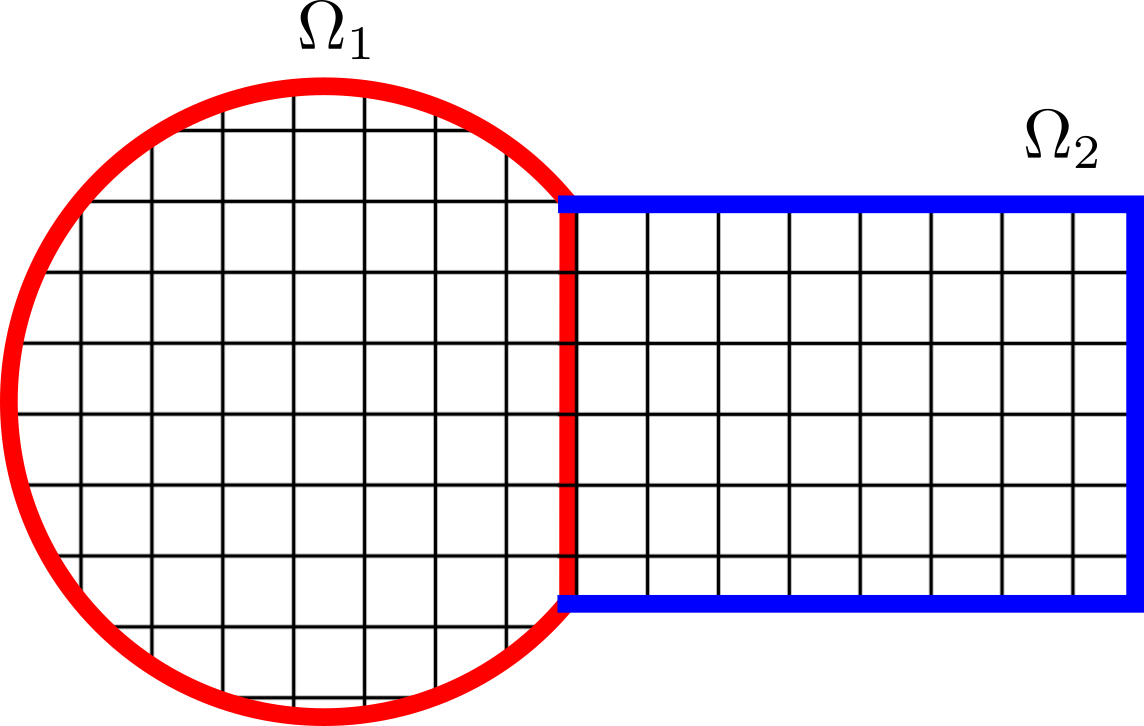}
        \caption{}
        \label{fig:conforming_mesh}
    \end{subfigure}
    \hspace{1em}
    \begin{subfigure}[b]{0.45\textwidth}
        \centering
        \includegraphics[width=0.5\linewidth]{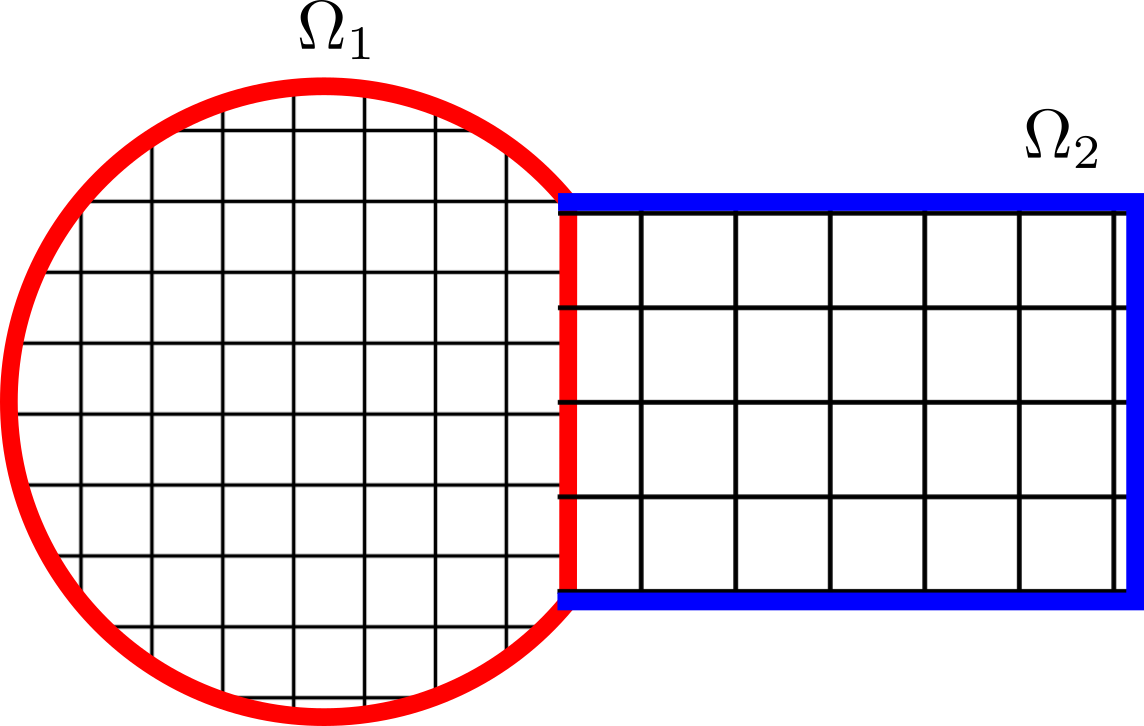}
        \caption{}
        \label{fig:nonconforming_mesh}
    \end{subfigure}
    \caption{(a) Conforming and (b) non-conforming meshes for adjacent subdomains.}
    \label{fig:conforming_nonconforming_mesh}
\end{figure}

The conforming and non-conforming high-resolution meshes along an interface are displayed in Fig. \ref{fig:conforming_nonconforming_mesh}. This characteristic significantly affects the algorithms for coupling subdomains. In some cases, additional mutual interpolation procedures are required for the reconstruction of global solutions \cite{zappon2023efficient}.

\subsubsection{Individual and generic decomposition}
\label{subsec:individual_generic_decomposition}

As indicated in Section \ref{sec:intro}, a global model can be decomposed into many independent subdomains or assembled by a few types of generic components (archetypes). The usage of the approaches depends on the geometric characteristics and the coupling algorithms. We will clarify the two strategies here.

\paragraph{Individual decomposition}
\label{subsubsec:individual_decomposition}

\emph{Individual decomposition} denotes a domain $\Omega$ composed of a finite number of open and bounded partitions, $\Omega = \cup_{m=1}^{N_\Omega} \Omega_m$. The set containing all interfaces is defined as $\Gamma$, and $\Gamma = \cup_{r=1}^{N_\Gamma} \Gamma_r$. The global boundary of $\Omega$ is assigned as $\partial \Omega$, and $\partial \Omega \cap \Gamma = \emptyset$. The local boundary of $\Omega_m$ is similarly noted as $\partial \Omega_m$,
and we have $\partial \Omega \cap \Gamma \neq \emptyset$. Indeed, the interfaces coincide with the specific boundary of partitions. Generally, a partition $ \Omega_m $ connects to several interfaces $\Gamma_r$. Thus, we define the set of interfaces connected to $ \Omega_m $ as $ \mathcal{I}_m = \{ \left. \Gamma_{r} \right| \Gamma_{r} \subset \partial \Omega_m\}$. The $i^\text{th}$ item in $\mathcal{I}_m$ is noted as $\Gamma_{m,i}$

\begin{figure}
    \centering
    \begin{subfigure}[b]{0.42\textwidth}
        \includegraphics[width=\linewidth]{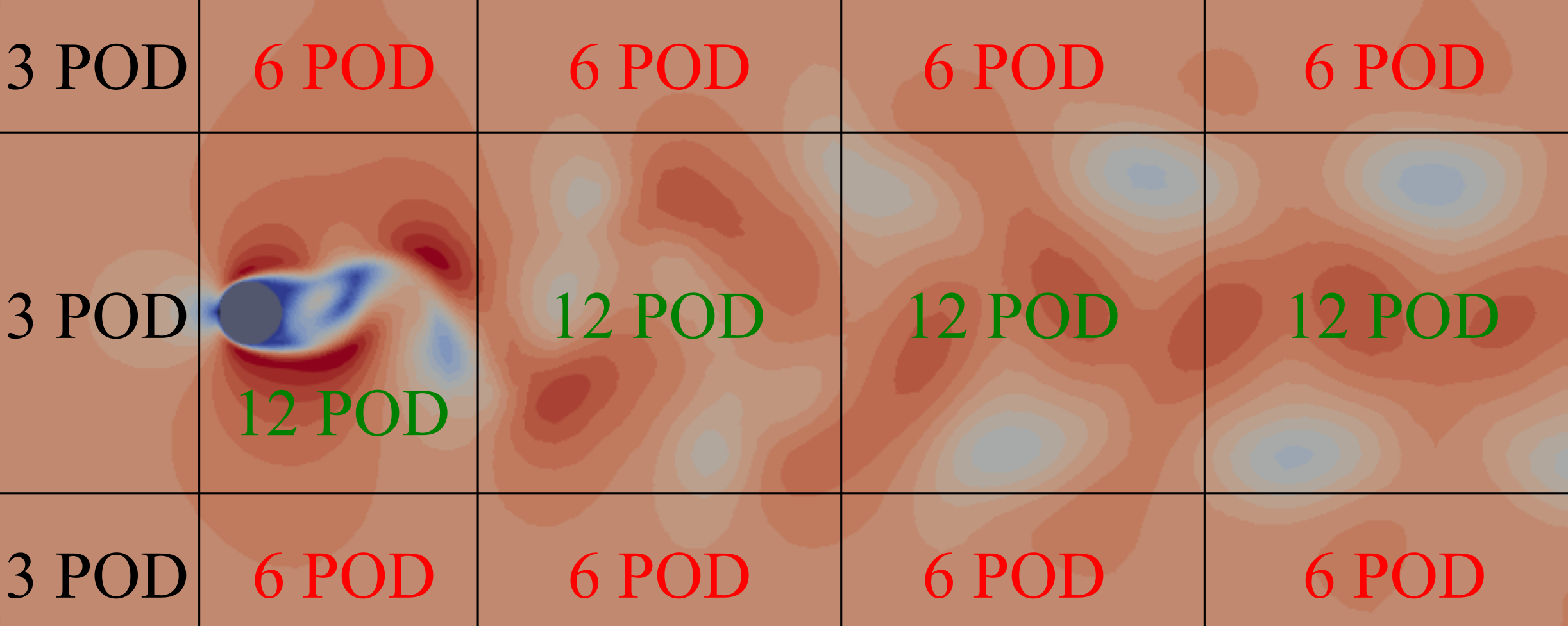}
        \caption{}
        \label{fig:individual_division_flow_cylinder}
    \end{subfigure}
    \hspace{1em}
    \begin{subfigure}[b]{0.25\textwidth}
        \includegraphics[width=\linewidth]{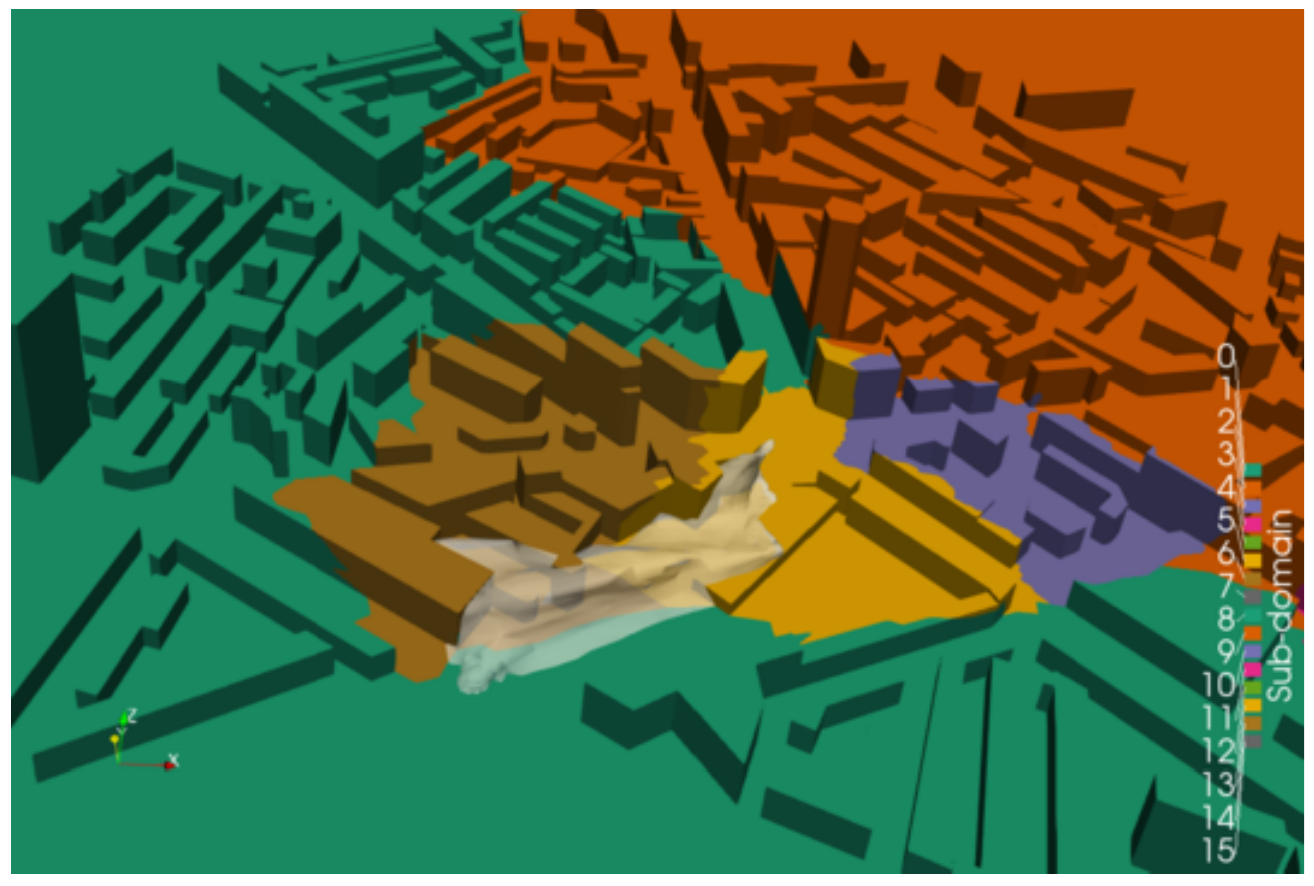}
        \caption{}
        \label{fig:individual_division_city_pollutant}
    \end{subfigure}
    \caption{Examples of domain divisions for problems without repeating geometric parts. (a) Kármán vortex street for flow around a cylinder, redrawn based on \cite{xiao2017domain}. (b) Pollutant transport in an urban environment, taken with permission from \cite{arcucci2020domain}, copyright owned by IOS Press. }
    \label{fig:individual_divisions_without_repeating}
\end{figure}

\begin{figure}[h]
    \centering
    \begin{subfigure}[b]{0.45\linewidth}
        \includegraphics[width=\linewidth]{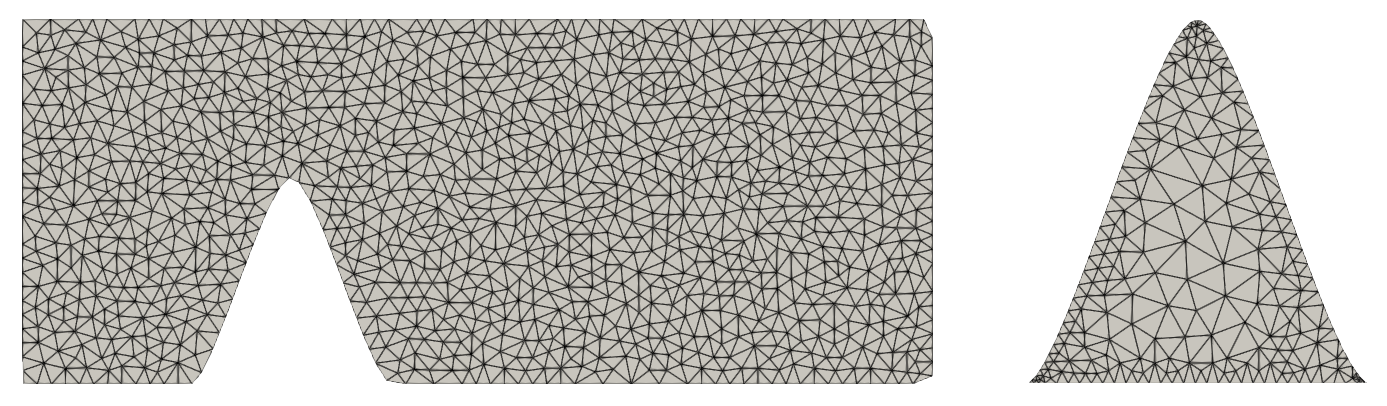}
        \caption{Stokes and linear elasticity \cite{discacciati2023localized}.}
    \end{subfigure}
    \hspace{1em}
    \begin{subfigure}[b]{0.4\linewidth}
        \includegraphics[width=\linewidth]{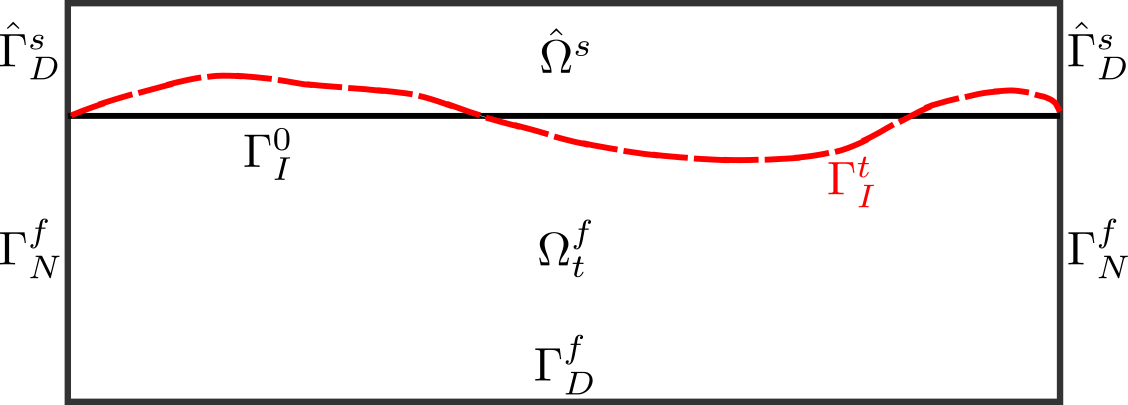}
        \caption{Navier--Stokes and linear elasticity \cite{prusak2023application}.}
    \end{subfigure}
    \caption{Decomposition for fluid-structure interaction problems. (a) Left: fluid; right: solid. Redrawn based on \cite{discacciati2023localized}. (b) $\hat{\Omega}_s$ and $\Omega_t^f$ are solid and fluid, respectively. Redrawn based on \cite{prusak2023application}.}
    \label{fig:individual_divisions_fsi}
\end{figure}

For domains without symmetry and periodicity, the splitting results in subdomains $\Omega_m$ of various shapes. Fig. \ref{fig:individual_divisions_without_repeating} displays two examples, i.e., laminar flow around a cylinder \cite{xiao2017domain} and pollutant transport within a city \cite{arcucci2020domain}. Furthermore, in case the model contains sub-regions governed by different PDEs, the entire geometry is divided considering the physical boundaries. The FSI problems are typical instances consisting of both fluid and solid partitions. Two FSI studies are presented in Fig. \ref{fig:individual_divisions_fsi}, which illustrate the interactions of a linear elastic structure and a fluid flow governed either by Stokes or incompressible Navier-Stokes equations.

\paragraph{Generic decomposition}
\label{subsubsec:Generic_decomposition}

The second strategy is called \emph{generic decomposition}. It takes its name to contrast with the previous \emph{individual} strategy. Assume we have $\hat{N}_\Omega$ reference subdomains, and the $k^{\text{th}}$ of them is denoted as $\hat{\Omega}^k$. We use the \emph{hat} symbol $\hat{\cdot}$ to characterize reference domains. Each archetype $\hat{\Omega}^k$ has $N_k$ instantiations, indexed by $m=1,\dots,N_k$, named as $\Omega^k_m$. These instances are obtained through a parametrized geometrical transformation $T^k (\boldsymbol{\xi}_m^k)$, namely,
\begin{equation}
\label{eq:geometric_transformation}
    \Omega^k_m \coloneqq T^k (\boldsymbol{\xi}^{k}_m) \cdot \hat{\Omega}^k.   
\end{equation}

We also note the total number of partitions with $N_\Omega$, and obviously, it holds that $N_\Omega = \sum_{k=1}^{\hat{N}_\Omega} N_k$. Consequently, the entire computational domain is expressed as
\begin{equation*}
    \Omega = \bigcup_{k=1}^{\hat{N}_\Omega} \bigcup_{m=1}^{N_k} \Omega^k_m.
\end{equation*}
Similarly, the interface set of $\Omega^k_m$ is given by $\mathcal{I}_m^k$.

\begin{figure}[h]
    \centering
    \begin{subfigure}[b]{0.58\linewidth}
        \includegraphics[width=\linewidth]{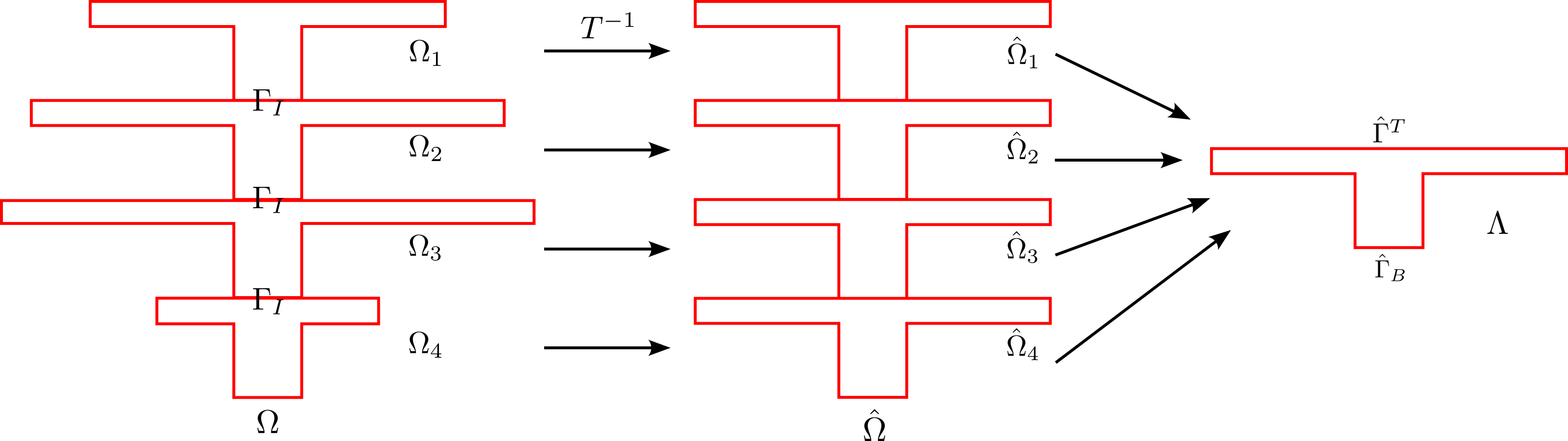}
        \caption{A multi-level thermal fin problem \cite{iapichino2016reduced}.}
        \label{fig:generic_division_thermal_fin}
    \end{subfigure}
    \hfill
    \begin{subfigure}[b]{0.4\linewidth}
        \includegraphics[width=\linewidth]{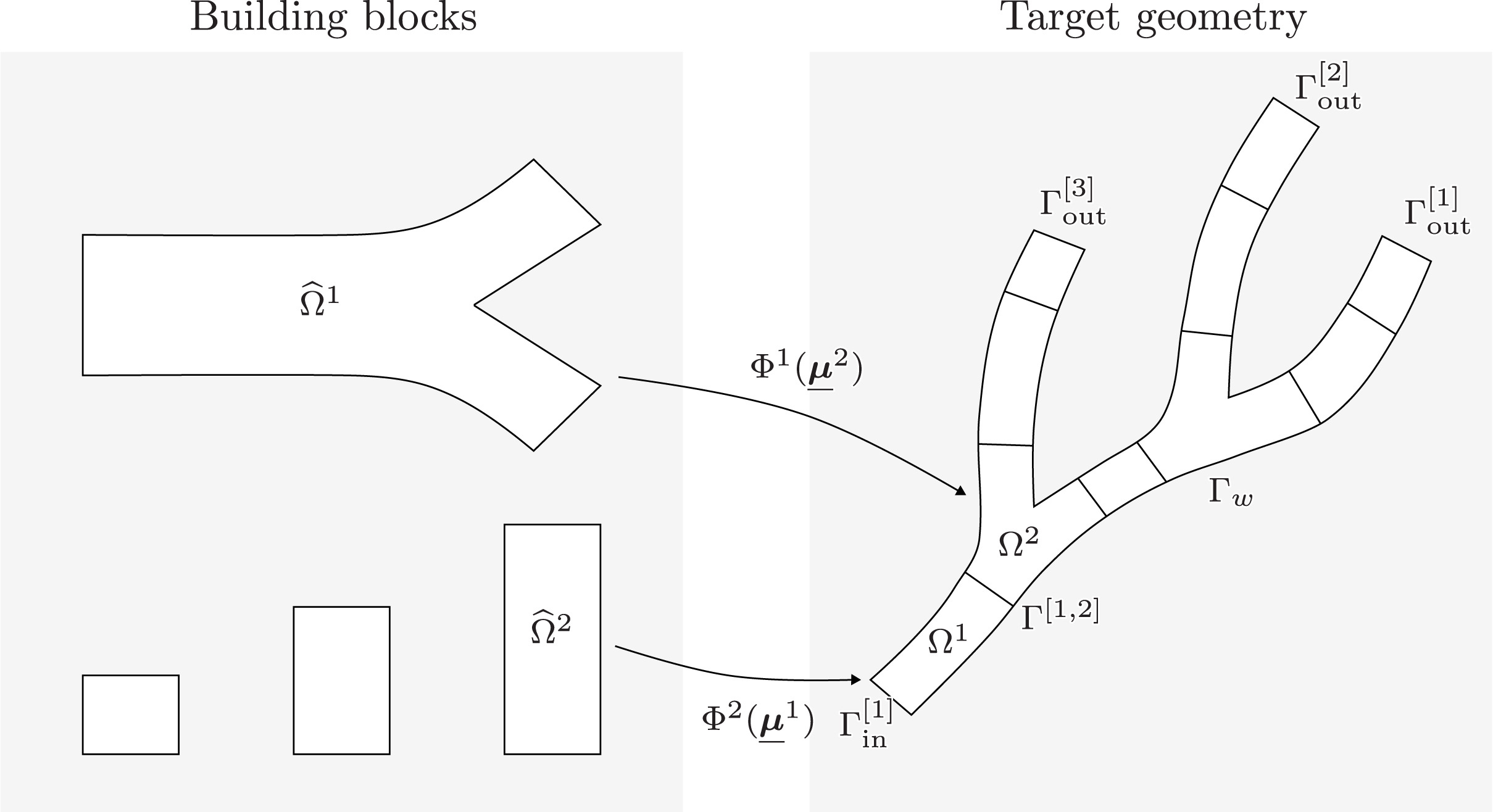}
        \caption{A blood vessel model \cite{pegolotti2021model}.}
        \label{fig:generic_division_blood_vessel}
    \end{subfigure}
    \caption{Several generic parts can be used to assemble entire models. The real blocks are obtained by geometric transformation from reference parts. (a) A multi-level thermal fin and the geometrical mapping to the reference block. Redrawn based on \cite{iapichino2016reduced}. (b) A blood vessel formed from blocks, taken with permission from \cite{pegolotti2021model}, copyright owned by Elsevier.}
    \label{fig:generic_devisions}
\end{figure}

Many realistic industrial facilities are designed to be periodic and symmetric, as these characteristics enhance process consistency, efficiency, and cost-effectiveness. Such structural regularity can be exploited to simplify simulations and significantly reduce computational costs at both FOM and ROM levels. The thermal fin pattern displayed in Fig. \ref{fig:generic_devisions} is a typical design for heat exchangers \cite{wang2016layer}. In \cite{iapichino2016reduced}, the authors analyzed the heat transfer in a multiple-level thermal fin, in which each block is transformed from a generic part. A similar concept appears in biological modeling, as presented in \cite{pegolotti2021model}, where the blood vessel structure is composed of four reference blocks. The applications demonstrate that a few archetypes can assemble even complex flow networks.

\subsection{parameterization techniques}
\label{sec:parameterization}

ROM is a real-time alternative to FOM to generate a dataset that can represent solutions within a specific range. Now, we will discuss three aspects regarding parameterization.

\subsubsection{Sampling strategies}
\label{subsec:Sampling_strategies}
Sampling consists of selecting a subset of the observed data, which can be used to estimate the features of the entire dataset \cite{thompson2012sampling}. It is important for applications like design optimization and uncertainty quantification. The theory of general sampling is extensive and beyond the scope of this work (see \cite{murthy1967sampling, thompson2012sampling} for details). Thus, we will only outline some studies that correspond to ROM, for the sake of completeness. Interested readers can also turn to  \cite{benner2015survey, quarteroni2015reduced, liu2024application} for a relatively concise summary of the topic.

\begin{figure}[h]
    \centering
    \includegraphics[width=0.25\linewidth]{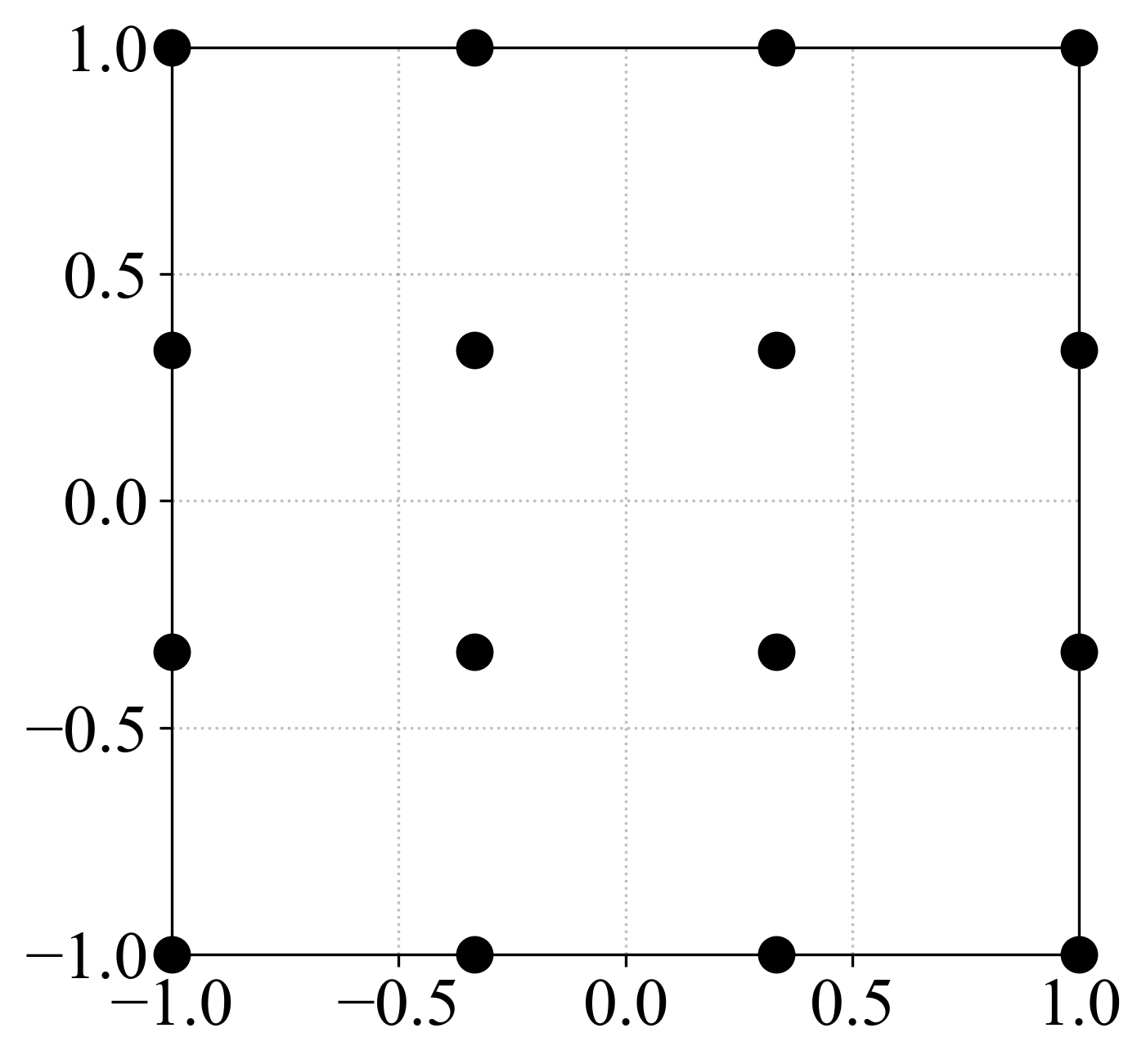}
    \hspace{1.5em}
    \includegraphics[width=0.25\linewidth]{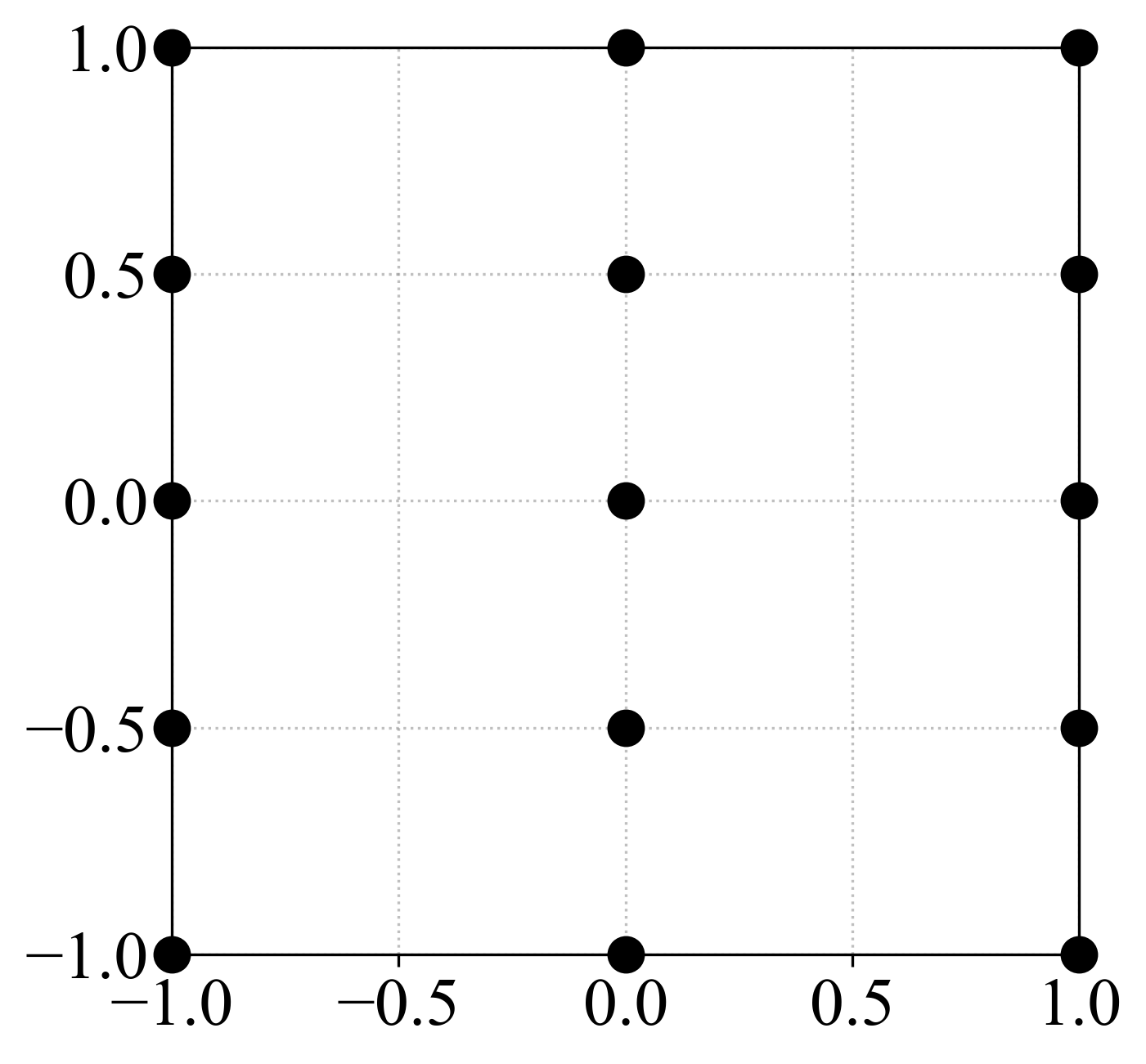}
    \hspace{1.5em}
    \includegraphics[width=0.25\linewidth]{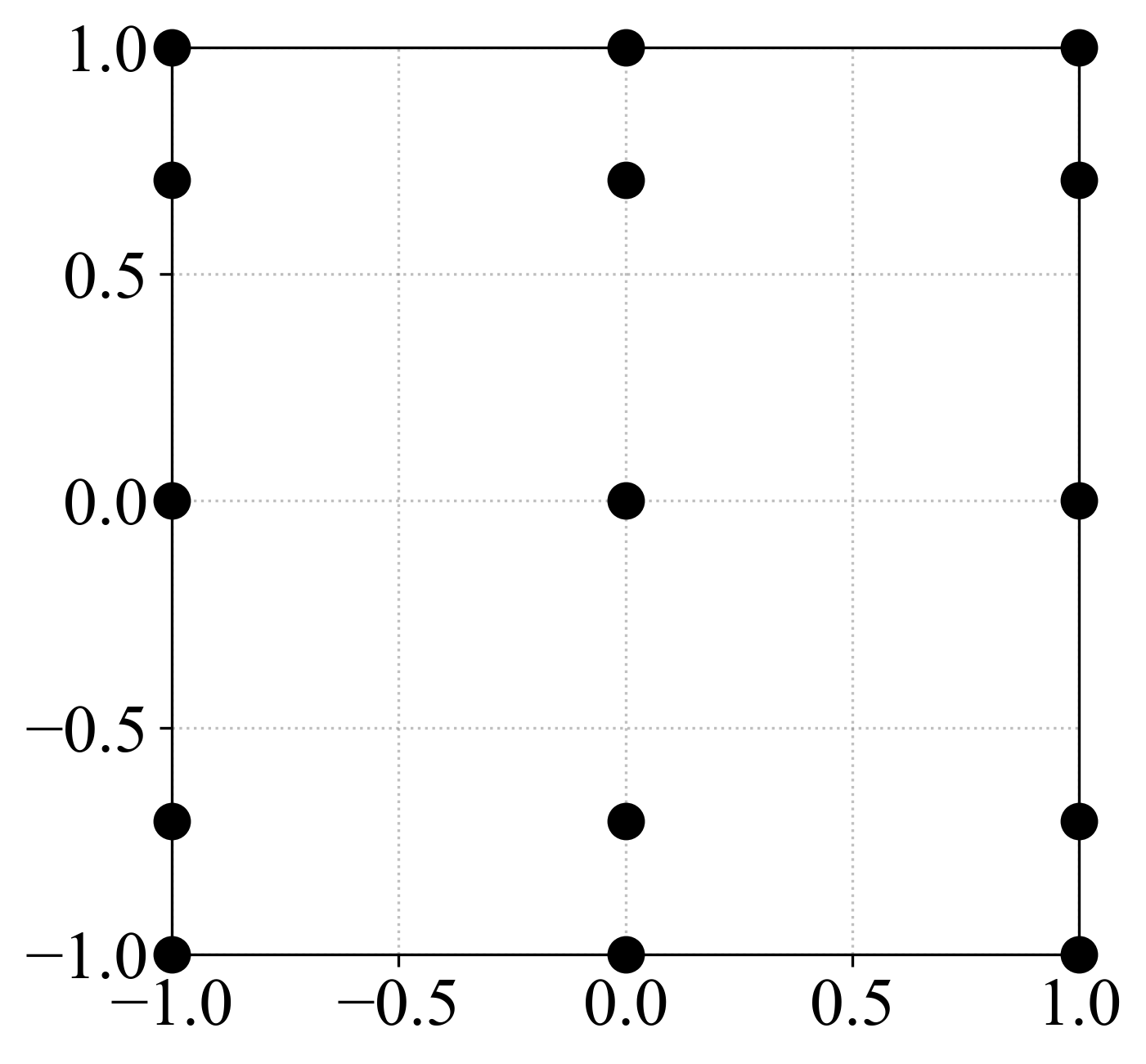}

    \includegraphics[width=0.25\linewidth]{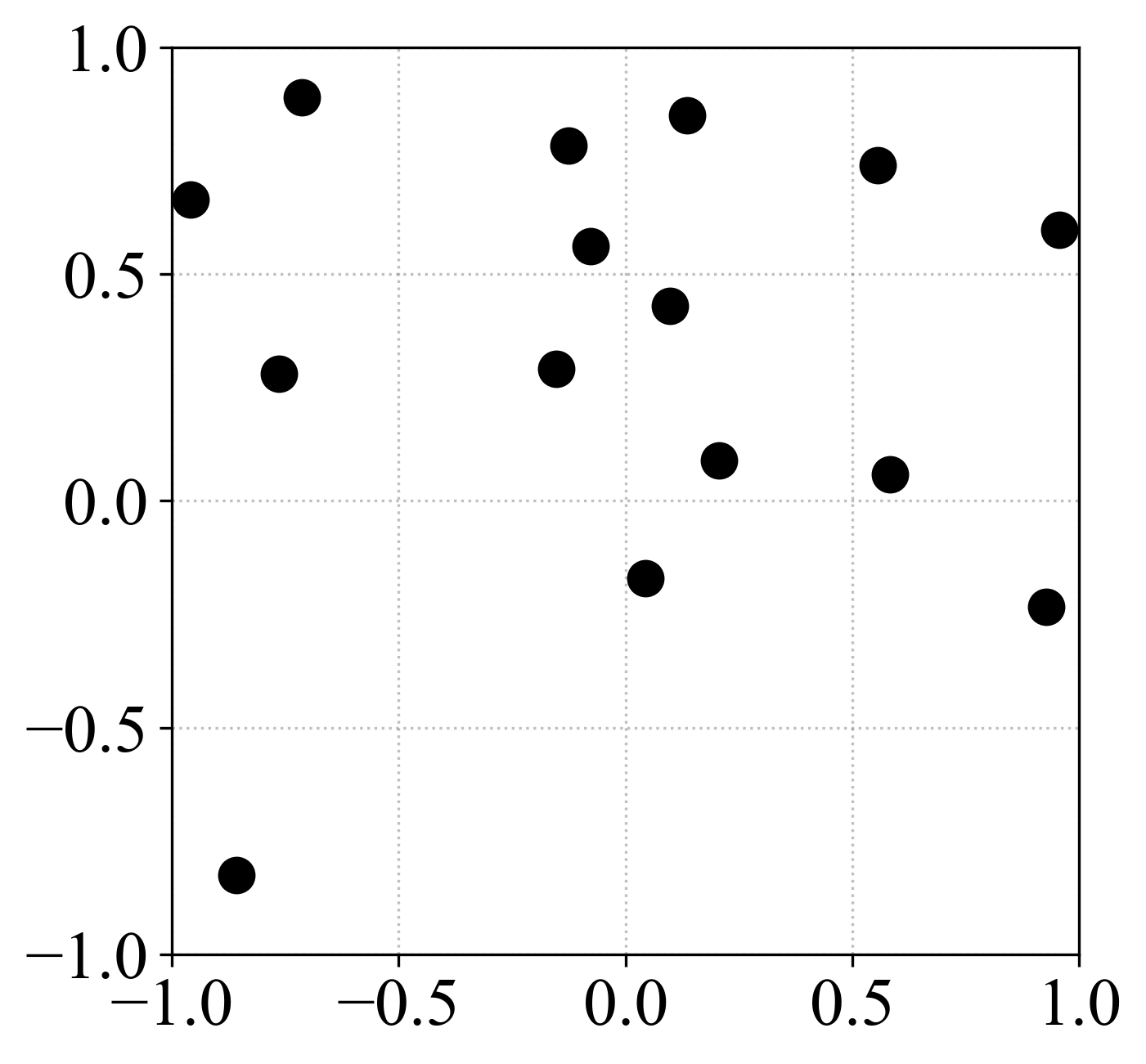}
    \hspace{1.5em}
    \includegraphics[width=0.25\linewidth]{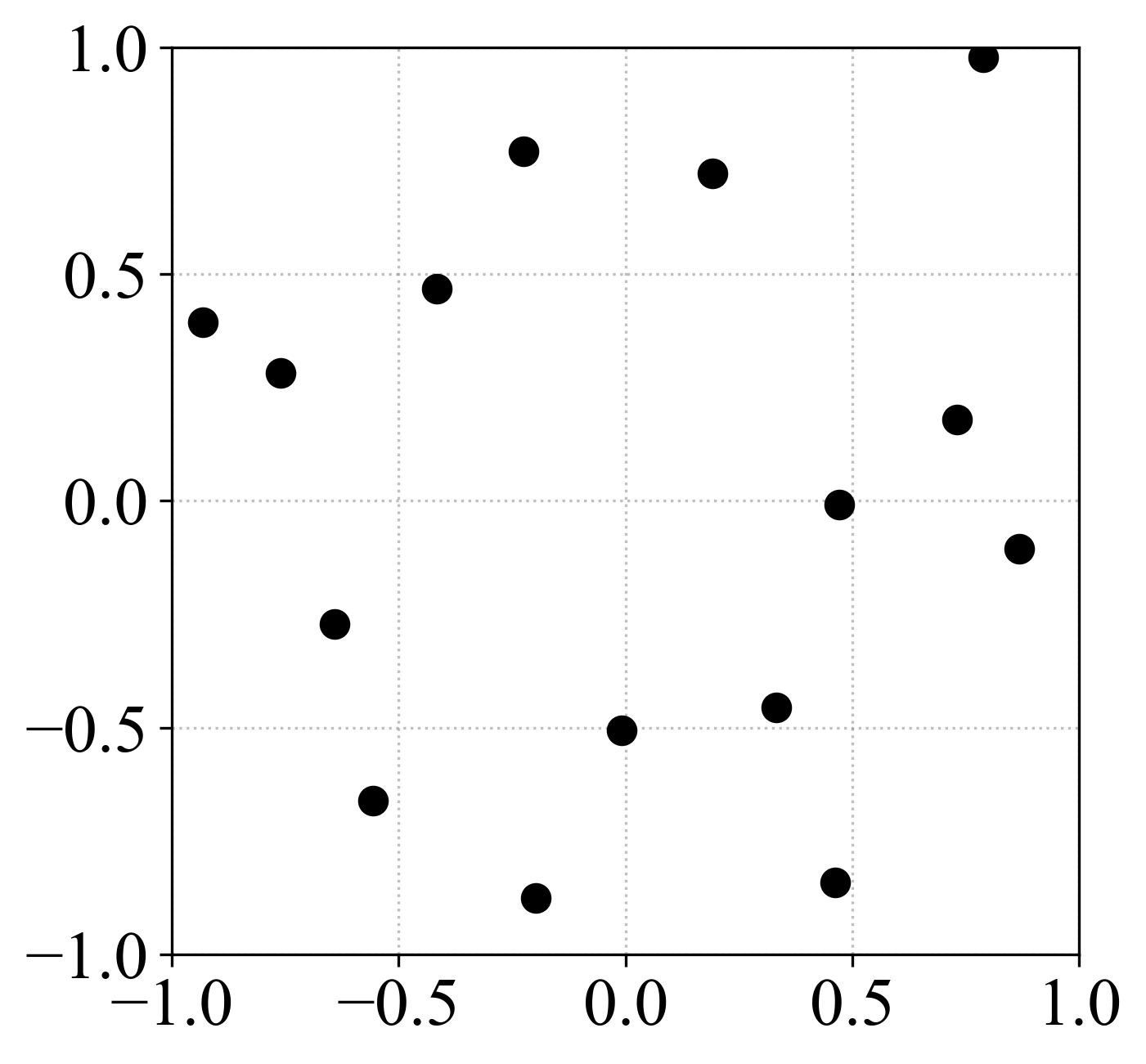}
    \hspace{1.5em}
    \includegraphics[width=0.25\linewidth]
    {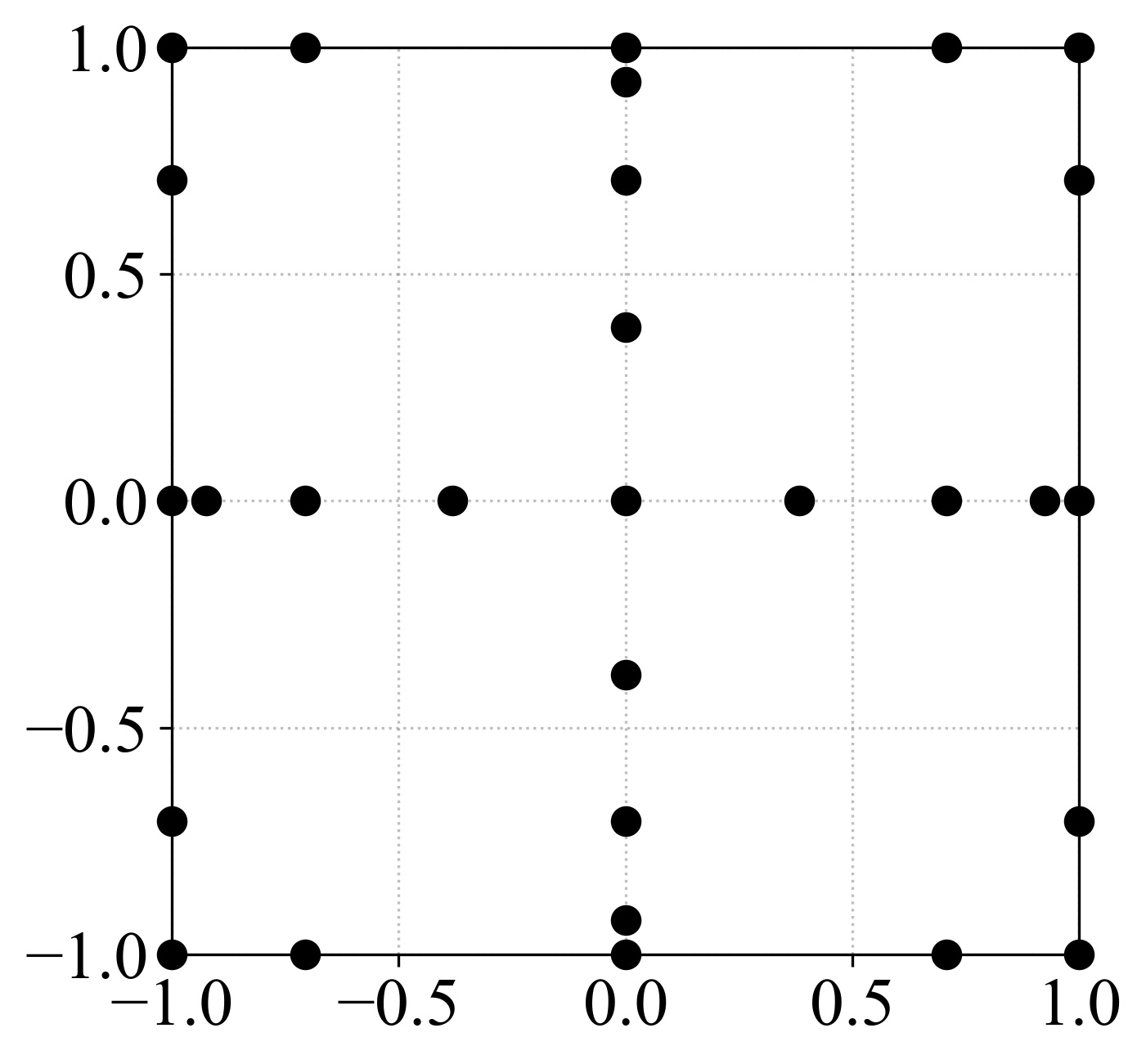}
    \caption{Sampling strategies in the space $\mathcal{P} = \left[ -1,1 \right]^2$. Grid-based (first row, from left to right): $4 \times 4$ uniform, $3 \times 5$ uniform, and $3 \times 5$ Clenshaw-Curtis. Statistical (second row, from left to right): random (15 points), Latin hypercube (15 points), and Smolyak sparse grid with Clenshaw-Curtis points. Figures redrawn based on \cite{quarteroni2015reduced}.}
    \label{fig:sampling_strategies}
\end{figure}

The choice of sampling strategies depends mostly on the dimension, $N_{\mathcal{P}}$, of the parameter space $\mathcal{P}$ \cite{benner2015survey}. \emph{Grid-based} approaches are the standard option for problems with low dimensions (typically $N_{\mathcal{P}} \leq 3$) \cite{bui-thanh2008model, quarteroni2015reduced}. For example, three grid-based techniques for $\mathcal{P} = \left[ -1,1 \right]^2$ are plotted in the first row of Fig. \ref{fig:sampling_strategies}. Non-uniform selections, such as \emph{Clenshaw-Curtis} points \cite{clenshaw1960method}, are an alternative. By employing these strategies, one can choose the same or different numbers of values in each dimension. 

When more parameters are considered, the number of grid points grows exponentially. Statistical algorithms like random (e.g., \emph{Monte Carlo}) and \emph{Latin hypercube} sampling \cite{mckay2000comparison}, as well as Smolyak sparse grid, are well-suited for these conditions \cite{bui-thanh2008model,quarteroni2015reduced} (see second row of Fig. \ref{fig:sampling_strategies}). As an additional alternative, we include the probabilistic sample method employing centroidal Voronoi tessellations \cite{ju2002probabilistic}.

The situation becomes more challenging when the dimension $N_{\mathcal{P}} > 10$. The common sampling techniques mentioned above are not feasible in these conditions. Namely, it is impossible to balance numerical cost with a good representation of $\mathcal{P}$. Therefore, more sophisticated problem-specific \emph{adaptive algorithms} should be utilized, as mentioned in \cite{bui-thanh2008model, benner2015survey,liu2024application}.

The adaptive workflow of these methods involves enriching the samples in sequence based on specific principles. The \emph{greedy method} is a widely used iterative approach for this task \cite{quarteroni2014reduced, hesthaven2016certified, benner2020snapshot}. It has been applied to various PDEs, such as elliptic \cite{veroy2003posteriori, rozza2008reduced}, parabolic \cite{grepl2005posteriori}, Stokes \cite{negri2015reduced}, and Navier-Stokes equations \cite{veroy2005certified, chen2018greedy, sleeman2025greedy}. There are also several alternative approaches, such as local sensitivity analysis \cite{bond2007piecewise}, the trust region algorithm \cite{fahl2000trust}, and a gradient-based method \cite{gong2018data}. The performance of different adaptive algorithms for the application of the ROM is compared in \cite{liu2024application, karcher2022adaptive}.

\subsubsection{Geometrical parameterization}
\label{subsec:Geometrical_parameterization}

Geometrical parameterization is important in industrial practice for design and optimization. See descriptions and examples in \cite{quarteroni2014reduced, quarteroni2015reduced, benner2020snapshot, benner2020applications, rozza2018advances, rozza2024real}. In terms of implementation, it is worth mentioning the Python package PyGem \cite{tezzele2021pygem}, which employs various techniques for generating geometrical samples. For the numerical analysis of propeller blades in particular, the geometrical parameterization and bottom-up construction can be performed using BladeX (Python Blade Morphing) \cite{gadalla19bladex}.

Mostly, geometrical parameterization aims to achieve two goals: (i) shape optimization; (ii) efficient simulation of multiple shapes. The former one is studied for many applications, including aircraft wings \cite{salmoiraghi2016advances}, marine propellers \cite{ivagnes2024shape}, and ship hulls \cite{demo2021hull}. 

In the other condition, solutions for numerous deformed geometries are required. For instance, in cardiovascular and blood system simulations, the shape of the computational domain is patient-specific \cite{siena2024accuracy, pegolotti2021model}.  

The procedure for geometrical parameterization is not unique. Different methods have been developed to meet specific requirements. For discussion, we identified two subcategories: \emph{Interpolation technique} and \emph{Free form deformation}.

\paragraph{Interpolation technique}

The general procedure for generating samples via interpolation can be summarized in three consecutive steps: (i) selecting parameters; (ii) deforming a reference geometry with regard to different parameters to generate a series of shapes; (iii) obtaining meshes from the reference one by interpolating the deformation of shapes. To clarify, we will provide three examples. 

Firstly, we regard the simple thermal fin problem displayed in Fig. \ref{fig:generic_division_thermal_fin}, in which the fins differ only in their length \cite{maday2004reduced, iapichino2016reduced}. 

\begin{figure}[h]
    \centering
    \includegraphics[width=0.35\linewidth]{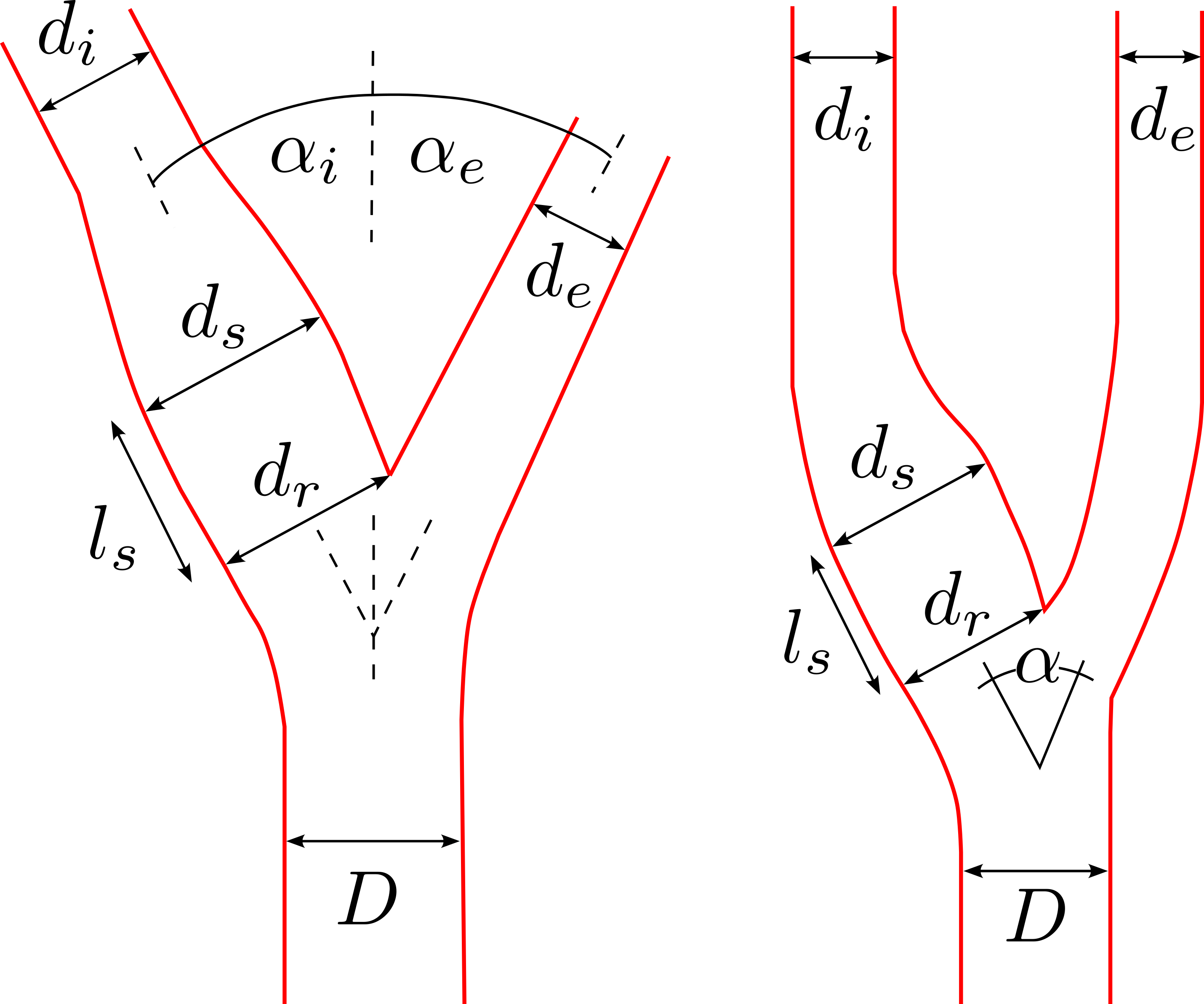}
    \caption{Parameters for a carotid artery bifurcation \cite{bressloff2007parametric, lee2008geometry, ding2001flow}. The shape depends on flow channel diameters, branch curvatures, and bifurcation angles. Figures redrawn based on \cite{manzoni2012model}.}
    \label{fig:carotid_artery_bifurcation}
\end{figure}

Secondly, we consider the complex model of Fig. \ref{fig:carotid_artery_bifurcation} that shows the parametric geometry of the human carotid artery bifurcation \cite{bressloff2007parametric, lee2008geometry, ding2001flow}. To characterize the domains, the authors defined several specific quantities, i.e., flow channel diameters, branch curvatures, and bifurcation angles. 

Lastly, a similar strategy is proposed by Anna Ivagnes et al \cite{ivagnes2024shape} for the optimization of marine propellers. The parameters are assigned at the level of a blade section, and thus, they determine the entire propeller.

\begin{figure}[h]
    \centering
    \includegraphics[width=\linewidth]{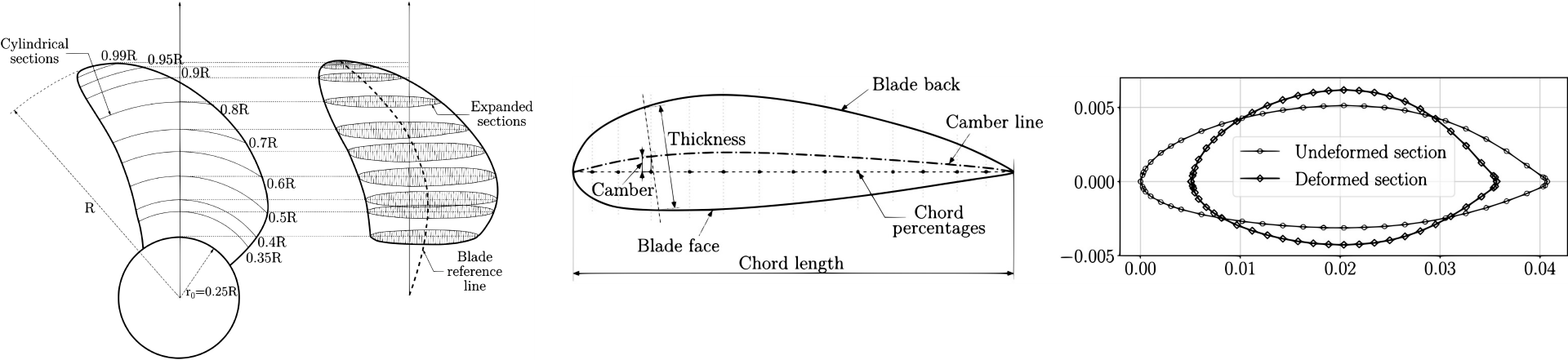}
    \caption{Definition of the parameters of a marine propeller's blade. The blade is represented by a series of sections, and each piece can be determined by several values. Figures taken with permission from \cite{ivagnes2024shape}, copyright owned by John Wiley \& Sons, Inc.}
    \label{fig:propeller_parameterization}
\end{figure}

Once obtain shape samples, we can compute the meshes via interpolation. \emph{Radial Basis Function} (RBF) interpolation is a widely used technique. Let us briefly describe it. Assume that input variables $x_i$ and results $y_i$, with $i = 1, \dots, N$ are known. Our objective is to approximate the function $f$, achieving $y_i = f(x_i)$. In this technique, the interpolation can be written as
\begin{equation}
\label{eq:RBFformula}
    y_i = \sum_{j=1}^{N} w_j \phi (\left\| x_i - x_j \right\|) \quad \text{ with } \quad j=1, \dots, N,
\end{equation}
where $\phi$ denotes a RBF that is weighted by coefficients $w_j$, $j = 1, \dots, N$. The unknowns of $w_j$ are computed by solving equation \eqref{eq:RBFformula} with $N$ pairs of $x_i$ and $y_i$. 

Once available, RBF interpolation can be applied to estimate $y$ corresponding to the numerous new $x$. The expression of $\phi$ is not unique. More details are shown in \cite{buhmann2000radial, schaback2007practical}. 

\begin{figure}[h]
    \centering
    
    \begin{subfigure}[b]{0.45\linewidth}
        \includegraphics[width=\linewidth]{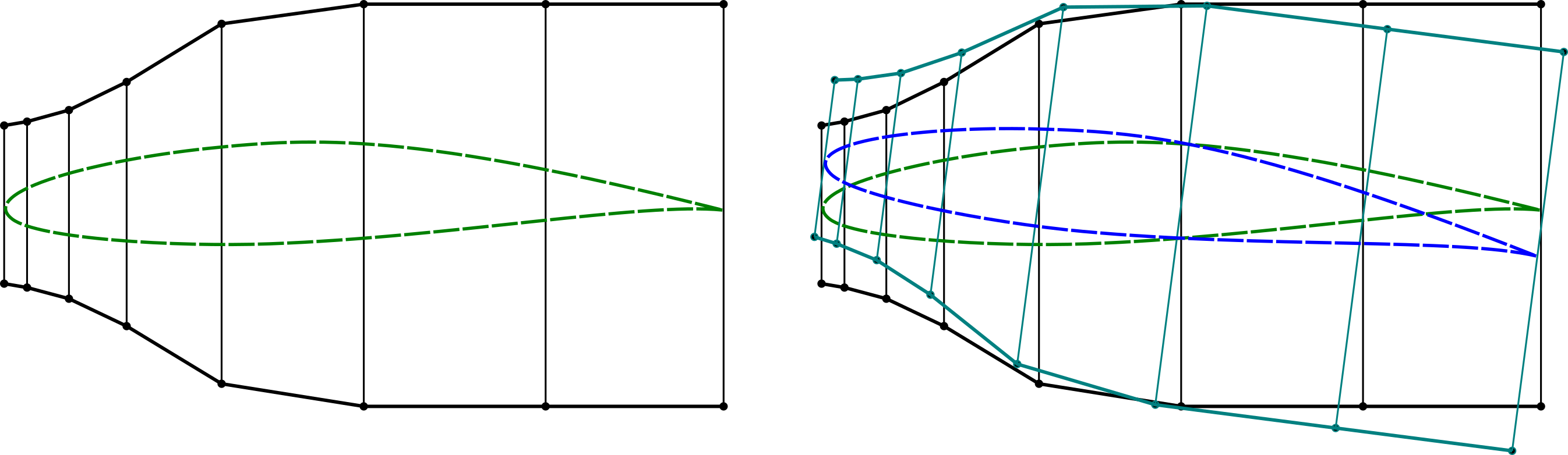}
        \caption{Deformation of an airfoil \cite{morris2008cfdbased}.}
        \label{fig:deformation_control_points_a}
    \end{subfigure}
    \hfill
    \begin{subfigure}[b]{0.45\linewidth}
        \includegraphics[width=\linewidth]{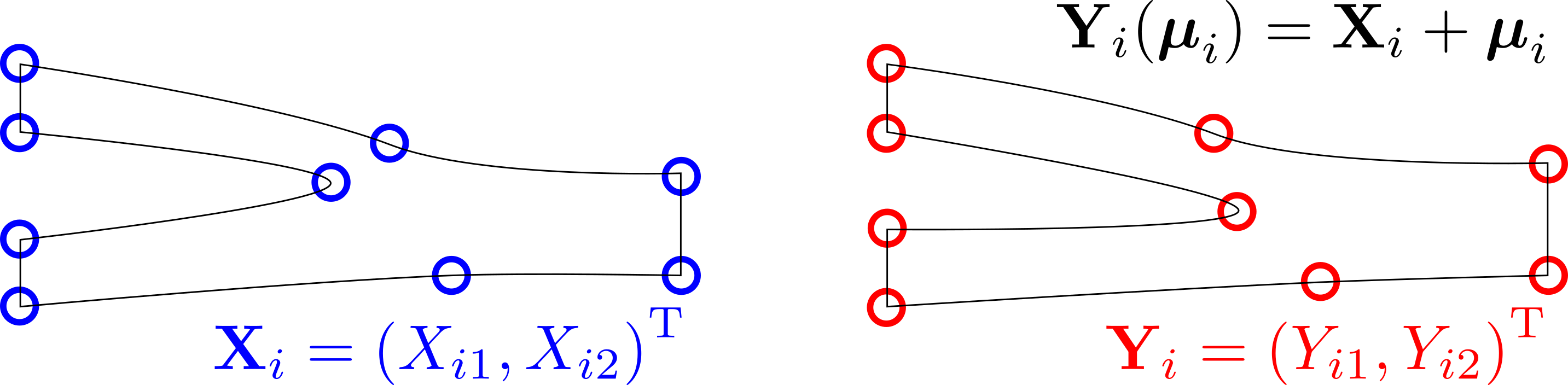}
        \caption{Deformation of a branch \cite{manzoni2012model}.}
        \label{fig:deformation_control_points_b}
    \end{subfigure}
    \caption{Generation of geometric samples  with control points. Applications for (a) an aircraft wing and (b) a bypass problem. (a) is redrawn based on \cite{morris2008cfdbased}; (b) is redrawn based on \cite{manzoni2012model}.}
    \label{fig:deformation_control_points}
\end{figure}

Note that the cases illustrated above employ problem-specific configurations, which are confined to individual applications. Now, we introduce a general strategy. This generates the geometrical shapes by the position of certain control points \cite{morris2008cfdbased, manzoni2012model}. The ideology is plotted in Fig. \ref{fig:deformation_control_points}. The set of control points is defined either surrounding or on the surfaces/edges of the model. The translations of points are regarded as parameters. The deformed outer face and mesh can be interpolated using RBF. 

Note that the geometric interpolation techniques can also be incorporated to approximate the mesh motion in FSI problems. For example, the usage of RBF is presented in \cite{forti2012Comparison, forti2014efficient}. \emph{Inverse Distance Weighting} (IDW) is an alternative to RBF. The results in \cite{witteveen2009explicit, witteveen2010Explicit} demonstrated that IDW is comparable in accuracy to RBF and faster in high fidelity CFD. The combination of ROM and IDW is presented in \cite{forti2012Comparison, damario2014Reducedorder}. However, the interpolation approaches become costly for numerous points. Thus, an adaptive selective RBF method is presented in \cite{forti2012Comparison}, and the \emph{Selective IDW} (SIDW) (also known as Reduced IDW) algorithm is proposed in \cite{damario2014Reducedorder, ballarin2019podselective}. The fundamental idea of the two variants is the same, in which only a subset of points is optimally identified for performing interpolation. See more details about these investigations in Section \ref{subsubsec:Monolithic_FSI}.

\paragraph{Free form deformation}
We will introduce another general and simple method that is well-suited for the geometrical transformation of 2D and 3D objects, namely the \emph{Free Form Deformation} (FFD). Note that the FFD also involves the movement of control points. We will now outline the procedure and indicate the differences compared to previous interpolation techniques.

\begin{figure}[h]
    \centering
    \includegraphics[width=0.4\linewidth]{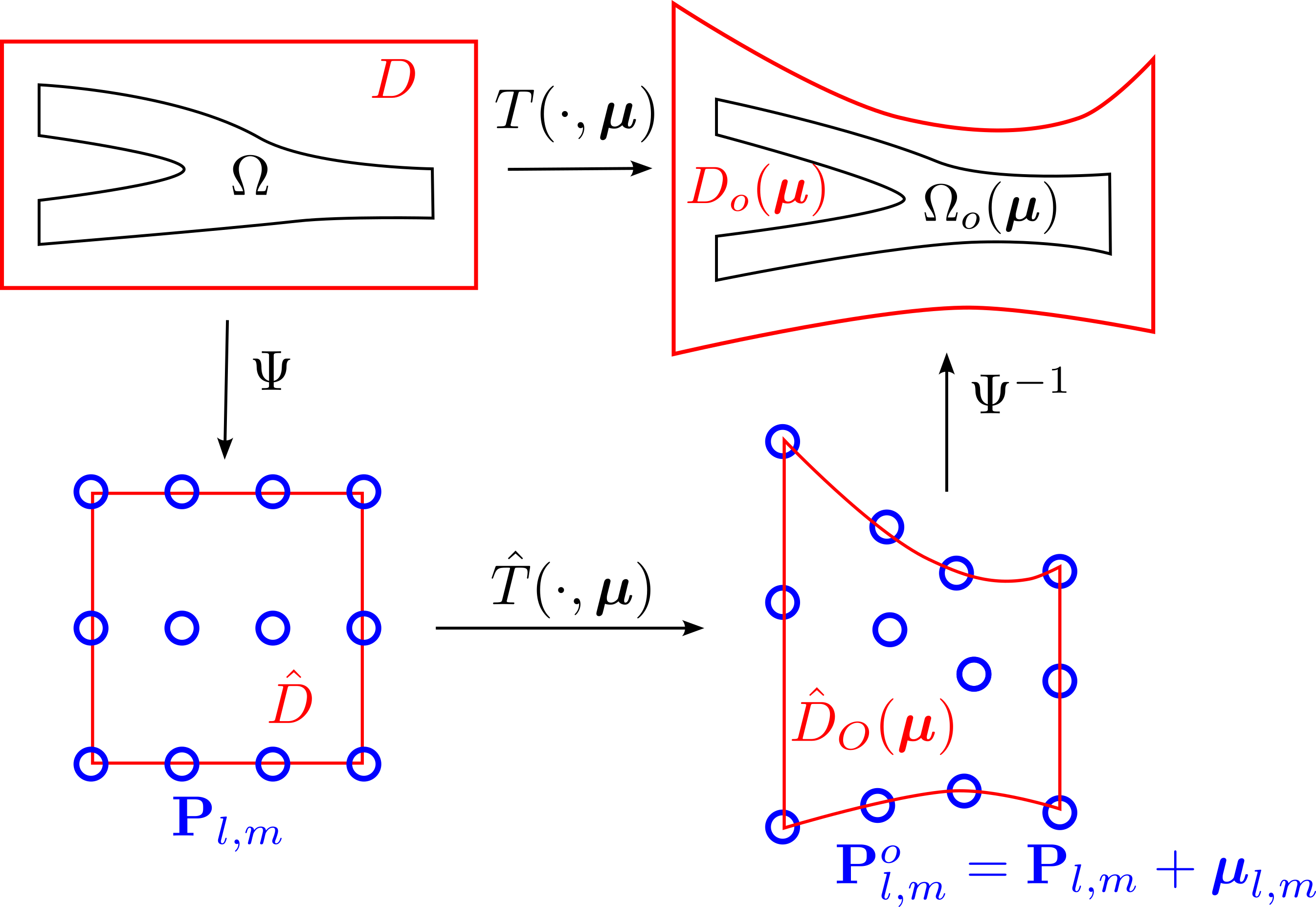}
    \caption{Two-dimensional FFD implementation for a bypass structure. FFD deforms a lattice (red) via moving the grid control points (blue). The new shape is obtained using mapping and interpolation. Figure redrawn based on \cite{manzoni2012shape}.}
    \label{fig:bypass_FFD_sketch}
\end{figure}

The procedure of FFD is sketched in Fig. \ref{fig:bypass_FFD_sketch}. The whole transformation $T(\cdot,\boldsymbol{\mu})$ consist on three steps $T(\cdot,\boldsymbol{\mu})=\psi^{-1}\circ\hat{T}(\cdot,\boldsymbol{\mu}) \circ \psi$. Firstly, a reference bypass structure, $\Omega$, is embedded in a rectangle, $D$. The latter is scaled, with transformation $\Psi$, into a unit cube $\hat{D}\equiv [0,1]^2$. 

Secondly, the transformation $\hat{T}$ is defined. We may consider $L$ and $M$ control points inside the cube $\hat{D}$. $L$ and $M$ are defined on the abscissa and ordinate axes, respectively. Each grid point is defined by its two coordinates and indexed by $\mathrm{P}_{l,m} = [l/L, m/M]$. The movement of each point is denoted $\boldsymbol{\mu}_{l,m}$. $\hat{T} (\hat{\mathbf{x}}; \boldsymbol{\mu})$ is expressed,
\begin{equation*}
    \hat{T} (\hat{\mathbf{x}}, \boldsymbol{\mu}) = \sum_{l=0}^L \sum_{m=1}^M b_{l,m}^{L,M} (\hat{\mathbf{x}}) \mathrm{P}^o_{l,m},
\end{equation*}
where is $\hat{\mathbf{x}}$ is any coordinate $\in \hat{D}\equiv [0,1]^2$, $b_{l,m}^{L,M}(\hat{\mathbf{x}})=b_{l}^{L}(s)b_{m}^{M}(t)$, $b_{l}^{L}(s)$ denotes tensor product of the 1-D Bernstein basis polynomials \cite{sederberg1986free} and $\mathrm{P}^o_{l,m} = \mathrm{P}_{l,m} + \boldsymbol{\mu}_{l,m}$ is the translated control point. By applying $\hat{T}(\boldsymbol{\mu})$, one can achieve $\hat{D}\to \hat{D}_0(\boldsymbol{\mu})$. Finally, $\Omega_0$ is obtained through an inverse mapping $\Psi$ from $\hat{D}_0$.

FFD becomes popular due to many valuable characteristics \cite{sederberg1986free, koshakji2013free}: (i) It is suitable for surfaces of any formulation or degree; (ii) Its implementation is independent of the domain and mesh; (iii) It can be applied locally or globally and preserves smoothness when incorporating deform only for subregions. 

There are many relevant studies about the combination of ROM and FFD discussing different phenomena: shape parametrized heat conduction problem \cite{garotta2020reduced}; the immersed bodies in Stokes flows \cite{ballarin2014shape}; PDE-Constrained Optimization Problems in Haemodynamics \cite{rozza2012reduction}; unsteady Navier-Stokes flow inside the Coronary artery bypass \cite{ballarin2015reduced, ballarin2016fast}; linear/nonlinear constrained design optimization of ship hull (governed by Navier-Stokes equations) \cite{padula2023generative}. Another highlighted investigation is addressed by Davide Forti \cite{forti2012Comparison}, in which he adopted RBF, IDW, and FFD for FSI scenarios.

\begin{figure}[h]
    \centering
    \includegraphics[width=0.6\linewidth]{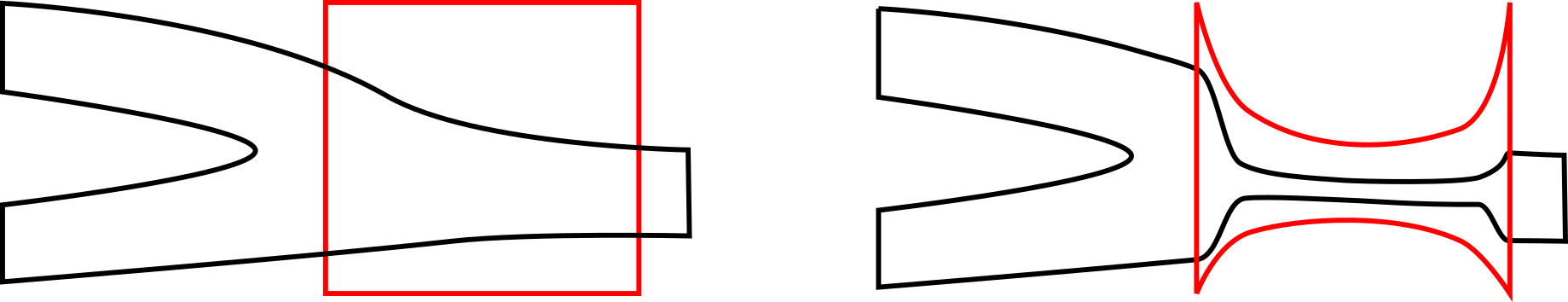}
    \caption{FFD for local deformation of a cardiovascular structure. e.g., a blood vessel. Figure redrawn based on \cite{siena2024accuracy}. }
    \label{fig:local_FFD}
\end{figure}

The FFD can also be applied locally. Fig. \ref{fig:local_FFD} shows a simulation of Navier-Stokes flow in a shrinking cardiovascular structure \cite{siena2024accuracy}. Another study for the aerodynamic shape optimization of a car's front end is presented in \cite{salmoiraghi2016advances}. \cite{forti2012Comparison} managed to combine the domain decomposition approach for FFD, in which an individual lattice is defined in each subdomain. This configuration allows local refinements of control points in a specific region.

\subsubsection{Parameter space reduction}
\label{subsec:Parameter_space_reduction}
The aforementioned techniques cover the general strategy for generating samples, including physical and geometrical. However, parameterization may become problematic when associated with a high-dimensional space. That implies difficulties in design optimization, sensitivity analysis, uncertainty quantification, etc. 

Some examples may be provided to clarify this issue. Suppose that shape optimization is performed employing FFD with a large amount of control points. The displacements of every point should be considered, that is, sampled. That is not practical for the design workflow of realistic applications. To deal with these difficulties, the \emph{parameter space reduction} technique arose. In summary, it involves moving points in a correlated manner, specifically by adjusting a few \emph{synthetic} parameters, and then obtaining samples based on this reduction.

Several studies have shown the key benefits of parametric reduction. Clearly, fewer inputs reduced the cost of ROM evaluation, which accelerates the analyses of sensitivity and uncertainty \cite{constantine2015active}. Additionally, it was shown that low-dimensional inputs also improve the accuracy of ROMs \cite{demo2019non-intrusive, tezzele2020enhancing, tezzele2019shape, padula2023generative, romor2024local}.

For completeness of the section, we will illustrate two popular techniques regarding parametric reduction in ROM: \emph{Active subspaces} and the \emph{AI-based generative model}. 

\paragraph{Active subspaces}
The \emph{Active Subspaces} (AS) method has recently been used for linear reduction in input spaces \cite{romor2024local}. Here, we will only present the steps of AS, whilst for more theoretical and mathematical details, as well as applications, we refer readers to \cite{constantine2015active}. 

We consider a multi-dimensional parameter space $\mathcal{P} \subset \mathbb{R}^{N_\mathcal{P}}$. We also regard the parametric function $f(\boldsymbol{\mu}): \mathbb{R}^{N_\mathcal{P}} \rightarrow \mathbb{R}$. The uncertainty in the input parameters is accounted for with a probability density function $\rho$. 

We utilize the gradient of $f$ with respect to $\boldsymbol{\mu}$, $\nabla_{\boldsymbol{\mu}} f(\boldsymbol{\mu}) \in \mathbb{R}^{N_\mathcal{P}}$. We form its non-centered covariance matrix $\Sigma$,
\begin{equation*}
    \Sigma = \mathbb{E} \left[ \nabla_{\boldsymbol{\mu}} f(\boldsymbol{\mu}) \nabla_{\boldsymbol{\mu}} f(\boldsymbol{\mu})^{\mathrm{T}} \right] = \int \left( \nabla_{\boldsymbol{\mu}} f(\boldsymbol{\mu}) \right) \left( \nabla_{\boldsymbol{\mu}} f (\boldsymbol{\mu}) \right)^\mathrm{T} \rho \ \mathrm{d}\boldsymbol{\mu},
\end{equation*}
where $\mathbb{E}$ is the expected value. Matrix $\Sigma$ is symmetric and positive semidefinite. Therefore, it has a real eigenvalue decomposition, 
\begin{equation*}
    \Sigma = \mathbf{W} \Lambda \mathbf{W}^{\mathrm{T}}, 
\end{equation*}
where  $\mathbf{W}$ is orthogonal and contains the eigenvectors in columns. $\Lambda = \text{diag} \{\lambda_1, \dots, \lambda_{N_\mathcal{P}}\}$ with the eigenvalues $\lambda_1 \geq \cdots \lambda_{N_\mathcal{P}} \geq 0$.

A key aspect of the method consists of realizing that eigenvalues and eigenvectors can be partitioned
\begin{equation*}
    \Lambda = \begin{bmatrix}
             \Lambda_1 &  \\
             & \Lambda_2  \\
            \end{bmatrix}, \quad \mathbf{W} = \left[ \mathbf{W}_1, \mathbf{W}_2 \right].
\end{equation*}
Note that $\mathbf{W}_1$ with a shape of $N_{\mathcal{P}} \times  M$ denotes the $M$-dimension \emph{active} subspace of $\mathbf{W}$. That is, the range of $\mathbf{W}_1$ captures the most of the variation of $\nabla_{\boldsymbol{\mu}} f (\boldsymbol{\mu})$. Then, the low-dimensional \emph{active} parameter can be defined as a linear combination of the original values, indeed, $\boldsymbol{\mu}_M = \mathbf{W}_1^\mathrm{T} \boldsymbol{\mu}$. With the active variables, we can approximate the function of interest $f$, 
\begin{equation*}
    f(\cdot, \boldsymbol{\mu}) \approx g(\cdot, \mathbf{W}_1^\mathrm{T} \boldsymbol{\mu}) = g(\cdot, \boldsymbol{\mu}_M).
\end{equation*}

The incorporation of ROM and AS has been tested in various scenarios. Nicola Demo et al. \cite{demo2019non-intrusive} revealed that their non-intrusive ROM involving AS outperforms the one disregarding this technique. The framework has also been utilized intensively for shape design and optimization. The applications include airfoil \cite{tezzele2020enhancing}, ship hull \cite{tezzele2019shape, tezzele2018dimension, padula2023generative}, and ship propeller blade \cite{mola2019efficient} optimizations. 

We remark on a recent publication \cite{romor2024local}, a slight variation of the method. There, the authors managed to partition the parameter space with clustering algorithms. Their method, the so-called \emph{Local active subspace} procedure, results in a more efficient and accurate dimension reduction than the standard one.

\paragraph{Generative model}

AI-based methods have been employed to generate 3D objects in computer vision. This suggests the feasibility of utilizing them to create geometrical samples. Given their relevance, we start introducing data-driven \emph{generative models} for parameter space reduction \cite{padula2023generative}. 

As mentioned before, Free Form Deformation (FFD) is a flexible and general technique for the generation of geometries. However, when considering linear/multi-linear constraints, its numerical cost can not be ignored. 

Thus, Guglielmo Padula et al. \cite{padula2023generative, padula2024generative, padula2025generative} proposed a strategy based on neural networks to generate efficiently constrained samples. It is composed of two parts: an Encoder $E: \mathbb{R}^N \rightarrow \mathbb{R}^n$ and a decoder $D: \mathbb{R}^n \rightarrow \mathbb{R}^N$, in which $n \ll N$. It achieves the approximation that $\boldsymbol{\mu} \approx D \circ E (\boldsymbol{\mu})$. The low-dimensional representation $E (\boldsymbol{\mu})$ is known as a \emph{latent space} of $\boldsymbol{\mu}$. 

The implementation of the generative model consists of three steps: (i) a few samples generated by constrained FFD are applied to train the network; (ii) instead of $\boldsymbol{\mu}$, the vector $E (\boldsymbol{\mu})$ is considered as parameter; (iii) new samples are produced utilizing decoder $D \circ E (\boldsymbol{\mu})$. The procedure can utilize a simple \emph{Autoencoder} (see section \ref{subsubsec:Neural_network}) or more complex frameworks, which were described and compared in \cite{padula2025generative}.  


This model has been validated utilizing several benchmark cases, demonstrating its effectiveness across different applications. These include the Poisson problem in Stanford Bunny geometry and Reynolds Averaged Navier-Stokes equations for the shape optimization of a naval hull. These tests have revealed that the method can represent constraints efficiently and result in a more accurate ROM \cite{padula2023generative, padula2024generative, padula2025generative}. Notably, both the \emph{active subspace method} and the \emph{generative model} confirm the significant benefits of exploiting parameter space reduction techniques.


\subsection{Snapshots and local reduced basis}
\label{sec:snapshots_local_RB}
Once the domain is properly partitioned and parameters are determined, the next step is to compute the local RB. The following paragraphs will briefly discuss three procedures for generating subdomain-level RBs. 

Note that the concept of RB is defined for intrusive methodologies, whereas for non-intrusive methods, the low-dimensional representation is a standard expression. We remind that \emph{RB} is still used hereafter as a common concept for both intrusive and non-intrusive frameworks.

\subsubsection{Localized global RB}
\label{subsubSec:localized_global_RB}
In the first approach, parametric numerical solutions for the entire geometry are collected $\mathbf{S}_{N_h \times N_\mu} = \left[ \mathbf{u}_1, \mathbf{u}_2, \dots, \mathbf{u}_{N_\mu} \right]$, where $N_\mu$ is the total amount of parameters and $N_h$ denotes the dimension of high fidelity solution. Each snapshot is given by $ \mathbf{u}_j \coloneqq \mathbf{u}_j (\mathbf{x}; \mu_j) \text{ for } j = 1, \dots, N_\mu$, and $\mu_j$ are either geometrical or physical parameters. A reduction technique, e.g., POD, is adopted to compute $N_{{\text{RB}}}$ global dominant modes. Thus, $\mathbf{u}_j (\mathbf{x}; \mu_j) \approx \sum_{i=1}^{N_{\text{RB}}} a_{i} (\mu_j) \mathbf{v}_i (\mathbf{x}), \quad \mathbf{v}_i \in \mathcal{V}$, where $\mathcal{V}$ is the global reduced space. Finally, the $i^\text{th}$ basis function of the local RB in subdomain $\Omega_m$ is defined as the restriction: $\mathbf{v}_{m,i} = \mathbf{v}_i \rvert_{\Omega_m}$. 

\subsubsection{Global solutions and local RB}
\label{subsubsec:global_solution_local_RB}
The second technique also requires global FOM results that are decomposed before any further treatment. The subdomain-level reduced subspaces are computed using the local values. There are two procedures for computing the RBs.

\begin{itemize}
    \item 
For the \emph{Individually Partition} strategy (Section \ref{subsubsec:individual_decomposition}), the local datasets are $\mathbf{S}_m = \left[ \mathbf{u}_{m,1}, \mathbf{u}_{m,2}, \cdots, \mathbf{u}_{m,N_\mu} \right]$, in which the local snapshots of $\Omega_m$ are obtained by the restriction $ \mathbf{u}_{m,j} \coloneqq \mathbf{u} (\mathbf{x}; \mu_j)\rvert_{\Omega_m} \text{ for } j = 1, \dots, N_\mu$. The RBs for each partition are calculated separately. The local approximations are $\mathbf{u}_{m,j} (\mathbf{x}; \mu_j) \approx \sum_{i=1}^{N_{\text{RB},m}} a_{m,i} (\mu_j) \mathbf{v}_{m,i} (\mathbf{x})$, in which $\mathbf{v}_{m,i} \in \mathcal{V}_m$.

This is shown in Fig. \ref{fig:individual_RB_backward_step}, in which the geometry of a backward-facing step problem is decomposed into two parts, i.e., $\Omega_1$ and $\Omega_2$ \cite{prusak2023optimisation}. Two different RBs are computed for each part, namely, $\mathcal{V}_1$ and $\mathcal{V}_2$. The two are \emph{glued} for visualization, as shown in Fig. \ref{fig:individual_RB_backward_step}.

\begin{figure}[h]
    \centering
    \begin{subfigure}[b]{0.45\linewidth}
        \includegraphics[width=\linewidth]{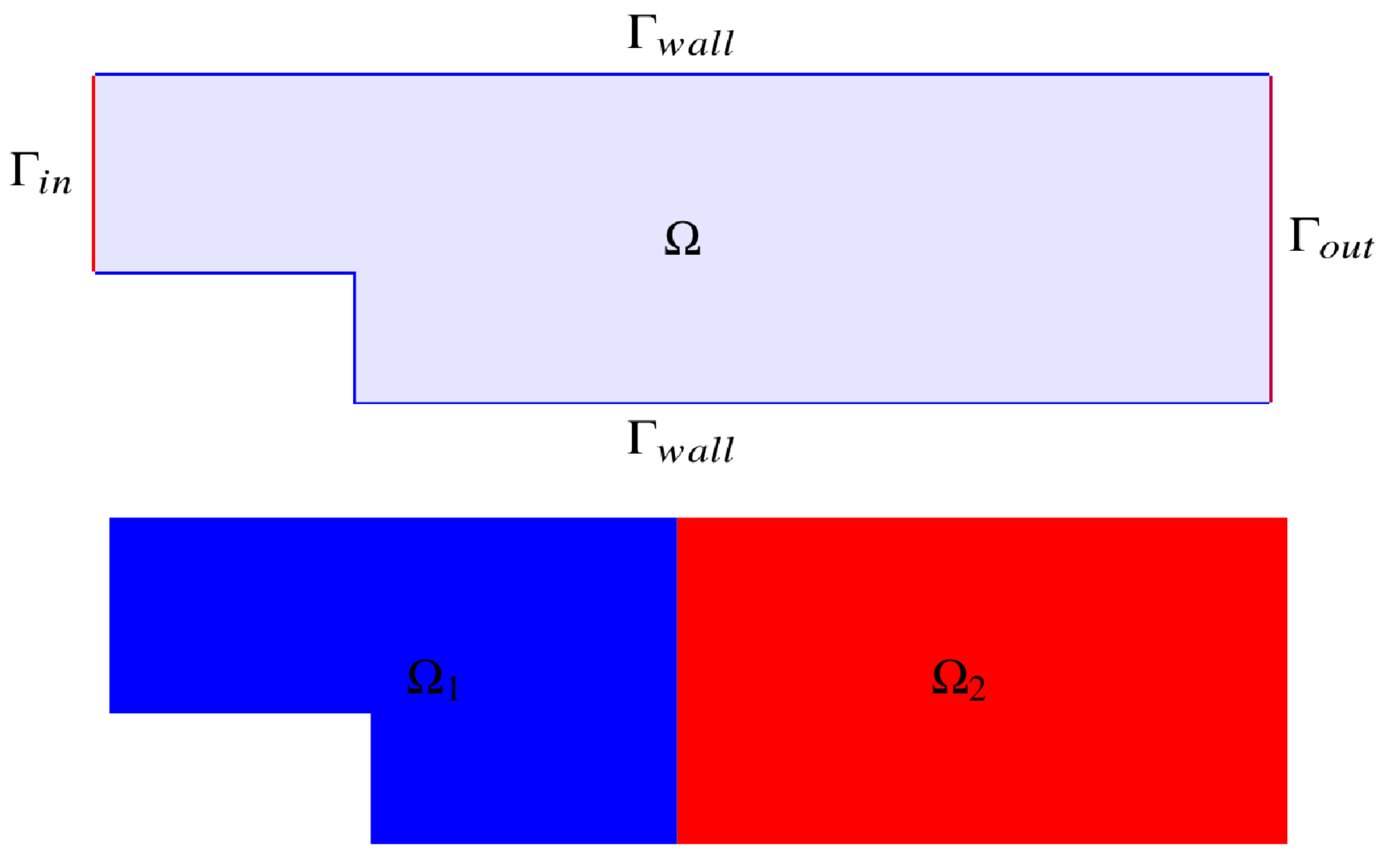}
        \caption{Geometry and subdomains.}
    \end{subfigure}
    \hspace{2em}
    \begin{subfigure}[b]{0.4\linewidth}
        \includegraphics[width=\linewidth]{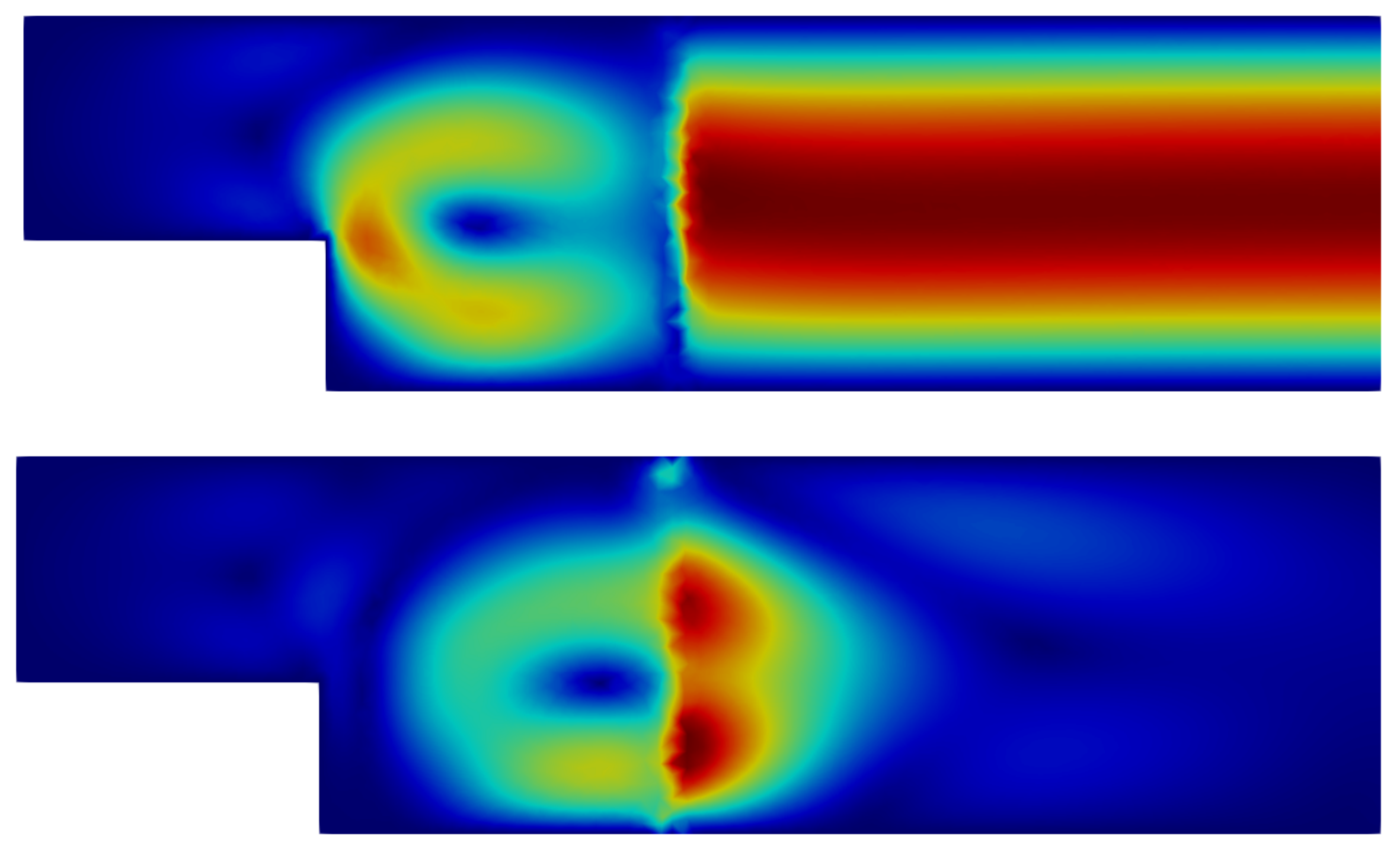}
        \caption{Velocity modes.}
        \label{fig:individual_vel_RB_backward_step}
    \end{subfigure}
    \caption{Individual domain decomposition and local RBs for a backward-facing step problem. The local velocity modes are \emph{glued} together for visualization. Taken with permission from \cite{prusak2023application},  copyright owned by the author.}
    \label{fig:individual_RB_backward_step}
\end{figure}

\begin{figure}[h]
    \centering
    \begin{subfigure}[b]{0.25\linewidth}
        \centering
        \includegraphics[width=\linewidth]{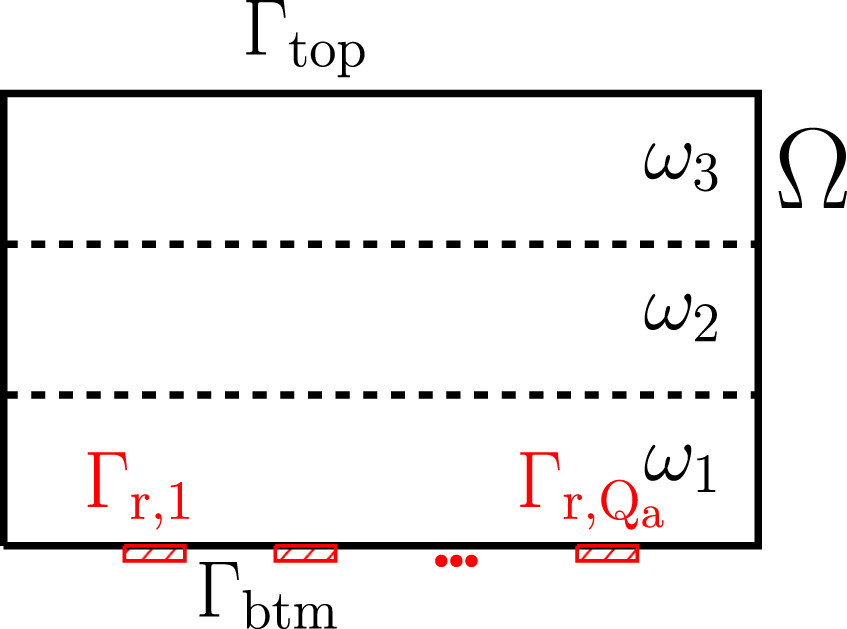}
        \caption{Global model.}
    \end{subfigure}
    \hspace{2em}
    \begin{subfigure}[b]{0.45\linewidth}
        \centering
        \includegraphics[width=\linewidth]{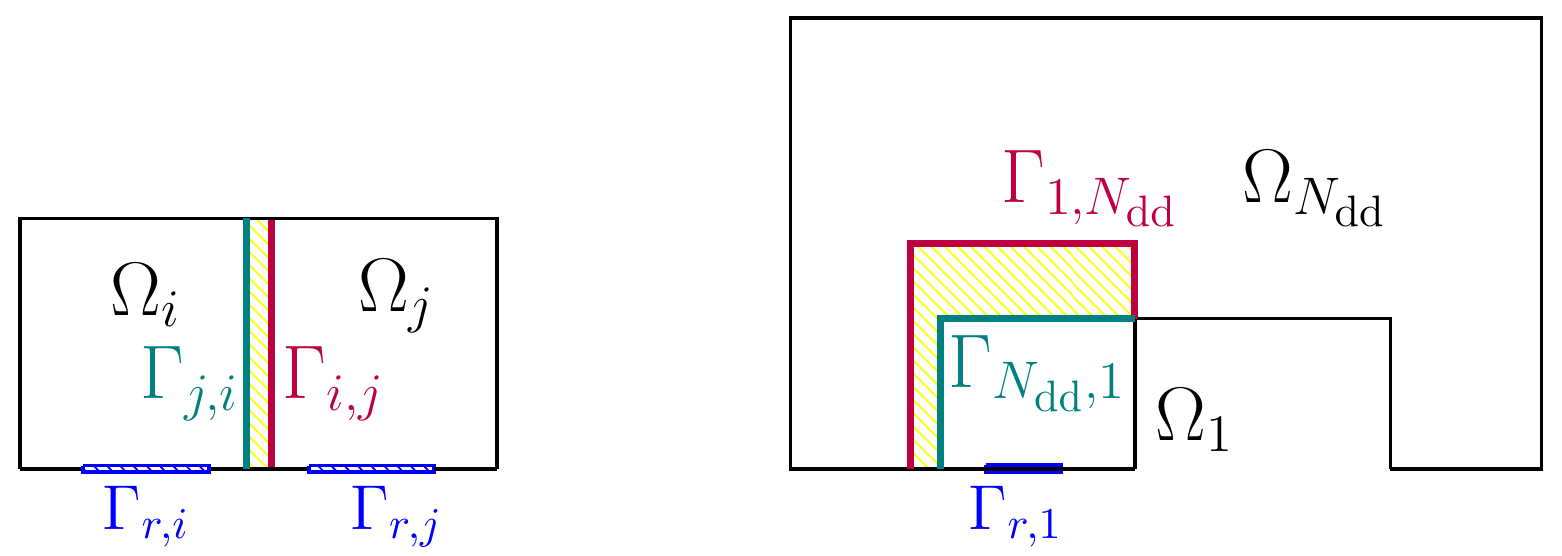}
        \caption{Assembly of the geometry.}
    \end{subfigure}
    \caption{Domain decomposition of a thermo-hydro-mechanical system. Many $\Omega_{int}^a$ and a single $\Omega_{ext}^a$ are composed to build the whole system. Taken with permission from \cite{sambataro2022component}, copyright owned by the author.}
    \label{fig:generic_RB_THM_system}
\end{figure}

    \item 
In the \emph{Generic Decomposition} approach (Section \ref{subsubsec:Generic_decomposition}), a local RB is constructed for each archetype. The RB of each instantiated block is transformed from its archetype RB. The number of local RBs is equal to the number of archetypes. Once $\hat{\Omega}^k$ and their instantiations $\Omega^k_m$ are defined, the local solution is obtained by restriction $\mathbf{u}^k_{m,j} \coloneqq \mathbf{u} \left. \left( \mathbf{x}; \mu_j \right) \right|_{\Omega^k_m}$. The local solutions in each instantiation of each archetype are collected as  $\mathbf{S}^k_m = \left[ \mathbf{u}^k_{m,1}, \cdots, \mathbf{u}^k_{m,N_\mu} \right]$. Then, they are transformed back to the generic division, e.g., $\hat{\mathbf{S}}^k_m = {T^k (\boldsymbol{\xi}^{k}_m)}^{-1} \cdot \mathbf{S}^k_m$ \footnote{If the subdomains have the same shape as their generic blocks, no transformation is need.}. They are further stacked into $\hat{\mathbf{S}}^k = \left[ \hat{\mathbf{S}}^k_1, \cdots, \hat{\mathbf{S}}^{k}_{N_k} \right]$. 

The reduction technique is applied at this stage to each $\hat{\mathbf{S}}^k$. For $\hat{N}_\Omega$ reference blocks, we can obtain in total $\hat{N}_\Omega$ RBs. The basis vectors are mapped to real blocks, which lead to $ \mathbf{v}^k_{m,i} = T^k (\boldsymbol{\xi}^{k}_m) \cdot \hat{\mathbf{v}}^k_{i}$. Then, the local approximations are 
\begin{equation}
\label{eq:generic_local_solution}
    \mathbf{u}^k_{m,j} (\mathbf{x}; \mu_j) \approx \sum_{i=1}^{N_{\text{RB},m}^k} a^k_{m,i} (\mu_j) \mathbf{v}^k_{m,i} (\mathbf{x}), \quad {\mathbf{v}}^k_{m,i} \in {\mathcal{V}}^k_m.
\end{equation}

Fig. \ref{fig:generic_RB_THM_system} shows the sketch of a thermo-hydro-mechanical system for radioactive waste disposal \cite{sambataro2022component,iollo2023one}, in which two archetype components, $\Omega_{int}^a$ and $\Omega_{ext}^a$, are adopted to build the facility. All parameterized solutions in $\Omega_{int}^a$ are stacked into a compact dataset, and POD is applied to compute the reduced subspace. As there is only a single $\Omega_{ext}^a$, its RB can be calculated directly. 

\end{itemize}

\subsubsection{Localized training and oversampling}
\label{subsubSec:Localized_training}
The third method is called \emph{localized training}, which supports geometries composed of generic parts. Instead of expensive high-fidelity modeling in the original large geometry, a set of parametric computations is carried out in several much smaller systems containing archetype blocks. In this frame, it is assumed that the dynamics that occur in small-scale networks can represent the physical behavior of the original large-scale systems. The approach can significantly accelerate the most expensive part of ROM, i.e., the offline stage. This has the potential to make this method become popular. 

The procedure to collect snapshots and compute RBs is similar to the one described above (Section \ref{subsubsec:global_solution_local_RB}), considering archetypes. Offline solutions restricted in each generic partition are gathered into separate matrices to build corresponding RBs. Let us consider the $k^{\text{th}}$ archetype $\hat{\Omega}^k$ with the set of parameters $\left\{\mu_j,\ j=1,\cdots, \hat{N}_\mu^k \right\}$ and the corresponding data $\hat{u}^k(\mu_j)$. As in Section \ref{subsubsec:Generic_decomposition}, the snapshot matrix $\hat{\mathbf{S}}^k = \left[ \hat{\mathbf{u}}^k_1, \cdots, \hat{\mathbf{u}}^k_{\hat{N}_\mu^k} \right]$ is formed, to which a reduction technique is applied. The reduced spaces $\hat{\mathcal{V}}^k = \text{span} \left\{ \hat{v}^k_i, i=1, \cdots, \hat{N}^k_{\text{RB}} \right\}$ are obtained, in which approximated solutions $\hat{u}^k(\mu_j) \approx \sum_{i=1}^{\hat{N}^k_{RB}} \hat{a}^k_i (\mu_j) \hat{v}^k_i$. Remark that the Gram-Schmidt procedure can be applied to orthogonalize the basis vectors. Finally, the archetypes RBs' may be transformed to each instantiated partition $\Omega_m^k$ utilizing the geometric transformation operator $T^k (\boldsymbol{\xi}_m^k)$ --as illustrated in equation \eqref{eq:geometric_transformation}--. We can obtain, $v_{m,i}^k = T^k (\boldsymbol{\xi}_m^k) \cdot \hat{v}^k_i$, and express a local result on $\Omega_m^k$ in the same form as equation \ref{eq:generic_local_solution}.

\begin{figure}[h]
    \centering
    \begin{subfigure}[b]{0.5\linewidth}
        \centering
        \includegraphics[height=4.5cm]{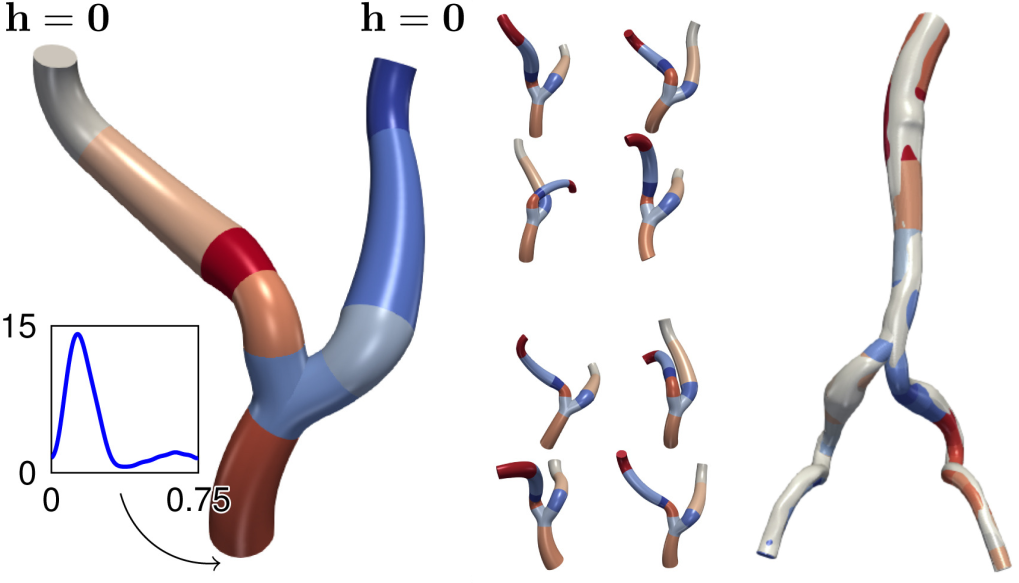}
        \caption{FOM geometry.}
        \label{fig:generic_division_artificial_blood_vessel_parameters}
    \end{subfigure}
    \hspace{2em}
    \begin{subfigure}[b]{0.4\linewidth}
        \centering
        \includegraphics[height=4.5cm]{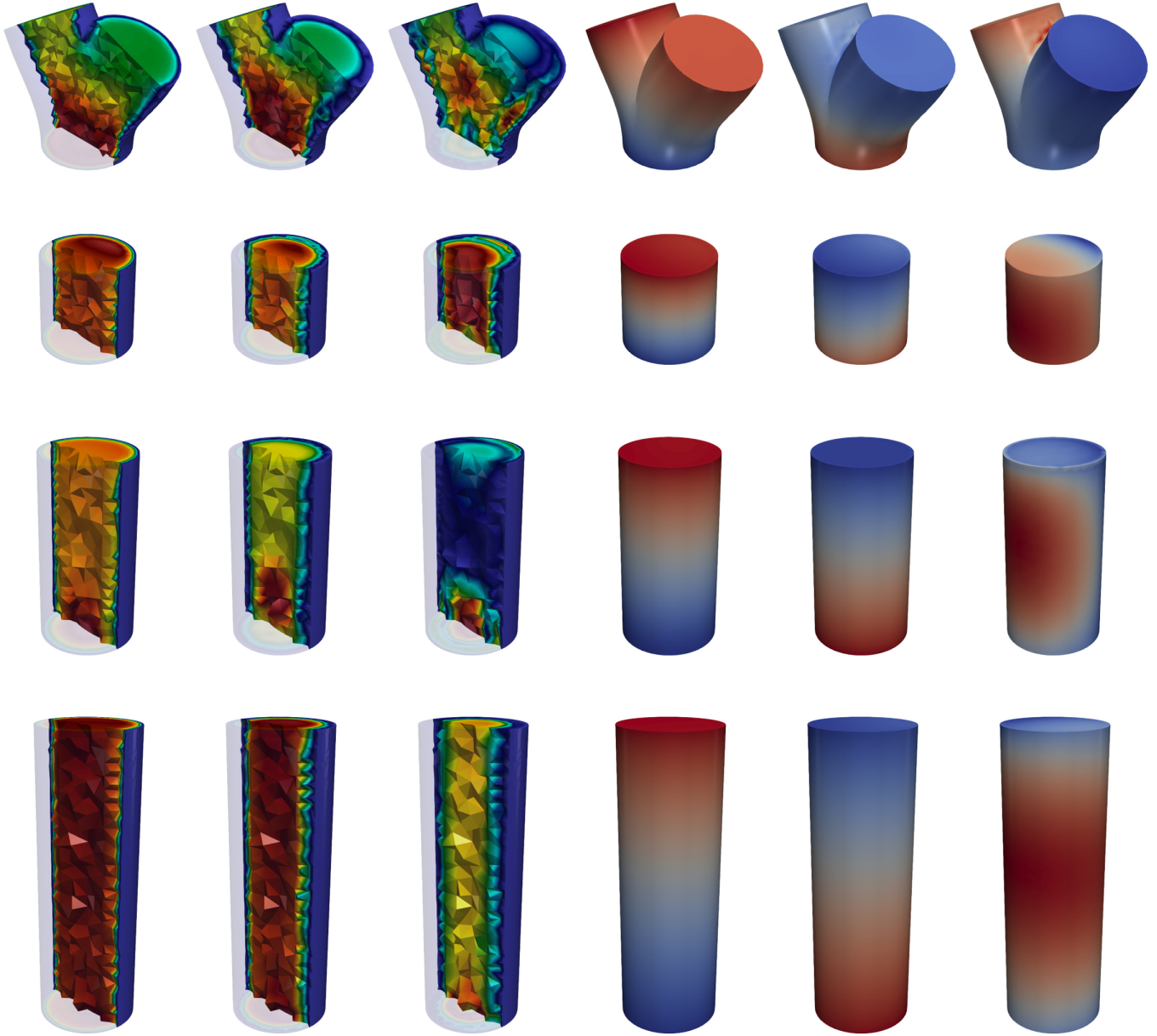}
        \caption{Velocity and pressure modes for archetypes.}
    \end{subfigure}
    \caption{The model for collecting high-dimensional data and dominant modes for archetypes. Figures taken with permission from \cite{pegolotti2021model}, copyright owned by Elsevier. (a) The \emph{artificial} model (left), geometrical parameterized samples for generating FOM solutions (middle), and a realistic global domain (right). (b) Subdomain-level modes for each archetype.}
    \label{fig:generic_geometry_RB_blood_vessel}
\end{figure}

Luca Pegolotti et al. \cite{pegolotti2021model} have applied this method to analyze the flow in blood vessels. The parameterized snapshots were collected from an \emph{artificial} (see Fig. \ref{fig:generic_division_artificial_blood_vessel_parameters}) geometry. It contains all four archetypes, but on a much smaller scale compared to real vessels. The local velocity and pressure POD modes for generic partitions are displayed in Fig. \ref{fig:generic_geometry_RB_blood_vessel}. 

\begin{figure}[h]
    \centering
    \includegraphics[width=0.6\linewidth]{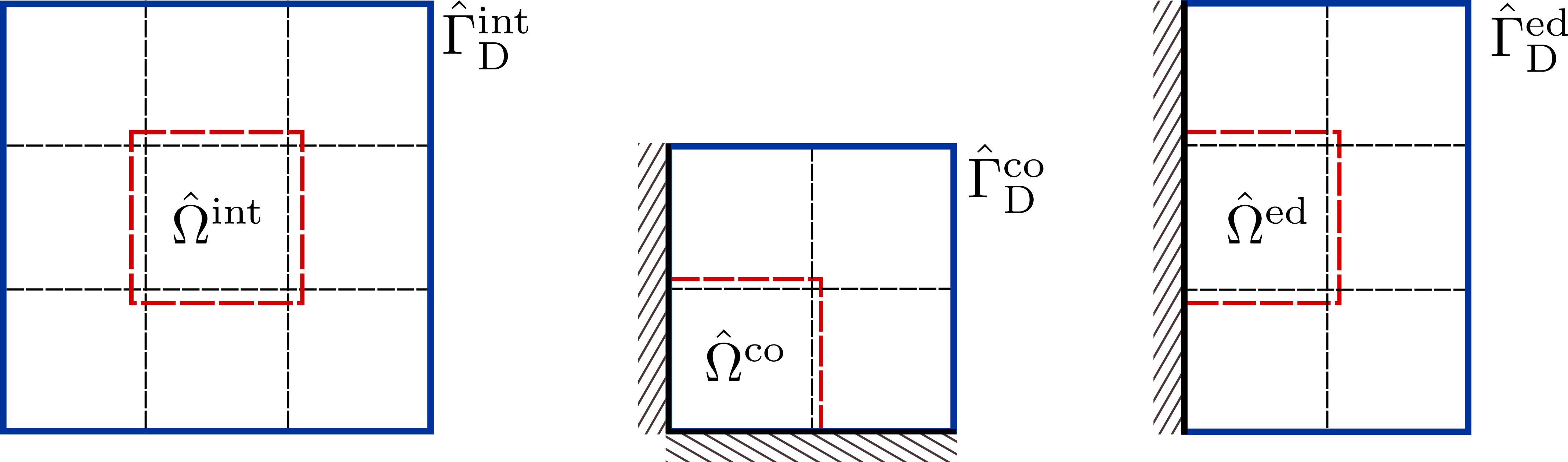}
    \caption{A sketch of the \emph{Oversampling} strategy for three archetype components, i.e., $\hat{\Omega}^{\text{int}}$, $\hat{\Omega}^{\text{co}}$ and  $\hat{\Omega}^{\text{ed}}$. Three small-scale models are constructed to generate high-fidelity solutions of each generic subdomain. Redrawn based on \cite{smetana2023localized}.}
    \label{fig:generic_subdomain_oversampling}
\end{figure}

In order to compute modes that are representative of the physical phenomena appearing in the global solution, an \emph{oversampling} strategy has been proposed \cite{hou1997multiscale}. The FOM simulations are performed in geometries that are slightly larger and comprise several archetypes. To the extended domains, some parameterized boundary conditions are applied, as shown in Fig. \ref{fig:generic_subdomain_oversampling}. The solutions within the reference subdomains are extracted by restriction to compute corresponding RBs. The oversampling approaches can identify low-dimensional structures and have been suggested and used extensively \cite{henning2013oversampling, buhr2020localized, babuvska2020multiscale, smetana2023localized}. Besides, the increasing computation cost for FOM simulations remains acceptable. 

\subsubsection{RB for interior, interface, and boundary.}
\label{sec:RBintintb}
We finalize this section by discussing the practical formation of the snapshot's matrix regarding the interior $\Omega$, interfaces $\Gamma$ (or the local boundary $\partial \Omega_m$), and global boundary $\partial \Omega$. Recall that the matrix will be treated with a reduction algorithm to obtain the RBs.

\begin{figure}[h]
    \centering
    \begin{subfigure}[b]{0.25\linewidth}
        \includegraphics[height=4cm]{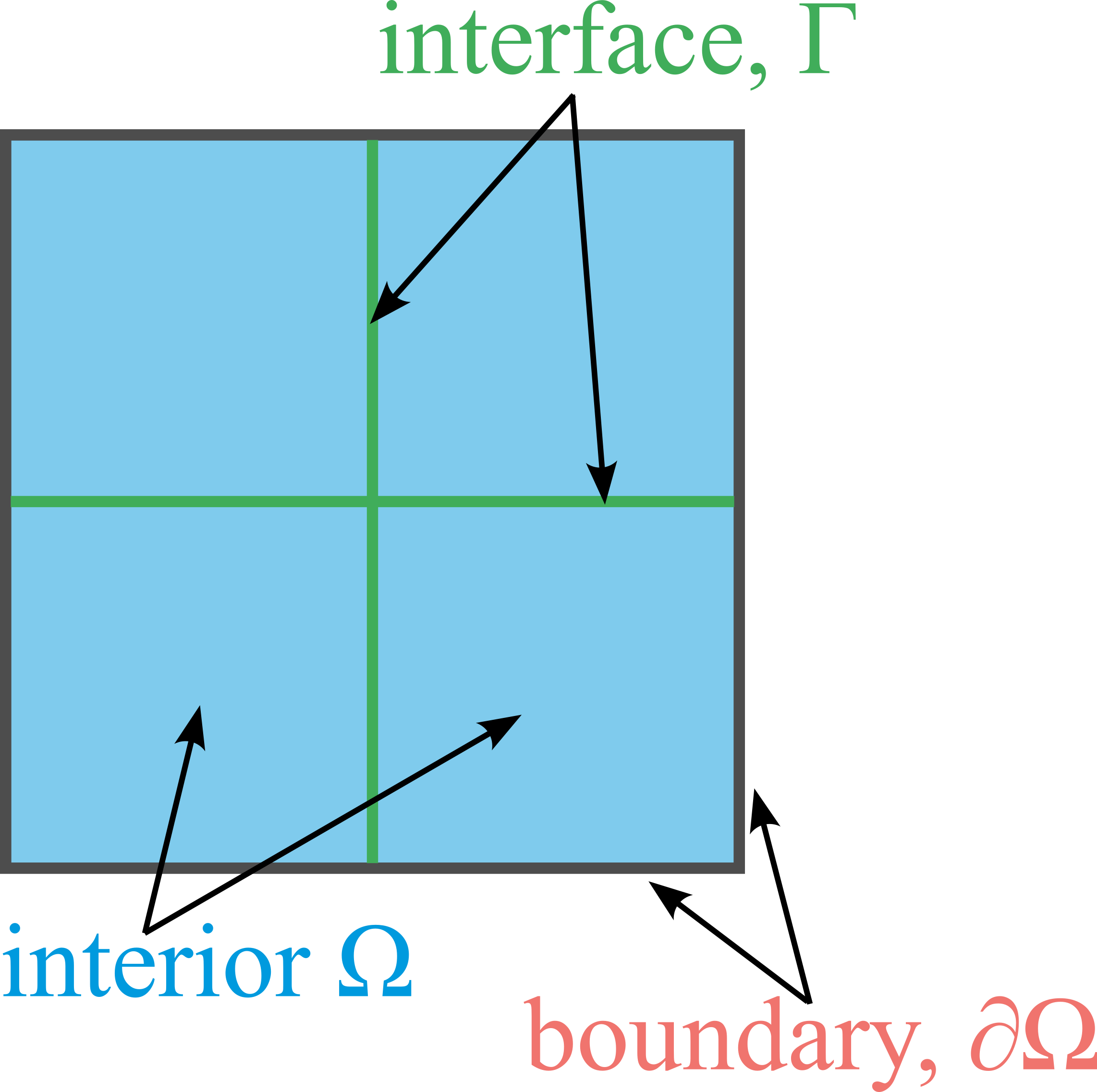}
        \caption{}
    \end{subfigure}
    \hspace{4em}
    \begin{subfigure}[b]{0.15\linewidth}
        \includegraphics[height=4cm]{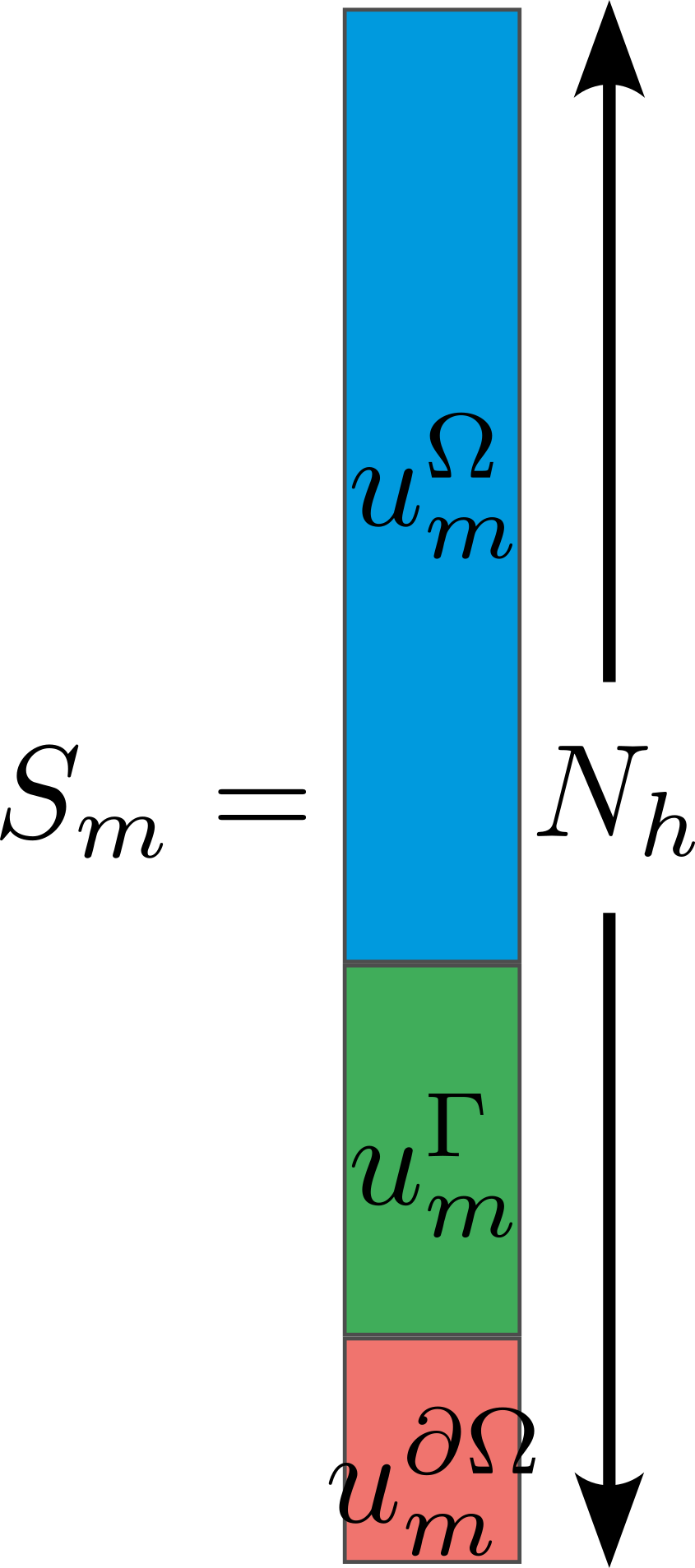}
        \caption{}
    \end{subfigure}
    \hspace{2em}
    \begin{subfigure}[b]{0.15\linewidth}
        \includegraphics[height=4cm]{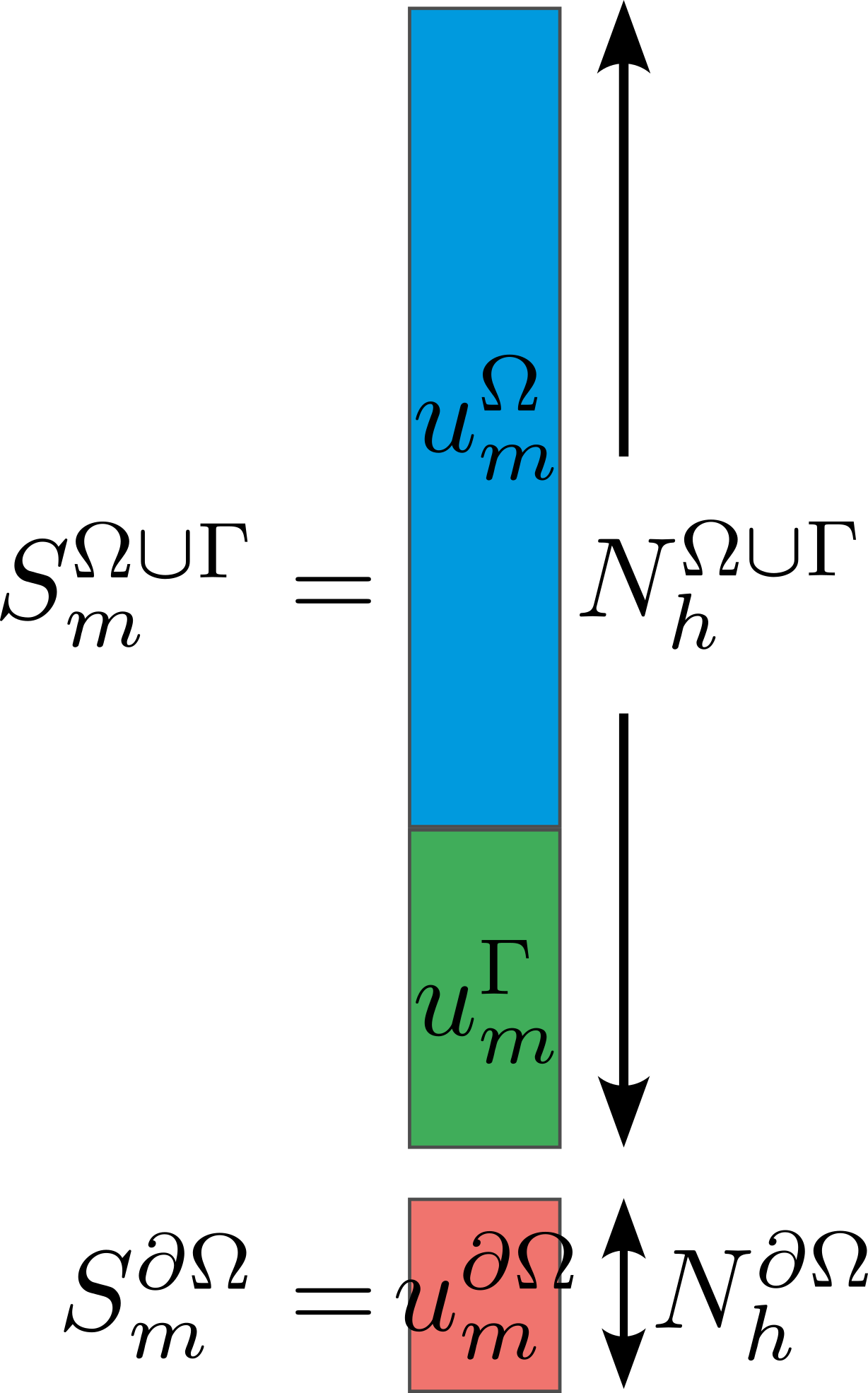}
        \caption{}
    \end{subfigure}
    \hspace{2em}
    \begin{subfigure}[b]{0.15\linewidth}
        \includegraphics[height=4cm]{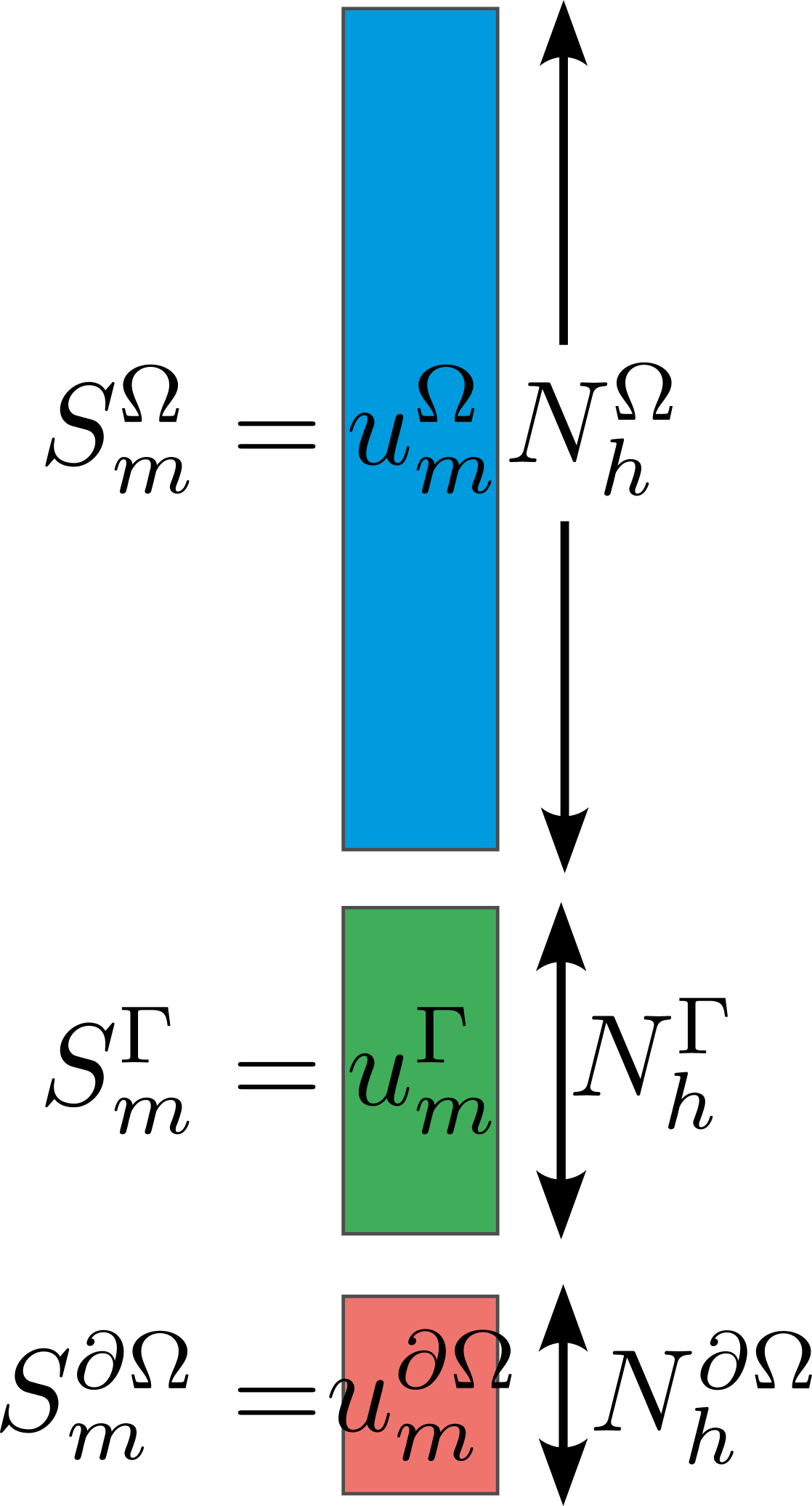}
        \caption{}
    \end{subfigure}
    \caption{Snapshots structures. (a) A domain and three regions: , $\Omega$; interface, $\Gamma$; and global boundary, $\partial\Omega$. (b) A compact subdomain snapshot. (c) Subdomain and boundary snapshots. (d) Separated snapshots.}
    \label{fig:regions_snapshots}
\end{figure}

Fig. \ref{fig:regions_snapshots} illustrates the three strategies to process offline solutions: (i) solutions of the three zones -interior, interfaces, and boundaries- are gathered, and a single compact reduced space is constructed; (ii) processing the interior and the interface together and the boundary separately; or (iii) each zone is post-processed separately.

The procedure to follow depends on the partition strategies described above and the coupling algorithm introduced in the following sections.

\section{Projection-based coupling algorithms}
\label{sec:projection_based_coupling_algorithms}

Galerkin projection is the fundamental procedure utilized to derive the weak formulation of PDEs. It is widely adopted to formulate and solve high-fidelity problems in the frame of many discretization algorithms, such as the Spectral Method (SM) \cite{canuto2006spectral} and Finite Element Method (FEM) \cite{ern2004theory}. 

\begin{figure}[h]
    \centering
    \includegraphics[width=0.7\linewidth]{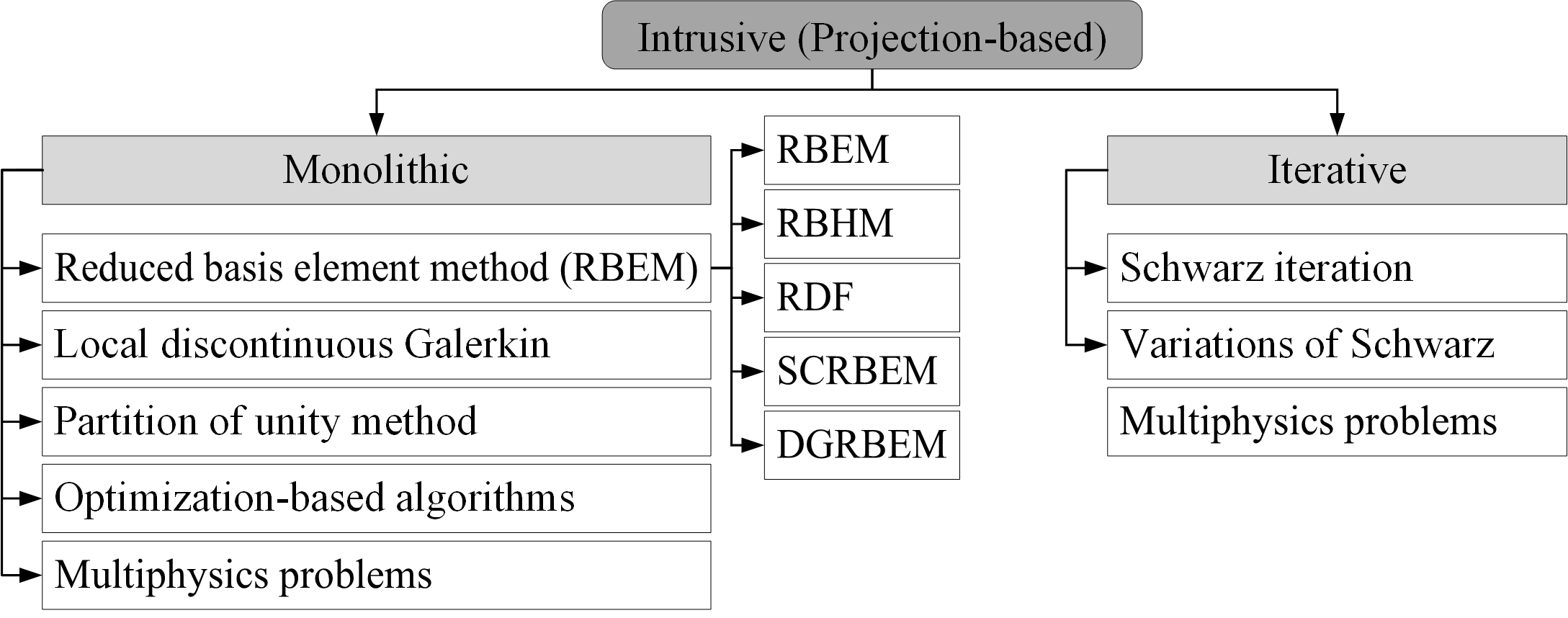}
    \caption{The categorization of intrusive techniques. The abbreviation: \emph{Reduced Basis Hybrid Method} (RBHM), The \emph{Reduced basis, Domain Decomposition, and Finite element} method (RDF), \emph{Discontinuous Galerkin Reduced Basis Element Method} (DGRBEM).}
    \label{fig:Intrusive_local_ROM}
\end{figure}

The technique is also commonly employed to derive intrusive reduced systems. Notably, it is utilized in our domain of interest: ROMs with domain decomposition. However, the usage of this approach is challenging. For intrusive local ROMs, the original PDEs have to be manipulated. In addition, the contributions of the interfaces and interactions between subdomains should also be included in the final formulations. In summary, intrusive ROMs assembling multiple partitions are effective but complex in terms of numerical and programming. Therefore, in this section, we provide a brief introduction to the various formulations of projection-based coupling algorithms, including both monolithic and iterative procedures. Fig. \ref{fig:Intrusive_local_ROM} displays the categories of the two groups \footnote{The abbreviation: {Reduced Basis Hybrid Method} (RBHM), The Reduced basis, Domain Decomposition, and Finite element method (RDF), {Discontinuous Galerkin Reduced Basis Element Method} (DGRBEM).}.

The differences and advantages of the two frameworks can be briefly summarized below. The monolithic methodology is more intrusive than the latter. It aims to create a large ROM involving all subdomains. Due to the compact reduced system, it is more favorable in terms of stability, convergence, and precision compared to iterative approaches \cite{wentland2024role}. However, the higher level of intrusiveness also results in extensive numerical and programming complexities. As shown in Fig. \ref{fig:Intrusive_local_ROM}, the monolithic category is more diverse.

In contrast, the iterative procedure aims to create multiple separate sub-problems corresponding to all partitions. The local ROMs are solved iteratively with updated boundary conditions from neighbors. The interface jump or error decreases iteratively, and the solver stops when it satisfies the predefined convergence criterion.

In the following sections, we will first give the basic formulation, namely, the strong and weak forms for PDEs. Then, the coupling techniques are explained. Note that a priori and a posteriori error analyses are usually included in the investigations, which are not illustrated in the review for brevity.

\subsection{Strong and weak form}

We first present the general formulation, which will serve as the foundation for the methods discussed in the following sections. The derivation of weak formulations from the strong form of a PDE has been extensively documented in the literature \cite{ern2004theory, canuto2006spectral, quarteroni2009numerical}. In this work, we will consider the strong form of a generic parametric PDE for the unknown $u \left ( x, t \right )$ variable, which is expressed as,
\begin{equation}
\label{eq:strong_form}
\begin{aligned}
    \dot{u} \left ( x,t \right ) + \mathcal{Q} \left ( u \left ( x, t \right ) ; \mu_j \right ) &= f \left ( x,t; \mu_j  \right ), & \quad \text{in} \ \Omega \times \left [ 0,t \right ], 
\end{aligned}
\end{equation}
with boundary conditions,
\begin{equation*}
\begin{aligned}
    \label{eq:bc_generic}
    u(x,t; \mu_j) &= g_D \left( x,t; \mu_j \right), & \quad \text{in} \ \Gamma_D \times \left [ 0,t \right ], \\
    \nabla u(x,t; \mu_j) \cdot \mathbf{n} &= \mathbf{g}_N \left( x,t; \mu_j \right), & \quad \text{in} \ \Gamma_N \times \left [ 0,t \right ] ,
\end{aligned}
\end{equation*} 
and initial values,
\begin{equation*}
\begin{aligned}
    u(x,0; \mu_j) &= u_0(x;\mu_j), & \quad \text{in} \ \Omega,
\end{aligned}
\end{equation*}
where $x$ and $t$ denote spatial and temporal coordinates, respectively. $\Omega \in \mathbb{R}^d$ is a $d$-dimensional domain, $d = 1,2,3$, with boundary $\partial \Omega$ decomposed into \emph{Dirichlet} $\Gamma_D$ and \emph{Neumann} $\Gamma_N$ conditions, $\partial \Omega = \Gamma_D \cup \Gamma_N$. In $\Gamma_D$ and $\Gamma_N$, conditions are assigned to $u_D$ and $\mathbf{g}_N$, respectively. The closure of the domain is written by $\overline{\Omega} \coloneqq \Omega \cup \partial \Omega$. We consider a parameter space $\mathcal{P}$, with elements $ \mu_j$. $u_0(x,0;\mu_j)$ is the initial condition and $f$ is a given source. Because we aim to construct the local ROM, regarding either linear or nonlinear physics, we consider $\mathcal{O}$ as an arbitrary differential operator, not necessarily linear. 

To simplify the notation, hereafter, we omit the dependent variables $\left( x, t\right) $ and parameters $ \mu_j$, e.g.,  $\dot{u}$ and $\mathcal{Q} \left ( u \right )$. The weak formulation is derived by means of the inner product with test functions $w \in \mathcal{W}$. The equations are required to hold only with respect to a certain test space $\mathcal{W}$, whose exact definition will vary among the formulations and will be pertinently detailed. Applying the integral-by-parts formula, volumetric integrals are divided into volumetric and surface parts. Indeed, $\mathcal{Q} (u) = \mathcal{Q}_\Omega (u) + \mathcal{Q}_\Gamma (u)$. $\Gamma$ denotes interfaces. Be aware that the $\mathcal{Q}$ is not necessarily equal to $\mathcal{Q}_\Omega$ and $\mathcal{Q}_\Gamma$. Then, by applying inner product with test function $w$, we rewrite equation \eqref{eq:strong_form} as the \emph{weak form},
\begin{equation}
\label{eq:weak_form_compact}
    \left(\dot{u}, w\right) + \mathcal{O} \left ( u, w \right ) = \mathcal{F}(w),  \quad \text{in} \ \Omega \times \left [ 0,t \right ],
\end{equation}
or
\begin{equation}
\label{eq:weak_form}
    \left(\dot{u}, w\right) + \mathcal{O}_\Omega \left ( u, w \right ) + \mathcal{O}_\Gamma \left ( u, w \right ) = \mathcal{F}(w), \quad \text{in} \ \Omega \times \left [ 0,t \right ],
\end{equation}
where $\left( \cdot, \cdot \right)$ denote the inner product operation, $\mathcal{O}_\Omega(u,w) = \left(\mathcal{Q}_\Omega(u), w\right)$, $\mathcal{O}_\Gamma(u,w) = \left(\mathcal{Q}_\Gamma(u), w\right)$, $\mathcal{F}(w) = \left. \left ( f, w \right ) \right|_{\overline{\Omega}}$. Remark that the global boundary contributions at $\Gamma_D$ and $\Gamma_N$ are involved in $\mathcal{F}(w)$ for simplification.

We may approximate $u$ as a function of trial functions $v_j \in \mathcal{V}$, with $\mathcal{V}$ a generic trial functional space. We get $u \approx \sum_{j=1}^n a_j v_j$ with $n$ the dimension of both $\mathcal{V}$ and $\mathcal{W}$ and $a_j$ the weighted coefficients of the basis vectors. To obtain a system containing $n$ equations we use inner product with test functions, $w_i\in\mathcal{W}, \ i=1,2, \cdots, n$, 
\begin{equation}
\label{eq:generic_ROM}
    \left(\sum_{j=1}^n \dot{a}_j v_j, w_i\right) + \mathcal{O}_\Omega \left ( \sum_{j=1}^n a_j v_j, w_i \right ) + \mathcal{O}_\Gamma \left ( \sum_{j=1}^n a_j v_j, w_i \right ) = \mathcal{F}(w_i) \quad v_j \in \mathcal{V}, w_i \in \mathcal{W}.
\end{equation}

Considering test/trial spaces with dimensions $N \ll n$, we obtain a ROM, where the $a_j$ are actually the unknowns of the problem. The fact that the number of these is several orders less than the cell/element number of a FOM simulation significantly promotes computational acceleration.

The algebraic form of the ROM \eqref{eq:generic_ROM} is 
\begin{equation}
\label{eq:generic_ROM_matrix}
    M \dot{\mathbf{a}} + C \mathbf{a} + D \mathbf{a}= F,
\end{equation}
where $\mathbf{a}$ is a vector containing coefficient of modes, and the element of matrices are addressed as $M_{ij} = (v_j, w_i)$, $C_{ij} = \mathcal{O}_\Omega (v_j, w_i)$, $D_{ij} = \mathcal{O}_\Gamma (v_j, w_i)$, and $F_{i} = \mathcal{F} (w_i) = (f, w_i)$.

In case the test and trial space are the same, the projection is called \emph{Bubnov--Galerkin} \cite{stenger1974convergence}, which is also referred to as \emph{standard Galerkin}. When the two sets of basis functions are different, it leads to the \emph{Petrov--Galerkin} or \emph{non--standard Galerkin} \cite{ern2004theory}. We kindly refer the reader to the books \cite{quarteroni2015reduced, hesthaven2016certified} for detailed instructions about the Galerkin projection-based ROM.

We will now start reviewing different methods available for the coupling task. 

\subsection{Monolithic}
\label{subsec:monolithic}
We begin our review from the first category of methods, known as \emph{Monolithic} techniques. These approaches rely on a set of local RBs, either a generic basis or subdomain-specific ones, to construct the global solutions. Interface contributions from adjacent partitions are incorporated into the original formulation either as additional terms or through a new set of equations enforcing the constraints. A noticeable advantage of these methods over Global ROMs is that each subdomain is governed by a smaller number of parameters, whether geometrical or physical, allowing for greater flexibility in parameterization. Moreover, for localized sensitivity and uncertainty analysis, the offline training effort is significantly reduced. Note that coupled systems for multiphysics are generally different from a single physics, so they are described in a separate section.

\subsubsection{Reduced Basis Element Method (RBEM)}
\label{subsubsec:RBEM}

\emph{Reduced Basis Element Method} (RBEM) was created by Yvon Maday and Einar M. Rønquist \cite{maday2002reduced} in 2002, based on the \emph{Mortar element method} \cite{bernardi1990new, belgacem1999mortar} that was already available and widely utilized for the decomposition of domains in high-resolution problems. As one of its main characteristics, the algorithm employs \emph{Lagrange multipliers} to enforce continuity across partition interfaces. As a result, the reduced system comprises (i) the original PDEs containing only volumetric terms, source, and boundary conditions (BCs); (ii) additional interface-integral terms corresponding to Lagrange multipliers; (iii) a new set of equations to constrain the jumps between adjacent subdomains \cite{maday2002reduced}. 

The decomposition into archetypes (see Section \ref{subsubsec:Generic_decomposition}) is utilized by RBEM to generate subdomains. The basis functions are constructed employing the localized training strategy described in Section \ref{subsubSec:Localized_training}. We recall that, for an instantiated component $\Omega_m^k$, the local data $u_m^k (\mu)$ can be represented as a linear combination of local basis functions $v_m^k$. At this stage, and for simplicity, we combine the indices of the archetype $k$ and the instantiation $m$ in a single one. Abusing of the notation, we denote them with $m$, namely $\Omega = \cup_{m=1}^{N_\Omega} \Omega_m$. That results in the local value
\begin{equation}
\label{eq:RBEM_solution}
    u_m (\mu) \approx \sum_{i=1}^{N_{\text{RB},m}} a_{m,i} v_{m,i}.   
\end{equation}

Now, we turn to a multiple-partition system and select one subdomain $\Omega_m$ as an example to explain the construction of a ROM based on RBEM. The interface connected to $\Omega_m$ and its neighbor $\Omega_n$ is denoted as $ \Gamma_{[mn]} \coloneqq \partial \Omega_m \cap \partial \Omega_n$. The subscript $[mn]$ is a combined index to designate a single face. As $\Omega_m$ might have several adjacent partitions, the indices of all neighbors are included in a set $N_\Gamma(m) = \left\{ n | \partial \Omega_m \cap \partial \Omega_n \neq \emptyset \right\}$.

Before illustrating the governing equations, we introduce notations for clarification. We start with the Hilbert space in $\Omega$, $H^1(\Omega)$. The test space belongs to $H^1(\Omega)$ and its trace is equal to $g_D$ on $\Gamma_D$ is $\mathcal{W}_{g_D} \coloneqq H^1_{g_D, \Gamma_D} (\Omega)$. Thus, for $\Omega_m$, $\mathcal{W}_{m, 0} \coloneqq H^1_{0, \Gamma_D} (\Omega_m)$ denotes functions in the test space that are homogeneous on the local $\Gamma_D \in \partial \Omega_m$. The trial space is $\mathcal{V}_m \coloneqq H^1 (\Omega_m)$. 

In $\Gamma_{[mn]}$, we also define the functional basis space $\mathbb{L}_{[mn]}$ \footnote{In the literature, $\mathbb{L}_{[mn]} \coloneqq H^{-1/2}_{00}(\Gamma_{[mn]})$, a definition of the latter exceeds the extent of this review, see \cite{quarteroni2009numerical} for further details.} and the basis vectors $\zeta_{[mn],i} \in \mathbb{L}_{[mn]}$. A possible choice of $\mathbb{L}_{[mn]}$ is as the space generated by a set of orthogonal polynomials on the interface $\Gamma_{[mn]}$, e.g., the Chebyshev polynomials. We refer interested readers to the references \cite{lovgren2006reduced, pegolotti2021model} for details. The Lagrange multiplier defined on $\Gamma_{[mn]}$ can be expressed in terms of the basis in the interface as $\lambda_{[mn]} \approx \sum_{i=1}^{N_{\lambda, [mn]}} q_{[mn],i} \zeta_{[mn],i}$, where $N_{\lambda, [mn]}$ is the dimension of the $\mathbb{L}_{[mn]}$ space.

The generic RBEM governing equation for $\Omega_m$ and its interfaces is written as
\begin{equation}
    \label{eq:generic_RBEM}
    \left(\dot{u}, w_{m,i}\right) + \mathcal{O}_\Omega \left ( u, w_{m,i} \right ) + \sum_{n \in N_{\Gamma} (m)} \int_{\Gamma_{[mn]}} \lambda_{[mn]} w_{m,i}  = \mathcal{F}(w_{m,i}), \quad w_{m,i} \in \mathcal{W}_{m,0},
\end{equation}
where The global boundary terms are included in $\mathcal{F}(w_{m,i})$ and $\mathcal{O}_\Gamma \left ( u, w_{m,i} \right ) = 0$ because $\left. w_{m,i} \right|$

For every interface, $\lambda_{[mn]} = - \lambda_{[n m]}$. Thus,
\begin{equation}
\label{eq:RBEM_interface}
    \int_{\Gamma_{[mn]}} \zeta_{[mn],i} \left( u_m - u_n \right) = 0 \quad \forall \zeta_{[mn],i} \in \mathbb{L}_{[mn]}.
\end{equation}
Note that the global Dirichlet boundary conditions are weakly imposed by the Lagrange multipliers in equation \eqref{eq:RBEM_interface}.

In case the operators $\mathcal{O}$ are linear, we may assemble the equation \eqref{eq:generic_RBEM} and \ref{eq:RBEM_interface} of each subdomain $m$ into a global algebraic form,
\begin{equation}
\label{eq:RBEM_matrix_local}
\begin{split}
    M_m\dot{\mathbf{a}} + C_m\mathbf{a} + D_m\mathbf{q} &= F_m \\
    Q_m\mathbf{a}&= 0,
\end{split}
\end{equation}
where $M_{m, ij} = (v_{m,j}, w_{m,i})$, $C_{m, ij} = \mathcal{O}_\Omega(v_{m,j}, w_{m,i})$ are diagonal matrices with a dimension of $N_{\text{RB},m}$. $D_m$ and $Q_m$ are a collection of many sub-blocks, i.e., $D_m = \left\{ \left. D_{[m n]} \right| n \in N_\Gamma (m) \right\}$, $Q_m = \left\{ \left. Q_{[m n]} \right| n \in N_\Gamma (m) \right\}$. Note that they are block matrices but not diagonal, and the structures depend on their neighbors. Each sub-block in $D_m$ contains elements $D_{[m n], ij} = \int_{\Gamma_{[m n]}}\zeta_{[m n], j}\ w_{m,i}$ of dimensions $N_{\text{RB},m} \times N_{\lambda, r}$, and those in $Q_m$ is $Q_{[m n], ij} = \int_{\Gamma_{[m n]}} v_{m,j} \zeta_{[m n], i} - \int_{\Gamma_{[m n]}} v_{n,j} \zeta_{[m n], i}$ of dimensions $ N_{\lambda, r} \times N_{\text{RB},m}$.

The total matrix can be constructed by stacking the various sub-matrices of index $m$,
\begin{equation}
\label{eq:RBEM_matrix_global}
\begin{split}
    M\dot{\mathbf{a}} + C\mathbf{a} + D\mathbf{q} &= F, \\
    Q\mathbf{a}&= 0,
\end{split}
\end{equation}
The $M$ is a diagonal block matrix, that is $M = \text{diag} \left( M_1, \cdots, M_m \right)$. Similarly, we have $C = \text{diag} \left( C_1, \cdots, C_m \right)$. The $M$ and $C$ are square matrices with a dimension of $\sum_{m=1}^{N_{\Omega}} N_{\text{RB},m}$. The $D$ and $Q$ can still be regarded as a collection of sub-matrices $\left\{D_1, \cdots, D_m \right\}$ and $\left\{Q_1, \cdots, Q_m \right\}$, however, not diagonal. We argue that by adopting extra linearization steps, the above procedures are also applicable for non-linear cases and build similar algebraic systems.

To the best of our knowledge, the \emph{Reduced Basis Element Method} (RBEM) was the first technique developed to combine domain decomposition and ROMs. Maday and Rønquist proposed this approach in \cite{maday2002reduced}, applying it to the Laplacian equation for a geometrically parameterized thermal fin problem (Fig. \ref{fig:generic_division_thermal_fin}) \cite{maday2004reduced}. Subsequent work with Lovgren extended RBEM to handle the steady Stokes and Navier-Stokes equations in a 2D blood vessel geometry \cite{lovgren2006stokes, lovgren2006reduced}.

Further developments of RBEM have since emerged. Chen et al. \cite{chen2011seamless} adapted the technique for the time-harmonic Maxwell's equation. A recent contribution utilizing RBEM to the unsteady Navier-Stokes flow in a 3D blood vessel can be found in \cite{pegolotti2021model}. We present their model in Section \ref{subsubSec:Localized_training} as an example for generic spatial decomposition (see Fig. \ref{fig:generic_division_blood_vessel}) and RBs for subdomains (see Fig. \ref{fig:generic_geometry_RB_blood_vessel}). 

Except for the incorporation with RBEM, the Lagrange multipliers can \emph{glue} different interfaces. Thus, the aforementioned formulation can be used to assemble several arbitrary geometries without algorithmic modifications. Charbel Farhat et al. followed this ideology to analyze the Helmholtz problems with plane waves in multiple subdomains \cite{farhat2009domain}. Recent studies by Amy de Castro and Irina Tezaur et al. \cite{de2022lagrange, de2023partitioned} have further demonstrated the method's capability. They managed to couple the FOM-ROM and ROM-ROM systems of an advection-diffusion equation. Note that a FOM-ROM system contains two subdomains, and they are simulated by a high-fidelity FOM and a ROM, respectively. Additionally, the ideology can be extended to a coupling of multiple FOMs and ROMs.

\paragraph{Reduced Basis Hybrid Method}

Based on RBEM, Iapichino et al. \cite{iapichino2012reduced} proposed a modification, the so-called \emph{Reduced Basis Hybrid Method} (RBHM). They intended to parametrize a steady Stokes problem for a blood vessel. The approach is almost the same as the RBEM except for the construction of RBs. To ensure continuity and consistency of the normal stress during the online stage, they used \emph{hybrid} RBs, which contain three components: (i) high-fidelity FE results collected from training models, $\hat{\mathbf{u}}^k (\mu_j)$; (ii) coarse finite element solutions $\mathbf{u}^{\text{H}} (\mu_p)$ of the whole problem considering parameters $\mu_p$ that one wants to predict with the ROM, and (iii) velocity supremizer basis vectors $\mathbf{s}_{m}^{\text{H}} (\mu_p)$. The last one entitles the basis to fulfill the inf-sup condition \cite{quarteroni2015reduced}, enhancing the numerical stability (details on \emph{supremizer enrichment} technique are provided by references \cite{rozza2005shape, rozza2007stability, hesthaven2016certified}). In each subdomain $\Omega_m$, the solution is expressed as
\begin{equation}
\begin{aligned}
    \mathbf{u}_m (\mu_p) & \approx T^k (\boldsymbol{\xi}_m^k) \cdot \left( \sum_{i=1}^{N_\mu} a_{m,i} \hat{\mathbf{u}}^k (\mu_i) + \sum_{i=1}^{N_\mu} c_{m,i} \hat{\mathbf{s}}^k (\mu_i) \right) + b \mathbf{u}^{\text{H}}_m (\mu_p) + d \mathbf{s}^{\text{H}}_m (\mu_p) \\
    & = \sum_{i=1}^{N_\mu} a_{m,i} \mathbf{v}_{m,i} + b \mathbf{u}_{m}^{\text{H}} (\mu_p) + \sum_{i=1}^{N_\mu} c_{m,i} \mathbf{s}_{m,i} + d \mathbf{s}_{m}^{\text{H}} (\mu_p).
\end{aligned}
\end{equation}
Note that this formulation substitutes the approximation given by equation \eqref{eq:RBEM_solution}. Afterwards, a global ROM approximation is obtained following the same steps of RBEM. 

\subsubsection{The Reduced basis, Domain Decomposition, and Finite element method (RDF)}
The \emph{Reduced basis, Domain Decomposition, and Finite element} method (RDF) procedure receives its name from the combination of the three approaches. It is designed for non-overlapping decomposition problems \cite{iapichino2012thesis, iapichino2016reduced}. 

The main difference with respect to RBEM and RBHM is that RDF does not need Lagrange multipliers that impose interface equality. Instead, the FE equations of the problem are projected to obtain ROM systems involving the interface basis functions that ensure the continuity of the solution. Indeed, each subdomain of the RDF can be regarded as a finite element in the original FEM-based FOM. The construction of RBs is explained as follows.

The internal RDF RBs are computed similarly to RBEM, using solutions of sub-problems defined in reference shapes. The sub-models are parameterized with random combinations \emph{Fourier series} or \emph{Lagrange polynomials} as Dirichlet BCs. Lifting functions are applied to obtain basis vectors vanishing on Dirichlet BCs. The lifted results are employed as RBs, which satisfy 
\begin{equation}
\label{eq:RDF_test_trial_space}
    \mathcal{V}_{m,0} = \mathcal{W}_{m,0} \coloneqq H^1_{0, \Gamma_D} (\Omega_m).
\end{equation}

The authors demonstrated that the additional interface parameters allow RDF to represent complex interface profiles better and thus capture variances of the final global solution more accurately. 

The FEM shape functions $\phi$ defined at the interfaces are employed to construct ROMs. Thus, the unknowns or Degrees of Freedom (DoFs)  \footnote{The unknowns of an approximation system can be named DoFs.} of the ROM are the same as FOM for interfaces.

Let us now consider $\hat{N}^h_{{\Gamma}}$ nodes on a reference internal interface $\hat{\Gamma}$. We have the basis
\begin{equation}
\label{eq:interface_space}
    \Phi^{\hat{\Gamma}} \coloneqq \text{span} \left\{ \phi^{\hat{\Gamma}}_l \in \left. \mathcal{V}^h \right|_{\hat{\Gamma}},\ l = 1, \cdots, N^h_{\hat{\Gamma}} \right\}.
\end{equation}
where $\phi^{\hat{\Gamma}}_l$ are basis functions of the FE nodes and $\mathcal{V}^h$ is the high dimensional FE basis space. 

It follows that the interfaces existing in the entire model can be transformed from $\hat{N}_{{\Gamma}}$ generic interfaces. The $e^{\text{th}}$ type of the generic interface is noted as $\hat{\Gamma}^e$. Each $\hat{\Gamma}^e$ them is instantiated $N_r^e$ times, in which the $r^{\text{th}}$ instantiation is $\Gamma^e_{r}$ and corresponding basis space $\Phi^e_r$. Thus, the total number of internal faces is $N_{\Gamma} = \sum_{e=1}^{\hat{N}_\Gamma} N_r^e$.

To simplify the notation in our explanation, we condense the archetype and realization index into a single one. We assume $\Omega = \cup_{m=1}^{N_\Omega} \Omega_m$ and $\Gamma = \cup_{r=1}^{N_\Gamma} \Gamma_r$. Therefore, the composed basis space for RDF systems is
\begin{equation}
\label{eq:RDF_basis_space}
    \mathcal{U} \coloneqq \bigcup_{m=1}^{N_\Omega} \mathcal{V}_{m,0} \oplus \bigcup_{r=1}^{N_\Gamma} \Phi_r,
\end{equation}
where $\Phi_r$ is defined on $\Gamma_{r}$. Take into account that the local spaces are defined only in subdomains, so we use a union operator instead of a summation to form the global space.

Remind that the interface set of $ \Omega_m $ is $\mathcal{I}_m$. Therefore, the local solution $u_m (\mu)$ can be represented by 
\begin{equation}
\label{eq:RDF_solution}
    {u}_{m} (\mu) \approx \sum_{i=1}^{N_{m}} a_{m,i} v_{m,i} + \sum_{\Gamma_r \in \mathcal{I}_m} \sum_{i=1}^{N_{r}} b_{r, i} \phi_{r,i},
\end{equation}
where $\phi_{r,i}$ are basis vector instantiated on $\Gamma_r$ that has in total $N_r$ nodes.

For a partition $\Omega_m$, we can project equation \eqref{eq:generic_ROM} onto $\mathcal{V}_{m,0}$ and obtain
\begin{equation}
\label{eq:RDF_inner}
    \left( \dot{u}_m, v_{m,i} \right) + \mathcal{O}_{\Omega} (u_m, v_{m,i}) + {\mathcal{O}}_{\Gamma} (u_m, v_{m,i}) = \mathcal{F}_m(v_{m,i}), \quad v_{m,i} \in \mathcal{V}_{m,0}.
\end{equation}

Suppose the two adjacent subdomains of $\Gamma_r$ are $\Omega_m$ and $\Omega_n$, by considering test space $\Phi_r$, equation \eqref{eq:generic_ROM} is reformulated as
\begin{equation}
\label{eq:RDF_face}
    \left. \left( \dot{u}_m, \phi_{r,l} \right)\right|_{\Gamma_r} + \left. \left( \dot{u}_n, \phi_{r,l} \right)\right|_{\Gamma_r} + {\mathcal{O}}_{\Gamma} (u_m, \phi_{r,l}) + {\mathcal{O}}_{\Gamma} (u_n, \phi_{r,l}) = \mathcal{F}_m(\phi_{r,l}) + \mathcal{F}_n(\phi_{r,l}) \quad \phi_{l} \in \Phi_{r},
\end{equation}
where we recall that $\Phi_{\Gamma}$ only defined in $\Gamma$, so the volumetric terms are eliminated.

For linear operators $\mathcal{O}$, the matrix form of a RDF system is constructed by composing all subdomains and interfaces, and we take
\begin{equation}
\label{eq:RDF_matrix}
\begin{aligned}
    M \dot{\mathbf{a}} + N \dot{\mathbf{b}} + C \mathbf{a} + D \mathbf{b} & = F_\Omega, \\
    K \dot{\mathbf{a}} + H \dot{\mathbf{b}} + Q \mathbf{a} + G \mathbf{b} & = F_\Gamma,
\end{aligned}
\end{equation}
where every matrix is a collection of many sub-matrices. The algebraic system also applies to nonlinear cases, in which linearized steps are employed to separate the inner and interface contributions.

The $M$, $N$, $K$, $H$ are block diagonal matrices, so we take $M = \text{diag} \left( M_1, \cdots, M_m \right)$, and $M_{m,ij} = (v_{m,j}, v_{m,i})$. The other three matrices, $N$, $K$, $H$, can be expressed similarly as the inner product of $v_{i}$ and/or $\phi_l$, but are not detailed here for compactness. Then, we have $C = \text{diag} \left( C_1, \cdots, C_m \right)$, with $C_{m,ij} = \mathcal{O}_\Omega \left( v_{m,j}, v_{m,i} \right)$. The structure of $Q$ is similar and the items are computed by $Q_{m, lj} = \mathcal{O}_\Gamma(v_{m,j}, \phi_{r,l})$. 

In contrast, as most of the partitions are multi-connected with others, $D$ and $G$ are constructed interface-wise, that is, $D = \left\{D_1, \cdots, D_r \right\}$, and $D_{r, il} = \left. \mathcal{O}_\Gamma(\phi_l, v_i) \right|_{\Gamma_r}$, $G_{r, lj} = \left. \mathcal{O}_\Gamma(\phi_j, \phi_l) \right|_{\Gamma_r}$. Note that here, we utilized the compact sub-index $r$ to designate the interface instead of the more descriptive $_{[mn]}$ sub-index.

The RDF provides a general framework to achieve global approximations with local RBs. Laura Iapichino et al. applied it for geometrical parameterization of an elliptic problem, the steady heat equation in the thermal fin \cite{iapichino2012thesis, iapichino2016reduced}. 

Immanuel Martini et al. \cite{martini2015reduced} derived a RDF-type reduced system while using a distinct strategy to construct reduced subspaces. The RBs are computed in two stages. The first step is to conduct a few global FOM simulations. Then, the interface modes are computed adopting values at the internal interface $\phi_{r,l}(\mu_j)$ of snapshots. Secondly, the interior RB is computed from a series of local problems with interface modes assigned as lifting functions at boundaries. Thus, the steady equation \eqref{eq:RDF_inner} defined in a subdomain is written by,
\begin{equation}
    \mathcal{O}_{{\Omega}} ({u}, {w}_i) = \mathcal{F}({w}_i) - \mathcal{O}_{{\Gamma}} (\phi_{r,l}, {w}_i),  \quad {w}_i \in {\mathcal{W}}_{0} \text{ and } \phi_{r,l} \in {\Phi}.
\end{equation}

Note that the last term in the right-hand side is formed from the standard lifting procedure and is thus known \footnote{That is equivalent to considering the second term in the left-hand side of equation \eqref{eq:RDF_solution} as known.}. The internal and boundary RBs, as well as the velocity supremizer basis, are gathered and orthonormalized by the \emph{Gram-Schmidt} process to produce the compact RBs. Finally, the coupled ROMs are solved in the frame of the RDF approach. The authors managed to model a steady-state problem (Stokes and Darcy equations) in a non-overlapping two-subdomain system, which considers geometric parameters. 

\subsubsection{Static Condensation Reduced Basis Element Method (SCRBEM)}
\label{subsubsec:SCRBEM}
\emph{Static Condensation Reduced Basis Element Method} (SCRBEM) was proposed for domains assembled by non-overlapping archetype components. Similar to the RBEM, generic decomposition (Section \ref{subsubsec:Generic_decomposition}) and localized training (Section \ref{subsubSec:Localized_training}) were employed in the methodology. As usual, the domain is decomposed $\Omega = \cup_{m=1}^{N_\Omega} \Omega_m$ with internal faces $\Gamma = \cup_{r=1}^{N_\Gamma} \Gamma_r$.

This approach is conceptually similar to the RDF method illustrated in the previous section. However, there are subtle differences, particularly in the definition and generation of the basis spaces, which will be explicitly highlighted later. The framework incorporates two key techniques: port reduction and static condensation. In the following sections, we will briefly discuss their fundamental principles and implementation.

\paragraph{Port Reduction RBEM}
\label{sec:portredrbem}
The space is generated using a \emph{Port Reduction} (PR) \cite{eftang2012adaptive, huynh2013static} technique. This aims to separate the space into two pieces for the surface and inner zones. Those take the names of \emph{port} and \emph{bubble} space, respectively.  

The port space $\Phi_r$ is defined only for interfaces $\Gamma_r$, and $\Phi_r \coloneqq H^1(\Gamma_r)$. The bubble space is the range of basis vectors that vanish on all boundaries of the local domain $\mathcal{V}_{m,0} \coloneqq \text{span} \{v \in H^1 (\Omega_m)\ |\ v|_{\partial \Omega_m} =0 \}$. Note that  contrary to the definition of \eqref{eq:RDF_test_trial_space}, now the vectors are homogeneous in the whole $\partial \Omega_m$. With these changes, like in equation \eqref{eq:RDF_basis_space}, the global test space is written as 
\begin{equation}
    \mathcal{W} \coloneqq \bigcup_{m=1}^{N_\Omega} \mathcal{V}_{m,0} \oplus \bigcup_{r=1}^{N_\Gamma} \Phi_r.
\end{equation}

We endeavor now strategies to span $\Phi_r$ and $\mathcal{V}_{m,0}$. The reduced port and bubble subspaces can be generated from two different problems, as presented in \cite{benaceur2022port}. 

\begin{figure}[h]
    \centering
    \begin{subfigure}[b]{0.5\linewidth}
        \centering
        \includegraphics[height=2cm]{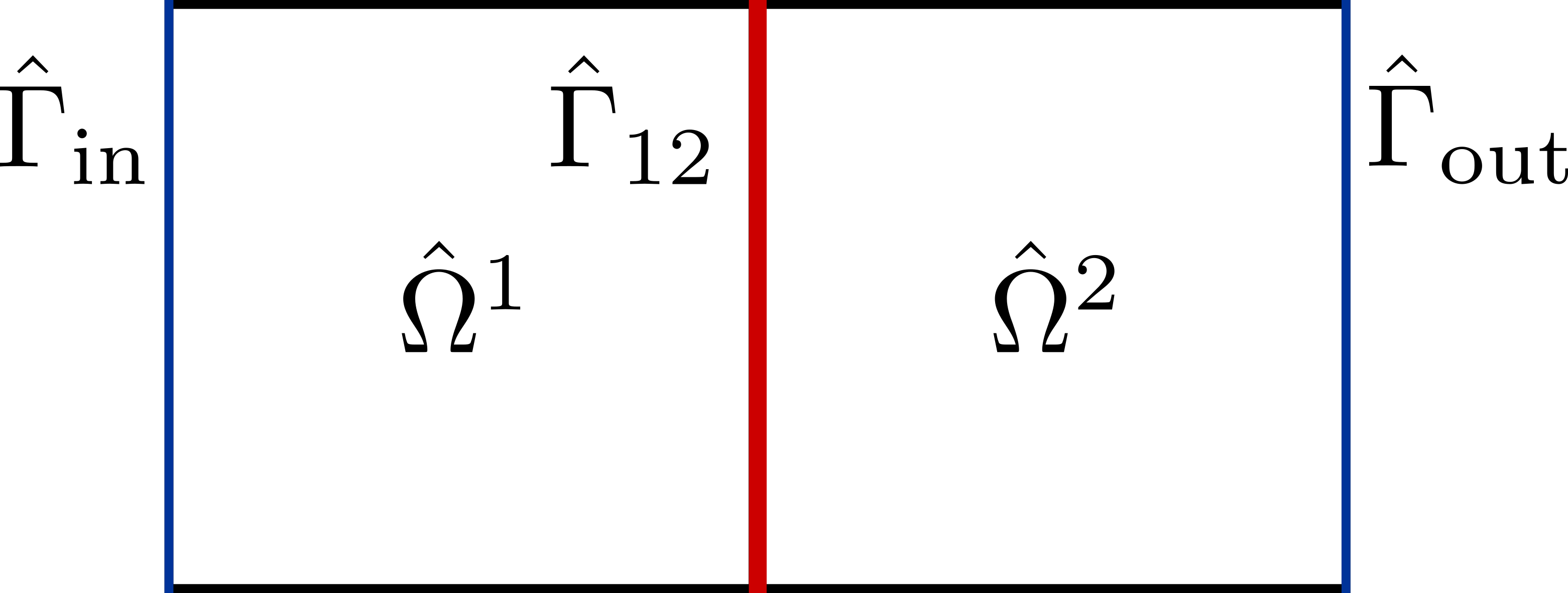}
        \caption{Port training problem.}
        \label{fig:port_training_problem}
    \end{subfigure}
    \hspace{2em}
    \begin{subfigure}[b]{0.3\linewidth}
        \centering
        \includegraphics[height=2cm]{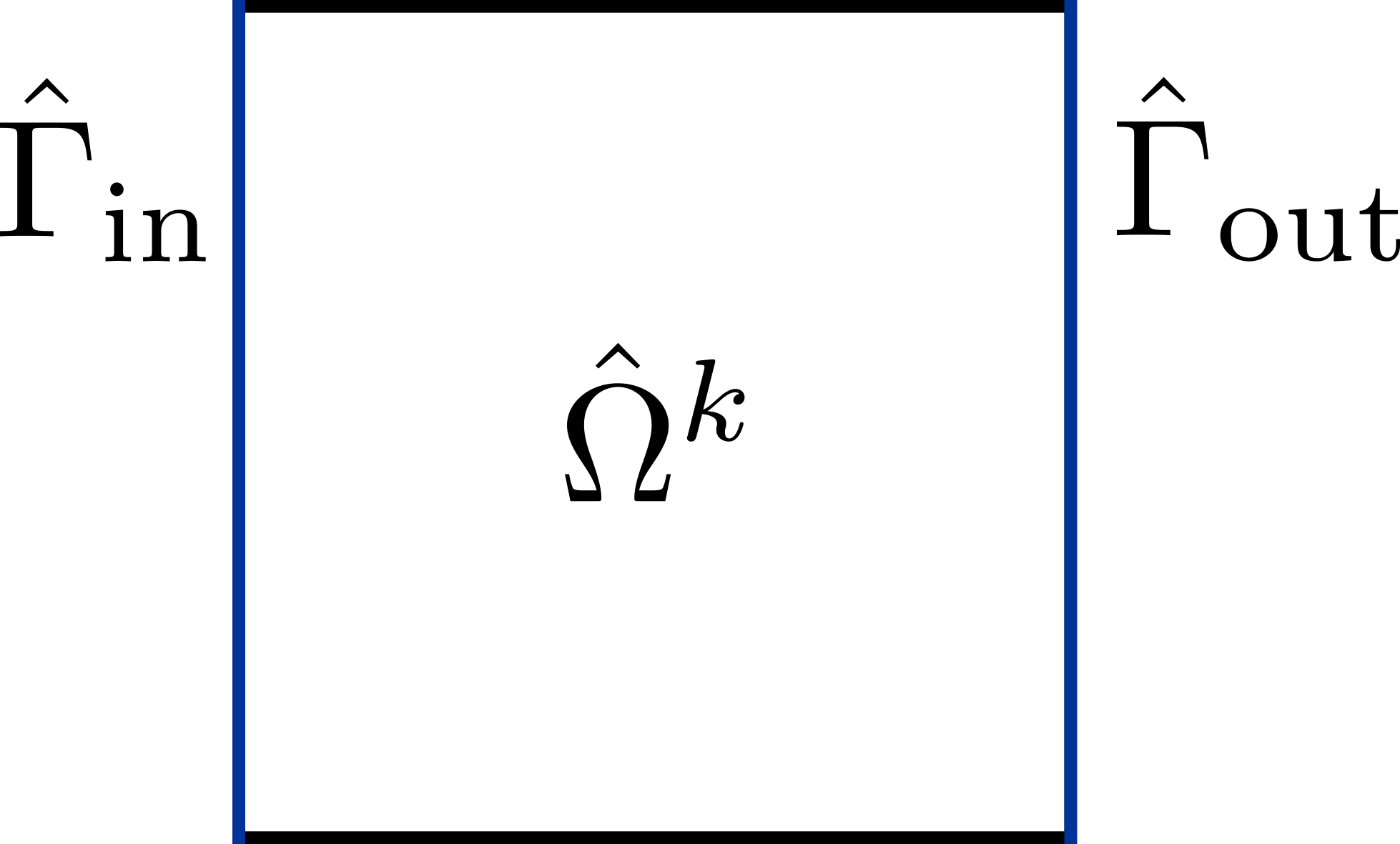}
        \caption{Bubble training problem.}
        \label{fig:bubble_training_problem}
    \end{subfigure}
    \caption{Port and Bubble training. Redraw based on \cite{benaceur2022port}. (a) The port training model for a generic $\hat{\Gamma}_{12}$ connects two reference partitions $\hat{\Omega^1}$ and $\hat{\Omega^2}$. Parametric Dirichlet BCs are assigned for $\Gamma_{\text{in}}$ and $\Gamma_{\text{out}}$. The results of interface $\Gamma_{mn}$ are collected as reduced Port space. (b) The bubble training for $\hat{\Omega}^k$ with parameterized BCs $\hat{\Gamma}_{\text{in}}$ and $\hat{\Gamma}_{\text{out}}$.}
    \label{fig:port_bubble_training_problem}
\end{figure}

The port training problem is sketched in Fig. \ref{fig:port_training_problem}. It solves a two-subdomain problem \footnote{The authors \cite{benaceur2022port} solve the original PDEs at the system. } with parameterized Dirichlet BCs on $\hat{\Gamma}_{\text{in}} $ and $\hat{\Gamma}_{\text{out}} $. Finally, the solution at the interface $\hat{\Gamma}_{\text{12}} $ is extracted to form a snapshot matrix. This will be post-processed with a dimensional reduction technique to obtain the basis $\Phi_r$. 

The bubble RB is obtained by solving the single domain problem of Fig. \ref{fig:bubble_training_problem}, in which the port basis functions are randomly weighted as boundary conditions $\hat{\Gamma}_{\text{in}}$, $\hat{\Gamma}_{\text{out}}$. Then, lifting functions \footnote{The authors \cite{benaceur2022port} adopt solutions of Stokes equations as lifting functions for Navier-Stokes problems.} are utilized to homogenize Dirichlet BCs, and POD is utilized to compute the reduced bubble space. More descriptive training steps are presented in \cite{benaceur2022port}.

In case the trial space $\mathcal{U}$ is the same as $\mathcal{W}$, the local approximations is expressed by ${u}_{m} (\mu) \approx \sum_{i=1}^{N_{m}} a_{m,i} v_{m,i} + \sum_{\Gamma_r \in \mathcal{I}_m} \sum_{i=1}^{N_{r}} b_{r, i} \phi_{r,i}$.

By applying Galerkin projection of equation \ref{eq:generic_ROM} to bubble $\mathcal{V}_{m,0}$ and port $\Phi_r$ spaces for $\Omega_m$, we aim to find $u_m$ that satisfy
\begin{equation}
\label{eq:SCRBEM_eq1}
\begin{aligned}
    \left( \dot{u}_m, v_i \right) + \mathcal{O}_{\Omega} (u_m, v_{m,i}) &= \mathcal{F}_m(v_{1,i}) &\quad & v_{m,i} \in \mathcal{V}_{m,0}, \\
    \left( \dot{u}_m, \phi_{r,l} \right)|_{\Gamma_r} + \left( \dot{u}_n, \phi_{r,l} \right)|_{\Gamma_r} + {\mathcal{O}}_{\Gamma} (u_m, \phi_{r,l}) + {\mathcal{O}}_{\Gamma} (u_n, \phi_{r,l}) &= \mathcal{F}_m(\phi_{r,l}) + \mathcal{F}_n(\phi_{r,l}) & \quad & \phi_{r,l} \in \Phi_{r},
\end{aligned}
\end{equation}
where $\Omega_m$ and $\Omega_n$ are connected by $\Gamma_r$. This system, which looks exactly equal to equations \eqref{eq:RDF_inner}-\eqref{eq:RDF_face} represents a different reality, hidden in its test-trial spaces.

The algebraic form of equation \eqref{eq:SCRBEM_eq1} is constructed by composing all subdomains and interfaces,
\begin{equation}
\label{eq:SCRBEM_matrix}
\begin{aligned}
    M \dot{\mathbf{a}} + C \mathbf{a} & = F_\Omega, \\
    H \dot{\mathbf{b}} + G \mathbf{b} & = F_\Gamma,
\end{aligned}
\end{equation}
where the structure of the matrices is analogous to the one of system \eqref{eq:RDF_solution} \footnote{That is the so-called port-reduction RBEM presented in \cite{benaceur2022port}.}.

Amina Benaceur and Anthony Patera applied the aforementioned procedure to model steady-state Navier-Stokes flow (including heat transport) for two-dimensional (2-D) rivers, considering two physical parameters are considered, i.e., Reynolds number and Prandtl number \cite{benaceur2022port}.

\paragraph{Static condensation}
\label{sec:sct}
We now illustrate the procedure of the so-called \emph{static condensation}. The basic idea of this method is to divide a domain into inner and boundary components and then make the DoFs in the interior enslaved by the boundary DoFs. Let us consider the matrix system of equations \eqref{eq:RDF_matrix} for steady-state conditions. The internal DoFs can be eliminated considering $\mathbf{a} = C^{-1}(F_\Omega-D\mathbf{b})$, to get the condensed system,
\begin{equation}
\label{eq:static_condenstaion_matrix}
    (G-QC^{-1}D) \mathbf{b} = F_{\Gamma}-QC^{-1}F_\Omega,
\end{equation}
where only DoFs in the boundary remain, and the internal DoFs are expressed in terms of the boundary ones.

Remark that the strategy of static condensation (i.e., eliminating interior DoFs) is not unique. Let us now extend this idea to more complex systems \cite{huynh2013complex}. Hereafter, we operate with the steady state of the general form \eqref{eq:generic_ROM} to simplify the analysis, but readers should keep in mind that the procedure also holds for time-dependent scenarios.

The framework of \emph{static condensation} involves port and bubble spaces. The solution is generated through a distinct strategy with three steps:
\begin{enumerate}
    \item 
Firstly, we intent to find $\hat{u}^k_{0} \in \hat{\mathcal{V}}^k_{0}$ as the solution of problem \eqref{eq:generic_ROM} defined in $\hat{\Omega}^k$,
\begin{equation}
\label{eq:SCRBEM_homo}
    \mathcal{O}_\Omega (\hat{u}^k_{0}, \hat{v}^k_{i}) = \mathcal{F}(\hat{v}^k_{i}) \quad \hat{v}^k_{i} \in \hat{\mathcal{V}}^k_{0},
\end{equation}
where $\hat{\mathcal{V}}^k_{0}$ is the bubble space defined in the archetype $\hat{\Omega}^k$.  $\hat{u}^k_{0} \in \hat{\mathcal{V}}^k_{0}$ satisfies the governing equation with homogeneous boundary conditions. Thus, $\mathcal{O}_{\Gamma}(\hat{u}^k_{0}, \hat{v}^k_{i})=0$.

\item
The next stage is to find port bases that account for local boundary contributions. The port modes of a given $\hat{\Gamma}^e$ are obtained by satisfying a boundary value problem of arbitrary choice. In \cite{huynh2013complex}, the Laplace equation is utilized. $\hat{\psi}^e_j$ can be computed by solving
\begin{equation*}
    \int_{\hat{\Omega}^k} \nabla \hat{\psi}^e_j \cdot \nabla \hat{v}^k_{0} = 0, \quad \hat{v}^k_{0} \in \hat{\mathcal{V}}^k_0
\end{equation*}
with Dirichlet BCs defined as,
\begin{equation*}
\begin{aligned}
    \hat{\psi}^e_j &= \hat{\chi}^e_j, & \quad & \text{on }\ \hat{\Gamma}^e, \\
    \hat{\psi}^e_j &= 0,  & \quad & \text{on }\ \partial \hat{\Omega}^k \setminus \hat{\Gamma}^e.
\end{aligned}
\end{equation*}
The $\hat{\psi}^e_j$ span a space on $\hat{\Gamma}^e$, $\hat{\Psi}^e = \left\{ \hat{\psi}^e_1, \cdots, \hat{\psi}^e_j \right\}$ of dimension $\hat{N}^e_\Gamma$. The $\hat{\chi}^e_j$ are eigenfunctions (see more in \cite{quarteroni2009numerical}) associated with the problem defined in $\hat{\Gamma}^e$, and they belong to the port space $\Phi_r$ defined above in Section \ref{sec:portredrbem}. 

\item
Then, we solve a series of subproblems for port modes on $\hat{\Gamma}^e$. Namely find $\hat{\varphi}^{e}_{0,j} \in \hat{\mathcal{V}}^k_0$,
\begin{equation}
    \mathcal{O}_\Omega (\hat{\varphi}^{e}_{0,j}, \hat{v}^k_{i}) = -\mathcal{O}_\Omega (\hat{\psi}^e_j, \hat{v}^k_{i}) \quad \hat{v}^k_{i} \in \hat{\mathcal{V}}_0^k.
\end{equation}
Here, the number of $\hat{\varphi}^{e}_{0,j}$ is equal to the number of eigenfunctions on each port space, i.e., $\hat{N}^e$.

We now create the port space $\hat{\Theta}_e \equiv \text{span} \left\{ \hat{\theta}_j = \hat{\varphi}^{e}_{0,j} + \hat{\psi}^e_j, j = 1, \cdots, \hat{N}^e \right\} $. 

\end{enumerate}

We transform the spaces, bases, and amplitudes to the instantiations $\hat{\Omega}_m\rightarrow\Omega_m$ and $\hat{\Gamma}_r\rightarrow\Gamma_r$ as usual. The final basis space can be constructed as
\begin{equation*}
    \mathcal{W} \coloneqq \bigcup_{m=1}^{N_\Omega} \mathcal{V}_{m,0} \oplus \bigcup_{r=1}^{N_\Gamma} {\Theta}_r.
\end{equation*}
Note that $\Theta_r \neq \Phi_r$. $\Phi_r$ is defined only on $\Gamma_r$ whilst $\Theta_r$ is built in both interior and surface by the procedure detailed above. 
Thus, the local solution can be represented as
\begin{equation*}
    u_m  \approx u_{m,0} + \sum_{\Gamma_r \in \mathcal{I}_m} \sum_{j=1}^{N^p_{r}} b_{r, j} \left( {\varphi}_{0, r, j} +  {\psi}_{r,j} \right),
\end{equation*}
where $\mathcal{I}_m$ is the interface set, $N^p_{r}$ is the number of retained eigenfunctions of $\Gamma_r$, and $b_{r, i}$ are the unknown amplitudes in $\theta_{r,j}$ modes. The first term $u_{m,0}$ satisfies the local PDEs with homogeneous interface conditions, and the second term represents the contribution of local BCs defined at interfaces.

Unlike the steady system of equation \ref{eq:SCRBEM_homo}. This time, we utilize a different test space, $\Theta_{r}$. Since the latter is not homogeneous in the interfaces, the term $\mathcal{O}_\Gamma ({u}_{m,0}, \theta_i)$ does not cancel. By applying Galerkin projection of the steady system \ref{eq:generic_ROM} onto the basis $\Theta_{r}$, we get a formulation with the only --port-- unknowns $b_{r, i}$,
\begin{equation}
\label{eq:SCRBEM_eq2}
    \mathcal{O}_\Omega ({u}_{m,0}+\sum_{j=1}^{N^p_{r}} b_{r, j} \left({\varphi}_{r,0,j} + {\psi}_{r,j}\right), \theta_{r,i}) 
    + 
    \mathcal{O}_\Gamma (\sum_{j=1}^{N^p_{r}} b_{r, j}  {\psi}_{r,j}, \theta_{r,i}) 
    =
    \mathcal{F}_m(\theta_{r,i}) 
    \quad  \theta_{r,i} \in \Theta_{r}.
\end{equation}

In case the $\mathcal{O}$ are linear or linearized, we have
\begin{equation}
\label{eq:SCRBEM_eq3}
\sum_{j=1}^{N^p_{r}}
\left(
    \mathcal{L}_\Omega ({\varphi}_{r,0,j} + {\psi}_{r,j}, \theta_{r,i}) 
    + 
    \mathcal{L}_\Gamma ( {\psi}_{r,j}, \theta_{r,i}) \right) b_{r, j} 
    =
    \mathcal{F}_m(\theta_{r,i}) 
    - \mathcal{L}_\Omega ({u}_{m,0}, \theta_{r,i}) 
    \quad  \theta_{r,i} \in \Theta_{r}.
\end{equation}
Or in matrix form,
\begin{equation*}
\label{eq:SCRBEM_matrix2}
    D \mathbf{b}  = F.
\end{equation*}
Since only port DoFs remain in the final system, the above steps fully accomplish the static condensation. 

Anthony T Patera et al. \cite{huynh2013complex} proposed the SCRBEM to simulate large-scale acoustic devices governed by the Helmholtz equation. Henceforth, a posteriori error estimation was studied in \cite{huynh2013static}. Their research also includes numerical results for the heat transfer \cite{eftang2013port} and a linear elasticity problem \cite{eftang2014port}. In collaboration with Sylvain Vallaghee \cite{vallaghe2014static}, the technique is utilized for analyses of a conjugate heat exchanger model containing 1D bulk fluid and solid heat conduction. 

The methodology of DoFs division and static condensation for ROMs is simultaneously adopted by Pierre Kerfriden et al. for the analyses of nonlinear fracture mechanics \cite{kerfriden2013partitioned}. A more recent publication by Lin Mu et al. employed it for linear steady-state convection-diffusion equations with random diffusivity and velocity \cite{mu2019domain}. 

\subsubsection{Discontinuous Galerkin-based local ROM methods}
\label{subsub:DD_DG_ROM}
\begin{figure}[h]
    \centering
    \includegraphics[width=0.75\linewidth]{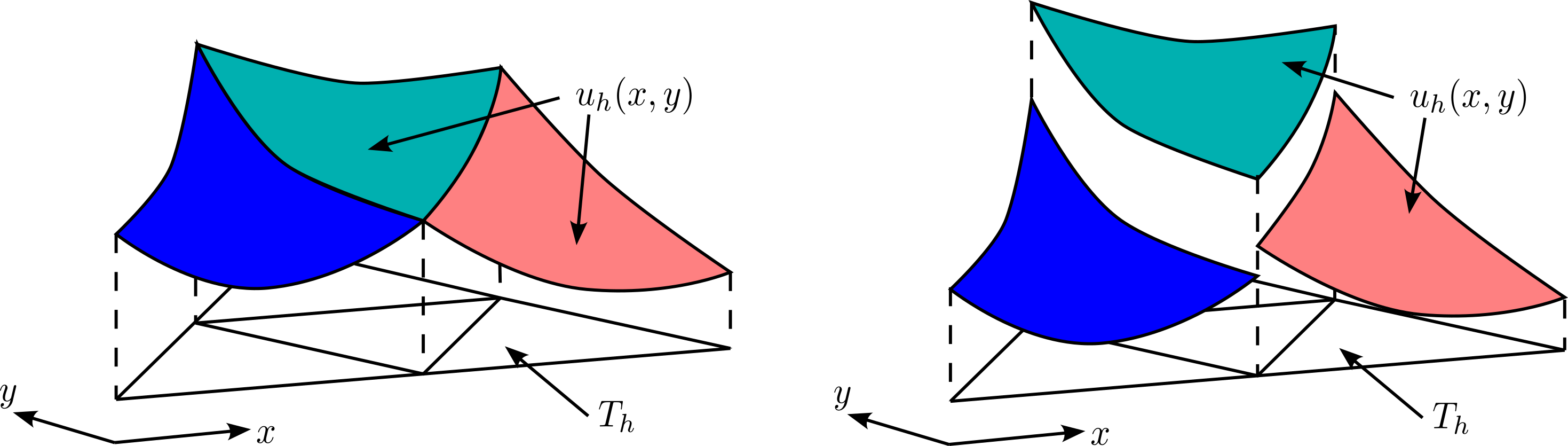}
    \caption{Schematic description of the basis functions employed for Continuous Galerkin (left) and Discontinuous Galerkin (right). Redrawn based on \cite{fidkowski2020output}.}
    \label{fig:CG_DG_scheme}
\end{figure}

\paragraph{Basic formulation of Discontinuous Galerkin}
The \emph{Discontinuous Galerkin} (DG) method was originally proposed as a high-order numerical technique intended for solving differential equations, combining features of the FE and finite volume frameworks \cite{lasaint1974finite, riviere2008discontinuous, cockburn2012discontinuous}. The standard Galerkin method (also known as Continuous Galerkin, CG) used in FE employs continuous basis functions across element interfaces. Unlike it,  DG utilizes discontinuous basis functions, as illustrated in Fig. \ref{fig:CG_DG_scheme}. To minimize variable jumps across interfaces, continuity and smoothness constraints are incorporated into the governing equations as penalization terms. The high-fidelity formulation also includes numerical flux terms defined on the interfaces. These also remain in the projection-based ROMs. Therefore, non-overlapping partitions can be \emph{glued} together directly in the online stage without Lagrange multipliers. 

The formulation depends on the exact equation to be approximated. Here, we use a Laplacian equation. Its weak formulation reads
\begin{equation}
\label{eq:laplacian_weak}
    - \int_\Omega \nabla u \cdot \nabla w + \int_{\partial \Omega} \nabla u w \cdot \mathbf{n} = 0
\end{equation}
where $u$ is the physical variable and $w$ is the test function, and $u, v \in H^1 (\Omega)$. To derive the DG formulation, we define on an interface $\Gamma_r = \Omega_m \cap \Omega_n$ the jump
\begin{equation*}
\label{eq:jump2}
\left. [w] \right |_{\Gamma_r} = \left. (w \cdot \mathbf{n}_m) \right |_{\Gamma_r} + \left. (w \cdot \mathbf{n}_n) \right|_{\Gamma_r},
\end{equation*}
and average operator,    
\begin{equation*}
\label{eq:jump_average1}
    \left. \{ w \} \right |_{\Gamma_r} = \frac{1}{2} \left( \left. w_m \right|_{\Gamma_r} + \left. w_n \right|_{\Gamma_r} \right), 
 \end{equation*}   
where $\mathbf{n}_m$ is the unitary normal vector pointing outwards $\Omega_m$ and $\mathbf{n}_m = -\mathbf{n}_n$.

As usual, we assume the entire domain is decomposed into $N_\Omega$ subdomains $\Omega = \cup_{m=1}^{N_\Omega} \Omega_m$, and the interface is collected into a set $\Gamma$. The integral over the whole model results in 
\begin{equation}
    - \sum_{\Omega_m \in \Omega} \int_{\Omega_m} \nabla u \cdot \nabla w + \sum_{\Gamma_r \in \Gamma } \int_{\Gamma_r} [\nabla u w] = \mathcal{F} (w),
\end{equation}
where $\mathcal{F} (w)$ involves global BCs. Remark that the term defined on the interface is called \emph{numerical flux} in the framework of the DG algorithm. Various schemes for the numerical flux have been discussed, and interested readers can turn to the references \cite{cockburn2004discontinuous, hesthaven2007nodal, riviere2008discontinuous} for more solid explanations. Here, we illustrate a typical strategy for dealing with it. 

By using the relationship 
\begin{equation*}
[uw] = [u] \{w\} + \{u\}[w],
\end{equation*}
we have
\begin{equation}
    - \sum_{\Omega_m \in \Omega} \int_{\Omega_m} \nabla u \cdot \nabla w + \sum_{\Gamma_r \in \Gamma } \int_{\Gamma_r} \left( [\nabla u ] \{w\} + \{\nabla u \} [w] \right) = \mathcal{F} (w).
\end{equation}
The solution $u$ and its gradient $\nabla u$ are continuous. Thus, 
\begin{equation*}
\label{eq:eqjumpgrad}
[u] = 0\; \text{and} \; [\nabla u] = 0,
\end{equation*}
to get
\begin{equation}
\label{eq:DGRBEM_1}
    - \sum_{\Omega_m \in \Omega} \int_{\Omega_m} \nabla u \cdot \nabla w + \sum_{\Gamma_r \in \Gamma } \int_{\Gamma_r} \{\nabla u \} [w] = \mathcal{F} (w).
\end{equation}

Two terms are added to the LHS of equation \eqref{eq:DGRBEM_1} to improve its numerical features: (i) $\frac{\gamma}{h} \int_{\Gamma} [ u ] [ w ]$ to penalize the jump of $u$, with $\gamma$ a factor and $h = \text{max} \left \{ \left | x - y \right |, \forall x, y \in \tau \right\} $ the largest diameter of the element in the mesh; (ii) A symmetric term $\varepsilon \int_{\Gamma} \{ \nabla w \} [u]$ where the constant $\varepsilon \in \{-1, 0, 1\}$, to created a formulation which is symmetric or not depending on the value of $\varepsilon$. The total results in 
\begin{equation}
\label{eq:DGRBEM}
    - \sum_{\Omega_m \in \Omega} \int_{\Omega_m} \nabla u \cdot \nabla w + \sum_{\Gamma_r \in \Gamma } \int_{\Gamma_r} \left( \{\nabla u \} [w] + \varepsilon \{ \nabla w \} [u] + \frac{\gamma}{h} \int_{\Gamma_r} [ u ] [ w ] \right) = \mathcal{F} (w).
\end{equation}
 
The above derivation can be extended to more complex problems, such as parabolic, Stokes, and Navier-Stokes equations. The major concern of the implementation is dealing with different terms, e.g., convection terms, and reformulating them as numerical flux regarding the jump and average conditions on interfaces. The detailed procedures and discussions are presented in \cite{hesthaven2007nodal, riviere2008discontinuous}. 

\paragraph{Discontinuous Galerkin RBEM}
The DG ROM approach integrating local RBs was proposed by Paola F. Antonietti et al. \cite{antonietti2016discontinuous} as a generalization and an improvement of both RBHM and RDF, called \emph{Discontinuous Galerkin Reduced Basis Element Method} (DGRBEM). It follows the philosophy of RBEM in domain decomposition. That is, it considers non-overlapping subdomains and a generic decomposition strategy (See Section \ref{subsubsec:Generic_decomposition}). The offline high-resolution DG simulations are performed in a generic partition, which is referred to as a local training strategy (discussed in Section \ref{subsubSec:Localized_training}). 

As RB functions are computed for each archetype, they are not matched on the interfaces satisfying the basic configuration of DG. Indeed, the local RBs, containing both internal and surface features, e.g., $\mathcal{V} \in H^1 (\Omega)$, can be directly adopted for equation \eqref{eq:DGRBEM}. Hence, every local solution can be approximated as a linear combination of local modes. We believe the derivation of a DGRBEM system and its matrix form is not complicated. Therefore, we don't clarify the details in this review.

Paola F. Antonietti et al. verified their algorithm with elliptic problems for multiple-division models, in which both physical properties and geometrical shapes are parameterized \cite{antonietti2016discontinuous}. Paolo Pacciarini et al. \cite{pacciarini2016spectral} extended the DGRBEM for parametrized Stokes problems, in which high-fidelity snapshots were collected by utilizing spectral element methods. A recent research from Seung Whan Chung et al. \cite{chung2024train} demonstrates the same methodology on three linear PDEs: the Poisson equation, the Stokes flow equation, and the advection-diffusion equation. 

We highly suggest the readers to an up-to-date investigation regarding \emph{Friedrichs' systems} \footnote{The Friedrichs' systems are regarded as a unified framework for various PDEs \cite{jensen2004discontinuous, romor2025friedrichs}: first and second order uniformly hyperbolic, elliptic and parabolic PDEs.} released by Francesco Romor et al. \cite{romor2025friedrichs}. The general technique is tested and validated by a series of PDEs, i.e., Maxwell equations, compressible linear elasticity, advection diﬀusion reaction, and incompressible Navier-Stokes equations. Note that the authors involved graph neural networks to efficiently infer the vanishing viscosity solution that is challenging for projection-based approaches.

Since the snapshots are collected from high-fidelity simulations at the subdomain level, the boundary conditions of local problems should be parameterized to approximate the global problems. Different types of boundaries can be imposed: non-homogeneous Neumann is adopted by Antonietti et al. \cite{antonietti2016discontinuous} and Paolo Pacciarini et al. \cite{pacciarini2016spectral}, while non-homogeneous Dirichlet is employed by Seung Whan Chung et al. \cite{chung2024train}.

\paragraph{Local discontinuous Galerkin ROM}
Note that the strategy for computing RBs in frames of DG is not unique. Andrea Ferrero et al. \cite{ferrero2018global} developed a \emph{local POD-DG} algorithm, whose formulations are the same as the DGRBEM except for domain decomposition and local RBs generation \cite{shah2020discontinuous, li2023simulation}. The approach considers individual decomposition of the whole problem and builds local POD RBs for each partition utilizing global solutions (as addressed in Section \ref{subsubsec:global_solution_local_RB}). Indeed, it can be regarded as a simplification of DGRBEM. As illustrated by Andrea Ferrero \cite{ferrero2018global}, the method results in good accuracy and significant acceleration for the 2-D Euler and RANS equations with the Spalart-Allmaras turbulence closure. 

We highlight that the DG-based systems for FOM and ROM are consistent except for basis vectors. This entitles DG to support the coupling of hybrid systems involving FOMs and ROMs, as presented by Sébastien Riffaud et al. \cite{riffaud2020reduced, riffaud2021dgdd}. They applied the method to model the isentropic vortex governed by 2D unsteady Euler equations, in which the computational domain is decomposed into non-overlapping individual subdomains. They also managed to couple ROMs and FOMs for unsteady transonic flows in a converging-diverging nozzle and over a NACA0012 airfoil in the presence of shocks \cite{riffaud2020reduced, riffaud2021dgdd}. They use a high-fidelity model to represent regions of complex nonlinear physical phenomena and a ROM to approximate the elsewhere efficiently, with which the computational complexity is significantly reduced and the accuracy is comparable to high-resolution solvers.

There is a series of articles published by Mario Ohlberger and Sven Kaulmann et al. proposed and investigated their \emph{Localized Reduced Basis Multiscale} (LRBMS) method. Although the algorithm is implemented for different equations and applications, the methodology and procedures are almost the same as Ferrero's method. The high-resolution DG is utilized to generate snapshots, and individual RBs are computed for each partition. It was originally introduced for general parametrized elliptic problems \cite{kaulmann2011new}. Then, it is improved to solve a two-phase flow in porous media \cite{Felix2012localized, kaulmann2015localized} and scalar parabolic evolution equations \cite{ohlberger2017true}. They illustrated a posteriori \emph{Error Estimations} of the LRBMS in detail for elliptic \cite{ohlberger2015error} and parabolic equations \cite{ohlberger2017true}. The online enrichment technique \footnote{This technique dynamically updates the reduced basis during the online phase. This improves capturing new parametric solutions that are deficiently represented by the current RB generated in the offline stage. Instead of relying solely on a precomputed reduced space, the ROM can identify regions where additional basis functions are needed and incorporate them adaptively. This helps achieve accuracy whilst remaining efficient.} is also included in their framework, and the details are discussed in \cite{ohlberger2015error,keil2023adaptive}.

\subsubsection{Partition of Unity method and ROM}
The \emph{Generalized Finite Element Method} (GFEM) is an extension of the classic FEM \cite{babuvska1994special}, which adopts a standard functional basis \footnote{It is typically a polynomial basis.} to approximate the common region of a problem and an enriched basis to capture specific local phenomena. The coupling is achieved following the  \emph{Partition of Unity Method} (PUM) \cite{melenk1996partition, babuvska1997partition}, which glues local bases defined in different regions to construct a global basis. 

The methodologies of GFEM and PUM can be incorporated within the framework of local ROM. They are utilized for problems consisting of overlapping divisions. The local RBs are computed using either local or global solutions. Obviously, they should not be continuous along interfaces. PUM is employed to weigh the local RBs and produce a continuous global space. Once the global RB is available, the reduced system \ref{eq:generic_ROM_matrix} can be easily constructed and solved. 

Let us suppose we have a series of local RBs, spanned by the components $\left\{ v_{m, i} \right\}_{m=1}^{N_{\Omega}}$. The $i^{\text{th}}$ vector of the PUM space can be defined as
\begin{equation*}
    v^{\text{PUM}}_i \coloneqq \text{span} \left\{ \bigcup_{m=1}^{N_\Omega} \psi_{m} v_{m,i}  \right\},
\end{equation*}
where $\psi_{m}$ is a PoU function defined for $\Omega_m$. Consequently, we can form a global reduced space 
\begin{equation*}
    \mathcal{V}^{\text{PUM}} \coloneqq \left\{ \sum_{i=1}^{N_{\text{RB}}}v^{\text{PUM}}_i \right\} \subset H^1_0(\Omega).
\end{equation*}

To finalize this description, we just need to define the $\psi_m$. We highlight that the formulation of the Partition of Unity (PoU) functions is not unique. Various examples are illustrated in \cite{melenk1996partition, babuvska1997partition}, as well as in the book \cite{bordas2023partition}. We will now give a basic description of the definition and construction of the PoU functions. Firstly, they satisfy
\begin{equation*}
\label{eq:PoU_function}
    0 \leq \left. \psi_{m} \right|_{\Omega_m} \leq 1, \quad \left. \psi_{m} \right|_{\Omega \setminus \Omega_m} = 0, \quad \sum_{m=1}^{N_\Omega} \psi_m (x) = 1, \quad \text{and} \quad \left\| \nabla \psi_{m} \right\|_{L^\infty (\Omega)} \leq C_m,
\end{equation*}
where $C_m$ is a constant upper bound of $\left\| \nabla \psi_{m} \right\|_{L^\infty (\Omega)}$. Here, the $\psi_{m}$ plays the role of the weights. 

\begin{figure}[h]
    \centering
    \includegraphics[width=0.5\linewidth]{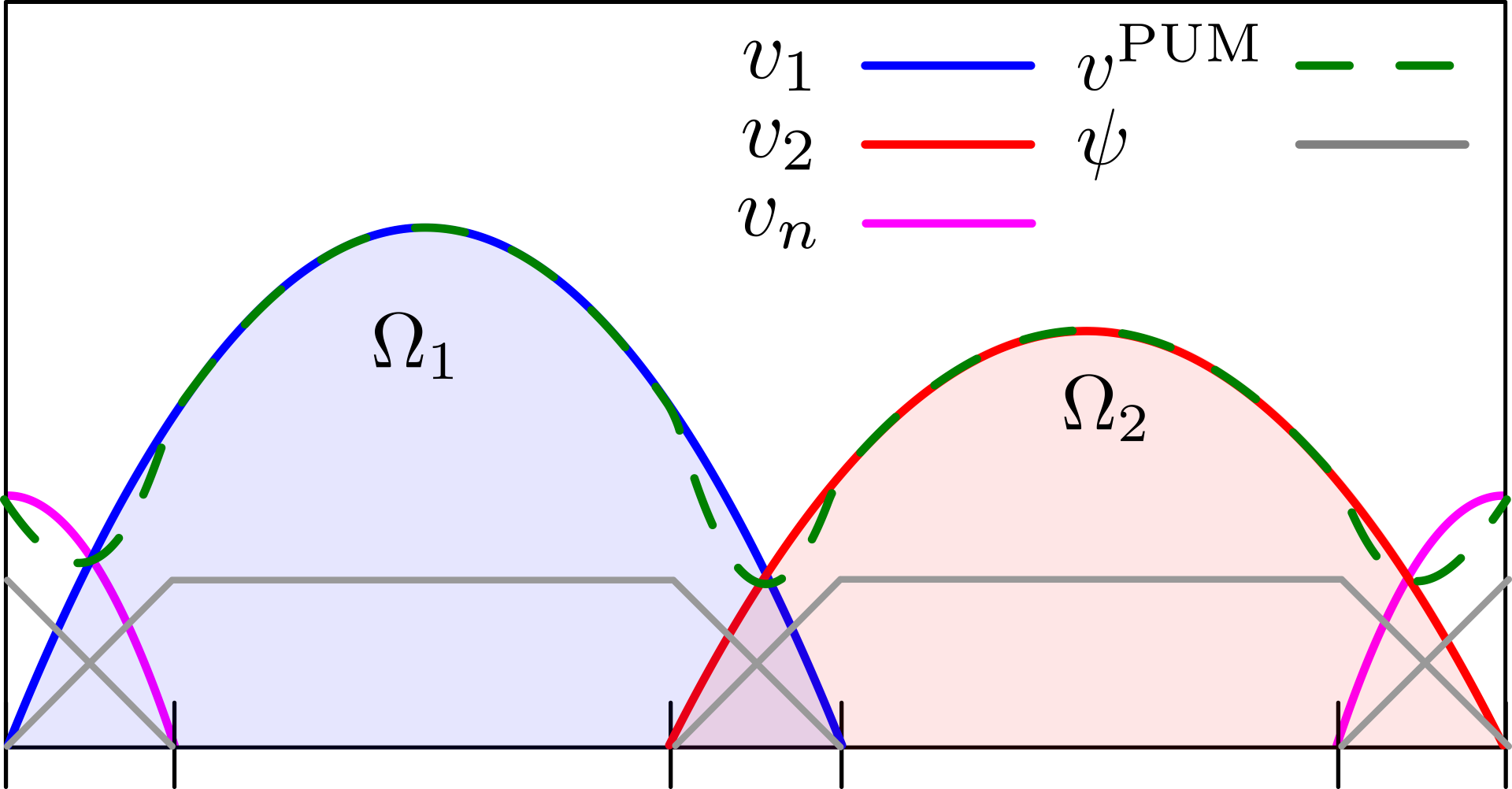}
    \caption{Sketch of PoU for a one-dimensional (1-D) case showing two subdomains $\Omega_1$ and $\Omega_2$. Their first basis functions $v_{1,1}$ and $v_{2,1}$ are marked in blue and red, respectively. The pink lines are vectors for segments out of the picture. Two PoU functions, $\phi_1$, $\phi_2$ are shown in gray lines. The vector resulting from the PUM procedure, $v^{\text{PUM}}_i$, is shown in green.}
    \label{fig:PUM_sketch}
\end{figure}

Here, we present a linear piecewise formulation for a 1-D problem as an example. The PoU functions and the resulting global basis are sketched in Fig. \ref{fig:PUM_sketch}. To satisfy the above requirements, the values of the PoU function in the local coordinate are given as
\begin{equation*}
    \psi_m (x) = \left\{ \begin{aligned}
    & \frac{x}{\delta_1}, & \quad & 0 \leq x \leq \delta_1 \\
    & 1, & \quad & \delta_1 < x < h-\delta_2, \\
    & -\frac{x}{\delta_2} + \frac{h}{\delta_2}, & \quad & h-\delta_2 \leq x \leq h,
    \end{aligned} \right.
\end{equation*}
where the local coordinate of a partition is assigned in the range of $x \in [0, h]$, and $0 \leq x \leq \delta_1$ and $h-\delta_2 \leq x \leq h$ are overlapping regions. Note that the same $\psi$ was defined here for all partitions $\Omega_m$. But this is not necessary, and different $\psi_m$ can be defined for each $\Omega_m$. Even more, multiple $\phi_{m,i}$ can be defined in each $\Omega_m$ and mode $v_{m,i}$ \cite{bordas2023partition}.  

The methodology of GFEM and its applications for MOR are presented in the references below. Ivo Babuška et al. \cite{babuska2011optimal} proposed the method to solve second-order elliptic PDEs. They also discussed the procedure to identify optimal local approximation spaces.
Additionally, they integrated into the method the \emph{oversampling} strategy to generate local solutions \cite{buhr2020localized, babuvska2020multiscale} and tested the approach for elliptic PDEs with random boundary conditions. Joan Baiges et al. adopted the strategy to couple bases for the construction of local ROMs and/or FOM-ROM systems to model incompressible Navier-Stokes flow around cylinders \cite{baiges2013domain} \footnote{The authors combined several local bases to build a global basis, and the methodology is similar to PUM, although they didn't explicitly mention the PUM in the text.}. They compared the performance of two spaces, i.e., localized global RBs (Section \ref{subsubSec:localized_global_RB}) and subdomain-level RBs (Section \ref{subsubsec:global_solution_local_RB}). Julia Schleu\ss{} et al. extended the procedure to analyze parabolic problems, proposing local space-time approximation spaces \cite{schleuss2022optimal}. 

A recent contribution from Kathrin Smetana et al. \cite{smetana2023localized} employed it to deal with nonlinear elliptic problems, considering \emph{generic decomposition} (Section \ref{subsubsec:Generic_decomposition}) with overlaps and \emph{oversampling training} (Section \ref{subsubSec:Localized_training}). The work also involved an adaptive basis enrichment algorithm when seeking global approximations in the online phase. They adopted an alternative formulation for the construction of the PoU functions 
\begin{equation*}
    \psi_{m,i} (x) = \frac{v_{m,i} (x)}{\sum_{n=1}^{N_\Omega} v_{n,i}(x)}, \quad x \in \Omega, m = 1, \cdots, N_\Omega,
\end{equation*}
where $v_{m, i}$ is a RB vector defined in $\Omega_m$ and $x$ denotes the global coordinate. Be aware that the subdomains are overlapping, thus $\sum_{n=1}^{N_\Omega} v_{n,i}(x) \neq v_{m,i} (x)$ in the shared regions.

\paragraph{Multiscale Finite Element Method}

A technique similar to the ideology of the previous approach is the so-called Multiscale Finite Element Method (MsFEM). The MsFEM employs a local basis to represent especially small-scale features inside of a large-scale solution \cite{hou1997multiscale}.  As clarified in \cite{efendiev2009multiscale}, MsFEM can be divided into two steps: (i) construct a subdomain-level basis that embeds essential multiscale information and (ii) couple the local bases to form a global formulation. 

From this description, one may conclude that  MsFEM and GFEM have similarities. Thus, it often happens that MsFEM and GFEM are discussed together in the literature \cite{babuska2011optimal, babuvska2020multiscale, buhr2020localized}. Resemblance appears in the following aspects: (i) the two are based on FEM and can be regarded as an extension of the method to enhance computational efficiency; (ii) the local basis is capable of capturing specific variations in subregions; (iii) the basis functions are continuous across the interfaces. Thus, one can resort to the general weak formulation to solve the system. 

In the early stage of MsFEM, researchers using this approach did not regard this technique as a ROM method \cite{hou1997multiscale}. It was later that it started to be considered as part of the local ROM approach \cite{efendiev2013generalized, chung2014generalized}. Note that the PUM is not incorporated in the MsFEM. The strategy for constructing the continuous basis will be briefly introduced.

Here, we adopt a linear elliptic equation for explanation and refer the readers for other examples, insight, and further discussion to  \cite{efendiev2009multiscale}. The problem adopted is
\begin{equation}
    \nabla \cdot (k(x) \nabla u ) = f, \quad \text{ in } \Omega,
\end{equation}
where $k(x)$ is the conductivity varying over multiple scales.

The domain is now spatially discretized into two scales, high-fidelity elements $e$ and subdomains $\Omega_m$ \footnote{The elements and blocks are not overlapping. We highlight that they are called fine/micro elements $e$ and coarse/macro blocks $\Omega_m$ in the frame of MsFEM.}. The aim is to compute basis functions in $\Omega_m$ and represent global solutions via the RB of $\Omega_m$ instead of the basis of $e$. The standard FEM basis vectors are denoted as $\phi^e_i$. As usual, $\phi_{m,i}$ is a vector for RB in $\Omega_m$. We compute the $\phi_{m,i}$ as the solutions of 
\begin{equation*}
    \nabla \cdot (k(x) \nabla \phi_{m,i}) = 0, \quad \text{ in }\Omega_m,
\end{equation*}
subjected to BCs
\begin{equation*}
    \phi_{m,i} = \phi^e_i, \quad \text{ on } \partial \Omega_m.
\end{equation*}

The condition $\phi_{m,i} = \phi^e_i$ induces a maximum limit in the dimension of the RB since the dimensions of the FEM basis are generally very reduced. In the terminology utilized in the literature for this method, one refers to the RB as \emph{multi-scale} basis\cite{efendiev2009multiscale}. 

The solution in each block is approximated with $N^{\text{RB}}_m$ vectors as $u_m = \sum_{i=1}^{N^{\text{RB}}_m} a_{m,i} \phi_{m,i}$. Also, $\mathcal{V}_m\equiv\mathcal{W}_m\equiv \text{span}\{\phi_{m,i}\}$. We recall that the BCs of $\phi_{m,i}$ are $\phi^e_i$. Since the $\phi^e_i$ are continuous as required by standard FEM, $\phi_{m,i}$ are conforming on the interfaces when assembling the entire domain. More precise, the global $\phi_i = \cup_m \phi_{i,m}$ are continuous in every interface. Indeed, $\mathcal{V} \in H^1_0 (\Omega) $

The MsFEM is to find $u \in \mathcal{V}_m$ such that 
\begin{equation}
    \sum_m \int_{\Omega_m} k(x) \nabla u \cdot \nabla \phi_{m,i} = \int_\Omega f \phi_{m,i}, \quad \forall \phi_{i} \in \mathcal{W}_m.
\end{equation}

The MsFEM was developed by Thomas Y. Hou et al., along with contributions from Yalchin Efendiev. It was applied to several problems, e.g., to linear and nonlinear elliptic problems, \cite{efendiev2004Multiscale, efendiev2007multiscale, efendiev2009multiscale}. 
There are alternative procedures to construct the RB. The oversampling technique (see Section \ref{subsubSec:Localized_training}) was initially utilized for solving local problems. Then, the framework was extended to develop the Generalized Multiscale Finite Element Methods (GMsFEM) for elliptic PDEs \cite{efendiev2013generalized, chung2014generalized, efendiev2014generalized}. This reduces the cost of constructing the RB when considering a complex input space. Patrick Henning et al. proposed an alternative method to generate the multiscale space $\mathcal{V}_m$, namely the Localized Orthogonal Decomposition (LOD) \cite{henning2014localizedSIAM, henning2014localizedESAIM}. Their approach was applied to elliptic PDEs with inhomogeneous Dirichlet and Neumann boundary conditions. The RB functions can be constructed adaptivity using error indicators to obtain larger accuracy for nonlinear problems \cite{chung2016adaptive}.

\subsubsection{Optimization-based algorithms}
\label{subsubsec:Optimization_based}
We now address the method based on a functional optimization strategy. Contrary to Section \ref{subsubsec:RBEM}, here we do not utilize the Lagrange multipliers procedure to deduce additional equations to add to the system. In the methods discussed below, an optimization problem is properly resolved.

\paragraph{Optimization-based Standard Galerkin}

The Optimization-based Domain Decomposition technique is accomplished by equalizing variables on both sides of the interfaces. At the same time, the governing equations and boundary conditions are satisfied. This methodology is intended for non-overlapping domains.

We follow the standard procedure to obtain $\Omega =\cup_{m=1}^{N_\Omega} \Omega_m$, $\Gamma = \cup_{r=1}^{N_\Gamma} \Gamma_r$, and the interfaces between $\Omega_m$ and $\Omega_n$ is noted as $\Gamma_{[mn]} = \left\{ \left. \Gamma_r \right| \Gamma_r = \partial \Omega_m \cap \partial \Omega_n, \Gamma_r \subset \partial \Omega_m \right\}$. For this methodology, we do not establish special restrictions on the local basis of the trial and test spaces. We refer the readers to any of the procedures described above for the derivations of the RBs spanning those spaces.  

We define a general target function $\mathcal{J} |_{\Gamma_{[mn]}}$ to be minimized to obtain a global solution, 
\begin{equation}
\label{eq:optimization_function}
    \mathcal{J} (u_m, u_n; g) \coloneqq \frac{1}{2} \int_{\Gamma_{[mn]}} \left | u_m - u_n \right|^2 + \frac{\gamma}{2} \int_{\Gamma_{[mn]}} \left | g \right|^2.
\end{equation}
Here, the first term denotes a jump of $u$ across the interface. The second term is known as the \emph{regularization} term. It arises as a way to impose the homogeneity of a secondary property in the interface and to ensure \emph{Well-Posedness}. Clearly, the parameter $\gamma > 0$ controls the relative importance of the terms of equation \eqref{eq:optimization_function}. We will give more details regarding $g$ below.

The optimization-based approach aims to find a solution that minimizes the target function. It should also satisfy the local governing PDEs equation \eqref{eq:generic_ROM} and pertinent BCs. We can also express the problem to fulfill the residual of control equations $\mathcal{R}_m (u_m, w_m) \coloneqq \left(\dot{u}_m, w_m \right) + \mathcal{O} \left ( u_m, w_m \right ) - \mathcal{F}(w_m), \quad (u_m, w_m) \in \left( \mathcal{V}_m, \mathcal{W}_m \right)$ and optimize the constraints. That is, to find $u$ such that 
\begin{equation}
\label{eq:optimization_based_ROM}
    \min \mathcal{J} (u_m;g) \quad \text{s.t.} 
    \left\{\begin{aligned}
        & \sum_{m=1}^{N_\Omega} \mathcal{R}_m (u, w) = 0, \quad \forall (u_m, w_m) \in \left( \mathcal{V}, \mathcal{W} \right), \\
        & u = g_D\  \text{ in }\ \Gamma_D,\ \text{ and }\ \partial u / \partial \mathbf{n} = \mathbf{g}_N\  \text{ in }\ \Gamma_N. \\
    \end{aligned} \right.
\end{equation}

Ivan Prusak et al. applied the methods for two incompressible Navier-Stokes benchmarks: the stationary backward-facing step and lid-driven cavity flow \cite{prusak2023application, prusak2023optimisation, prusak2024optimisation}. The authors utilized eq \eqref{eq:optimization_function} with a regularization term $\mathbf{g}$ that accounts for viscous stress tensor on the surface, following the ideas of Gunzburger \cite{gunzburger1999optimization}. It is defined as $\mathbf{g} = \left. \left( \nu \frac{\partial \mathbf{u}_m}{\partial \mathbf{n}_n} - p_m \mathbf{n}_m \right) \right|_{\Gamma_{[mn]}}$. Global high-fidelity solutions were collected to generate separated local RBs for different non-overlapping partitions. Both FOM and ROM problems were solved by a \emph{gradient-based optimization algorithm}. 
The exact optimality system is reformulated in terms of a Lagrangian functional, which is detailed in the aforementioned references, as well as in \cite{gunzburger1999optimization}. 

Tommaso Taddei et al. followed a similar strategy to analyze flow dynamics (Navier-Stokes equations) in blood vessel shape systems that are assembled by two archetype components \cite{taddei2024non}. They suggested penalizing the discontinuity of both velocity and pressure fields. They utilized a modified regularization term. This results in the following optimization function
\begin{equation}
    \mathcal{J} \coloneqq \frac{1}{2} \int_{\Gamma_{[mn]}} \left | \mathbf{u}_m - \mathbf{u}_n \right|^2 + \frac{1}{2} \int_{\Gamma_{[mn]}} \left | p_m - p_n \right|^2 + \frac{\gamma}{2} \int_{\Gamma_{[mn]}} \left( \left | \nabla_{\Gamma_{[mn]}} \mathbf{g} \right|^2 + \left | \mathbf{g} \right|^2 + \left | h \right|^2 \right),
\end{equation}
where $\mathbf{g}$ also denotes interface viscous stress, $\left. h \right|_{\Gamma_{[mn]}} = \left. \frac{\partial p}{\partial \mathbf{n}} \right|_{\Gamma_{[mn]}}$, and $\nabla_{\Gamma_{[mn]}} \mathbf{g}$ denotes the gradient of $\mathbf{g}$ in the tangential direction of ${\Gamma_{[mn]}}$.

To generate RBs for interface control, they proposed a \emph{pairwise-training} approach. They performed high-fidelity simulations for systems containing two partitions. Since each archetype might have several inlets and outlets, they construct various small systems to represent all possible connections between the two archetypes. The two-subdomain systems were resolved considering random Dirichlet and Neumann BCs. Also, they were parameterized with a range of Reynolds numbers and shapes. The internal RBs are constructed considering localized training (See Section \ref{subsubSec:Localized_training}). Adaptive enrichment procedure \footnote{This approach continuously evaluates the error of ROM and enriches the reduced basis when necessary. New FOM simulations are carried out for the worst predicted parametric data points, and new modes are added to improve approximation accuracy.} was also integrated into the framework. 

Finally, we summarize here a method that can be classified to this typology, but that is a combination of several of the methodologies described in this document. Giulia Sambataro et al. developed an innovative scheme named \emph{One-Shot Overlapping Schwarz} (OS2) method \cite{iollo2023one, sambataro2022component}. They reformulated the standard iterative Schwarz procedure (see Section \ref{subsec:iterative}) to be an optimization strategy that minimizes the lump sum of interface jumps \footnote{Although the authors use the name \emph{Schwarz} that is regarded as an iterative approach, we categorize their method into the group after analyzing their methodology.}. The subdomains are represented through \emph{port} and \emph{bubble} spaces (see Section \ref{subsubsec:SCRBEM}). Those are defined to represent interiors and internal boundaries. Then, the \emph{static condensation} procedure is utilized to eliminate the interior DoF (see Section \ref{sec:sct}). Therefore, the global approximation is reconstructed by a reduced system involving only interface DoFs.

\paragraph{Least squares Petrov-Galerkin}
\label{para:Least_squares_Petrov-Galerkin}

The second optimization-based procedure we discuss is based on the \emph{Least Squares Galerkin} (LSG) method\cite{quarteroni2009numerical, quarteroni2015reduced, carlberg2017galerkin}. A LSG system is formulated as an \emph{Optimal Control} problem, aiming to minimize the norm of the residual in a least-squares sense. Consequently, the ROM for a domain decomposition problem can be written as a system with constraints at the interfaces. 

It must be emphasized that the governing equations and enforced constraints are the inverse of those defined in the Optimization-based Standard Galerkin. The possibility of this swap often occurs in constrained minimization problems. Therefore, with the same conditions as the previous section, the optimality system is written as 
\begin{equation}
\label{eq:least_square_galerkin}
    \min \frac{1}{2} \sum_{m=1}^{N_\Omega} \left \| \mathcal{R}_m \left( u^\Omega_m \cup u^\Gamma_m, w_m \right) \right\|^2_2,
\end{equation}
subject to the restrictions,
\begin{equation*}
    \sum_{\Gamma \in \partial \Omega_m \cap \partial \Omega_n} \int_{\Gamma} \left | u_m^\Gamma - u_n^\Gamma \right|  = 0,
\end{equation*}
where $u^\Omega_m$ and $u_m^\Gamma$ denote interior and local surface values, respectively. The constraints are utilized to impose equality at interfaces -- or overlapping regions if they exist -- of two neighbouring subdomains $\Omega_m$ and $\Omega_n$. For a detailed strategy for solving the above minimization problem, we refer interested readers to the original publication.

Chi Hoang et al. have utilized the technique for parameterized Laplacian and Burgers' equations in a non-overlapping partitioned computational domain \cite{hoang2021domain}. This study has also incorporated several topics thoroughly. Notably, it discusses the utilization of strong or weak equality constraints on the interfaces. The performance of different strategies to build RBs was also investigated. Specifically, the behavior of RBs constructed with separated interior/boundary or full-subdomains, see Section \ref{sec:RBintintb}. The benefits of hyper-reduction are also revealed in their analysis. The authors adopted two training procedures to generate snapshots, one that requires full-system solutions (Section \ref{subsubsec:global_solution_local_RB}) and a second that only demands data from a single subdomain (see Section \ref{subsubSec:localized_global_RB}). 

Note that this methodology is prone to be integrated into a framework of Neural Networks to solve the reduced system, as investigated in \cite{diaz2023nonlinear, diaz2024fast}.

\subsubsection{Multiphysics problems}
\label{subsubsec:Monolithic_FSI}
Multiphysics models are typically governed by two or more PDEs, which might necessitate partitioning to account for several subregions. ROM is extensively applied to simulate these phenomena. However, the coupling of several physics is different from a single continuum. Considering the distinct characteristics, a separate description of the topic is presented below.

\paragraph{Fluid-structure interaction}
Among multiphysics, FSI is especially relevant in many practical applications. There are two widely used approaches for solving high fidelity FSI: fully-coupled/monolithic methods and partitioned/iterative methods \cite{keyes2007domain}. More details can be seen from \cite{quarteroni2014reduced, rozza2022advanced}. We will present the former here, and the latter will be described in Section \ref{subsec:iterative}.

According to the observed references, \emph{Arbitrary Lagrangian-Eulerian} (ALE) method is the most frequently applied for solving FSI-FOMs \cite{nonino2020application}. The formulation incorporates the fluid mesh displacement as an additional variable, enabling the fluid domain to follow the deformation of structures. Subsequently, solutions from the two regions are collected as snapshots, including displacement and velocity for both fields, as well as fluid pressure. Readers can refer to \cite{rozza2022advanced, nonino2020application} for more detailed descriptions of FSI-ROM.

Before applying Galerkin projection to construct the intrusive ROM, RB should be computed. Although the FSI problems consist of two physics, one can compute a global reduced basis for the entire model. For example, given the displacement fields for solid $s$ and fluid $f$ as $d^s = [d^s_1, \dots, d^s_{N^s_h}]^\mathbf{T}$ and $d^f= [d^f_1, \dots, d^f_{N^f_h}]^\mathbf{T}$, each global snapshot can be defined as $S_d = [d^f_1, \dots, d^f_{N^f_h}, d^s_1, \dots, d^s_{N^s_h}]^\mathbf{T}$. The same strategy can be used for other variables. Thus, due to the continuity of high-fidelity solutions, the resulting basis functions ensure continuity at the interface. The RBs can construct a monolithic reduced system. Erwan Liberge et al. applied the method for a transient flow around an oscillating cylinder \cite{liberge2010reduced}. Francesco Ballarin and Gianluigi Rozza extended it for more general parameterization cases for a coupling of Navier-Stokes flow and linear elastic structures \cite{ballarin2016podgalerkin}. 

Note that the above investigations employing a global RB for two systems are not strictly regarded as local ROM approaches. We begin discussing studies that consider two separate local RBs for fluid and solid, respectively. An example is presented to clarify the coupling. 

Given two sets of governing equations for a simplified FSI formulated by ALE as 
\begin{equation}
\label{eq:FSI}
\begin{aligned}
    \left(\dot{x}^f, w^f\right) + \mathcal{O}^f_{\Omega^f} \left ( x^f, w^f \right ) + \mathcal{O}^f_{\Gamma} \left ( x^f, w^f \right ) & = \mathcal{F}^f (w^f), \\
    \left(\dot{x}^s, w^s\right) + \mathcal{O}^s_{\Omega^s} \left ( x^s, w^s \right ) + \mathcal{O}^s_{\Gamma} \left ( x^s, w^s \right ) & = \mathcal{F}^s (w^s),
\end{aligned}
\end{equation}
where $x$ contains all field values (geometry displacement $d$, velocity $t$, fluid pressure $p$), $w$ is the test function, the superscripts $s$ and $f$ denote solid and fluid, respectively. 

Assume a constraint $ \left. x^f \right|_\Gamma = \left. x^s \right|_\Gamma$, and the interface terms in \eqref{eq:FSI} can be replaced by $\mathcal{O}^f_{\Gamma} \left ( x^f, w^f \right ) = \mathcal{O}^f_{\Gamma} \left ( x^s, w^f \right )$ and $\mathcal{O}^s_{\Gamma} \left ( x^f, w^s \right ) = \mathcal{O}^s_{\Gamma} \left ( x^f, w^s \right )$, that denotes
\begin{equation}
\label{eq:FSI_coupled}
\begin{aligned}
    \left(\dot{x}^f, w^f\right) + \mathcal{O}^f_{\Omega^f} \left ( x^f, w^f \right ) + \mathcal{O}^f_{\Gamma} \left ( \underline{x^s}, w^f \right ) & = \mathcal{F}^f (w^f), \\
    \left(\dot{x}^s, w^s\right) + \mathcal{O}^s_{\Omega^s} \left ( x^s, w^s \right ) + \mathcal{O}^s_{\Gamma} \left ( \underline{x^f}, w^s \right ) & = \mathcal{F}^s (w^s),
\end{aligned}
\end{equation}
where the structural variable $x^s$ acts as a term in the fluid ROM, and vice versa for $x^f$. 

The fluid and solid ROMs are assembled into a coupled system. The same treatments can be applied to other constraints if necessary. Readers can turn to the articles for applications in aeroelasticity (i.e., Euler equations and linear/nonlinear elastic structure) \cite{lieu2006reduced-order, colciago2014comparisons} and the Cardiovascular System (i.e., incompressible Navier-Stokes flow and linear elasticity) \cite{colciago2014comparisons}. 

We highlight several recent contributions from Monica Nonino et al. \cite{nonino2020application, nonino2021monolithic, nonino2023reduced}. The authors successfully addressed a crucial issue in transport-dominated nonlinear FSI-ROMs, namely, the slow decay of the Kolmogorov $n$-width. That was achieved by a method called \emph{transport maps}. They solved FSI systems adopting an ALE formulation like \eqref{eq:FSI_coupled}. Moreover, interface constraints can be weakly enforced by Lagrange multipliers (similar to equation \eqref{eq:RBEM_interface}). Assume the geometry displacements $d$ are equal on an interface $d^f = d^s$. Then the additional equations are expressed as, 
\begin{equation*}
    \int_{\Gamma} \zeta_{i} \left( d^f - d^s \right) = 0 \quad \forall \zeta_{i} \in \mathbb{L},
\end{equation*}
where $\mathbb{L}$ is a functional space for the Lagrange multipliers $\zeta_{i}$.

Be aware that \emph{bifurcating phenomena} may occur in FSI problems, which indicate non-unique stable solutions as parameters are varied \cite{hale2012dynamics}. Moaad Khamlich et al. studied a special bifurcating behavior (known as the Coandă effect) in a FSI problem comprised of incompressible Navier-Stokes flow and a linear elastic solid \cite{khamlich2022model, rozza2022advanced}. In their simulations, the FOM coupling is achieved through the ALE framework, and interface constraints are enforced using Lagrange multipliers (similar to those employed by Monica Nonino et al.). More descriptions about model order reduction for bifurcation problems can be found in \cite{pitton2017computational, pichi2019reduced, pichi2020reduced, pichi2022driving}. 

\paragraph{Mesh motion in FSI using interpolation}
The aforementioned methods (i.e., the ALE formulations) require solving a subsystem of equations to estimate the motion of the fluid mesh in cases of unsteady FSI phenomena. According to the references observed, geometric reduction of mesh morphing can also be achieved using interpolation. Control points in the structural interfaces can be used as input to predict the fluid grids at each time step. 

Davide Forti and Gianluigi Rozza \cite{forti2014efficient} adopted RBF to create the interpolants for FSI-ROMs. To further reduce the cost of interpolation, they proposed a greedy procedure to minimize the number of control points. The basic idea is to select surface points that can be used to accurately approximate the eigenmodes of the solid. They tested the approach in two cases: (i) adaptive selection in a commercial aircraft; (ii) external viscous fluid flow past a deformable rectangular obstacle.

Apart from RBF, the \emph{Inverse Distance Weighting} (IDW) formulation can also be adopted. Jeroen A.S. Witteveen et al. exploited IDW for high-fidelity FSI with a 2D airfoil and a 3D wing \cite{witteveen2009explicit}. They demonstrated that IDW is comparable in accuracy to RBF, while reducing computational costs by around 100 times \cite{witteveen2009explicit}. 

Alessandro D'Amario incorporated IDW and ROM in his thesis to simulate typical FSI phenomena: the fluid flow around an aircraft wing and a ship hull \cite{damario2014Reducedorder}. However, the IDW approach is not practical for complex geometries consisting of a large number of points. To overcome this limitation, D'Amario and Francesco Ballarin et al. proposed the Selective IDW (SIDW) method, also called Reduced IDW \cite{damario2014Reducedorder, ballarin2019podselective}. This algorithm optimally extracts a subset of control points based on a geometric criterion, thereby reducing the cost of interpolation. Furthermore, additional constraints can be enforced to enhance SIDW, resulting in the so-called enriched SIDW (ESIDW). The authors exploited the ESIDW technique to analyze three benchmark cases, including the structural deformation of a wing, and fluid mesh motion around a wing and a rotating hull.

Lastly, we highly recommend the comprehensive investigation presented in Davide Forti's thesis \cite{forti2012Comparison}, in which RBF, IDW, and FFD are employed to compute mesh motion when constructing FSI-ROMs. Additionally, the author proposed two novel strategies for improving RBF and FFD. The adaptive section of RBF points is integrated as a hyper-reduction. Domain decomposition is incorporated with FFD to enable different control point refinement for different locations. The performance of these techniques is compared with respect to FSI phenomena governed by incompressible Navier-Stokes equations and linear elasticity: (i) an external laminar fluid flow past an elastic obstacle and (ii) an internal flow in a cylindrical deformable vessel. In conclusion, he emphasized the significant potential of RBF in terms of computational cost and accuracy.

\paragraph{Embedded boundary method for FSI}
Geometric parameterization analysis typically involves applying transformations to reference domains or remeshing. This section introduces a so-called \emph{Embedded Boundary Method} (EBM) as an alternative approach. It incorporates the concept of domain decomposition to parameterize complex shapes and FSI problems. An important step of the framework is "cutting" finite elements from a mesh, which can be considered a type of spatial partition and is therefore presented here.

The introduction consists of three parts: decomposition, FOM formulation, and snapshots and ROM. See more in \cite{wang2011Algorithms, main2018Shifteda, main2018Shiftedb}.

\begin{figure}[h]
    \centering
    \includegraphics[width=\linewidth]{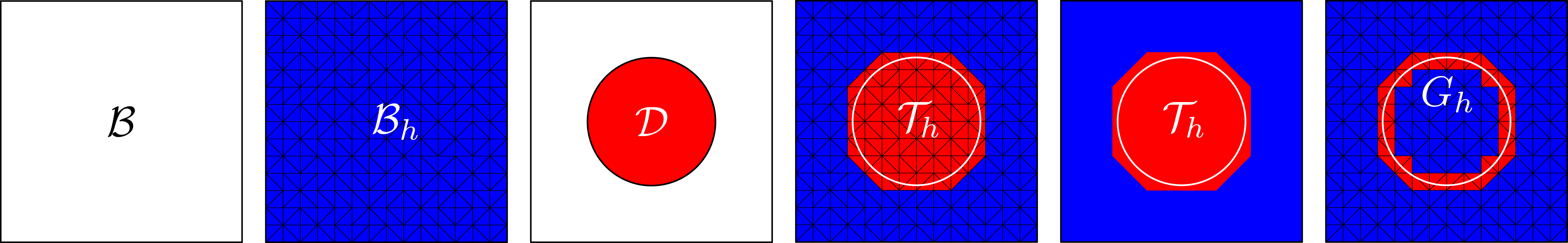}
    \caption{Sketch and definition of the embedded boundary method. From left to right: the background geometry $\mathcal{B}$, the background mesh $\mathcal{B}_h$, the solid region $\mathcal{D}$, the \emph{embedded} mesh $\mathcal{T}_h$, the \emph{surrogate} geometry $\mathcal{D}_{\mathcal{T}}$, and the "cut" elements $G_h$. The location and shape of $\mathcal{D}$ can be parameterized. Figures redrawn based on \cite{karatzas2020Projectionbased}.}
    \label{fig:embedded_boundary_method_concept}
\end{figure}

The decomposition and notation of the EBM can be summarized as follows (see also in Fig. \ref{fig:embedded_boundary_method_concept}.): (i) discretize a domain $\mathcal{B}$ that embeds the solid $\mathcal{D}$; (ii) generate a background mesh $\mathcal{B}_h$; (iii) "cut" the approximate \emph{surrogate} geometry $\mathcal{D}_{\mathcal{T}}$ according to the embedded mesh $\mathcal{T}_h$ and $\mathcal{D}$. Note that the subdomain $\mathcal{D}$ can be deformed and relocated in cases of parameterization. 

After decomposition, the EBM weak formulation of the concerned problem is imposed in $\mathcal{D}_\mathcal{T}$ or $\mathcal{T}_h$ as 
\begin{equation*}
    \left(\dot{u}_d, w \right) + \mathcal{O} \left ( u_d, w \right ) = \mathcal{F}(w) \quad \text{ in } \mathcal{D}_\mathcal{T} (\text{or } \mathcal{T}_h),
\end{equation*}
where $\mathcal{O}$ is the EBM operator consisting of penalty terms that are not extended here.

The values in the rest region, $u_b \text{ in } \mathcal{B} \setminus \mathcal{D}_\mathcal{T}$, can be obtained from two strategies, namely
\begin{equation*}
(1)\ u_b = 0 \text{ in } \mathcal{B} \setminus \mathcal{D}_\mathcal{T} \quad \text{and} \quad  (2) \left \{
\begin{aligned}
 &-\Delta u_b = 0, &  \text{ in } & \mathcal{B} \setminus \mathcal{D}_\mathcal{T}, \\
 &u_b = \left. u \right|_{\partial \mathcal{D}_\mathcal{T}}, & \text{ on } & \partial \mathcal{D}_\mathcal{T}, \\
 & u_b = 0, & \text{ in } & \partial \mathcal{B}, \\
\end{aligned}
\right.
\end{equation*}
where (1) is called a \emph{natural smooth extension} and (2) is a \emph{harmonic extension}.

If the shape of $\mathcal{D}$ is parameterized, a set of FOM simulations is performed for each sample. It is obvious that all calculations are performed in the same background mesh with different $\mathcal{D}_\mathcal{T}$. Then, the solutions in $\mathcal{B}$ are collected as snapshots, indeed, $\mathbf{S} = [u_1, \dots, u_{N_\mu}]$, in which $u_i = u_d(\mu_i) \cup u_b(\mu_i)$. The RB can be computed via POD for the entire computational domain. Consequently, Galerkin/least-squares projection formulations can be applied to construct a ROM.

The high-fidelity EBM has been used for analyzing fluid dynamics and FSI \cite{wang2011Algorithms}. An early work of EBM-based ROM was carried out by Maciej Balajewicz and Charbel Farhat \cite{balajewicz2014Reduction}. The authors formulated ROM employing the least-squares projection and tested it with three FSI problems: (i) one-dimensional FSI based on the viscous Burgers equation; (ii) a compressible viscous flow around a heaving rigid cylindrical body; (iii) a viscous flow inside a square cavity with a rotating ellipsoidal body.

Efthymios N. Karatzas et al. have published a sequence of studies employing two innovative embedded frameworks, e.g., the \emph{Shifted Boundary} (SB) method \cite{main2018Shifteda, main2018Shiftedb} and \emph{cut element FEM} (CutFEM) \cite{pasquariello2016Cutcell}, and the POD-Galerkin strategy ROM. Note that the novelty of the SB method lies in its weak formulation, whereas CutFEM involves additional integral treatment for cutting elements intersecting with the interface (e.g., $G_h$ in Fig. \ref{fig:embedded_boundary_method_concept}).

The SB-based ROMs are utilized to analyze various phenomena, including a heat conduction \cite{karatzas2018Reduced}, a Stokes flow around an embedded circular cylinder and more complex obstacles \cite{karatzas2019Reduceda}, a steady incompressible Navier-Stokes flow past an embedded rectangular domain \cite{karatzas2020ReducedOrder}, a fourth-order nonlinear geometrical PDE, e.g., Cahn-Hilliard system \cite{karatzas2021Reduced}, and shallow water hyperbolic equations \cite{zeng2022Embedded}. The CutFEM incorporating with ROM for the Darcy flow pressure model and the Steady Stokes problem can be found in \cite{karatzas2020Projectionbased}.

Although EBM-ROMs are not widely adopted for solving domain decomposition problems, we believe they have potential in this area due to their inherent characteristics of spatial partition.

\paragraph{Complex multiphysics}
The coupled formulations \eqref{eq:FSI_coupled} can be extended to multiphysics problems. Andrea Manzoni et al. constructed ROMs for cardiac electrophysiology and mechanics \cite{manzoni2018reduced, cicci2023projection}. We recommend \cite{quarteroni2014reduced} for a general description and \cite{benner2020applications} for applications about coupled multiphysics ROMs. 

Lagrange multipliers are also applied to couple multiphysics phenomena. Alberto Corigliano et al. employed the technique to analyze nonlinear elastic-plastic structural dynamics and the electro-mechanical system \cite{corigliano2014combined, corigliano2015model, dossi2015combined}. 

Note that the investigations mentioned above require a two-way exchange between the models. However, in some scenarios, the interaction is one-way. For example, given a two-equation problem, one has a fixed interface condition (maybe homogeneous), and the other subproblem needs interface quantities as inputs. Elena Zappon et al. presented both steady and unsteady one-way coupling for second-order elliptic problems \cite{zappon2023staggered, zappon2023efficient}. The authors succeeded in coupling non-conformal interfaces mesh at both the FOM and ROM stages, using a so-called \emph{INTERNODES} method \cite{deparis2016internodes, gervasio2018analysis}.

\subsection{Iterative}
\label{subsec:iterative}
We will now continue to address the second family of coupling methodologies: \emph{iterative} procedures. The methods have been widely used in high-fidelity numerical analysis. They can couple sub-problems based on any discretization techniques (e.g., finite differences, finite element, finite volume, and spectral element methods) \cite{quarteroni2009numerical}. The framework can also be adapted to an assembly of local ROMs regarding either overlapping or non-overlapping partitions.

The structure of the methods can be divided into two stages: (i) in each iteration, each sub-problem is solved independently, imposing data transferred from adjacent partitions as BCs or, alternatively, known global boundaries; (ii) iterations are performed considering the whole system until discontinuities or residuals satisfy a predefined threshold. 

The implementations of the most widely used iterative algorithms, the Schwarz method and its variations, are described in sequence. Besides, compared to a single continuum, the coupling of multiphysics phenomena is more complex and requires more interface constraints. A separate section on the topic is presented at the end.

\subsubsection{Schwarz method}
\label{subsubsec:Schwarz_method}
The \emph{Schwarz method} (also referred to as \emph{Schwarz alternating method}) is a classic iterative domain decomposition method formulated by H. A. Schwarz around 1870. It can be applied to spatial partitions with overlaps. Denoting $\Omega = \cup_{m=1}^{N_\Omega} \Omega_m$, for two overlapping parts $\Omega_m$ and $\Omega_n$, their interfaces are defined as $\Gamma_{[mn]} = \partial \Omega_m \cap \Omega_n$ and $\Gamma_{[n m]} = \partial \Omega_n \cap \Omega_m$. The two adjoining sub-problems and their boundary and initial conditions can be formulated as
\begin{equation}
\label{eq:schwarz_problem}
\begin{split}
& \left \{
\begin{aligned}
    \left(\dot{u}^{(t)}_m, w_m\right) + \mathcal{O} \left ( u^{(t)}_m, w_m \right ) &= \mathcal{F}(w_m) &\quad &\text{ in } \Omega_m, \\
    {u}^{(t)}_m &= {u}^{(t-1)}_n & \quad & \text { on } \Gamma_{[mn]}, \\
    {u}^{(t)}_m &= g_D & \quad & \text { on } \partial \Omega_m \cap \Gamma_D, \\
    \frac{\partial {u}^{(t)}_m}{\partial \mathbf{n}} &= \mathbf{g}_N & \quad & \text { on } \partial \Omega_m \cap \Gamma_N, \\
\end{aligned}
\right. \\[1em]
& \left\{
\begin{aligned}
    \left(\dot{u}^{(t)}_n, w_n\right) + \mathcal{O} \left ( u^{(t)}_n, w_n \right ) &= \mathcal{F}(w_n) &\quad &\text{ in } \Omega_n, \\
    {u}^{(t)}_n &= \left\{ \begin{aligned} &{u}^{(t-1)}_m \text{ or } \\ &{u}^{(t)}_m \end{aligned} \right. & \quad & \text { on } \Gamma_{[mn]}, \\
    {u}^{(t)}_n &= g_D & \quad & \text { on } \partial \Omega_n \cap \Gamma_D, \\
    \frac{\partial {u}^{(t)}_n}{\partial \mathbf{n}} &= \mathbf{g}_N & \quad & \text { on } \partial \Omega_n \cap \Gamma_N, \\
\end{aligned}
\right.
\end{split}
\end{equation}
where the superscript $(t)$ indicates iteration steps, $\Gamma_D$ and $\Gamma_N$ are global boundaries. Depending on the iteration BCs in sub-problem \eqref{eq:schwarz_problem}, the method is called \emph{multiplicative Schwarz} when choosing $u^{(t)}_m$, or \emph{additive Schwarz} when the selection is $u^{(t-1)}_m$.

It has been proven that the Schwarz method always converges for second-order elliptic equations, with a rate that increases as the measure of the overlapping region increases \cite{quarteroni2009numerical}. Additionally, relaxation can be utilized to improve the convergence, which modifies the interface condition to be $u_m^{(t)} = \beta u_n^{(t-1)} + \left( 1-\beta \right) u_m^{(t-1)}$, with $\beta$ as the relaxation factor.

Each sub-problem is approximated separately. Thus, local ROMs can be directly constructed using Galerkin projection, resulting in a formulation expressed in terms of equation \eqref{eq:generic_ROM}.

The early stage of the iterative scheme is addressed as a hybrid technique to couple FOMs and ROMs. Marcelo Buffoni et al. \cite{buffoni2009iterative} applied the multiplicative Schwarz method to combine a ROM and a FOM for two non-linear problems, including a Laplace equation with non-linear boundary conditions and compressible Euler equations. Davide Cinquegrana et al. \cite{cinquegrana2011hybrid} analyzed the aerodynamic flow field around a 2-D airfoil geometry governed by the Reynolds Average Navier-Stokes equations with the multiplicative Schwarz procedure. The region surrounding the airfoil is resolved by a high-fidelity model, and the rest is approximated with ROM. Junpeng Song and Hongxing Rui \cite{song2024reduced} utilized it to analyze the convection-dominated diffusion equation for a geometry divided with overlaps. Localized global RBs were computed to approximate sub-problems.

In general, the Schwarz algorithm is utilized for divisions with overlaps. Nevertheless, a recent publication from Irina Tezaur et al. \cite{wentland2024role} demonstrates the possibility of obtaining a stable and accurate coupled model for non-overlapping cases. Their approach is evaluated by three challenging 2-D nonlinear hyperbolic problems: the shallow water equations, Burgers' equation, and the compressible Euler equations. They employed a hyper-reduction technique, which samples a few elements on the interfaces to exchange local BCs. Thus, they also assessed the effect of the boundary sampling approach. 

A novelty in the frame of Schwarz iteration can be seen from the study of Marco Discacciati et al \cite{discacciati2024overlapping}. They incorporate the overlapping multiplicative Schwarz method for simulations of parametric elliptic problems via proper generalized decomposition (PGD). Instead of reduction by POD, the PGD is employed to approximate high-fidelity solutions.

The basic ideology of PGD is that the solution can be approximated as a sum of separable functions, indeed,
\begin{equation*}
    u(\mathbf{x};\mu_1,\mu_2) = \sum_{i=1}^{N^{\text{PGD}}} X_{i}(\mathbf{x})M_{i}(\mu_1)K_{i}(\mu_2),
\end{equation*}
where $\mu_1$ and $\mu_2$ are two different parameter sets. $X$, $M$ and $K$ are modes representing $\mathbf{x}$, $\mu_1$ and $\mu_2$, respectively. 

The implementation of PGD in \cite{discacciati2024overlapping} is briefly illustrated hereafter. The authors divided the local solution $u_m(\mu)$ into two parts
\begin{equation*}
    u_m(\mathbf{x};\mu,\mu_{\text{BC}}) = v_{m,0}(\mathbf{x};\mu) + \tilde{v}_{m}(\mathbf{x};\mu,\mu_{\text{BC}}),
\end{equation*}
where $v_{m,0}(\mu)$ denotes the local solution satisfying homogeneous interface conditions, global boundary conditions, and source terms. $\tilde{v}_{m}(\mathbf{x};\mu,\mu_{\text{BC}})$ only represents the interface contribution without any sources and boundary conditions, and the interfaces are parametrized with $\mu_{\text{BC}}$. In the article, this represents the artificial boundary conditions imposed during the localized training (i.e., Dirichlet BCs at the interfaces) to generate snapshots.

The interface terms are further regarded as a summation of two items, namely,
\begin{equation*}
    \tilde{v}_{m}(\mathbf{x};\mu,\mu_{\text{BC}}) = v_{m}(\mathbf{x};\mu,\mu_{\text{BC}}) + \sum_{\Gamma_r \in \mathcal{I}_m} \sum_{q=1}^{N_{\Gamma_r}}\Lambda_{q}^{\Gamma_r} (\mu_{\text{BC}}) \eta _{q}^{\Gamma_r} (\mathbf{x}),
\end{equation*}
where $\mathcal{I}_m$ is the set of interfaces associated with the $\Omega_m$. $\eta _{q}^{\Gamma_r}$ and $\Lambda_{q}^{\Gamma_r}$ are suitable basis functions and corresponding coefficients that are capable of approximating interface quantities.

The two parts $v_{m,0}$ and $v_{m}$ are approximated separately following the strategy of PGD. This yields two PGD expansions
\begin{equation*}
\begin{split}
    v_{m,0}(\mathbf{x};\mu) &= \sum_{i=1}^{N_{m,0}^{\text{PGD}}} X_{m,0,i}(\mathbf{x})M_{m,0,i}(\mu), \\
    v_{m}(\mathbf{x};\mu,\mu_{\text{BC}}) &= \sum_{i=1}^{N_r^{\text{PGD}}} X_{r,i}(\mathbf{x})M_{r,i}(\mu)K_{r,i}(\mu_{\text{BC}}).
\end{split}
\end{equation*}

Finally, the complete approximation is written by 
\begin{equation}
    u_m(\mathbf{x};\mu,\mu_{\text{BC}}) = \sum_{i=1}^{N_{m,0}^{\text{PGD}}} X_{m,0,i}(\mathbf{x})M_{m,0,i}(\mu) + \sum_{i=1}^{N_m^{\text{PGD}}} X_{m,i}(\mathbf{x})M_{m,i}(\mu)K_{m,i}(\mu_{\text{BC}}) + \sum_{\Gamma_r \in \mathcal{I}_m} \sum_{q=1}^{N_{\Gamma_r}}\Lambda_{q}^{\Gamma_r} (\mu_{\text{BC}}) \eta _{q}^{\Gamma_r} (\mathbf{x}).
\end{equation}

Note that the algorithms for computing PGD modes and boundary basis functions fall outside the scope of this review. Readers may consult the original publications. We also recommend a book \cite{croft2015proper} for a comprehensive explanation of PGD.

Once the PGD expansions are available, they can be integrated into the weak formulation to solve local problems. And each local solution is iteratively updated employing the Schwarz approach. It is easy to check that, following the previous steps, the unknowns are condensed to the interfaces, i.e., $\Lambda^{\Gamma_r}_q (\mu_{\text{BC}})$. Thus, it can be regarded as a static condensation procedure.

The method is tested in three problems considering overlapping subdomains: a two-domain parametric Poisson equation, a two-domain convection-diffusion equation, and a multi-domain thermal problem with discontinuous conductivity. Note that, for the multi-partition cases, the snapshots are collected through localized training (see Section \ref{subsubSec:Localized_training}) with arbitrary Dirichlet boundary conditions at the interfaces. 

\subsubsection{Variations of the Schwarz method}
\label{subsubsec:Variations_Schwarz_method}
The classic Schwarz method accomplishes \emph{Dirichlet-Dirichlet} coupling between adjacent subdomains but can be extended to Neumann and Robin conditions in a non-overlapping framework. For \emph{Dirichlet-Neumann}, the interface condition is formulated as,
\begin{equation}
\label{eq:D-N_conditions}
\left.
\begin{aligned}
    {u}^{(t)}_m &= {u}^{(t-1)}_n \\
    \frac{\partial {u}^{(t)}_n}{\partial \mathbf{n}}  &= \frac{\partial {u}^{(t)}_m}{\partial \mathbf{n}} \\
\end{aligned}
\right\}
\quad \text { on } \Gamma_{[mn]}.
\end{equation}
Whilst the \emph{Robin-Robin} iteration is written as 
\begin{equation}
\label{eq:R-R_conditions}
\left.
\begin{aligned}
    \frac{\partial {u}^{(t)}_m}{\partial \mathbf{n}} + \gamma_m {u}^{(t)}_m &= \frac{\partial {u}^{(t-1)}_n}{\partial \mathbf{n}} + \gamma_m {u}^{(t-1)}_n  \\
    \frac{\partial {u}^{(t)}_n}{\partial \mathbf{n}} + \gamma_n {u}^{(t)}_n  &= \frac{\partial {u}^{(t)}_m}{\partial \mathbf{n}} + \gamma_n {u}^{(t)}_m \\
\end{aligned}
\right\}
\quad \text { on } \Gamma_{[mn]}. 
\end{equation}
$\gamma_m$ and $\gamma_n$ are non-negative parameters to control the importance of each term. The choice of these parameters will affect the convergence. The values should be selected from a sensitivity study \cite{gander2012best}. As indicated in \cite{discacciati2023localized}, the Robin parameter $\gamma$ is prior fixed to ensure the well-posedness of the governing equation, which typically requires $\gamma > 0$.

We start the discussion with studies that employ the Dirichlet-Neumann method. Immanuel Maier and Bernard Haasdonk \cite{maier2011iterative, maier2014dirichlet} implemented it for an elliptic PDE (i.e., static heat equation with homogeneous Dirichlet boundary conditions). Separated RBs were computed for interfaces and inner regions. They managed to glue RBs belonging to two non-overlapping subdomains of distinct shapes. A posteriori error estimation concludes that the residuals of the ROM were comparable to the FOM. 

The framework was then developed for incompressible Navier-Stokes equations by Ricardo Reyes \cite{reyes2020element}. The authors adopted a non-overlapping decomposition. Additionally, they revealed that the performances of the Dirichlet-Neumann are comparable to the local DG ROM (see Section \ref{subsub:DD_DG_ROM}).

Irina Tezaur et al. \cite{barnett2022schwarz} further tested its possibility for FOM-FOM, FOM-ROM, and ROM-ROM coupling. In their study, snapshots were collected from a sequence of FOM-FOM coupled simulations. Individual RBs were computed for each subdomain. The performance of different coupling systems (FOM-ROM, ROM-ROM, etc) was evaluated on a 1-D nonlinear wave propagation problem. The authors also revealed that the performance of Dirichlet-Neumann and Dirichlet-Dirichlet is comparable \cite{tezaur2022schwarz}. We also recommend readers two presentations addressed by Irina Tezaur et al. \cite{mota2021schwarz, tezaur2022schwarz}, in which they intended to employ the Schwarz alternating scheme to achieve the coupling of multi-scale and multi-physics problems, e.g., FSI approximated by a FOMs-ROMs system. Furthermore, they also discussed the possibility of creating a general framework to couple FOMs with intrusive and/or non-intrusive ROMs for more complex and realistic applications.

Several up-to-date publications \cite{mota2023fundamentally, koliesnikova2023dirichlet} have demonstrated that the Dirichlet-Neumann iteration scheme is suitable for solid mechanics, e.g., contact problems in elastodynamics. 

A recent innovation by Elena Zappon et al. \cite{zappon2024reduced} is utilizing the Dirichlet-Neumann to couple problems discretized with non-conforming high-fidelity mesh for different subdomains. The interpolation between the non-conforming meshes at the interface is achieved by \emph{INTERNODES} \cite{deparis2016internodes, gervasio2018analysis} method. Furthermore, hyper-reduction based on a \emph{discrete empirical interpolation method} is applied \cite{chaturantabut2010nonlinear}, which allows interpolating the parametric Dirichlet-Neumann interface conditions with a small amount of \emph{magic points}. Their approach has been tested in several second-order elliptic/parabolic problems that consider two partitions.

Except for Dirichlet-Neumann, Niccolò Discacciati, and Jan S. Hesthaven integrated the Robin-Robin approach to analyze a series of parametrized nonlinear problems, including Advection-diffusion-reaction coupling, time-dependent diffusion equation in multiple domains, etc \cite{discacciati2023localized}. The authors also proposed a localized training procedure (see Section \ref{subsubSec:Localized_training}) by parametrizing Robin conditions. Snapshots and RBs are constrained in archetype components. 

The \emph{Neumann-Neumann} iterative algorithm differs slightly from the methods above. It requires solving additional correction problems. Due to the lack of applications for local ROM coupling (to the best of our knowledge), it is not illustrated here. Its formulation is presented in \cite{quarteroni2009numerical}.

\subsubsection{Multiphysics problems}
\label{subsubsec:Iterative_FSI}

For the studies mentioned above, a pair of conditions is exchanged among adjacent subdomains. However, for multiphysics problems, more constraints must be enforced at the interfaces.

Since FSI is the most typical multiphysics problem, we describe it first. For a problem combining incompressible Navier-Stokes flow and elasticity, three interface conditions are required to guarantee geometric consistency, velocity continuity, and the balance of normal forces \cite{ballarin2016podgalerkin}. They are defined on $\Gamma$ as 
\begin{equation*}
\begin{aligned}
    d^f - d^s &= 0, \\
    u^f - u^f &= 0, \\
    \sigma^f \cdot \mathbf{n}^f - \sigma^s \cdot \mathbf{n}^s &= 0,
\end{aligned}
\end{equation*}
where superscripts $f$ and $s$ indicate fluid and solid respectively, $d$ is the displacement of the geometry, $u$ is the velocity ($u^s = \frac{\partial d^s}{\partial t}$), $\mathbf{n}$ is the normal vector and $\sigma \cdot \mathbf{n}$ is the normal force on the interfaces. 

Note that although more variables are considered in multi-physics scenarios, iterative schemes still support the coupling. The stepping interface conditions are formulated as those in equations \eqref{eq:schwarz_problem}, \eqref{eq:D-N_conditions}, and \eqref{eq:R-R_conditions}.  This coupling procedure can be applied to both FOM and ROM levels. Some applications are presented in \cite{vierendeels2007implicit, aletti2017reduced, nonino2023projection}. 

For other multiphysics problems, different interface constraints should be enforced, but the iterative algorithm remains the same. In this context, we found two studies particularly interesting regarding the implementation of multiphysics: (i) linear thermo-poroelasticity \cite{ballarin2024projection-based}; (ii) nuclear reactor thermal-hydraulics (Navier-Stokes flow with energy balance) coupled with neutronics (scalar diffusion equation) \cite{riva2024multi-physics}.


\section{Data-driven Techniques}
\label{sec:Data_driven_Techniques}
The rapid development of \emph{Machine Learning} (ML) techniques and computational power has recently encouraged the research community to incorporate pure data-driven methods to construct ROMs. These approaches are fundamentally non-intrusive, as they operate directly on the high-fidelity solutions without requiring any modification of the underlying governing equations. A review published in 2021 \cite{heinlein2021combining} discussed the combination of ML and domain decomposition. Due to the growing interest in this topic and various recent developments, we present an overview of current advances in this area. Also, several relevant references reviewed in \cite{heinlein2021combining} are also included for completeness.

In the following sections, we will provide the basic preliminaries for constructing data-driven models. Then, the strategies for combining domain decomposition and ROM will be explained with respect to each technique. In our opinion, the classification of the non-intrusive approaches is not as straightforward as the intrusive ones. That is due to the variety of ML frameworks and the flexibility of incorporating different techniques to couple the local ROMs. Thus, considering their fundamental principles, we have created four categories to roughly classify them. Those are: (i) Schwarz-based iteration; (ii) interpolation algorithm; (iii) optimization-based technique; (iv) \emph{Physical Informed Neural Network} (PINN) based methods. 

\begin{figure}[h]
    \centering
    \includegraphics[width=0.8\linewidth]{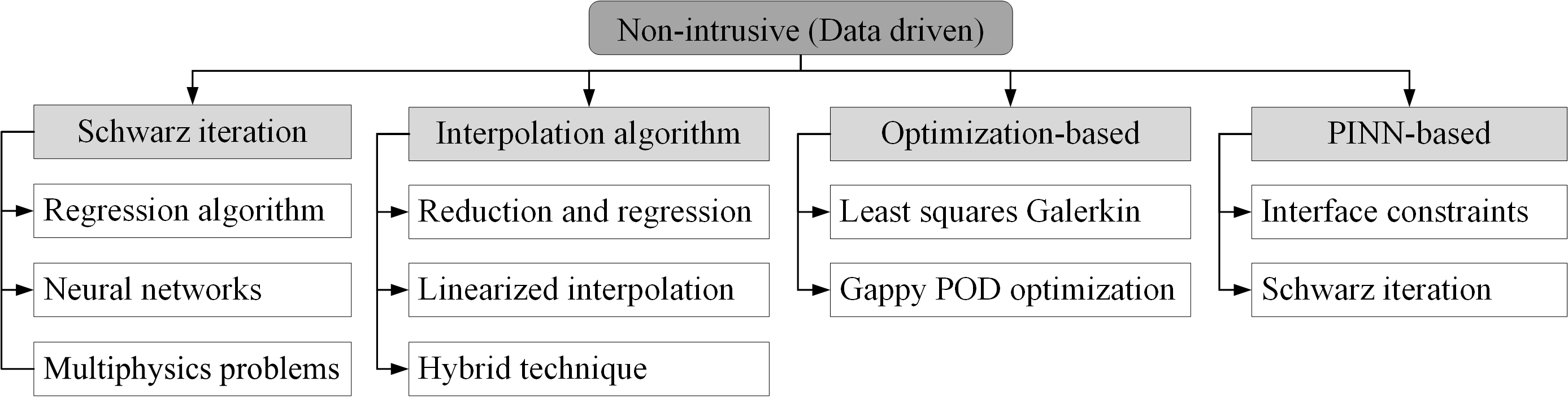}
    \caption{The classification of data-driven techniques. The abbreviation: \emph{Physical Informed Neural Network} (PINN).}
    \label{fig:data_driven_local_ROM}
\end{figure}

The classification of the data-driven techniques is illustrated in Fig. \ref{fig:data_driven_local_ROM}. 

It is worth noting that the domain decomposition and basis construction strategies illustrated in Section \ref{sec:dd_local} also support the narrations hereafter. Additionally, the coupling strategies of the aforementioned intrusive procedures may be incorporated into the non-intrusive systems.

\subsection{Preliminaries}
\label{subsec:data_driven_preliminaries}

The fundamental ideology of the data-driven strategy is building a surrogate model that correlates inputs and output results. The model is trained and validated on sampled parameters and corresponding results. Once available, it can make predictions of unknown parameters with a cost several orders of magnitude lower than FOM calculations. The surrogate model can be built using various techniques, including conventional regression algorithms and Neural Networks, especially the more advanced PINNs. They will be briefly explained in the paragraphs that follow.

\subsubsection{Regression algorithms}

Regression processes in the frame of machine learning are diverse. Least squares or linear regression can be utilized for simple scenarios. For more complex conditions, several advanced methods are adopted, such as RBF interpolation \cite{buhmann2000radial} and Gaussian Process Regression (GPR) \cite{williams2006gaussian}, as well as Artificial Neural Network (ANN).

Suppose $\mathbf{u} (\mu)$ is a set of quantities that need to be predicted, and $(\mu)$ is the parameter. The basic procedure for constructing a regression model $\mathcal{Z}$ consists of \emph{training} and \emph{testing} stages. $\mathcal{Z}$ is trained using observed parameters and corresponding solutions, i.e., $\mathbf{u} (\mu) \approx \mathcal{Z} (\mu)$, and then, for unknown parameters ${\mu}'$, the predictions is approximated as ${\mathbf{u}'} (\mu) \approx \mathcal{Z} ({\mu}')$. 

Moreover, dimensionality reduction techniques like POD (i.e., $\mathbf{u} \approx \sum_{i=1}^{N_\text{RB}} \alpha_i \mathbf{v}_i$) can be involved in the procedure. Instead of regressing whole snapshots, $\mathcal{Z}$ can also be used to predict POD coefficients $\boldsymbol{\alpha}$. Indeed, $\boldsymbol{\alpha}' (\mu) \approx \mathcal{Z} ({\mu}')$, and ${\mathbf{u}'} (\mu) \approx \sum_{i=1}^{N_\text{RB}} \boldsymbol{\alpha}' (\mu) \mathbf{v}_i$ \cite{rozza2022advanced, padula2024brief}. This approach is called \emph{POD interpolation} (PODI).

The dimensionality reduction can also be accomplished by applying neural networks, like the \emph{autoencoder} (AE) \cite{brunton2024data}. The AE can generate a vector to represent high-dimensional data. The vector is referred to as \emph{latent variables} (also known as a latent space). More explanations about AE can be seen in Section \ref{subsubsec:Neural_network}.

Note that, due to the variety of regression processes, their formulations cannot be explained in detail here. We suggest that readers consult the books \cite{buhmann2000radial, williams2006gaussian, bonaccorso2018machine, brunton2024data} for more detailed explanations.

\subsubsection{Neural network}
\label{subsubsec:Neural_network}
\begin{figure}[h]
    \centering
    \includegraphics[width=0.45\linewidth]{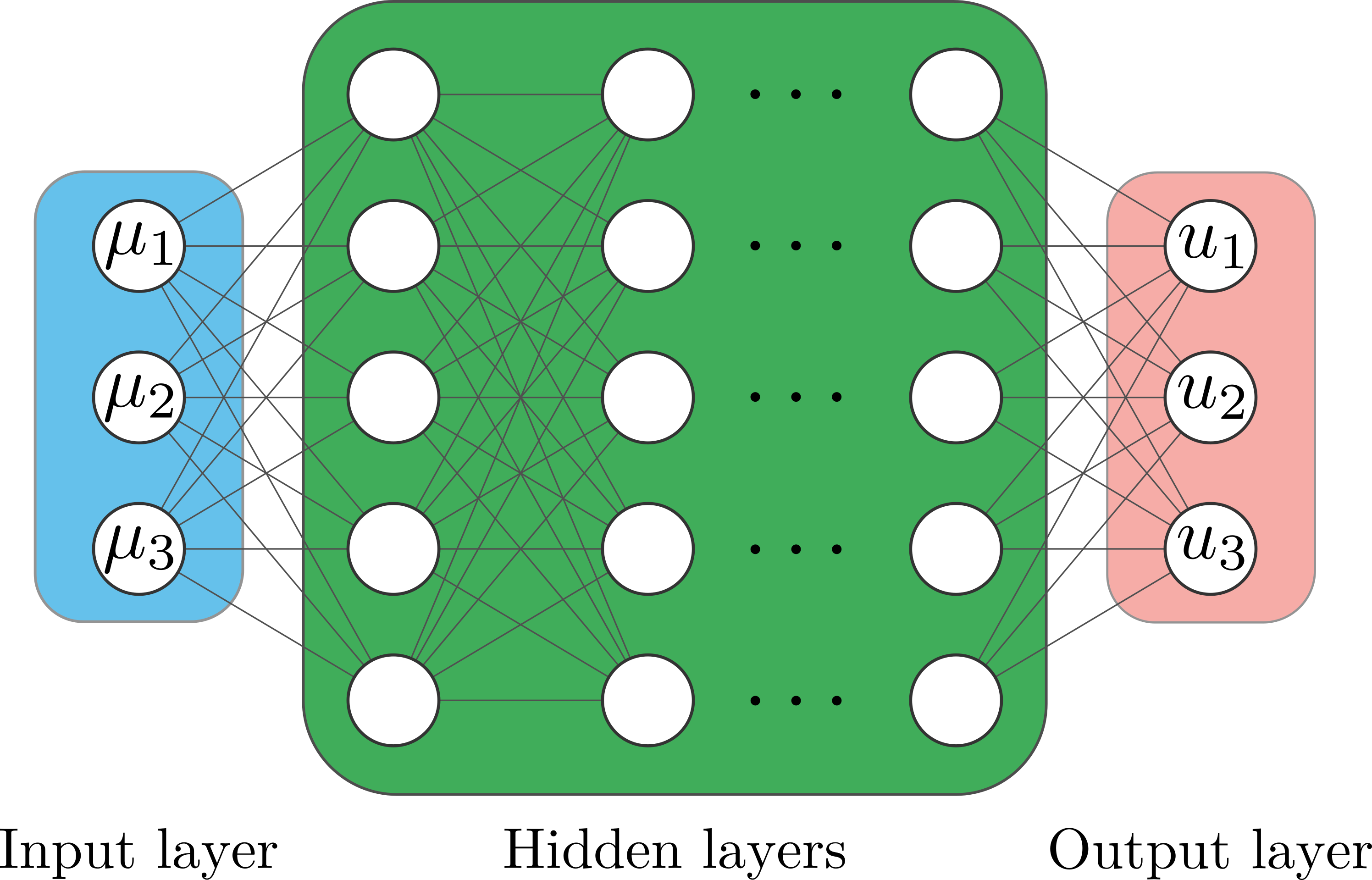}
    \caption{A sketch of neural network architecture.}
    \label{fig:Neural_networks}
\end{figure}

Recently, ANNs have been widely incorporated to construct regression frameworks for complex and nonlinear situations \cite{padula2024brief}. Considering its relevance to the following narration, we briefly discuss the basic idea of the ANN. Fig. \ref{fig:Neural_networks} illustrates a simple example of a network architecture. Parameters $\mu$ are regarded as inputs, namely \emph{input layer}. Target results are extracted from the \emph{output layer}. The quantities pass through several \emph{hidden layers} to build a surrogate for inputs and outputs.

The aim of a neural network in the training stage is to reduce the \emph{loss function}, indeed, the error between its predictions $\mathcal{Z} (\mu)$ and the actual target values $\mathbf{u} (\mu)$. For example, a standard loss function, \emph{Mean Squared Error }(MSE), is defined as
\begin{equation}
\label{eq:mse}
    \mathcal{E} = \frac{1}{N_{\mathbf{u}}} \lVert \mathbf{u} (\mu) - \mathcal{Z} (\mu) \rVert_2^2,
\end{equation}
where $N_{\mathbf{u}}$ denotes the amount of $\mathbf{u} (\mu)$.

The minization is achieved through the value transition process in all neurons. For each neuron, the input $x$ and output $y$ are related as
\begin{equation}
\label{eq:nn_in_out}
    y = \sigma(wx + b),
\end{equation}
where $\sigma$ is the \emph{activation function} that is selected before training, $w$ is the \emph{weight} and $b$ is the \emph{bias}. $w$ and $b$ are tuned to minimize the loss, i.e., the training stage. More in-depth discussions regarding the implementations are not presented here, but can be found in \cite{nielsen2015neural, aggarwal2018neural}.

\begin{figure}[h]
    \centering
    \includegraphics[width=0.35\linewidth]{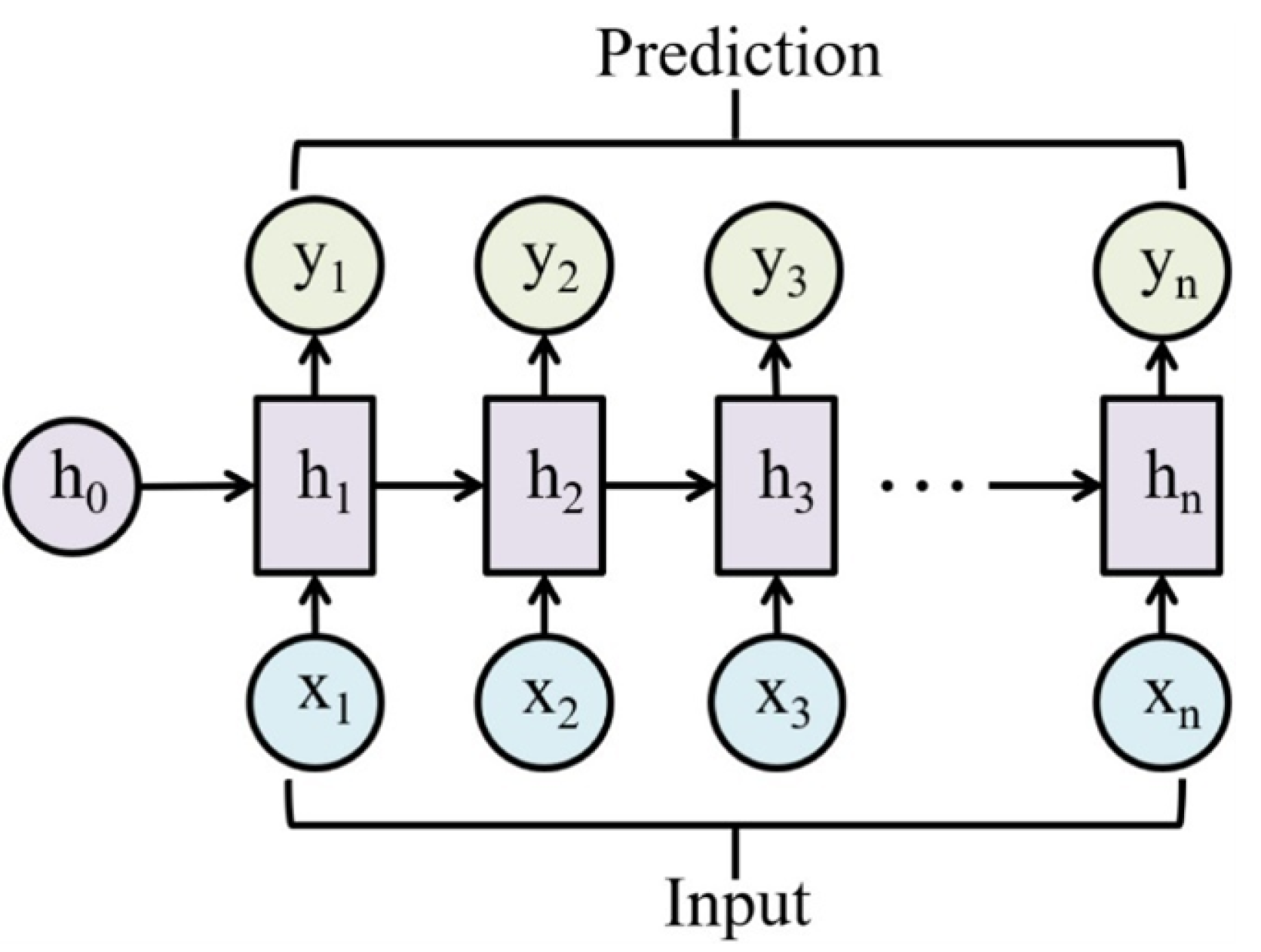}
    \hfill
    \includegraphics[width=0.5\linewidth]{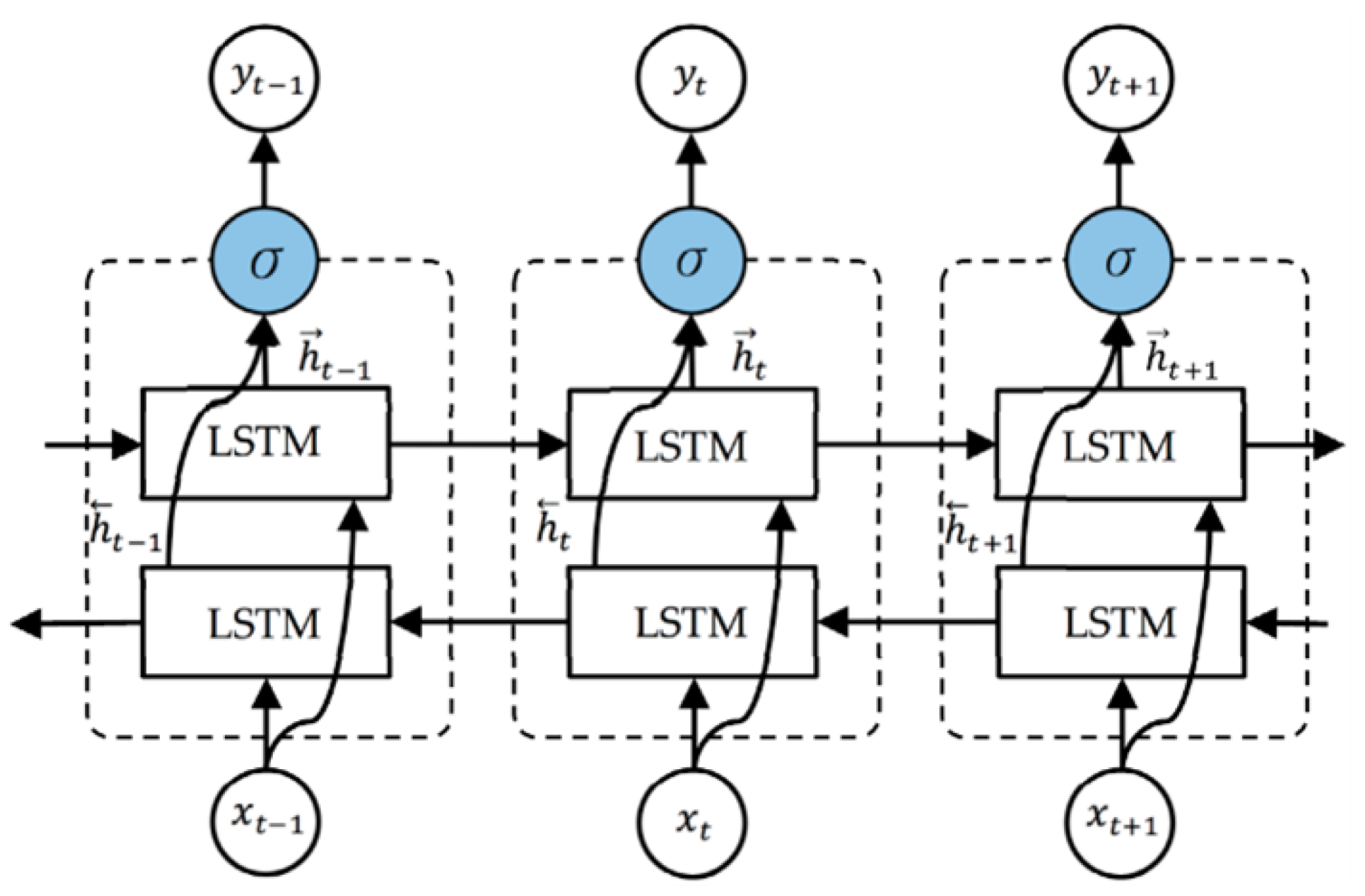}
    \caption{Architectures of recurrent neural networks for a sequential dataset. \textit{Left}: Standard RNN with initial input layer $h_0$ and middle hidden nodes $h_1, \dots, h_n$. Multiple inputs $x_1, \dots, x_n$ and corresponding outputs $y_1, \dots, y_n$. \textit{Right}: Long Short-Term Memory networks. An additional path of nodes for feedback from the next sequence. Figures are taken with permission from \cite{mienye2024recurrent}, copyright owned by MDPI.}
    \label{fig:RNN_LSTM}
\end{figure}

The neural network can be further classified considering the exchange between neurons. For the ANN plotted in Fig. \ref{fig:Neural_networks}, the inputs are only transferred forward to reach the final outputs. However, for sequential data and transient phenomena, it requires results at different temporal points. This can be achieved by the so-called \emph{Recurrent Neural Networks}, in which the output of nodes can feed back into the model \cite{brunton2024data, mienye2024recurrent}. It is also worth mentioning an advanced variant, \emph {Long Short-Term Memory}, which can transfer information both forwards and backwards. The two are displayed in Fig. \ref{fig:RNN_LSTM} for clarification.

\begin{figure}[h]
    \centering
    \includegraphics[width=0.75\linewidth]{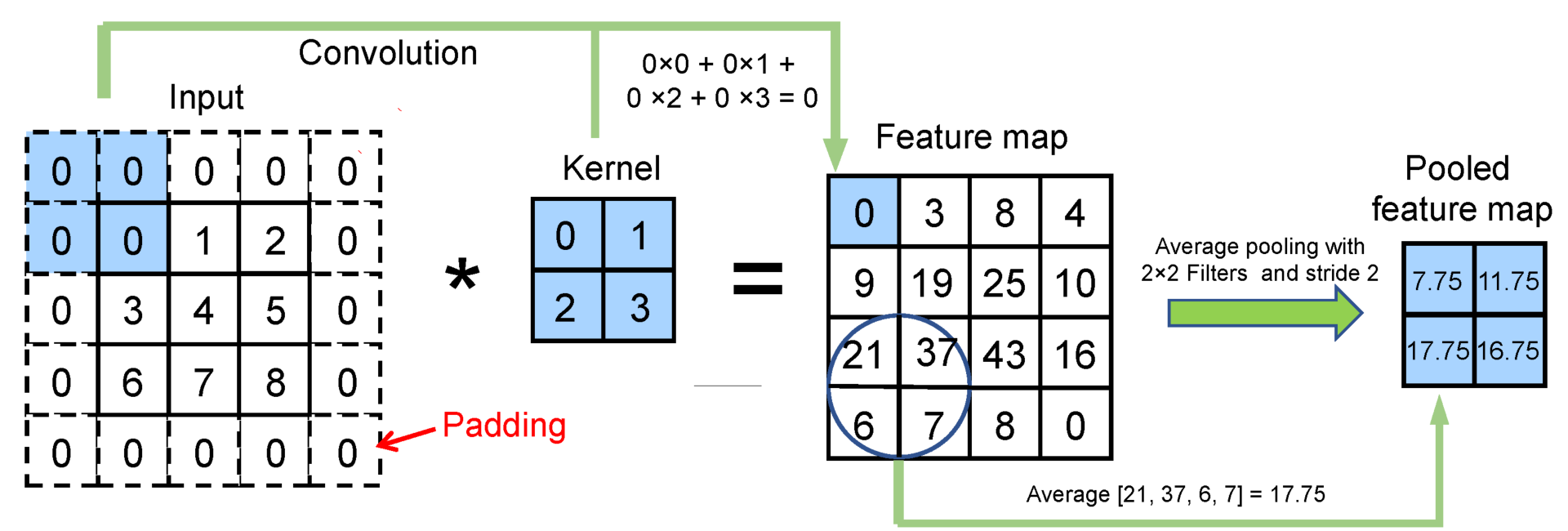}
    \caption{Example of convolution and pooling layers in the convolutional neural network. Figure taken with permission from \cite{li2022smart}, copyright owned by MDPI.}
    \label{fig:CNN}
\end{figure}

Another widely used architecture is called \emph{Convolutional Neural Network}, which consists of special additional layers to extract data in sub-regions. These denote the convolution and pooling process shown in Fig. \ref{fig:CNN}.

Other frameworks, such as \emph{Graph Neural Networks} \cite{pichi2024graph, romor2025friedrichs} and the neural operator \cite{chen2025tensor, diab2025temporal}, are gaining popularity nowadays. They are mentioned here for completeness and without further explanation. 

\begin{figure}[h]
    \centering
    \includegraphics[width=0.5\linewidth]{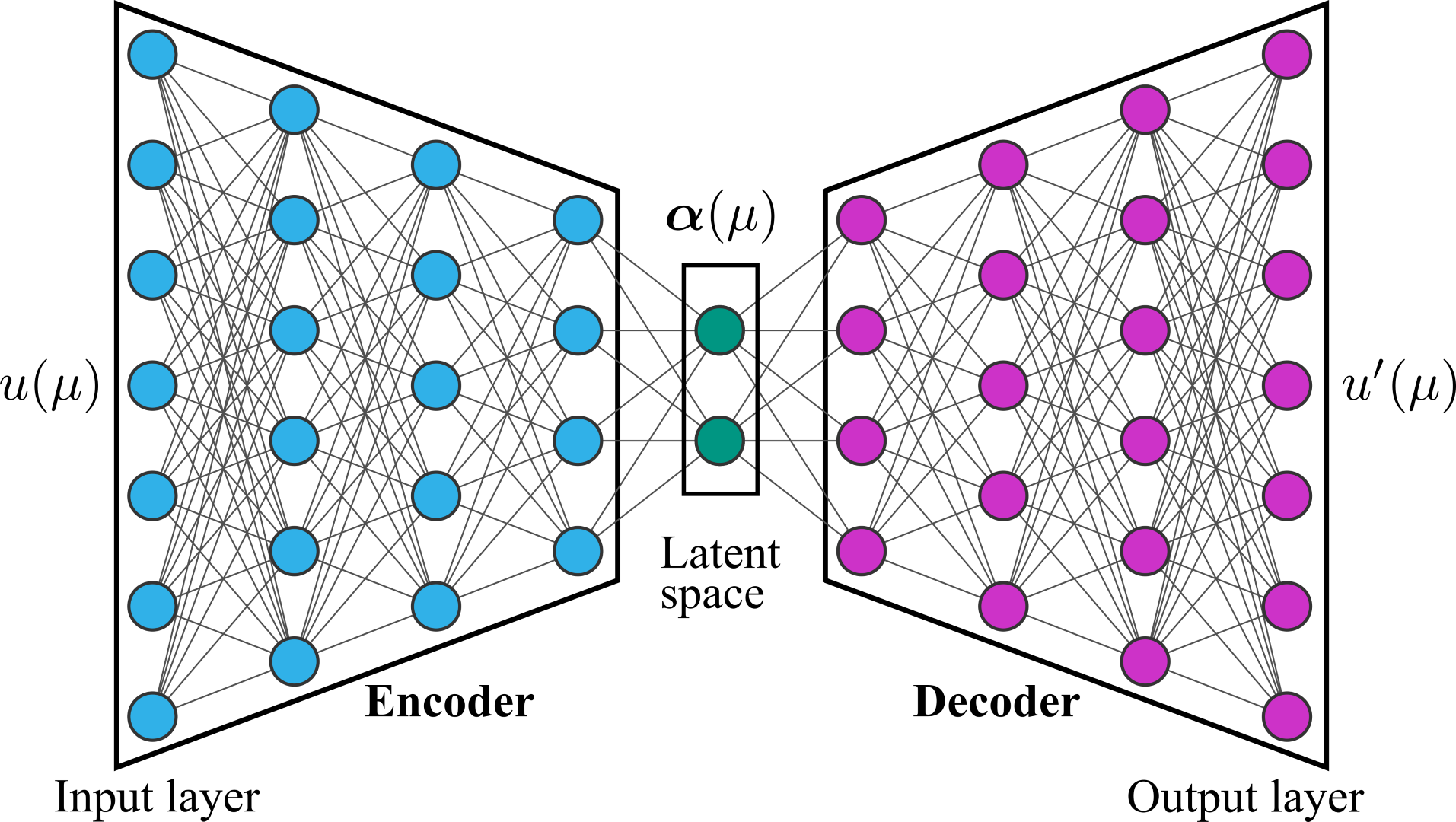}
    \caption{A sketch of an autoencoder. High-fidelity snapshots $\mathbf{u}(\mu)$ are compressed by the \emph{Encoder} into a low-dimensional \emph{latent space}  $\boldsymbol{\alpha} (\mu)$. The \emph{Decoder} is applied to reconstruct fields $\mathbf{u}'(\mu)$ from $\boldsymbol{\alpha} (\mu)$. }
    \label{fig:autoencoder}
\end{figure}

A widely used dimensionality reduction technique, \emph{autoencoder}, is explained as follows. Its architecture is sketched in Fig. \ref{fig:autoencoder}. The inputs are snapshots $\mathbf{u}(\mu)$, and the dimensionality is reduced by the \emph{Encoder} into a low-dimensional \emph{latent space} $\boldsymbol{\alpha} (\mu)$. By analogy with POD, $\boldsymbol{\alpha} (\mu)$ can be regarded as coefficients of a very low dimension. The reconstructed fields $\mathbf{u}'(\mu)$ are obtained via the \emph{Decoder} operation. The framework is built to minimize the difference between $\mathbf{u}(\mu)$ and $\mathbf{u}'(\mu)$, and the loss function is given by $ \mathcal{E} = \frac{1}{N_u} \lVert \mathbf{u}(\mu) - \mathbf{u}'(\mu) \rVert_2^2$.

\subsubsection{Physical Informed Neural Network (PINN)}
\label{subsubsec:pinn}
PINN is a type of regression model that involves knowledge of physical phenomena (i.e., PDEs) within the framework of neural networks. It was proposed by George Em Karniadakis et al. in 2017 \cite{raissi2017physics, raissi2019physics}, and has recently been adapted for domain decomposition problems. The following paragraphs provide a brief introduction to PINN. Its implementations for local ROMs will be presented in Section \ref{subsec:pinn_localROM} and \ref{subsec:Schwarz_based_iterative_scheme}.

In contrast to the typical networks that aim to match input and output (see equation \eqref{eq:mse}), the loss function of a PINN is intrinsically designed to integrate governing PDEs, boundary, and initial conditions. Namely, the MSE of a PINN is given by
\begin{equation}
\label{eq:pinn_loss}
    \mathcal{E} = w_f \mathcal{E}_f + w_b \mathcal{E}_b,
\end{equation}
where $\mathcal{E}_f$ and $\mathcal{E}_b$ are the PDE loss using predicted solutions and the error results from the boundary/initial conditions, respectively. The two are weighted by factors $w_f$ and $w_b$ to comprise the total MSE. Remark that $\mathcal{E}_f$ can be either residuals of the governing equations or an optimal formulation of the model problem. Properties of the MSE loss (equation \eqref{eq:pinn_loss}) guarantee that the PINN's predictions satisfy both the governing equations and physical constraints. 

\begin{figure}[h!]
    \centering
    \includegraphics[width=0.6\linewidth]{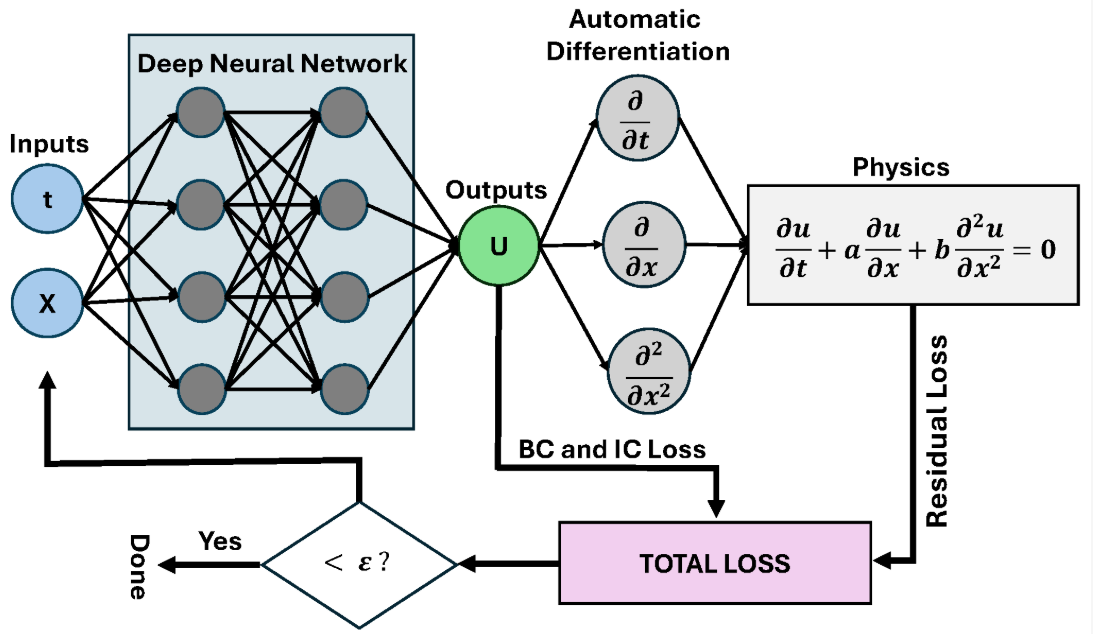}
    \caption{Architecture of a PINN \cite{trahan2024Quantum}. The inputs are coordinates $x$ in the computational domain and time $t$. The outputs are physical variables $\mathbf{u}(x, t)$ over the spatial and temporal domains. The MSE consists of two parts: the PDEs loss and the error with respect to the constraints (initial and boundary conditions). The parameters of the \emph{Deep Neural Network}, namely $w$ and $b$, are trained to reduce the combined loss. Figure taken from \cite{trahan2024Quantum}, copyright owned by MDPI.}
    \label{fig:PINNs}
\end{figure}

The architecture of a PINN is depicted in Fig. \ref{fig:PINNs}, which consists of a distinct structure compared to standard networks. The inputs are spatial coordinate $x$ and temporal coordinate $t$, and the output is the solution $\mathbf{u}(x, t)$. Suppose a PDE $\mathcal{D} \left ( \mathbf{u}(x, t) \right) = 0$ and a PINN prediction $\mathbf{u}'$, the residual loss is written as $ \mathcal{E}_f = \frac{1}{N_f} \lVert \mathcal{D} \left ( \mathbf{u}' (x,t) \right) \rVert_{2, \Omega}^2 $. Similarly, the constraint error is obtained by $\mathcal{E}_b = \frac{1}{N_b} \lVert \mathbf{u}(x, t) - \mathbf{u}'(x,t) \rVert^2_{2, \partial \Omega}$. In the two expressions, $N_f$ and $N_b$ denote the size of the corresponding data. $\Omega$ and $\partial \Omega$ indicate the interior and boundary regions for computing them, as discussed below.

\begin{figure}[h!]
    \centering
    \includegraphics[width=0.4\linewidth]{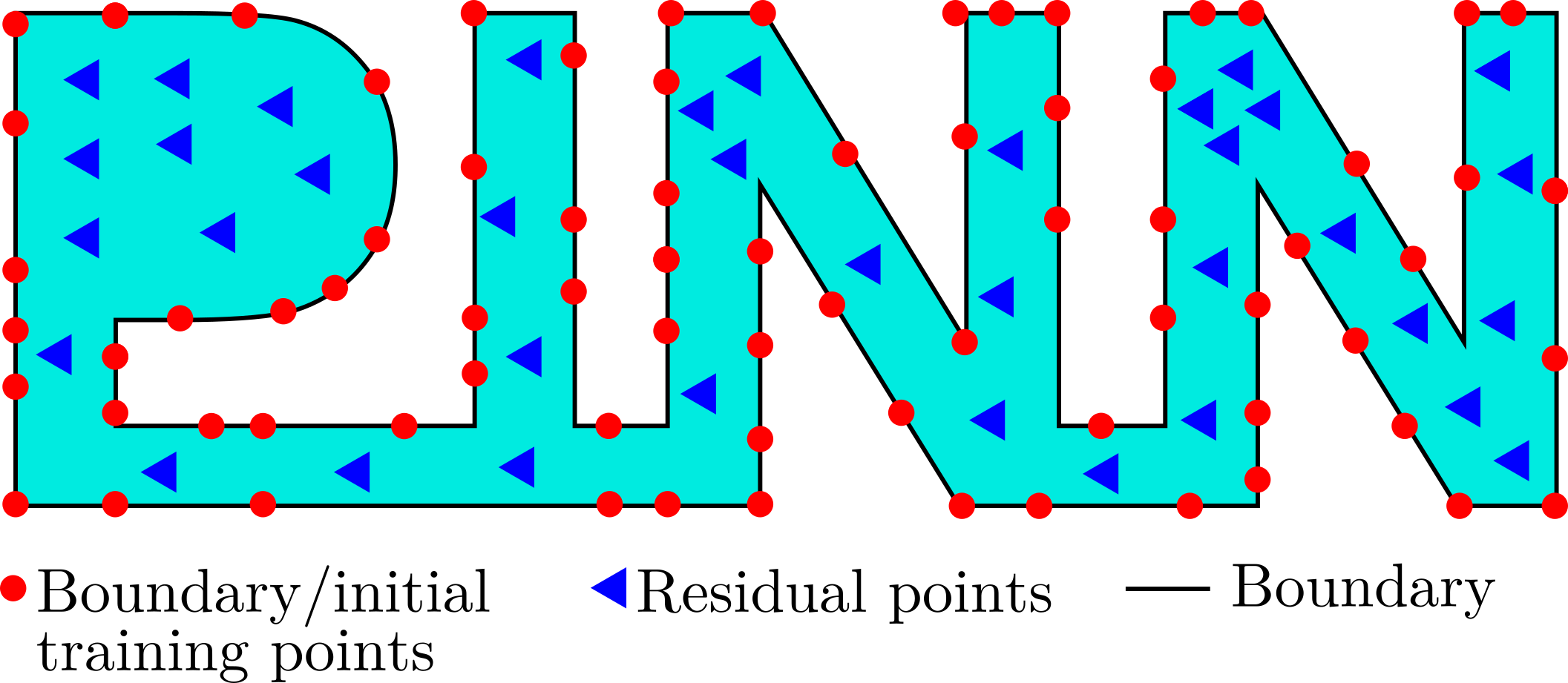}
    \caption{PINN training data points for a domain to compute loss functions. Figure redrawn based on \cite{jagtap2020conservative}.}
    \label{fig:pinns_training_points}
\end{figure}

The inputs $x$ and outputs $u(x)$ of the PINN framework are confined in a set of pre-selected data points as shown in Fig. \ref{fig:pinns_training_points}. The points are divided into groups for the computation of boundary/initial losses $\mathcal{E}_b$ and the residual loss $\mathcal{E}_f$. The PINN is trained only using the points during the offline stage. Once converged, it can be utilized to predict solutions at any location and time step within the computational domain. Thus, the selections of training points significantly affect the accuracy of a PINN. Details about those aspects can be seen in the specific review articles \cite{cuomo2022scientific, zhao2024comprehensive}.

Note that the inputs and outputs of a PINN are interconnected through a series of neurons as defined in equation \eqref{eq:nn_in_out}. This architecture enables the computation of derivatives using the chain rule, for example, $\partial \mathbf{u}(x, t)/\partial x$ and $\partial \mathbf{u}(x, t)/\partial t$, as well as higher-order derivatives. This capability forms the mathematical foundation for evaluating the residual loss term $\mathcal{E}_f$. This approach is known as \emph{automatic differentiation} in the framework of computational algebra \cite{griewank2008evaluating}.

Unlike traditional ROMs, PINN can still solve PDEs without requiring data. The approximations can be obtained purely by satisfying governing equations and constraints, as shown in Fig. \ref{fig:PINNs}. However, when FOM results are available at specific locations, the standard data MSE (equation \eqref{eq:mse}) can be integrated into the total PINN loss (equation \eqref{eq:pinn_loss}). That indicates $\mathcal{E} = w_f \mathcal{E}_f + w_b \mathcal{E}_b + w_d \mathcal{E}_d$, where $\mathcal{E}_d = \frac{1}{N_u} \lVert \mathbf{u}(x, t) - \mathbf{u}'(x, t) \rVert_2^2$. 

The inputs and outputs of PINNs depend on the problem's configuration. Recently, PINNs have been intensively employed to solve both forward and inverse problems \cite{jagtap2020conservative, dwivedi2021distributed, kharazmi2021hp, moseley2023finite}. In forward problems, PINNs approximate solutions of PDEs with fully specified parameters, boundary conditions, and initial conditions. The network inputs typically consist of spatial coordinates $x$, temporal coordinates $t$, and parameters $\mu$, while the output is the solution field $u(x,t;\mu)$ over the computational domain. In contrast, inverse systems aim to identify unknown parameters or conditions from observed data. Hence, a data loss term that enforces consistency with observed solutions is included. The inverse PINN can adopt the same inputs as the forward, while also providing extra quantities as outputs for inverse problems. More references for their implementations are presented in Section \ref{subsec:pinn_localROM}.

It is important to emphasize that while PINNs involve knowledge of the governing PDEs, their principle differs significantly from traditional intrusive methods. Recall that intrusive techniques numerically solve reduced algebraic systems derived from the PDEs. In contrast, the PINN framework replaces conventional solvers with neural networks that learn the solution operator. Since PINNs do not explicitly solve the PDE system through algebraic manipulation, they are appropriately classified as non-intrusive ROMs.

\subsection{Schwarz-based techniques}
\label{subsec:Schwarz_based_iterative_scheme}
We now start reviewing the first techniques, Schwarz-based techniques, in the category of data-driven methods. We will give a general description and then review the implementations in the references. Note that not only Dirichlet but also Neumann conditions are discussed in this subgroup.

The approaches that will be presented are almost the same as the Schwarz method and its variations. Instead of solving intrusive ROMs to update solutions in subdomains, non-intrusive surrogate models are built to solve each subproblem employing interface conditions and iteratively update local approximations until convergence is achieved. 

\subsubsection{Formulations}
We adopt multiplicative Schwarz (equation \eqref{eq:schwarz_problem}) for the explanation, but the algorithms are also suitable for other conditions presented in Section \ref{subsec:iterative}. 

Assume two overlapped subdomain $\Omega_m$ and $\Omega_n$, and their interfaces $\Gamma_{[mn]} = \partial \Omega_m \cap \Omega_n $ and $\Gamma_{[nm]} = \partial \Omega_n \cap \Omega_m $.  $\mathcal{S}_m$ denotes an operator that obtains the local solution $u_m$ regarding the local boundary conditions $\left. u_m \right|_{\Gamma}$. Namely $u_m = \mathcal{S}_m \left. u_m \right|_{\Gamma}$. 

The condition for the Schwarz algorithm denotes $\left. u_m^{(t)}\right|_{\Gamma_{[mn]}} = \left. u_n^{(t-1)} \right|_{\Gamma_{[mn]}}$ in equation \eqref{eq:schwarz_problem}. We assume other boundaries, e.g., $\partial \Omega_m \setminus  \Gamma_{[mn]}$, are fixed. Then, the local solution $u_m$ at step $t$ is given by 
\begin{equation*}
    u_m^{(t)} = \mathcal{S}_m \left. u_m^{(t)}\right|_{\Gamma_{[mn]}} = \mathcal{S}_m \left. u_n^{(t-1)} \right|_{\Gamma_{[mn]}}.
\end{equation*}

Note that $\Gamma_{[nm]}$ is an inner face of $\Omega_m$, and $\mathcal{T}_m$ is defined to be a restriction operator that extracts the interface values. That is $\left. u_m^{(t)} \right|_{\Gamma_{[nm]}} = \mathcal{T}_m u_m^{(t)} $. Combining the above formulation leads to 
\begin{equation*}
    \left. u_m^{(t)} \right|_{\Gamma_{[nm]}} = \mathcal{T}_m \circ \mathcal{S}_m \left. u_n^{(t-1)} \right|_{\Gamma_{[mn]}}, 
    \quad \text{ and } \quad
    \left. u_n^{(t)} \right|_{\Gamma_{[mn]}} = \mathcal{T}_n \circ\mathcal{S}_n \left. u_m^{(t)} \right|_{\Gamma_{[nm]}}.
\end{equation*}

Now, we combine $\mathcal{T}_m$ and $ \mathcal{S}_m$ as a single operator $\mathcal{Z}_m = \mathcal{T}_m \circ \mathcal{S}_m$, and the same to $\mathcal{Z}_n = \mathcal{T}_n \circ \mathcal{S}_n$. Then, the above formulation is rewritten as 
\begin{equation}
\label{eq:non_intrusive_schwarz}
    \left. u_m^{(t)} \right|_{\Gamma_{[nm]}} = \mathcal{Z}_m \left. u_n^{(t-1)} \right|_{\Gamma_{[mn]}}, 
    \quad \text{ and } \quad
    \left. u_n^{(t)} \right|_{\Gamma_{[mn]}} = \mathcal{Z}_n \left. u_m^{(t)} \right|_{\Gamma_{[nm]}}.
\end{equation}

Moreover, the relaxation strategy can be adopted to improve the convergence of the Schwarz method \cite{quarteroni1999domain}. Thus, equation \eqref{eq:non_intrusive_schwarz} is expressed as
\begin{equation}
\label{eq:non_intrusive_schwarz_relax}
\begin{aligned}
    \left. u_m^{(t)} \right|_{\Gamma_{[nm]}} &= \mathcal{Z}_m \left. u_n^{(t-1)} \right|_{\Gamma_{[mn]}},
    \quad \text{ and } \quad \\
    \left. u_n^{(t)} \right|_{\Gamma_{[mn]}} &=  \mathcal{Z}_n \left. \tilde{u}_m^{(t)} \right|_{\Gamma_{[nm]}} \hspace{0.5em} \text{ with } \hspace{0.5em} \left. \tilde{u}_m^{(t)} \right|_{\Gamma_{[nm]}} = \omega \left. u_m^{(t)} \right|_{\Gamma_{[nm]}} + (1-\omega) \left. u_m^{(t-1)} \right|_{\Gamma_{[nm]}},
\end{aligned}
\end{equation}
where $\omega$ is the relaxation factor.

Remark that the dimensionality reduction techniques like POD and neural networks are generally employed beforehand, and then $\mathcal{Z}$ is constructed for the regression for low-dimensional variables.

\subsubsection{Applications}
Now, the problems turn to creating a model to approximate the operator $\mathcal{Z}$ for each subdomain. It is clear that the construction of the surrogate model $\mathcal{Z}$ is not unique. The following paragraphs will discuss the implementation of both regression algorithms and neural networks to construct local ROMs.

Shi Cheng et al. used the method to study 2D nonlinear elliptic equations \cite{chen2021reduced}. They collected snapshots utilizing localized training and oversampling (see Section \ref{subsubSec:Localized_training}). The $\mathcal{Z}$, indeed a ROM, is constructed by \emph{two-layer Neural Networks}. The authors also showed several relevant prepressing and postpressing, and comparisons of different configurations. 

The technique can be extended to Dirichlet-Neumann conditions for non-overlapping divisions. In such a case, two surrogate models $\mathcal{Z}^D$ and $\mathcal{Z}^N$ are constructed to exchange Dirichlet and Neumann conditions on a shared interface. The coupling has the form of 
$\left. u_m^{(t)} \right|_{\Gamma_{[mn]}} = \mathcal{Z}^D_n \left. \frac{\partial u_n^{(t-1)}}{\partial \mathbf{n}} \right|_{\Gamma_{[mn]}}$
and 
$\left. \frac{\partial u_n^{(t)}}{\partial \mathbf{n}} \right|_{\Gamma_{[mn]}} = \mathcal{Z}^N_m \left. u_m^{(t)} \right|_{\Gamma_{[mn]}}$ for a two subdomain system with $\Omega = \Omega_m \cup \Omega_n$ and $\Gamma_{[mn]} = \partial \Omega_m \cap \partial \Omega_n$.

Niccolò Discacciati et al. adopted the Dirichlet-Neumann iteration (named \emph{boundary-to-boundary mapping} approach) to analyze a series of linear and nonlinear multi-physics scenarios \cite{discacciati2024model}, including steady and transient nonlinear heat equations, as well as a FSI problem. The localized training (see Section \ref{subsubSec:Localized_training}) is utilized to generate datasets. POD is adopted to reduce the dimensionality of snapshots. Regression-based ROMs are built using two techniques, ANN and \emph{Vectorial Kernel Orthogonal Greedy Algorithm}. The technical aspects of the two models are beyond the scope of this review. They are described in the original document. The method is also applied to both conforming and non-conforming high-fidelity meshes. Note that the authors also published an intrusive version of the framework, which is presented in Section \ref{subsubsec:Variations_Schwarz_method} \cite{discacciati2023localized}.

\subsubsection{Multiphysics phenomena}
The procedure can also be applied to multiphysics problems, especially FSI phenomena. 

Dunhui Xiao et al. \cite{xiao2016non-intrusive} utilized POD and radial basis functions to create ROM for several test cases: (i) flow past a cylinder; (ii) 2D free-falling square in water; (iii) vortex-induced vibrations of an elastic beam. The authors coupled Navier-Stokes flow and linear elastic structural motion. 

M. Barzegar Gerdroodbary et al. \cite{barzegar2025predictive} used POD and \emph{Convolutional Neural Network} to extract dominant features. Subsequently, they utilized \emph{Long Short-Term Memory} neural network to build a surrogate model. They simulated an Abdominal Aortic Aneurysm governed by unsteady Reynolds average Navier-Stokes and elastodynamics equations. A significant conclusion can be drawn from this study: the POD-based ROM is more efficient, while the pure neural network model is more accurate. 

Xinshuai Zhang et al. \cite{zhang2022data-driven} incorporated an \emph{Autoencoder} for data reduction and a so-called \emph{sparse identification of the nonlinear dynamics} algorithm to identify governing equations with respect to low-dimensional variables. The innovative framework was tested in a 2D case, flow past a cylinder. 

Moreover, the coupling of ROM and FOM can also follow the procedure. Renkun Han et al. \cite{han2022deep} proposed a method that couples a deep neural network-based fluid ROM and a high-fidelity structural dynamic solver. 

For investigations for multiphysics applications, beyond FSI, we suggest two references: (i) \cite{gobat2023reduced} that uses deep learning-based approaches for micro-electro-mechanical-systems; (ii) \cite{riva2024multi-physics} that couples neutronics and thermal-hydraulics in nuclear reactor cores.

\subsection{Interpolation algorithm}
The second group is for methods that employ interpolations. We will present the general principles of the interpolation-based approach, following studies that implement it in different ways.

The interpolation scheme employs iterations to couple local ROMs, and generally, the dimensionality reduction (e.g., POD) is applied to pre-process the high-fidelity data. The steps for constructing an interpolation-based ROM for domain decomposition problems are outlined in the following paragraphs. 

Suppose the POD is used to compute the reduced basis of each subdomain, $u_m=\sum_i^{N_\text{RB,m}} \alpha_{m,i} v_{m,i}$. Recall that the indices of all neighbors of $\Omega_m$ are included in a set $N_\Gamma(m) = \left\{ n | \partial \Omega_m \cap \partial \Omega_n \neq \emptyset \right\}$. Then, a interpolation model $\mathcal{Z}_{m}$ for the POD coefficient $\boldsymbol{\alpha}_m$ of $\Omega_m$ can be expressed as

\begin{equation}
\label{eq:interpolation_based_DD_ROM}
    \boldsymbol{\alpha}_{m}^{(t)} = \mathcal{Z}_{m} \left( \boldsymbol{\alpha}_{m}^{(t-1)}, \boldsymbol{\alpha}_{n}^{(t-1)}, \mu_m \right) \quad n \in N_\Gamma(m),
\end{equation}
where $\boldsymbol{\alpha}$ is a vector contains all coefficients, $\mathcal{Z}_{m}$ is the interpolation model, $\mu_m$ is the parameter and $t$ is the step number of iterations. The equation indicates that the subdomain approximation at iteration step $(t)$ depends on several inputs, including its own results in $(t-1)$ and the neighbours' results in $(t-1)$, as well as parameters $\mu_m$. 

\begin{figure}[h]
    \centering
    \includegraphics[width=0.3\textwidth]{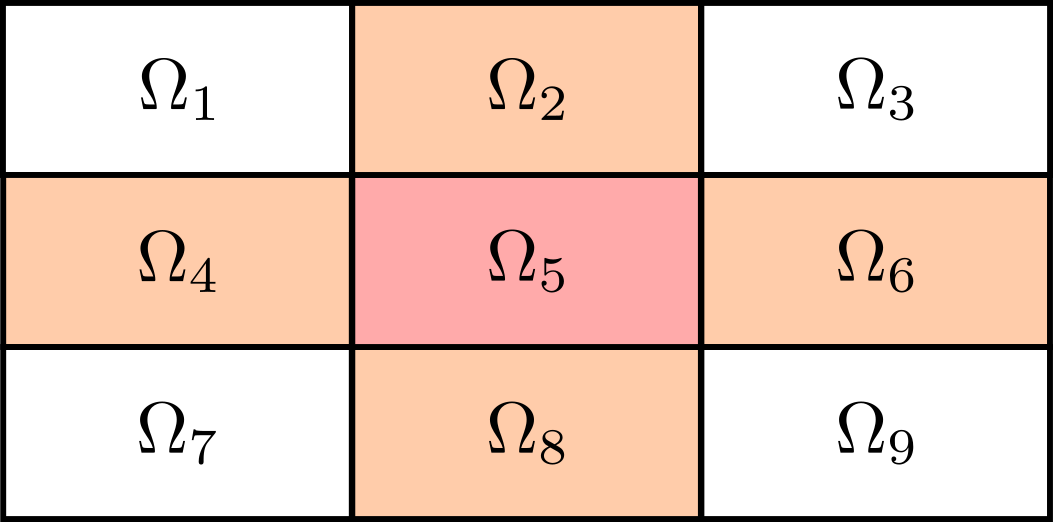}
    \caption{A subdomain and its neighbours. Redraw based on \cite{xiao2019non}.}
    \label{fig:subdomains_neighbours}
\end{figure}

To be clearer, an example is depicted in Fig. \ref{fig:subdomains_neighbours}, in which the target division is $\Omega_5$ and its neighbors are $\Omega_2,\ \Omega_4,\ \Omega_6,\text{ and } \Omega_8$. Therefore, we have $\boldsymbol{\alpha}_{5}^{(t)} = \mathcal{Z}_{5} \left( \boldsymbol{\alpha}_{5}^{(t-1)}, \boldsymbol{\alpha}_{2}^{(t-1)}, \boldsymbol{\alpha}_{4}^{(t-1)}, \boldsymbol{\alpha}_{6}^{(t-1)}, \boldsymbol{\alpha}_{8}^{(t-1)}, \mu_5 \right)$. Also, in case a subdomain is connected to the global boundary, the boundary conditions can be involved in the interpolation. 

Be aware that the above approach can incorporate other dimensionality reduction algorithms, and various techniques can be applied to build the regression model $\mathcal{Z}$, e.g., those introduced in Section \ref{subsec:data_driven_preliminaries}. Moreover, the way of incorporating $\mathcal{Z}$ to build ROMs is diverse. Hence, after reviewing the references, we identify three interpolation strategies. We will provide their ideologies and applications for different equations. 

\subsubsection{Standard interpolation}
The first method we presented comprises the two steps above, namely, dimensionality reduction and interpolation. The studies share the same procedure, while they differ in the techniques employed in the two stages.

Dunhui Xiao et al. have published several investigations exploiting the pure regression strategy \cite{xiao2017domain, xiao2019domainfluid, xiao2019domaintur}. They studied several large-scale problems governed by transient incompressible Navier-Stokes equations: (i) 2D flow past a cylinder; (ii) 2D and 3D urban street canyon test case; and (iii) 3D air flow around London South Bank University. 

Although the three applications are distinct, their procedures are nearly the same. The entire geometry is decomposed similarly into non-overlapping subdomains (as Fig. \ref{fig:subdomains_neighbours}). The global snapshots are computed and extracted for the partitions. POD is applied to calculate a set of separated dominant modes for each of them (see Section \ref{subsubsec:global_solution_local_RB}). The ROM is constructed through PODI with RBF, and a different number of POD modes is adopted for different regions to approximate the global fields. The results of non-intrusive domain decomposition ROM are comparable to FOM simulations and global ROM. In addition, fewer modes are required to represent "less important" areas, and more are used to capture dominant subdomain-level variances. To conclude, the number of modes in the local ROM is fewer than in the global procedure, and thus, the strategy can balance cost and accuracy. 

Besides POD, Claire E. Heaney et al. also incorporate AE to investigate two incompressible Navier-Stokes problems, a 2D single-phase flow past a cylinder and 3D two-phase in-pipe flow \cite{heaney2022ai}. In these cases, there are no overlaps among the divisions. Four dimensionality reduction methods are used and compared in their analysis, namely POD, a \emph{convolutional AE}, an \emph{Adversarial Autoencoder}, and a \emph{hybrid SVD AE}. A \emph{predictive neural network} is constructed to achieve the interpolation.

The flow past a cylinder case is modeled considering global snapshots and local RB (see Section \ref{subsubsec:global_solution_local_RB}). For the multiphase case, the tube is split axially into partitions of the same shape. Then, a set of generic local POD modes or AE latent space is utilized to reconstruct the global solutions. This is the generic decomposition and local RB construction indicated in Section \ref{subsubsec:global_solution_local_RB}. Thus, the authors managed to extend a ROM built from ten subdomains to approximate a longer pipe that consists of 100 subdomains. They concluded that the AE-based reduced techniques outperform POD, and the best among the four techniques is the convolutional AE.

\subsubsection{Linearized interpolation}

A modification of the above direct interpolation algorithm is proposed by Cong Xiao et al \cite{xiao2019non, xiao2019subdomain, xiao2021efficient} to model a large-scale reservoir, namely \emph{Subdomain POD-Trajectory-Piecewise-Linearization (TPWL)}. Given the prediction of the last step $\boldsymbol{\alpha}_{m}^{t-1}$ and parameter $\mu_m$, the surrogate model $\mathcal{Z}_m$ can predict a  new solution using a first-order expansion around the training data point,
\begin{equation}
    \boldsymbol{\alpha}^{(t)}_m = \boldsymbol{\alpha}_{tr}^{(t)} + \mathbf{E}_{m,\boldsymbol{\alpha}_\text{tr}}^{(t)} \left( \boldsymbol{\alpha}^{(t-1)}_m - \boldsymbol{\alpha}_{tr}^{(t-1)} \right) + \mathbf{G}_{m,\mu_\text{tr}}^{(t)}
    \left( \mu_m - \mu_{tr} \right),
\end{equation}
with
\begin{equation*}   
    \mathbf{E}_{m,\boldsymbol{\alpha}_\text{tr}}^{(t)} = \frac{\partial \mathcal{Z}_{m}^{(t)} }{\partial \boldsymbol{\alpha}_{tr}^{(t-1)}}, \quad
    \mathbf{G}_{m,\boldsymbol{\alpha}_\text{tr}}^{(t)} = \frac{\partial \mathcal{Z}_{m}^{(t)} }{\partial \mu_{tr}},
\end{equation*}
where $(\boldsymbol{\alpha}_{tr}^{(t)}, \boldsymbol{\alpha}_{tr}^{(t-1)}, \mu_{tr})$ belong to the \emph{closest} training set to $\boldsymbol{\alpha}^{(t)}_m$. The selection algorithm for closeness is presented in \cite{he2015constraint}. Note that $\mathcal{Z}_m$ is created based on RBF, so the expression of $\mathbf{E}_{m,\boldsymbol{\alpha}_\text{tr}}^{(t)}$ and $\mathbf{G}_{m,\boldsymbol{\alpha}_\text{tr}}^{(t)}$ can be obtained analytically after offline training.

The reservoir is partitioned into non-overlapped subdomains, and local RBs are obtained for each subdomain via POD (as indicated in Section \ref{subsubsec:global_solution_local_RB}). The results of global POD and local POD are compared. It can be concluded that, compared to global modes, fewer local modes are needed in each subdomain for the same level of reconstruction accuracy. Note that the local ROM results are discontinuous along the subdomain interface. Thus, the fields are smoothed by solving an additional minimization problem.

Furthermore, Cong Xiao et al. aim to solve an inverse problem. They supposed that the governing equations and solutions in a few locations of the computational domain are known. Still, the parameters that exist in the PDEs over the domain are unknown. Therefore, they aimed to compute parameter fields exploiting known solutions. That can be achieved via an optimization formulation that minimizes the differences between observations and computations. 

The details and formulations of the smoothing algorithm and the inverse problem are beyond the focus of this review, and they are well explained in \cite{xiao2021efficient}.

\subsubsection{Hybrid technique}
A hybrid interpolation technique for inverse problems with non-overlapped decomposition is proposed by Rossella Arcucci et al. \cite{arcucci2020domain, arcucci2020adaptive}. The authors simplify equation \eqref{eq:interpolation_based_DD_ROM} and interpolate POD coefficients as 
\begin{equation*}
    \boldsymbol{\alpha}_{m}^{(t)} = \mathcal{Z}_{m} \left( \boldsymbol{\alpha}_{m}^{(t-1)}, \mu_m \right).
\end{equation*}

Then, the authors proposed a so-called Domain Decomposition Reduced Order Data Assimilation (DD-RODA) process to solve the inverse problem. The method results in an unconstrained least squares formulation
\begin{equation}
    u_m^{\text{RODA}} = \arg\min_{u_m} \left\{ \lVert u_m - u_m^{(t)} \rVert + \lVert d_m^{(t)} - M \left(u_m^{(t)}\right) \rVert
    \right\},
\end{equation}
where $u_m^{(t)}$ is the ROM prediction in time $(t)$, $d_m^{(t)}$ is the observed value of a finite number of locations in time $(t)$, and $M \left(u_m^{(t)}\right)$ is an operator extract observed points from the prediction fields.

The method is applied to solve the air pollution problem of London South Bank University, modeled by 3D incompressible Navier-Stokes equations. Known measurements in several locations, $d_M^{(t)}$, are utilized to calibrate local predictions. The results reveal that DD-RODA performs better than standard local ROM. Note that no extra procedure is applied to smooth the local predictions, which results in discontinuities at interfaces for global fields.

\subsection{Optimization-based technique}
We turn now to the third category, the so-called Optimization-based methodologies. The technique aims to formulate optimal systems to smooth the discontinuities among the subdomain interfaces. 

Two methods, based on the least squares Galerkin projection and Gappy POD, are included in the group. Due to their diversity, details can be observed from the explanation below rather than a general description here.

\subsubsection{Least squares Galerkin}

Youngsoo Choi et al. adopted the least squares Galerkin method to construct the optimal systems of domain decomposition problems \cite{diaz2023nonlinear, diaz2024fast}. The numerical setups are almost the same as those presented in their previous publication \cite{hoang2021domain} in Section \ref{para:Least_squares_Petrov-Galerkin}. Similar objective functions are constructed, while solutions are approximated using data-driven techniques instead of Galerkin projection-based formulations.

In equation \eqref{eq:least_square_galerkin}, the residuals are formulated via a Galerkin projection. In contrast, in the non-intrusive procedure shown in \cite{diaz2023nonlinear, diaz2024fast}, the authors use an AE (defined as $\mathcal{A}$) to construct low-dimensional representations. Suppose a parameter set $\{ \mu_j \}_{j=1}^{N_\mu}$ and FOM solutions defined for interiors $\mathbf{u}_m^\Omega (\mu_j) \coloneqq \left. \mathbf{u} (\mu_j) \right |_{\Omega_m} $ and interfaces $\mathbf{u}_m^\Gamma (\mu_j) \coloneqq \left. \mathbf{u} (\mu_j) \right |_{\Gamma \in \partial \Omega_m} $. The autoencoders aim to minimize the respective MSE losses
\begin{equation*}
    \mathcal{E}_m^\Omega = \frac{1}{N_\mu} \sum_{j=1}^{N_\mu} \left\| \mathbf{u}_m^\Omega (\mu_j) - \mathcal{A}^\Omega_m (\mu_j) \right\|_2^2, \quad
    \mathcal{E}_m^\Gamma = \frac{1}{N_\mu} \sum_{j=1}^{N_\mu} \left\| \mathbf{u}_m^\Gamma (\mu_j) - \mathcal{A}^\Gamma_m (\mu_j) \right\|_2^2
\end{equation*}
for internal $\mathbf{u}_m^\Omega$ and surface $\mathbf{u}_m^\Gamma$ values of each subdomain $m = 1, \cdots, N_\Omega$. 

A general formulation for approximations of any new parameter $\mu'$, i.e., ${\mathbf{u}'}_m^\Omega (\mu') $, can be expressed as ${\mathbf{u}'}_m^\Omega (\mu') \approx \left. \mathcal{A}(\mu') \right|_{\Omega_m} = \mathcal{A}^\Omega_m (\mu')$ and ${\mathbf{u}'}_m^\Gamma (\mu') \approx \left. \mathcal{A}(\mu') \right|_{\Gamma\in \partial \Omega_m} = \mathcal{A}^\Omega_m (\mu')$. Similar to equation \eqref{eq:least_square_galerkin}, ROM for a geometry $\Omega = \sum_{m=1}^{N_\Omega} \Omega_m$ is now formulated as an optimality system, 
\begin{equation}
    \min_{\Omega_m,\ m =1, \cdots, N_{\Omega}} \frac{1}{2} \sum_{m=1}^{N_\Omega} \left \| \mathcal{R}_m \left( \mathcal{A}^\Omega_m (\mu') \cup \mathcal{A}^\Gamma_m (\mu') \right) \right\|^2_2,
\end{equation}
subject to the constraints,
\begin{equation*}
    \sum_{\Gamma \in \partial \Omega_m \cap 
    \partial \Omega_n}^{N_\Gamma} \int_{\Gamma} \left | \mathcal{A}^\Gamma_m (\mu') - \mathcal{A}^\Gamma_n (\mu') \right|  = 0,
\end{equation*}
where $\Omega_m$ and $\Omega_n$ are two adjacent partitions. This nonlinear constrained minimization problem is solved via Sequential Quadratic Programming (SQP). While a detailed discussion of SQP exceeds the scope of this review, a comprehensive treatment can be found in \cite{diaz2024fast}.

The method is applied to a 2D steady-state Burgers' equation for overlapped subdomains. The authors compared the pure data-driven approach with their previous LS-ROM (see Section \ref{para:Least_squares_Petrov-Galerkin}). The results indicate that LS-ROM achieves higher speedup, and the advantages are more evident when integrating hyper-reduction techniques (achieved via sampling of the mesh elements). Nevertheless, this non-intrusive method performs better in terms of accuracy.

\subsubsection{Gappy POD optimization}

Nikhil Iyengar et al. developed a different optimization framework using Gappy-POD \cite{iyengar2022nonlinear, iyengar2024domain}. This method reconstructs optimal global solutions by combining local and global POD modes. We will first outline Gappy-POD \cite{everson1995karhunen}, and then explain its implementation for approximating smooth fields in multiple subdomains.

Gappy POD is an approach to finding the optimal approximation in a domain when POD modes and values are available at a few locations. The process for achieving this is discussed as follows.

Assume a field $\mathbf{u}$ and its POD modes $\mathbf{v}_i$, we have $\mathbf{u} \approx \sum_{i=1}^{N_\text{RB}} \alpha_i \mathbf{v}_i$. A solution with some unknowns in the domain is defined as $\mathbf{u}^* = \mathbf{m} \odot \mathbf{u}$, 
where $\odot$ denotes element-wise product and $\mathbf{m} = [m_1, \cdots, m_n]$ is a mask vector. $\mathbf{m}$ have the same dimension as $\mathbf{u}$ and it contains zeros in the positions where $\mathbf{u}$ is unknown and ones for the known locations. The goal of Gappy-POD is to find a set of coefficients $\tilde{\alpha}$ to reconstruct the approximation as $\tilde{\mathbf{u}} \approx \sum_{i=1}^{N_\text{RB}} \tilde{\alpha}_i \mathbf{v}_i$. 

Indeed, the procedure can be addressed as: $\mathbf{v}_i$ and $\mathbf{u}^*$ are known, and the solution in the whole domain can be obtained once $\tilde{\alpha}_i$ is determined. Remark that $\tilde{\alpha}$ can be computed by minimizing the error $\mathcal{E} = \left\| \mathbf{u}^* - \mathbf{m} \odot \left( \sum_{i=1}^{N_\text{RB}} \tilde{\alpha}_i \mathbf{v}_i \right) \right\|_2$. The Least-squares algorithm can be employed to solve the optimization problem, yielding a linear algebraic system
\begin{equation*}
    M \tilde{\boldsymbol{\alpha}} = \mathbf{f},
\end{equation*}
where $M_{ij} = (\mathbf{v}_i, \mathbf{v}_j)$ and $f_i = (\mathbf{v}_i, \mathbf{u}^*)$.

Then, the approximation obtained via the Gappy-POD, $\mathbf{u}'$, is expressed as 
\begin{equation}
\label{eq:gappy_POD_reconstruction}
    u' = \left\{ \begin{matrix} \sum_{i=1}^{N_\text{RB}} \tilde{\alpha}_i \mathbf{v}_i, &\quad &\mathbf{m}=0 \\ \mathbf{u}^*, &\quad &\mathbf{m}=1 \end{matrix} \right. \ .
\end{equation}

\begin{figure}
    \centering
    \includegraphics[height=3cm]{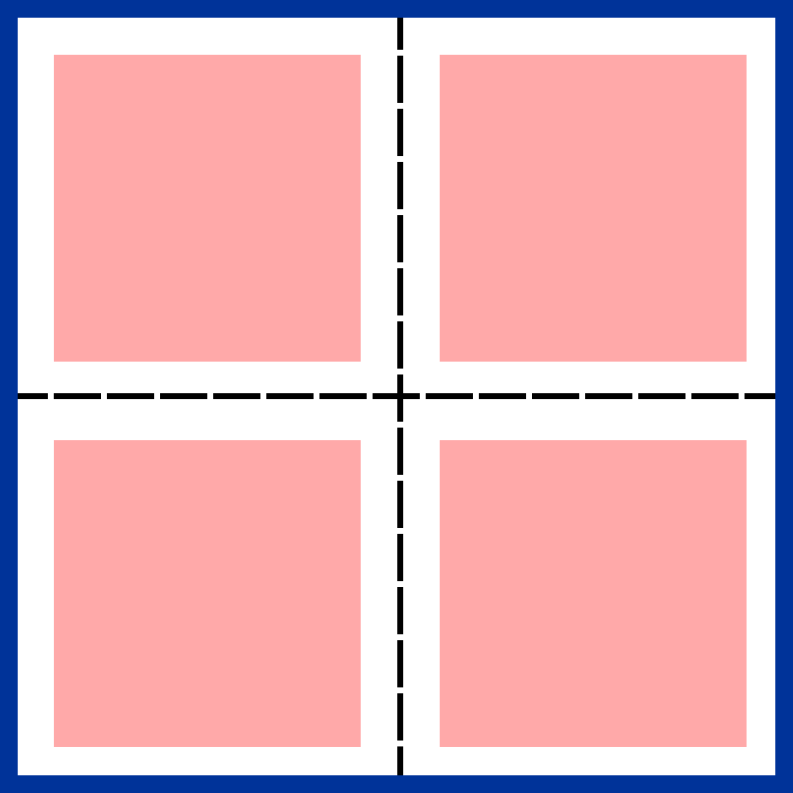}
    \caption{Global reconstruction using Gappy-POD. The entire domain is decomposed into non-overlapping subdomains (split by dashed lines). Local ROM is applied to predict results in partitions (the red region) that serve as $\mathbf{u}^*$, i.e., areas with $m_i\neq 0$. The remaining white regions are approximated using Gappy-POD.}
    \label{fig:gappy_POD}
\end{figure}

The implementation of Gappy-POD for domain decomposition problems proceeds as follows. First, full-domain simulations are performed to generate global solutions. The global POD modes $\mathbf{v}_i$ are computed. Then, the domain is decomposed into non-overlapping subdomains, and solutions in each subdomain are extracted as local snapshots. Local POD modes are calculated for each subdomain. PODI is used to construct ROMs for each partition. As illustrated in Fig. \ref{fig:gappy_POD}, $\mathbf{u}^*$ is non-zero in the red areas of each division, while it is zero in the rest of the domain. Consequently, the global POD coefficients $\tilde{\alpha}$ are obtained via Gappy-POD. Finally, the complete solution is reconstructed as \eqref{fig:gappy_POD}. Note that the local ROMs alone produce non-smooth predictions across subdomain interfaces. Thus, the authors employ Gappy-POD and global POD modes to ensure solution smoothness throughout the entire domain.

The authors validated the method using three cases: (i) quasi-1D Flow through a Converging-Diverging Nozzle; (ii) flow over a Wedge with Shocks; (iii) transonic Flow over an Airfoil (i.e., Reynolds-Averaged Navier-Stokes equations). The comparison with the pure global ROM reveals the benefits of adopting the aforementioned local ROM technique.

\subsection{Physical Informed Neural Network}
\label{subsec:pinn_localROM}
The fourth technique introduced within the data-driven category is PINN. The procedure of PINNs to solve PDEs in a single domain is already explained in Section \ref{subsubsec:pinn}. Note that the equation and boundary/initial condition losses also apply to domain decomposition problems. Additionally, the existence of interfaces among subdomains requires extra treatment at these interfaces. That denotes the continuity of solutions over the global domain.

Since the partitions can be either overlapping or non-overlapping, this also results in different strategies for dealing with the interfaces. Thus, we categorize the two groups based on the available references. If no overlaps among divisions, the interface continuity conditions denote the equality of the variable and its flux on both sides of the face \cite{canuto2006spectral}. If subdomains overlap, the Schwarz iteration with the exchange of Dirichlet conditions can be employed \cite{quarteroni1999domain}. The two scenarios are presented separately as follows.

\subsubsection{Interface constraints}
\label{subsubsec:pinn_interface_constraints}
For two adjoining domains with a single interface, to ensure the assembly of local approximations equal to the global approximation, both the variables and their flux should be equal on both sides of an interface \cite{canuto2006spectral}.
Consequently, in the frame of PINN, two losses are defined to satisfy the equality. Assume a subdomain $\Omega_m$ and its interface set $\mathcal{I}_m = \left \{ \Gamma_r | \Gamma_r \in \partial \Omega_m \right\} $, the two terms are written by 
\begin{equation}
\label{eq:PINN_interface_terms}
\begin{split}
    \mathcal{E}_{\Gamma_m,\text{flux}} &= \sum_{\Gamma_r \in \mathcal{I}_m}  \frac{1}{N_{\Gamma_r}} \lVert q (\mathbf{u})|_{\Gamma_r^+} \cdot \mathbf{n} - q (\mathbf{u})|_{\Gamma_r^-} \cdot \mathbf{n} \rVert_2^2, \\
    \mathcal{E}_{\Gamma_m,\mathbf{u}} &= \sum_{\Gamma_r \in \mathcal{I}_m}  \frac{1}{N_{\Gamma_r}} \lVert \mathbf{u}|_{\Gamma_r^+} - \{\mathbf{u}\}|_{\Gamma_r} \rVert_2^2, \\
\end{split}    
\end{equation}
where $q$ denotes flux, the two sides of a interface are noted with $^+$ and $^-$, $\mathbf{n}$ is the normal vector of the face, $N_{\Gamma_r}$ are training points assigned at $\Gamma_r$, and $\{ \mathbf{u} \} = \frac{1}{2} \left( \mathbf{u}|_{\Gamma_r^+} + \mathbf{u}|_{\Gamma_r^-} \right) $ is the average operator. Thus, the local loss for $\Omega_m$ can be formulated as 
\begin{equation}
\label{eq:local_pinn_loss}
    \mathcal{E}_m = w_{f,m} \mathcal{E}_{f,m} + w_b \mathcal{E}_b + w_{\Gamma_m,\text{flux}} \mathcal{E}_{\Gamma_m,\text{flux}} + w_{\Gamma_m,\mathbf{u}} \mathcal{E}_{\Gamma_m,u},
\end{equation}
where all terms are weighted by different factors and $\mathcal{E}_b$ appears only when $\Omega_m$ connected to global boundaries.

Be aware that a common feature of the following studies is that interface constraints are included as additional loss terms. However, the strategy of including the constraints is not unique. Over the past five years, various PINN-based architectures have been proposed. One can construct local PINNs for each subdomain or a single PINN whose loss function is the sum of all partition-level losses. Moreover, FOM solutions and POD can be incorporated to reduce the complexity of the PINN structure. These aspects are described below.

\paragraph{Network for physical fields}
A conservative PINN (cPINN) is proposed by Ameya D. Jagtap et al. \cite{jagtap2020conservative} to handle non-overlapping domain decomposition problems, which comprise local neural networks for each subdomain. That yields same local error formulation (see equation \eqref{eq:local_pinn_loss}) for each subdomain as $\mathcal{E}_m = w_{f,m} \mathcal{E}_{f,m} + w_{b,m} \mathcal{E}_{b,m} + w_{\Gamma_m} \left( \mathcal{E}_{\Gamma_m,\text{flux}} + \mathcal{E}_{\Gamma_m,u} \right)$.

\begin{figure}[h]
    \centering
    \includegraphics[width=0.5\linewidth]{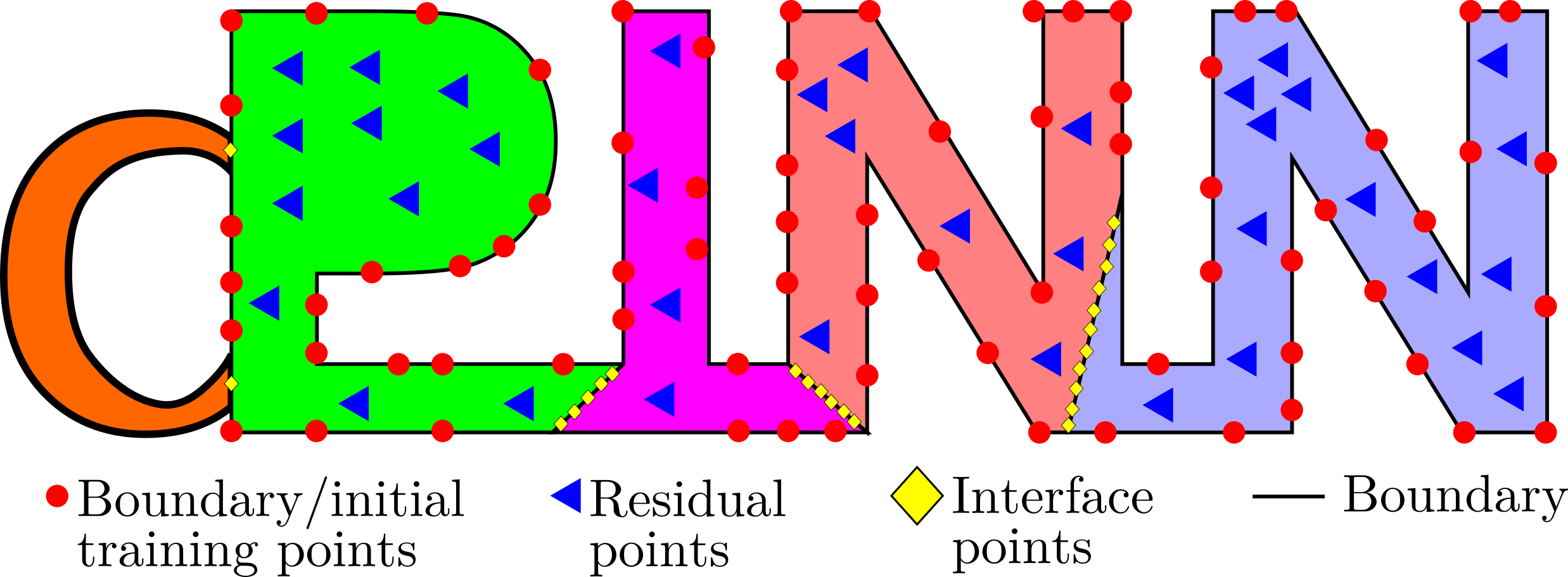}
    \caption{Computational domain and cPINN training data points. Figure redrawn based on \cite{jagtap2020conservative}.}
    \label{fig:PINNs_points_subdomains}
\end{figure}

The cPINN was used to solve both forward and inverse problems of various nonlinear PDEs, including Burgers, Korteweg-De Vries, incompressible Navier-Stokes, and compressible Euler equations. The computational domain and the location of training points are displayed in Fig. \ref{fig:PINNs_points_subdomains}. The whole geometry is decomposed into several subregions, and the two additional errors are computed at the interface points.

The authors tested and compared the effect of PINN structures, including the number of hidden layers, the number of neurons per layer, the activation functions, and the size of training points. Their numerical experiments for the 1D Burgers' equation demonstrate that increasing the number of hidden layers and neurons can improve accuracy. A key advantage of cPINN is that different network architectures can be defined for each subdomain to better capture local physical phenomena, thereby balancing the training cost and accuracy. 

Ameya D. Jagtap et al. further developed the cPINN and proposed the \emph{eXtended PINN} (XPINN), which can be extended to any PDEs (both forward and inverse problems) with non-overlapping partitions \cite{jagtap2020extended, hu2022extended}. The key modification of XPINNs is that a residual continuity loss is added on common interfaces, which is given by 
\begin{equation*}
    \mathcal{E}_{\Gamma_m,f} = \sum_{\Gamma_r \in \mathcal{I}_m}  \frac{1}{N_{\Gamma_r}} \left\| \left. \mathcal{D} \left(\mathbf{u}' (x,t)\right) \right|_{\Gamma_r^+} - \left. \mathcal{D} \left(\mathbf{u}' (x,t)\right) \right|_{\Gamma_r^-} \right\|_2^2.
\end{equation*}

The combined local MSE for is now termed by $\mathcal{E}_m = w_{f,m} \mathcal{E}_{f,m} + w_{b,m} \mathcal{E}_{b,m} + w_{\Gamma_m} \left( \mathcal{E}_{\Gamma_m,f} + \mathcal{E}_{\Gamma_m,u} \right)$. Note that the equality and residual continuity terms can sufficiently enforce the XPINN approximation to satisfy the global governing PDE. Thus, the flux matching constraints $\mathcal{E}_{\Gamma_m,\text{flux}}$ are not necessary and can be optionally imposed depending on the PDEs. 

The cPINN and XPINN were compared in terms of parallel efficiency for 1D Burgers' and 2D steady-state incompressible Navier-Stokes equations \cite{shukla2021parallel}. The results show that cPINNs are more efficient for spatial decomposition problems. In contrast, XPINNs are more flexible for temporal domain decomposition and are additionally well-suited for arbitrarily shaped and complex subdomains.

Vikas Dwivedi et al. proposed a different framework for solving forward and inverse problems of several PDEs among non-overlapping subdomains \cite{dwivedi2021distributed}, which is named \emph{Distributed PINN} (DPINN). Their method integrates temporal decomposition and trains separate networks for each time step. That denotes an additional temporal continuity loss $\mathcal{E}_{\text{t}}$ to match two sequential time steps $(t-1)$ and $(t)$, which is given by
\begin{equation*}
    \mathcal{E}_{\text{t}} = \frac{1}{N_t} \lVert \mathbf{u}^{(t)} - \mathbf{u}^{(t-1)} \rVert_2^2.
\end{equation*} 

Additionally, the total loss function of DPINN is distinct from previous architectures. The authors also define the individual network for each subdomain. However, DPINN doesn't define $\mathcal{E}_m$ for each subdomain. Instead, the total loss is expressed as a summation of all $\mathcal{E}_m$. That is 
\begin{equation*}
    \mathcal{E} = \sum_{m=1}^{N_\Omega} w_{f,m} \mathcal{E}_{f,m} + w_b \mathcal{E}_b + \sum_{r=1}^{N_\Gamma} \left( w_{\Gamma_r,\mathbf{u}} \mathcal{E}_{\Gamma_r,\mathbf{u}} +w_{\Gamma_r,\text{flux}} \mathcal{E}_{\Gamma_r,\text{flux}} \right) + \sum_{t=1}^{N_t}  w_{t} \mathcal{E}_{t},
\end{equation*}
where $N_\Omega$, $N_\Gamma$, and $N_t$ are the number of subdomains, interfaces, and time steps, respectively.

The DPINN is tested in several 2D benchmark problems, including steady heat conduction, Laplacian, Burgers', advection, and steady Navier-Stokes equations.

The PINNs discussed above utilize the strong form for loss formulation. Given the prevalence of weak formulations based on Galerkin methods in the intrusive framework, Ehsan Kharazmi et al. developed the \emph{hp}-variational PINN (\emph{hp}-VPINN) \cite{kharazmi2021hp} for nonoverlapping divisions. In this approach, the trial space is represented globally by a single neural network across the computational domain. The piecewise polynomials (i.e., Legendre polynomials, as mentioned in the work) are applied separately as test functions for every partition. In this setting, the VPINN achieves \emph{h}-refinement through increased spatial resolution and \emph{p}-refinement via higher-order polynomial test functions. One can analogize this framework to FEM employing the Petrov-Galerkin projection. Each subdomain can be regarded as an element of FEM, and different test and trial spaces are adopted.

Instead of multiple local PINNs, a single PINN is built and trained to approximate the global solution. In such a case, the total loss function of \emph{hp}-PINN, $\mathcal{E}^{hp}$, is a summation of all local losses $\mathcal{E}_m$ as shown in equation \eqref{eq:local_pinn_loss}. Thus,
\begin{equation*}
\begin{aligned}
    \mathcal{E}^{hp} &= (\mathcal{E}, \mathbf{w}), \\
    \mathcal{E} &= \sum_{m=1}^{N_\Omega} w_{f,m} \mathcal{E}_{f,m} + w_b \mathcal{E}_b,
\end{aligned}
\end{equation*}
where $(\cdot, \cdot)$ denotes the inner product operation, $\mathbf{w}$ is the test function, and $N_\Omega$ is the number of subdomains. Note that only a single PINN is employed to approximate the entire solution, and the polynomials test functions are continuous at interfaces. Thus, no interface constraints are required in $\mathcal{E}^{hp}$. 

\begin{figure}[h]
    \centering
    \includegraphics[width=\linewidth]{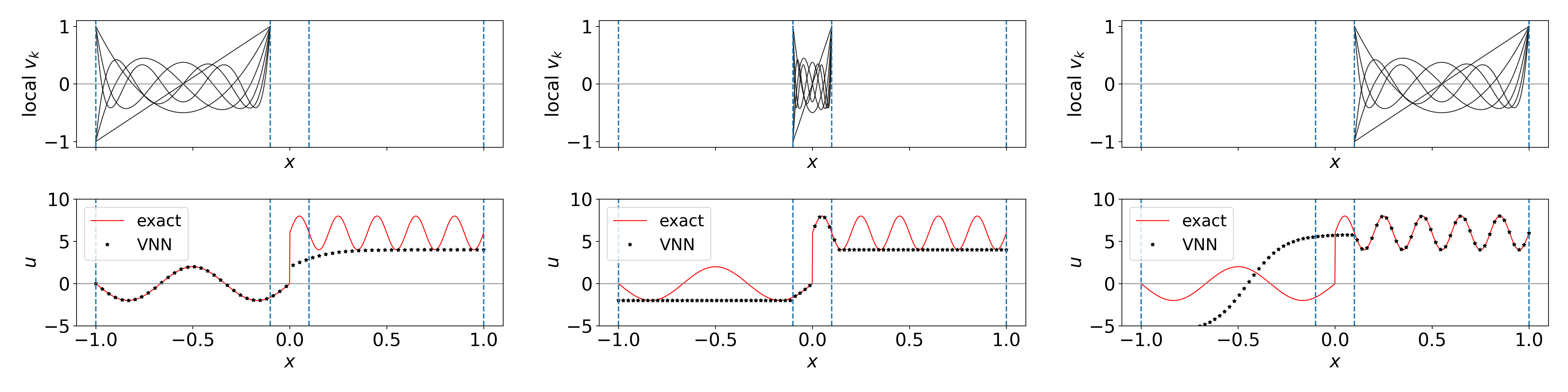}
    \caption{Local test functions of VPINN over each subdomain (first row). Results of each subdomain using VPINN (second row). The dashed blue lines are the sub-domain boundaries. Figures redrawn based on \cite{kharazmi2021hp}.}
    \label{fig:VPINN_local_test}
\end{figure}

An example of local test functions and the solution is illustrated in Fig. \ref{fig:VPINN_local_test}. The test functions are only defined inside the corresponding subdomains. All local solutions are \emph{glued} together to form a global result.

The efficiency and accuracy of \emph{hp}-VPINN were validated and compared with PINN in Poisson's equation and the advection-diffusion equation. The results of 2D Poisson's equation demonstrate that the \emph{hp}-VPINN with domain decomposition performs better than standard PINN, and the \emph{h}-refinement (increasing the number of subdomains in specific regions) yields a better local approximation.

\paragraph{Network employing POD modes}
Rather than directly approximating solutions, Xinyu Pan and Dunhui Xiao proposed a novel approach that employs neural networks to obtain POD coefficients, called the \emph{Domain Decomposition Physics-Data Combined Neural Network} (DD-PDCNN) \cite{pan2024domain}. The DD-PDCNN loss function is similar to cPINN, that is
\begin{equation*}
    \mathcal{E} = w_f \mathcal{E}_f + w_b \mathcal{E}_b + w_I \mathcal{E}_I + w_{\Gamma, \text{flux}} \mathcal{E}_{\Gamma,\text{flux}} + w_{\Gamma, \mathbf{u}} \mathcal{E}_{\Gamma,\mathbf{u}},
\end{equation*}
where $\mathcal{E}_I$ and $\mathcal{E}_b$ represents the initial and boundary condition terms, respectively. These terms are similar to those of cPINN, but they are formulated in terms of POD modes rather than solutions. 

Note that the usage of POD modes $\mathbf{u} = \sum_{i=1}^{N_\text{RB}} \alpha_i \mathbf{v}_i$ yields substantial computational savings in computing loss functions for the PINN. For example, assume a residual $\mathcal{E}_f = \Delta \mathbf{u}$, the error can be rewritten as $\Delta \mathbf{u} = \sum_{i=1}^{N_\text{RB}} \alpha_i \Delta \mathbf{v}_i$. Hence, $\Delta \mathbf{v}_i$ can be pre-computed and fixed for the training. The losses in each training step can be obtained by multiplying them by the POD coefficients.

The technique was applied to three cases: the Korteweg-de Vries equation, Kovasznay flow, and the incompressible Navier-Stokes equations. In comparison to the global PDCNN, the DD-PDCNN comprised fewer hidden layers and neurons, and it is far more accurate. That is consistent with the aforementioned PINN-based architectures. A disadvantage of this approach is that it requires high-fidelity solutions to compute POD modes. 

\subsubsection{Schwarz iteration}
The last sub-category of the data-driven techniques combines the PINN and the Schwarz method. According to the literature, the combination is not unique. It can be applied to the training stage, in which the values at overlapping regions are exchanged between adjacent domains. Also, one can construct PINN-based ROMs and then couple them iteratively via the Schwarz algorithm. The details of the two procedures are presented as follows.

\paragraph{Schwarz for training PINNs}
The Schwarz procedure is employed here to iteratively assign local boundary values of multiple subdomain-level PINNs. That means, for $\Omega_m$, the loss function is given by 
\begin{equation*}
    \mathcal{E}_m = w_{f,m} \mathcal{E}_{f,m} + w_{b,m} \mathcal{E}_{b,m},
\end{equation*}
where $\mathcal{E}_{b,m}$ is error of local boundaries. Suppose $\Gamma_{[mn]} = \Omega_n \cap \partial \Omega_m$, we can obtain $\mathcal{E}_{b,m}$ via 
\begin{equation*}
    \mathcal{E}_{b,m} = \frac{1}{N_{\Gamma_{[mn]}}} \left \| \left. \mathbf{u}_m^{(t)} \right|_{\Gamma_{[nm]}} - \left. \mathbf{u}_n^{(t-1)} \right|_{\Gamma_{[nm]}} \right\|,
\end{equation*} 
where $(t)$ denotes the step of training. 
Consequently, in the $(t)$ iteration, the BC of $\Omega_m$ on $\Gamma_{[mn]}$ is updated using the approximation $\mathbf{u}_n^{(t-1)}$. 

Ke Li et al. \cite{li2019d3m} applied this procedure to model two elliptic equations - a classical Poisson equation and a steady-state Schrödinger equation - considering an overlapping decomposition. They initiated separate PINNs for each subdomain. Then, in each iteration of the training, values of overlapping regions are shared with neighbours. 

This approach is known as the \emph{Deep Domain Decomposition Method} (D3M or DDM). It is worth noting that, in their study, the PDEs loss $\mathcal{E}_f$ is formulated concerning an optimization system of the problem, instead of the residual presented in Section \ref{subsubsec:pinn}.

A nearly identical framework to DDM was developed by Wuyang Li et al. for elliptic problems \cite{li2020deep} with overlaps. They investigated the influence of several factors on the performance of the solver. The accuracy of different PINN structures is compared, including the number of hidden layers, the number of neurons per layer, and the size of the training points. Different sizes of overlap and partitions are also compared. Within the range of their tests, more complex neural networks and data points, as well as more subdomains, yield more accurate ROM approximations. The accuracy increases with the overlapping region, which is already proven for intrusive methods \cite{quarteroni2009numerical}.

Ben Moseley et al. developed a so-called \emph{Finite Basis PINN} (FBPINN) that enables higher flexibility for the training stage. \cite{moseley2023finite}. The workflow of FBPINN can be summarized into three steps. Firstly, multiple small neural networks are instantiated for each subdomain. Each local solution is weighted by a so-called window function, which is smooth, differentiable, and vanishes outside the specific subregion. Since partitions overlap, as well as window functions, the solutions of neighbouring parts are incorporated into the MSE of each subproblem (i.e., Schwarz method). Then, all local PINNs are trained simultaneously. Finally, in the online stage, their predictions are weighted by the window functions to assemble a global solution. 

\begin{figure}[h]
    \centering
    \begin{subfigure}[b]{0.45\linewidth}
        \centering
        \includegraphics[width=\linewidth]{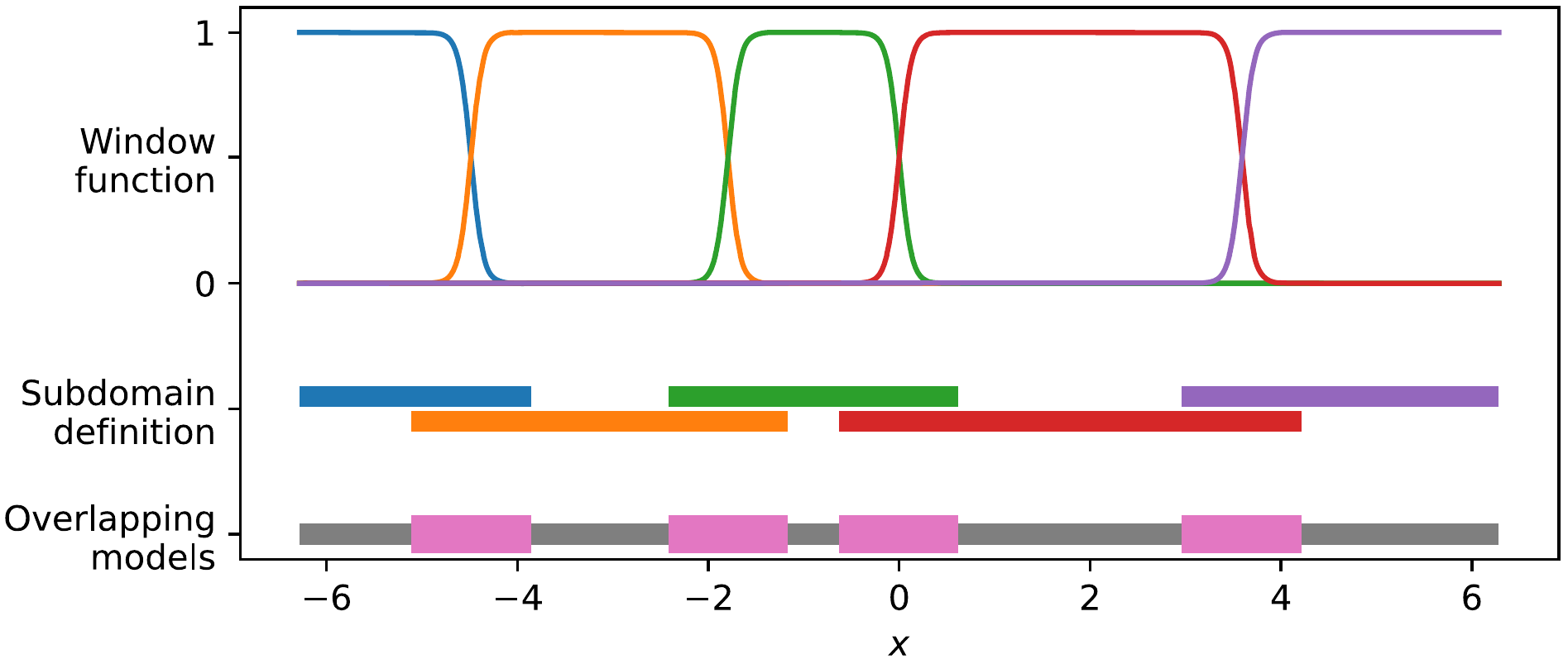}
        \caption{}
    \end{subfigure}
    \hspace{0.5em}
    \begin{subfigure}[b]{0.45\linewidth}
        \centering
        \includegraphics[width=\linewidth]{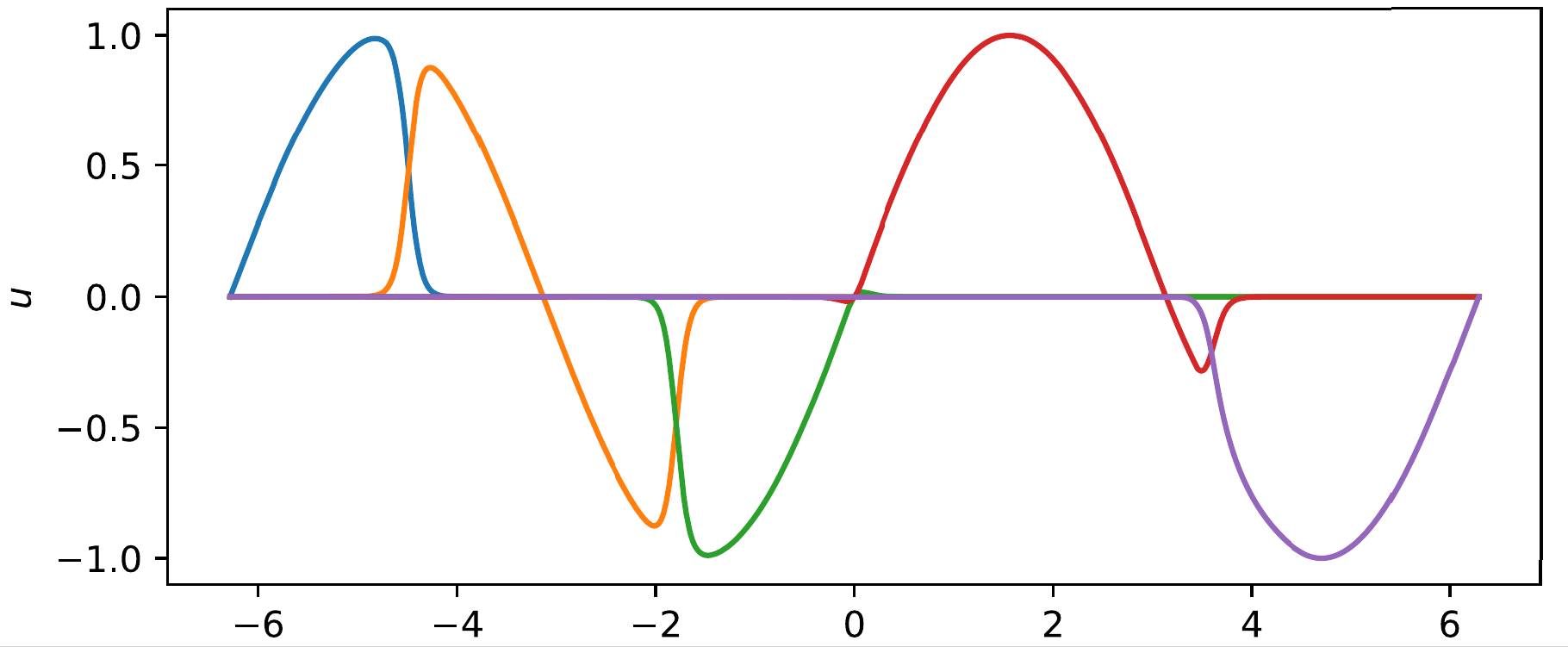}
        \caption{}
    \end{subfigure}
    \caption{One dimensional example of the FBPINN solver. (a) Domain decomposition and window functions for each subdomain. (b) Local solutions. Figures taken from \cite{moseley2023finite}, copyright owned by Springer.}
    \label{fig:FBPINN}
\end{figure}

An 1D example of the FBPINN solver is shown in Fig. \ref{fig:FBPINN}. The entire model is decomposed into overlapping subdomains. Separated window functions are defined for each partition. Local results are weighted by corresponding window functions and summed to form total solutions.

The authors also present a \emph{flexible} training schedule, where only a subset of local FBPINNs is trained at each iteration, rather than all of them. This strategy allocates computational resources by focusing on complex locations while reducing efforts for simpler regions, thereby achieving an effective balance between accuracy and computational efficiency.

The FBPINN is applied to approximate the first- and second-order differential functions, and then utilized to model Burgers' and wave equations. A comparison with the standard PINN (see Section \ref{subsubsec:pinn}) reveals that FBPINN converges faster and results in lower final error. Additionally, each local FBPINN has fewer layers and neurons compared to the global PINN, which can significantly reduce training cost.

A \emph{multilevel FBPINN} (MFBPINN) was proposed by Victorita Dolean et al. \cite{dolean2022finite, hrebenshchykova2023multilevel, dolean2024multilevel}, combining the FBPINN framework \cite{moseley2023finite} with the Schwarz iteration. They weighted all local losses using window functions and united them as a global PINN loss. At each iteration, only a subset of subdomains is activated and trained. Then, the updated values on the overlapping regions are shared with neighbouring subdomains, as in the Schwarz method. The novelty of this approach is that more computational resources are allocated to train portions with higher losses. The strategy achieves a balance between cost and convergence, thereby improving overall computational efficiency. 

\begin{figure}[h]
    \centering
    \includegraphics[width=0.6\linewidth]{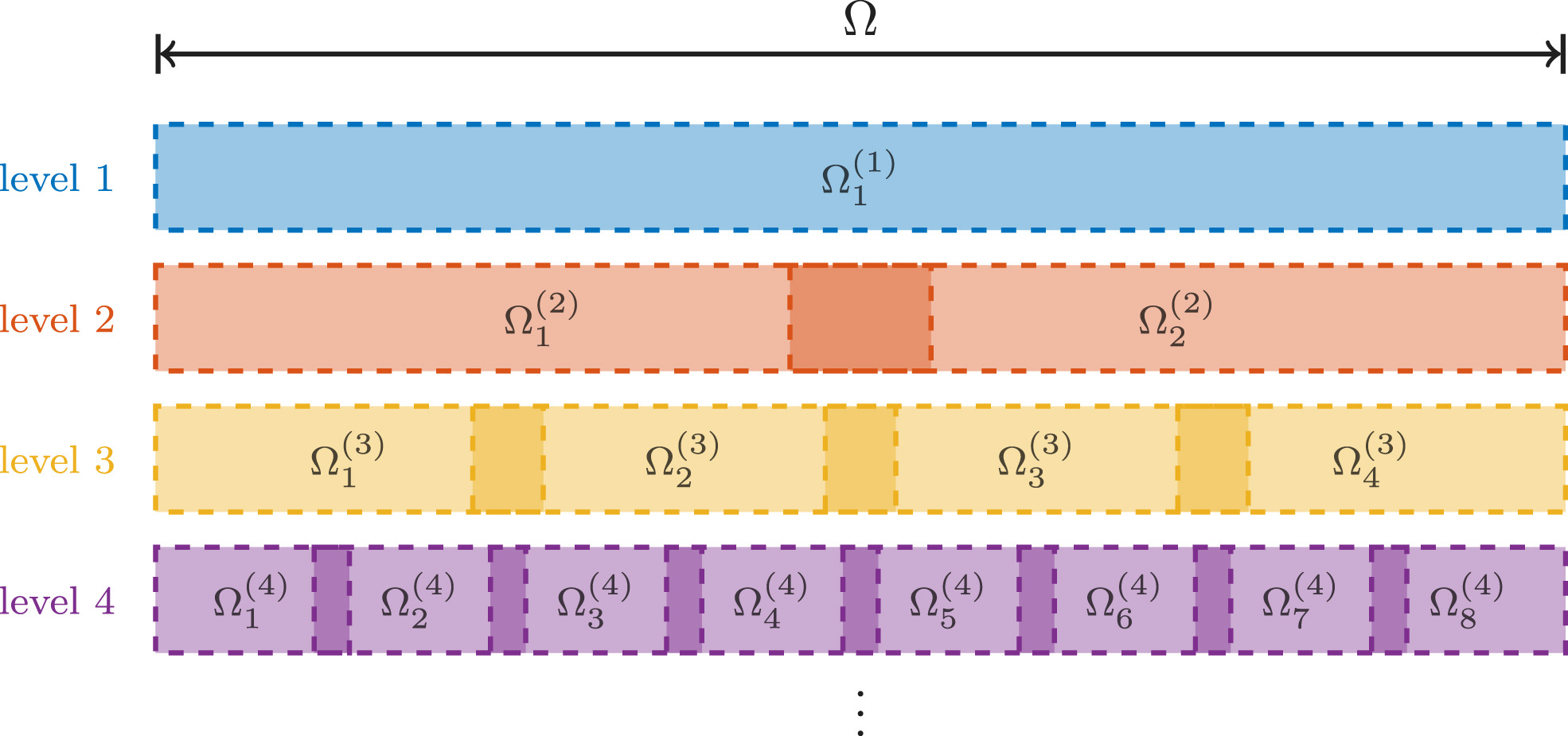}
    \caption{Hierarchy of levels used in the multilevel FBPINN. Figure taken from \cite{dolean2024multilevel}, copyright owned by Elsevier.}
    \label{fig:MFBPINN}
\end{figure}

An innovative aspect of this method compared to previous PINNs is that the geometry is divided into different numbers of subdomains. That is referred to as four decomposition \emph{levels} in the work, as shown in Fig. \ref{fig:MFBPINN}. 

The MSE of partitions at all levels is weighted by window functions and then summed as the global loss function. The authors indicated that employing the configuration can improve the numerical scalability and convergence properties.

The algorithm is applied for the simulation of several problems: 1D ordinary differential equations, Laplacian problems, and the Helmholtz Equation. The results demonstrate that very small overlapping sizes affect the convergence. A MFBPINN with fewer hidden layers and neurons can achieve an error several orders of magnitude lower than the standard PINN. In particular, the multilevel strategy can improve the accuracy.

For the PINN architectures mentioned above, both local and global boundary conditions are constrained via loss functions. William Snyder et al. published an investigation that focuses on parameterizing the Dirichlet BCs \cite{snyder2023domain}. 

In this approach, the solution on a subdomain $\Omega_m$ is represented as
\begin{equation*}
    u'_m (x;\mu) = \gamma_m(x) \boldsymbol{\theta}_m (x;\mu) + \varrho_m (x) g_{D,m} (x),
\end{equation*}
where $\boldsymbol{\theta}_m (x;\mu)$ is the neural network approximating the solution over $\Omega_m$, $\gamma_m(x)$ is a smooth function (can be referred as a test function) that vanishes on $\partial \Omega_m$, $g_{D,m} (x)$ represents the solution at the overlapping regions transferred from neighboring subdomains and $\varrho_m(x)$ serves as the weighting function defined at the overlaps. The authors indicated that this formulation can \emph{strongly enforce} Dirichlet BCs. In previous PINNs, Dirichlet BCs are added as penalty terms in the loss function, which is regarded as a \emph{weak enforcement}. 

The authors compared three configurations: (i) weak imposition of both global BCs and interfaces; (ii) strongly constrained conditions for all Dirichlet boundaries; (iii) mixed approach combining strong enforcement of global BCs and weak imposition at interfaces. The results of a 1D advection-diffusion problem reveal that the mixed scheme (iii) converged faster and provided a more accurate solution compared to the other two strategies. Moreover, they successfully coupled PINN-ROM and FOM within the framework.

\paragraph{Schwarz for coupling PINN ROMs}
We remind you that all PINN + Schwarz architectures reviewed above can only predict solutions for the specific domain in which they are trained. To address this shortcoming, Arthur Feeney et al. \cite{feeney2023breaking} developed a framework that utilizes PINN and the Schwarz algorithm to couple multiple local ROMs in the online stage. 

The approach is implemented as follows. In the offline stage, PINN is constructed to solve \emph{Boundary Value Problems} (BVPs) in a small subdomain with arbitrary Dirichlet BCs. Then, in the online stage, local ROMs are iteratively assembled into a global approximation using the Schwarz method. We remark that the procedure is similar to the \emph{boundary-to-boundary mapping} presented in Section \ref{subsec:Schwarz_based_iterative_scheme}.

\begin{figure}[h]
    \centering
    \begin{subfigure}[b]{0.18\textwidth}
        \centering
        \includegraphics[height=3cm]{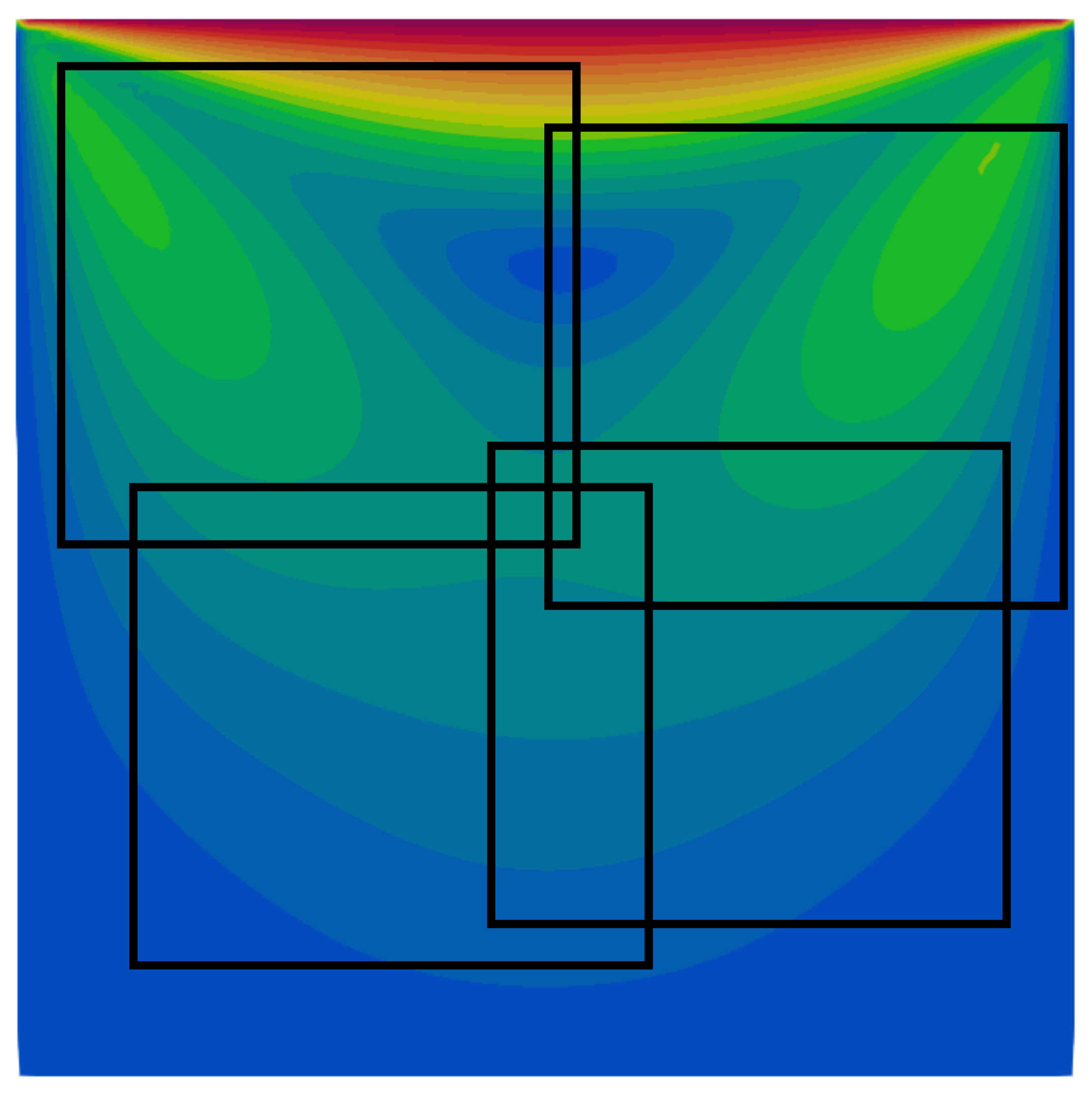}
        \caption{Local snapshots}
        \label{fig:MFP_snapshots}
    \end{subfigure}
    \begin{subfigure}[b]{0.8\textwidth}
        \centering
        \includegraphics[height=3.5cm]{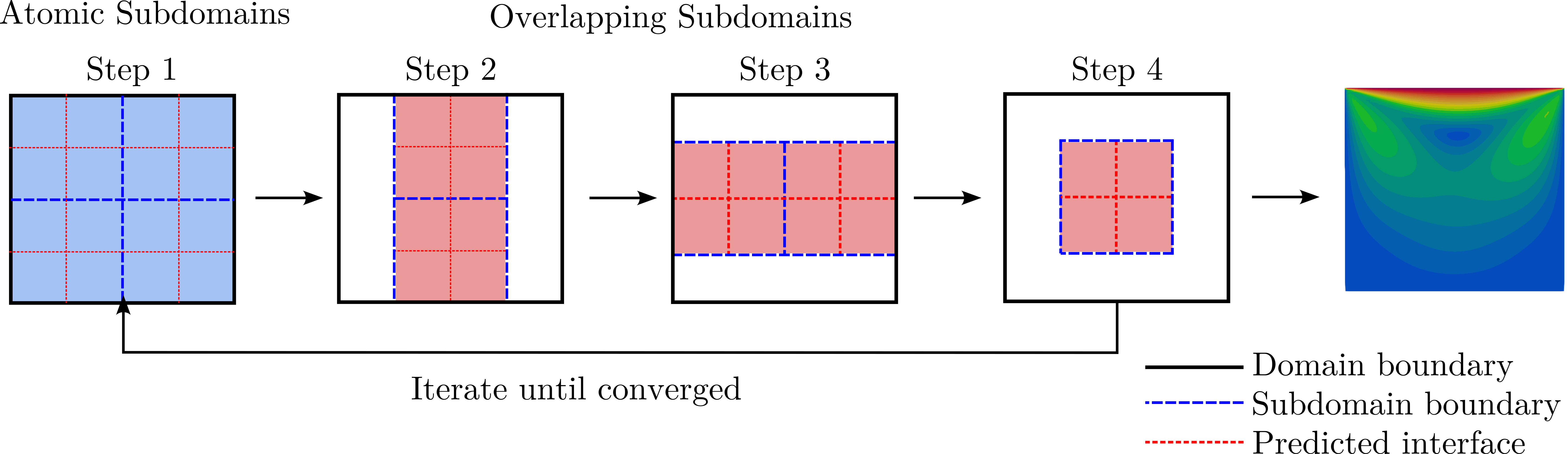}
        \caption{Mosaic Flow predictor}
        \label{fig:Mosaic_Flow_predictor}
    \end{subfigure}
    \caption{A sketch of the \emph{Distributed Mosaic Flow Predictor}. (a) FOM simulations are performed in a slightly larger model, and local snapshots are extracted from the solutions, which is similar to the oversampling technique explained in Section \ref{subsubSec:Localized_training}. (b) \emph{Mosaic Flow Predictor} contains four sub-steps in each iteration. Both atomic and overlapping subdomains overlap with their neighbours, and the solutions are communicated with the Schwarz iteration. Figures redrawn based on \cite{wang2022mosaic}.}
    \label{fig:Dist-MF}
\end{figure}

The authors proposed two tools to achieve the two steps: (i) a PINN-based \emph{SubDomain solver Network} to solve BVPs; (ii) a \emph{Mosaic Flow Predictor} to couple local ROM. The combined architecture is referred to as the \emph{Distributed Mosaic Flow Predictor}, and the two stages are illustrated in Fig. \ref{fig:Dist-MF}. The oversampling strategy (see Section \ref{subsubSec:Localized_training}) is applied to generate parametric high-fidelity solutions. Then, various local snapshots are collected and used for training. To achieve BC parameterization, the input layer of the \emph{SubDomain solver Network} includes not only the spatial and temporal variables ($x$ and $t$) but also BC parameters. 

The \emph{Mosaic Flow Predictor} is referred to as the iterative scheme shown in Fig. \ref{fig:Mosaic_Flow_predictor}. Each iteration contains four sub-steps. The atomic subproblems are solved to obtain a global approximation. Then, the solutions of the so-called \emph{overlapping subdomains} are computed to accelerate the propagation of global boundary information into the interior regions, i.e., speed up the iterative process. The two types of subdomains overlap with their neighbours, and their local BCs are exchanged following the Schwarz iteration.

The authors apply the method to analyze a 2D Laplace equation. Efforts are devoted to testing the accuracy and parallelization performance using different graphics cards. They have proven that it is capable of predicting a very large domain comprising up to 32 $\times$ 32 partitions.

It is worth noting that, in \cite{wang2022mosaic}, a \emph{Genomic Flow Network} based on deep neural networks was created as an alternative to replace \emph{SubDomain solver Network}. The coupling of ROMs was also achieved via \emph{Mosaic Flow Predictor}.

The numerical experiments indicate that the ROM trained on a small piece can be applied to predict significantly larger domains. For Laplace and Navier-Stokes equations, the size of the total domain is 1200 times and 12 times that of a subdomain, respectively. The training of the \emph{Genomic Flow Network} is 1-3 orders of magnitude faster compared to the standard PINN. Moreover, its accuracy is significantly better compared to XPINN.

We highlight that applying the Schwarz algorithm for coupling many subdomains can reduce the cost of training PINNs. In contrast, using the Schwarz iteration for PINN training does not avoid the computational expense over the global domain.

Lastly, we would like to note that PINNs have also become popular for solving multiphysics problems. However, in the papers we noticed, monolithic networks are constructed to predict flow and structural dynamics, which can not be completely considered as local ROMs. As a hint, we list some references for different applications: (i) FSI problems \cite{lee2024parametric, farea2025learning}; (ii) multiphysics in chemical engineering, including Navier-Stokes, energy conservation, mass transport and chemical kinematic equations \cite{wu2023application, ryu2024multiphysics}; (iii) electro-thermal coupling (two Laplace equations) \cite{ma2022preliminary}.


\section{Summary and Conclusions}
\label{sec:summary_conclusions}

In this paper, we have tried to present an overview of the \emph{local ROM methodology}: a challenging and rapidly developing field. We aim to address an informative and brief description, avoiding the complexities of a strict mathematical derivation that can be found anyway in the original literature. In this manner, which may be too shallow for some readers, has helped us create a narrative that makes this paper possible, considering the broad scope of the methods.

Considering the sources consulted, we realize that
\emph{Domain Decomposition} in frames of \emph{Model Order Reduction} techniques, is dominated by a fundamental ideology: minimize the jumps across the subdomain interfaces, meanwhile satisfying the governing PDEs and the boundary conditions.

From our literature review, we may conclude that the available techniques are based on a few concepts. However, such a general overview may be shadowed by an extensive amount of individually readapted techniques that are proposed and consequently modified to fulfill the authors' necessities. In summary, we may deduce that the available methods can be adequately classified into two categories: (i) \emph{intrusive}; (ii) \emph{non-intrusive}. The former involves manipulating the original PDEs to derive reduced-dimensional systems. In contrast, the second group is constructed purely based on FOM solutions. 

We remind that the strategy for domain decomposition and the computation of local RBs are prerequisites of any technique. Thus, their overview is presented first. Then, considering the distinct principles between the two groups of coupling techniques, we will summarize their research status separately.

\subsection{Common aspects} 

We may summarize that the selection of \emph{Overlapping} or \emph{Non-overlapping} of subdomains and the usage of  \emph{Conforming} or \emph{Non-conforming} mesh is a choice incorporated in the specific coupling methodology. The choice depends on the following three criteria: (i) overlapping decomposition has better numerical characteristics; (ii) non-overlapping subdomains appear to be more flexible for more complex and realistic problems. It also reduces the total amount of interfaces; (iii) non-conforming meshes are suitable for problems in which spatial resolution differences are more pronounced. 

Regarding the procedure utilized for the decomposition, we may formulate a couple of conclusions. \emph{Individual decomposition} is a very flexible technique, appealing to a wide range of scenarios. In spite of this, we believe that the other category, the so-called \emph{generic decomposition}, will be more attractive to realistic applications in the next years. This is due to its capabilities to represent geometries with repeating divisions.

Before performing a series of high-fidelity simulations, the parameters should be determined. \emph{Sampling strategies} are critical for the efficient modeling of parametric problems. Refer to the literature, different algorithms can be applied, depending on the dimension of the parameter space. The investigations for a range of physical variables are well developed. Meanwhile, the shape design and optimization are also important in industrial applications. Thus, \emph{geometrical parameterization} is proposed. The interpolation-based techniques and Free-Form Deformation are widely adopted for generating geometrical samples. Additionally, performing FOM simulations becomes a challenge with respect to a high-dimensional parameter space. That can be solved via \emph{parameter space reduction}, which has demonstrated improvements in terms of efficiency and accuracy. The active subspace method is a broadly employed procedure for obtaining low-dimensional synthetic representations of the original space. Recently, the AI-based generative model has become a powerful alternative. Various studies are being carried out to find its potential.

After parameterization, it turns to the stage of \emph{Snapshots collection} and \emph{Local RB construction}. We may conclude from the engineering point of view that: (i) \emph{Localized global RB} may not be attractive in frames of ROM due to their expense in the offline stage; (ii) \emph{Global solutions and local RB} may be applied for some --few-- specific problems due to the ability of local basis vectors to capture small-scale phenomena. Nevertheless, it will show limited growth due to the extensive cost of training. (iii) \emph{Localized training and oversampling technique} is a very promising technique for large-scale facilities due to its inexpensive offline stage cost and capabilities to represent subdomain-level variance via the local RBs. 

However, the advantages of the localized training may be offset by several difficulties: (i) Complexities will grow rapidly when trying to couple numerous subdomains, e.g., thousands of divisions; (ii) The local boundaries (interfaces) must be parametrized to well approximate global models.

\subsection{Coupling techniques}
\paragraph{Projection-based methods}
Now, we turn our attention to the intrusive coupling approaches. They denote reduced algebraic systems explicitly derived from the original PDEs, and then solved numerically. The available intrusive methods can be categorized into two sub-groups (see Table \ref{tab:summary_local_ROMs}): (i) \emph{monolithic}; (ii) \emph{iterative}. 

In turn, these two are characterized by a few conceptions. Monolithic techniques are based on one of the following three ideas: (i) use specific interface basis --test and/or trial-- to enforce continuity; (ii) create a scheme to minimize by penalty or constrain the discontinuity of local basis vectors; (iii) construct a continuous --global-- basis using some strategy to impose a coherent superposition. Even simpler, \textit{iterative} algorithms follow a single conception: create a set of local problems, resolve them iteratively by exchanging information through interfaces, and continue until jump and error thresholds are satisfied.

From our literature review, we have observed that techniques belonging to the monolithic group can be arranged into ten categories, according to their ideology. Iterative techniques can be further classified according to the type of interface condition. Table \ref{tab:summary_local_ROMs} summarizes the properties of each technique according to the classification addressed in Section \ref{sec:classappr}. Readers can easily notice the variety of problems to which each of the methods has been applied. Most of the techniques have been utilized in elliptic and parabolic equations. This is probably because these basic problems are used to test and validate the methods.

Considering the monolithic techniques, we highlight that the method that inherits the concept of RBEM (i.e., RBEM, RDF, SCRBEM, DGRBEM) is more prospective for application. They share two significant advantages: (i) A localized training strategy is incorporated to reduce the offline cost; (ii) Entire models are assembled by instantiations transformed from archetype subdomains, and geometrical parameterization can be exploited during the transformation. 

Moreover, we believe \emph{DGRBEM} is the most promising for engineering problems. This is due to its advantages in the following four aspects: (i) It utilizes localized training, resulting in an efficient offline stage; (ii) Although DG is quite novel in the field of assembling local ROMs, we consider that it is mature. It has been developed for many years in the context of high-fidelity solutions; (iii) It incorporates a more compact RB, including interior and domain interfaces. The discontinuities of the basis vectors are involved in the DG framework. Thus, no extra procedures are needed to obtain separated interface spaces and ensure continuity; (iv) FOM and ROM formulations are almost consistent except RBs. This facilitates the coupling of FOMs and ROMs in complicated, practical cases.

The three categories organized under the item of iterative algorithms are inherently and, in comparison, not so distinct. In spite of this, a preferred method can be decided considering the conditions at the boundary. Dirichlet-Dirichlet iterations generally require overlapping subdomains. This introduces additional difficulties for complex geometries. Therefore, we forecast a trend to avoid Dirichlet-Dirichlet coupling and choose Dirichlet-Neumann or Robin-Robin configurations for engineering applications. We believe that the most important challenge in the near future will be to propose efficient schemes for the communication between sub-problems for assembling numerous local FOMs and ROMs.

Regarding the comparison of monolithic and iterative techniques, we can forecast the following trends. \textit{Monolithic} methods are, in general, more complex and intrusive. This complexity arises when deriving reduced systems with better numerical features. Pitifully, this also results in additional difficulties in its implementation. They must be extensively redeveloped in accordance with the specific demands of the problem. In contrast, \textit{iterative} methods enjoy a more general formulation. They are also more prone to assembling numerous high-resolution systems and ROMs, as well as multi-physics scenarios. Therefore, this may allow us to foretell that the implementation of iterative techniques may occur sooner in commercial codes and industrial applications. In contrast, \textit{monolithic} may remain a long-term path for the simulation of engineering problems. A common trend of both approaches is that they should be further developed to be suitable for a large amount of subdomains.

\paragraph{Non-intrusive methods}
Recently, pure data-driven methodologies have been increasingly utilized to construct ROM, as well as for local ROMs. Based on our observations, we can conclude that the principle of non-intrusive ROMs aims to establish parameter-solution mappings. This framework considers two essential aspects: (i) the algorithms employed for constructing surrogate models; (ii) input and output data. Apart from the parameters (e.g., physical properties and boundary conditions), the interface values between subdomains are also involved accordingly as inputs and/or outputs of the local ROMs.

Considering the procedure of coupling, the non-intrusive local ROMs can be categorized into four groups. A summary regarding their applications is listed in Table \ref{tab:summary_local_ROMs_data_driven}.

There are two subcategories that incorporate interpolations to achieve global approximations. They follow a similar procedure. Independent ROMs are constructed for each partition, and their fields are updated iteratively with inputs from adjacent subdomains. The difference is the quantities exchanged between subproblems. For Schwarz-based approaches, interface values from neighbours are considered as inputs. In contrast, for interpolation-based methods, the solutions of adjacent subdomains are adopted to update each local ROM. 

The optimization-based techniques seem more \emph{monolithic}. They intend to minimize interface discontinuities among all subdomains. The optimal systems are solved using data-driven approaches.

PINN-based ROM frameworks demonstrate greater diversity. Those imposing extra interface constraints can be regarded as \emph{monolithic}, in which networks are built to diminish residuals and interface jumps simultaneously. Alternatively, the Schwarz iteration can be incorporated to train PINNs, computing local residuals along the interfaces. Additionally, PINNs can be used to create boundary value problem solvers (local ROMs) and then couple subdomain-level problems.

Regarding the techniques, we believe the Schwarz-based and iteration approaches are more practical for industrial practices. Our opinion is based on the following three points: (i) They are very flexible in employing different dimensionality reduction algorithms and surrogate models; (ii) They are less expensive than PINNs (still time-consuming nowadays) in the offline stage; (iii) According to the observed studies, they support both overlapping and non-overlapping decompositions. 

Now, let us discuss the PINN framework. It is still not practical for large-scale applications due to the prohibitively high computational cost for training neural networks. However, it contains physical knowledge and demonstrates its potential for solving inverse problems. As the computational hardware is growing recently, it might become popular in the near future. 

\paragraph{Complex systems} 
The coupling of multiple local ROMs is also well-suited for \emph{multiphysics} problems, where each physics can be modeled with an individual ROM. Three types of complex systems are described in the review: (i) the \emph{fluid-structure interaction} is one of the most common scenarios; (ii) \emph{multiphysics} phenomena comprising fluid dynamics, structural mechanics, electronics, neutronics, etc; (iii) \emph{bifurcations} in fluid flows. Especially for FSI problems, two innovative techniques are included: the reduction of mesh motion using interpolation algorithms and the embedded boundary method. 

Both intrusive and non-intrusive frameworks have been employed for simulating those conditions. For multiphysics processes, the variables are generally confined only to specific subregions. Thus, the aspects of generic spatial decomposition, localized training strategy, and RBs defined in reference domains do not always apply to them. Nevertheless, the same algorithms can still be utilised to assemble the ROMs, such as \emph{Lagrange multiplier} and \emph{iterative} schemes. Besides, in some situations, homogeneous conditions are defined at the interfaces, and then the subdomains can be coupled straightforwardly.

\subsection{Conclusions}

According to the literature review and aforementioned summaries, we may conclude that:

\begin{itemize}
    \item Recent developments in local ROMs have indicated a shift from intrusive to non-intrusive methodologies. The former have been developed for around two decades. In contrast, most machine learning-based studies appeared within the last five years. However, due to their advantages and flexibility, they are becoming popular, especially those utilizing neural networks.

    \item Iterative algorithms are well-suited for coupling local ROMs, in both intrusive and non-intrusive frameworks. Since each subproblem can be constructed independently, it is capable of coupling different sub-models, including different ROMs, FOMs with ROMs, and multiphysics. 

    \item Many large-scale engineering devices are assembled by numerous repeating subdomains. Thus, the generic domain decomposition and localized training procedures are very promising. They have the potential to reduce the computational cost during the offline stage and accurately reconstruct global high-fidelity solutions. 
    
    \item Regarding parameterization, various strategies have been investigated for different situations. Geometrical parameterization and shape optimization have also been well studied in the past decades. RBF, IDW, and FFD are feasible tools for accomplishing this task. Recently, advanced parameter space reduction employing AI and neural networks is full of potential to improve the accuracy and efficiency of ROMs.

    \item The local ROM framework has been widely applied to single-physics problems for decades. It is also well-suited to complex phenomena, such as fluid-structure interaction and bifurcations in fluid flow. Additionally, two innovative techniques, e.g., embedded method and interpolation algorithms, have been proposed for the geometric reduction of FSI problems. In the near future, the trend will be toward more challenging \emph{multiphysics} phenomena involving fluid dynamics, structural mechanics, electronics, neutronics, etc.
    
    \item Projection-based approaches contain more comprehensive mathematical foundations, e.g., a priori and a posteriori error estimations, convergence, and well-posedness analyses. They are also more mature in applications, such as geometrical parameterization, FOMs-ROMs coupling, and multiphysics problems.

    \item For non-intrusive methods, the procedure is similar for different scenarios, thus they can be easily adapted to solve various problems. Unlike intrusive techniques, which require manipulation of governing equations, they have to be redeveloped with respect to PDEs and boundary conditions.
    
\end{itemize}

Lastly, we consider it may be adequate to mention that the aforementioned opinions and prognoses are mainly based on an engineering point of view. This might tend to underestimate very important mathematical aspects --well-posedness, error estimations-- emphasizing feasibility or implementation ease. We hope for the indulgence of those readers regarding our conclusions as controversial.

\newpage

\begin{landscape}
\begin{table}[h]
\caption{Intrusive local ROMs. Classifications with respect to domain decomposition, parameterization, RB construction and applications. Abbreviations: \emph{DD} denotes Domain Decomposition, \emph{Optim.} refers to Optimization-based, \emph{Diff.} is different, \emph{Arche.} means archetype, and \emph{CV-Diff.} infers convection-diffusion or advection-diffusion.}
\label{tab:summary_local_ROMs}
\resizebox{\columnwidth}{!}{
\begin{tabular}{lllccccccccccccccc}
\toprule
\multicolumn{3}{l}{\multirow{3}{*}{}}                                         & \multicolumn{11}{c}{Monolithic}                                                                                                                                                                              & \multicolumn{3}{c}{Iterative}                                                                 \\ \cmidrule(r){4-14} \cmidrule(l){15-18}
\multicolumn{3}{l}{}                                                          & \multirow{2}{*}{RBEM} & \multirow{2}{*}{RDF} & \multicolumn{2}{c}{SCRBEM} & \multicolumn{2}{c}{DG-based local ROM} & \multicolumn{2}{c}{GFEM/PUM} & \multicolumn{2}{c}{Optimization} & \multirow{2}{*}{\MP}  & \multirow{2}{*}{D-D}     & \multirow{2}{*}{D-N} & \multirow{2}{*}{R-R} & \multirow{2}{*}{\MP} \\ \cmidrule(r){6-7} \cmidrule(lr){8-9} \cmidrule(lr){10-11} \cmidrule(l){12-13}
\multicolumn{3}{l}{}                                                          &                       &                      & PR           & SC          & DGRBEM            & Local DG           & PUM          & MsFEM         & Optim.           & LSG           &                       &                          &                      &                      &  $\bullet$           \\ \cline{1-18}
\multirow{6}{*}{DD}         & \multirow{2}{*}{Adjacent}   & Overlapping       &                       &                      &              &             &                   &                    & $\bullet$    &               & $\bullet$        & $\bullet$     & $\bullet$             & $\bullet$                &                      &                      &  $\bullet$           \\
                            &                             & Non-overlapping   & $\bullet$             & $\bullet$            & $\bullet$    & $\bullet$   & $\bullet$         & $\bullet$          &              & $\bullet$     &                  &               & $\bullet$             & $\bullet$                & $\bullet$            & $\bullet$            &  $\bullet$           \\
                            & \multirow{2}{*}{Mesh}       & Conforming        & $\bullet$             & $\bullet$            & $\bullet$    & $\bullet$   & $\bullet$         & $\bullet$          & $\bullet$    & $\bullet$     & $\bullet$        & $\bullet$     & $\bullet$             & $\bullet$                & $\bullet$            & $\bullet$            &  $\bullet$           \\
                            &                             & Non-conforming    & $\bullet$             &                      &              &             & $\bullet$         &                    &              &               &                  &               & $\bullet$             & $\bullet$                & $\bullet$            &                      &  $\bullet$           \\
                            & \multirow{2}{*}{Subdomain}  & Individual        &                       &                      &              &             &                   & $\bullet$          &              &               & $\bullet$        & $\bullet$     & $\bullet$             & $\bullet$                & $\bullet$            & $\bullet$            &  $\bullet$           \\
                            &                             & Archetypes        & $\bullet$             & $\bullet$            & $\bullet$    & $\bullet$   & $\bullet$         &                    & $\bullet$    & $\bullet$     & $\bullet$        & $\bullet$     &                       &                          &                      & $\bullet$            &                      \\
\multirow{2}{*}{Parameters} &                             & Physical          & $\bullet$             & $\bullet$            & $\bullet$    & $\bullet$   & $\bullet$         & $\bullet$          & $\bullet$    & $\bullet$     & $\bullet$        & $\bullet$     & $\bullet$             & $\bullet$                & $\bullet$            & $\bullet$            &  $\bullet$           \\
                            &                             & Geometrical       & $\bullet$             & $\bullet$            & $\bullet$    & $\bullet$   & $\bullet$         & $\bullet$          &              &               & $\bullet$        &               & $\bullet$             & $\bullet$                &                      &                      &  $\bullet$           \\
\multirow{4}{*}{\RBsnap}    & \multirow{3}{*}{Global}     & Localized RB      &                       &                      &              &             &                   & $\bullet$          &              &               &                  &               &                       &                          &                      &                      &                      \\
                            &                             & Diff. local RB    &                       &                      &              &             &                   & $\bullet$          &              &               & $\bullet$        & $\bullet$     & $\bullet$             & $\bullet$                & $\bullet$            & $\bullet$            &  $\bullet$           \\
                            &                             & Arche. local RB   &                       &                      &              &             &                   &                    &              &               & $\bullet$        &               &                       &                          &                      &                      &                      \\
                            & Local                       & Arche. RB         & $\bullet$             & $\bullet$            & $\bullet$    & $\bullet$   & $\bullet$         &                    & $\bullet$    & $\bullet$     &                  &               &                       &                          &                      & $\bullet$            &                      \\ \cline{1-18}
\multirow{11}{*}{Equations} & \multirow{2}{*}{Elliptic}   & Second order      & $\bullet$             & $\bullet$            &              & $\bullet$   & $\bullet$         & $\bullet$          & $\bullet$    & $\bullet$     & $\bullet$        & $\bullet$     & $\bullet$             & $\bullet$                & $\bullet$            & $\bullet$            &  $\bullet$           \\
                            &                             & Helmholtz         & $\bullet$             &                      &              & $\bullet$   &                   &                    &              &               &                  &               &                       &                          &                      &                      &                      \\
                            & \multirow{2}{*}{Parabolic}  & Second-order      &                       & $\bullet$            &              &             &                   & $\bullet$          & $\bullet$    & $\bullet$     & $\bullet$        &               & $\bullet$             &                          & $\bullet$            & $\bullet$            &  $\bullet$           \\
                            &                             & CV-Diff.          &                       &                      &              & $\bullet$   & $\bullet$         &                    &              &               &                  &               &                       & $\bullet$                & $\bullet$            & $\bullet$            &                      \\
                            & \multirow{2}{*}{Hyperbolic} & Wave              &                       &                      &              &             &                   &                    &              &               &                  &               &                       & $\bullet$                &                      &                      &                      \\
                            &                             & Transport         &                       &                      &              &             &                   &                    &              & $\bullet$     &                  &               &                       &                          &                      &                      &                      \\                            
                            & Other                       & Maxwell's         & $\bullet$             &                      &              &             & $\bullet$         &                    &              &               &                  &               &                       &                          &                      &                      &                      \\
                            & \multirow{4}{*}{Fluid}      & Burgers'          &                       &                      &              &             &                   &                    &              &               &                  & $\bullet$     &                       & $\bullet$                &                      &                      &                      \\
                            &                             & Euler             &                       &                      &              &             &                   & $\bullet$          &              &               &                  &               & $\bullet$             & $\bullet$                &                      &                      &                      \\
                            &                             & Stokes            & $\bullet$             & $\bullet$            &              &             & $\bullet$         &                    &              & $\bullet$     &                  &               &                       &                          &                      & $\bullet$            &                      \\
                            &                             & Navier stokes     & $\bullet$             &                      & $\bullet$    &             & $\bullet$         &                    & $\bullet$    & $\bullet$     & $\bullet$        &               & $\bullet$             & $\bullet$                &                      &                      & $\bullet$            \\ \cline{1-18}
\multirow{2}{*}{Coupling}   &                             & ROM-ROM           & $\bullet$             & $\bullet$            & $\bullet$    & $\bullet$   & $\bullet$         & $\bullet$          & $\bullet$    & $\bullet$     & $\bullet$        & $\bullet$     & $\bullet$             & $\bullet$                & $\bullet$            & $\bullet$            & $\bullet$            \\ 
                            &                             & FOM-ROM           &                       &                      &              &             &                   & $\bullet$          &              &               &                  &               &                       & $\bullet$                & $\bullet$            &                      &                      \\
\bottomrule
\end{tabular}
}
\end{table}
\end{landscape}

\newpage

\begin{table}[h]
\caption{Non-intrusive local ROMs. Methods, classification of domain decomposition, RB construction and applications. Abbreviations: \emph{DD} denotes Domain Decomposition}
\label{tab:summary_local_ROMs_data_driven}
\begin{tabular}{lllccccc}
\toprule
\multicolumn{3}{l}{}                                                          & \multirow{2}{*}{Schwarz}  & \multirow{2}{*}{Interpolation} & \multirow{2}{*}{Optimization} & \multicolumn{2}{c}{PINN}      \\ \cmidrule(l){7-8}
\multicolumn{3}{l}{}                                                          &                             &                                &                               & Loss term      & Iterative    \\ \cline{1-8}
\multirow{6}{*}{DD}         & \multirow{2}{*}{Adjacent}   & Overlapping       & $\bullet$                   & $\bullet$                      & $\bullet$                     & $\bullet$      & $\bullet$    \\ 
                            &                             & Non-overlapping   & $\bullet$                   & $\bullet$                      & $\bullet$                     & $\bullet$      & $\bullet$    \\
                            & \multirow{2}{*}{Mesh}       & Conforming        & $\bullet$                   & $\bullet$                      & $\bullet$                     &                &              \\
                            &                             & Non-conforming    & $\bullet$                   &                                &                               &                &              \\
                            & \multirow{2}{*}{Subdomain}  & Individual        & $\bullet$                   & $\bullet$                      & $\bullet$                     & $\bullet$      & $\bullet$    \\
                            &                             & Archetypes        & $\bullet$                   & $\bullet$                      &                               &                & $\bullet$    \\
\multirow{2}{*}{Parameters} &                             & Physical          & $\bullet$                   & $\bullet$                      & $\bullet$                     & $\bullet$      & $\bullet$    \\
                            &                             & Geometrical       &                             &                                & $\bullet$                     &                &              \\
\multirow{3}{*}{\RBsnap}    & \multirow{2}{*}{Global}     & Diff. local RB    & $\bullet$                   & $\bullet$                      & $\bullet$                     & $\bullet$      & $\bullet$    \\
                            &                             & Arche. local RB   &                             & $\bullet$                      &                               &                &              \\
                            & Local                       & Arche. RB         & $\bullet$                   &                                &                               &                & $\bullet$    \\ \cline{1-8}
\multirow{7}{*}{Equations}  &                             & Elliptic          & $\bullet$                   & $\bullet$                      &                               & $\bullet$      & $\bullet$    \\
                            &                             & Parabolic         & $\bullet$                   & $\bullet$                      &                               & $\bullet$      & $\bullet$    \\
                            &                             & Hyperbolic        &                             &                                &                               &                &              \\
                            & \multirow{4}{*}{Fluid}      & Burgers'          & $\bullet$                   &                                & $\bullet$                     & $\bullet$      &              \\
                            &                             & Euler             &                             &                                &                               & $\bullet$      &              \\
                            &                             & Stokes            &                             &                                &                               &                &              \\
                            &                             & Navier stokes     & $\bullet$                   & $\bullet$                      & $\bullet$                     & $\bullet$      & $\bullet$    \\ \cline{1-8}
\multirow{2}{*}{Coupling}   &                             & ROM-ROM           & $\bullet$                   & $\bullet$                      & $\bullet$                     & $\bullet$      & $\bullet$    \\ 
                            &                             & FOM-ROM           & $\bullet$                   &                                & $\bullet$                     &                &              \\
\bottomrule
\end{tabular}
\end{table}

\section*{Acknowledgements}
\addcontentsline{toc}{section}{Acknowledgements}
The authors acknowledge Dr. Jorge Yanez (KIT, the Karlsruhe Institute of Technology, Germany) for his help in improving the language used in this paper. They also thank the Institute for Thermal Energy Technology and Safety (ITES) and the Karlsruhe House of Young Scientists (KHYS) at KIT for funding Mr Shenhui Ruan's visit to the International School for Advanced Studies (SISSA), thereby facilitating collaboration.

\bibliographystyle{abbrvnat}
\bibliography{ref}

@article{roelofs2019towards,
    title     = {Towards validated prediction with {RANS CFD} of flow and heat transport in a wire-wrap fuel assembly},
    author    = {Roelofs, Ferry and Uitslag-Doolaard, Heleen and Dovizio, Daniele and Mikuz, Blaz and Shams, Afaque and Bertocchi, Fulvio and Rohde, Martin and Pacio, Julio and Di Piazza, Ivan and Kennedy, Graham and others},
    journal   = {Nuclear Engineering and Design},
    volume    = {353},
    pages     = {110273},
    year      = {2019},
    publisher = {Elsevier}
}

@article{courtessole2024experimental,
    title     = {Experimental investigation of {MHD} flows in a {WCLL TBM} mock-up},
    author    = {Courtessole, C and Brinkmann, H-J and B{\"u}hler, L and Roth, J},
    journal   = {Fusion Engineering and Design},
    volume    = {202},
    pages     = {114306},
    year      = {2024},
    publisher = {Elsevier}
}

@article{chen2017experimental,
    title     = {Experimental investigation of a novel multi-tank thermal energy storage system for solar-powered air conditioning},
    author    = {Chen, Liang and Jin, Sumin and Bu, Guangfeng},
    journal   = {Applied Thermal Engineering},
    volume    = {123},
    pages     = {953--962},
    year      = {2017},
    publisher = {Elsevier}
}

@article{batta2017cfd,
    title     = {{CFD} analysis of pressure drop across grid spacers in rod bundles compared to correlations and heavy liquid metal experimental data},
    author    = {Batta, A and Class, AG},
    journal   = {Nuclear Engineering and Design},
    volume    = {312},
    pages     = {121--127},
    year      = {2017},
    publisher = {Elsevier}
}

@book{moukalled2016finite,
    title     = {The finite volume method},
    author    = {Moukalled, Fadl and Mangani, Luca and Darwish, Marwan and Moukalled, F and Mangani, L and Darwish, M},
    year      = {2016},
    publisher = {Springer}
}

@book{ern2004theory,
    title     = {Theory and practice of finite elements},
    author    = {Ern, Alexandre and Guermond, Jean-Luc},
    volume    = {159},
    year      = {2004},
    publisher = {Springer}
}

@book{canuto2007spectral,
    title     = {Spectral methods: evolution to complex geometries and applications to fluid dynamics},
    author    = {Canuto, Claudio and Hussaini, M Yousuff and Quarteroni, Alfio and Zang, Thomas A},
    year      = {2007},
    publisher = {Springer Science \& Business Media}
}

@article{zhong2024stochastic,
    title     = {A stochastic {Galerkin} lattice {Boltzmann} method for incompressible fluid flows with uncertainties},
    author    = {Zhong, Mingliang and Xiao, Tianbai and Krause, Mathias J and Frank, Martin and Simonis, Stephan},
    journal   = {Journal of Computational Physics},
    volume    = {517},
    pages     = {113344},
    year      = {2024},
    publisher = {Elsevier}
}

@article{batta2024cfd,
    title     = {{CFD} Validation of Forced and Natural Convection for the Open Phase of {IAEA} Benchmark {CRP-I31038}},
    author    = {Batta, Abdalla and Class, Andreas G},
    journal   = {Arabian Journal for Science and Engineering},
    pages     = {1--7},
    year      = {2024},
    publisher = {Springer}
}

@article{buhler2024geometric,
    title     = {Geometric optimization of electrically coupled liquid metal manifolds for {WCLL} blankets},
    author    = {B{\"u}hler, Leo and Mistrangelo, Chiara},
    journal   = {IEEE Transactions on Plasma Science},
    year      = {2024},
    publisher = {IEEE}
}

@book{bergman2011fundamentals,
    title     = {Fundamentals of heat and mass transfer},
    author    = {Bergman, Theodore L},
    year      = {2011},
    publisher = {John Wiley \& Sons}
}

@article{fan2021recent,
    title     = {Recent development of hydrogen and fuel cell technologies: {A} review},
    author    = {Fan, Lixin and Tu, Zhengkai and Chan, Siew Hwa},
    journal   = {Energy Reports},
    volume    = {7},
    pages     = {8421--8446},
    year      = {2021},
    publisher = {Elsevier}
}

@article{roos2021thermocline,
    title     = {Thermocline control through multi-tank thermal-energy storage systems},
    author    = {Roos, Philipp and Haselbacher, Andreas},
    journal   = {Applied Energy},
    volume    = {281},
    pages     = {115971},
    year      = {2021},
    publisher = {Elsevier}
}

@inproceedings{Class2010Coarse,
    author       = {{Class}, Andreas G. and {Viellieber}, Mathias O. and {Himmel}, Steffen R.},
    title        = {{Coarse Grid CFD} for underresolved simulation},
    booktitle    = {APS Division of Fluid Dynamics Meeting Abstracts},
    series       = {APS Meeting Abstracts},
    volume       = {63},
    year         = {2010},
    month        = {nov},
    organization = {}
}

@article{viellieber2015coarse,
    title     = {{Coarse-Grid-CFD} for the Thermal Hydraulic Investigation of Rod-Bundles},
    author    = {Viellieber, Mathias and Class, Andreas},
    journal   = {Pamm},
    volume    = {15},
    number    = {1},
    pages     = {497--498},
    year      = {2015},
    publisher = {Wiley Online Library}
}

@article{class2014two,
    title     = {Two way coupled micro/meso scale method for wind farms},
    author    = {Class, Andreas and Viellieber, Mathias and Moussiopoulos, Nicolas and Barmpas, Fotios},
    journal   = {PAMM},
    volume    = {14},
    number    = {1},
    pages     = {585--586},
    year      = {2014},
    publisher = {Wiley Online Library}
}

@article{du2019thermal,
    title     = {Thermal-hydraulics analysis of flow blockage events for fuel assembly in a sodium-cooled fast reactor},
    author    = {Du, Peng and Shan, Jianqiang and Zhang, Bo and Leung, Laurence KH},
    journal   = {International Journal of Heat and Mass Transfer},
    volume    = {138},
    pages     = {496--507},
    year      = {2019},
    publisher = {Elsevier}
}

@book{todreas2021nuclear2,
    title     = {Nuclear systems Volume {II}: {Elements} of thermal hydraulic design},
    author    = {Todreas, Neil E and Kazimi, Mujid S and Massoud, Mahmoud},
    year      = {2021},
    publisher = {CRC Press}
}

@article{prusak2024time,
    title   = {A time-adaptive algorithm for pressure dominated flows: a heuristic estimator},
    author  = {Prusak, Ivan and Torlo, Davide and Nonino, Monica and Rozza, Gianluigi},
    journal = {arXiv preprint arXiv:2407.00428},
    year    = {2024}
}

@book{rozza2022advanced,
    title     = {Advanced Reduced Order Methods and Applications in Computational Fluid Dynamics},
    author    = {Rozza, Gianluigi and Stabile, Giovanni and Ballarin, Francesco},
    year      = {2022},
    publisher = {SIAM}
}

@book{quarteroni2015reduced,
    title     = {Reduced basis methods for partial differential equations: an introduction},
    author    = {Quarteroni, Alfio and Manzoni, Andrea and Negri, Federico},
    volume    = {92},
    year      = {2015},
    publisher = {Springer}
}

@book{benner2020snapshot,
    title     = {Model Order Reduction: {Volume} 2: {Snapshot}-Based Methods and Algorithms},
    author    = {Benner, Peter and Schilders, Wil and Grivet-Talocia, Stefano and Quarteroni, Alfio and Rozza, Gianluigi and Miguel Silveira, Lu{\'\i}s},
    year      = {2020},
    publisher = {De Gruyter}
}

@book{benner2020applications,
    title     = {Model order reduction: volume 3 applications},
    author    = {Benner, Peter and Schilders, Wil and Grivet-Talocia, Stefano and Quarteroni, Alfio and Rozza, Gianluigi and Miguel Silveira, Lu{\'\i}s},
    year      = {2020},
    publisher = {De Gruyter}
}

@book{quarteroni1999domain,
    title     = {Domain decomposition methods for partial differential equations},
    author    = {Quarteroni, Alfio and Valli, Alberto},
    year      = {1999},
    publisher = {Oxford University Press}
}

@article{amsallem2010towards,
    title   = {Towards real-time computational-fluid-dynamics-based aeroelastic computations using a database of reduced-order information},
    author  = {Amsallem, David and Cortial, Julien and Farhat, Charbel},
    journal = {AIAA journal},
    volume  = {48},
    number  = {9},
    pages   = {2029--2037},
    year    = {2010}
}

@inproceedings{washabaugh2012nonlinear,
    title        = {Nonlinear model reduction for {CFD} problems using local reduced-order bases},
    author       = {Washabaugh, Kyle and Amsallem, David and Zahr, Matthew and Farhat, Charbel},
    booktitle    = {42nd AIAA Fluid Dynamics Conference and Exhibit},
    pages        = {2686},
    year         = {2012},
    organization = {}
}

@article{amsallem2012nonlinear,
    title     = {Nonlinear model order reduction based on local reduced-order bases},
    author    = {Amsallem, David and Zahr, Matthew J and Farhat, Charbel},
    journal   = {International Journal for Numerical Methods in Engineering},
    volume    = {92},
    number    = {10},
    pages     = {891--916},
    year      = {2012},
    publisher = {Wiley Online Library}
}

@article{diaz2024fast,
    title     = {A fast and accurate domain decomposition nonlinear manifold reduced order model},
    author    = {Diaz, Alejandro N and Choi, Youngsoo and Heinkenschloss, Matthias},
    journal   = {Computer Methods in Applied Mechanics and Engineering},
    volume    = {425},
    pages     = {116943},
    year      = {2024},
    publisher = {Elsevier}
}

@article{discacciati2023localized,
    title     = {Localized model order reduction and domain decomposition methods for coupled heterogeneous systems},
    author    = {Discacciati, Niccol{\`o} and Hesthaven, Jan S},
    journal   = {International Journal for Numerical Methods in Engineering},
    volume    = {124},
    number    = {18},
    pages     = {3964--3996},
    year      = {2023},
    publisher = {Wiley Online Library}
}

@article{riffaud2021dgdd,
    title     = {The {DGDD} method for reduced-order modeling of conservation laws},
    author    = {Riffaud, S{\'e}bastien and Bergmann, Michel and Farhat, Charbel and Grimberg, Sebastian and Iollo, Angelo},
    journal   = {Journal of Computational Physics},
    volume    = {437},
    pages     = {110336},
    year      = {2021},
    publisher = {Elsevier}
}

@article{xiao2017domain,
    title   = {A domain decomposition method for the non-intrusive reduced order modelling of fluid flow},
    journal = {Computer Methods in Applied Mechanics and Engineering},
    volume  = {354},
    pages   = {307-330},
    year    = {2019},
    issn    = {0045-7825},
    author  = {D. Xiao and F. Fang and C.E. Heaney and I.M. Navon and C.C. Pain}
}

@incollection{arcucci2020domain,
    title     = {A domain decomposition reduced order model with data assimilation ({DD-RODA})},
    author    = {Arcucci, Rossella and Casas, C{\'e}sar Quilodr{\'a}n and Xiao, Dunhui and Mottet, Laetitia and Fang, Fangxin and Wu, Pin and Pain, Christopher and Guo, Yi-Ke},
    booktitle = {Parallel Computing: {Technology} Trends},
    pages     = {189--198},
    year      = {2020},
    publisher = {IOS Press}
}

@phdthesis{prusak2023application,
    title   = {Application of optimisation-based domain--decomposition reduced order models to parameter-dependent fluid dynamics and multiphysics problems},
    author  = {Prusak, Ivan},
    year    = {2023},
    school  = {SISSA},
    address = {}
}

@article{pegolotti2021model,
    title     = {Model order reduction of flow based on a modular geometrical approximation of blood vessels},
    author    = {Pegolotti, Luca and Pfaller, Martin R and Marsden, Alison L and Deparis, Simone},
    journal   = {Computer methods in applied mechanics and engineering},
    volume    = {380},
    pages     = {113762},
    year      = {2021},
    publisher = {Elsevier}
}

@article{iapichino2016reduced,
    title     = {Reduced basis method and domain decomposition for elliptic problems in networks and complex parametrized geometries},
    author    = {Iapichino, Laura and Quarteroni, Alfio and Rozza, Gianluigi},
    journal   = {Computers \& Mathematics with Applications},
    volume    = {71},
    number    = {1},
    pages     = {408--430},
    year      = {2016},
    publisher = {Elsevier}
}

@article{wang2016layer,
    title     = {Layer pattern thermal design and optimization for multistream plate-fin heat exchangers—A review},
    author    = {Wang, Zhe and Li, Yanzhong},
    journal   = {Renewable and Sustainable Energy Reviews},
    volume    = {53},
    pages     = {500--514},
    year      = {2016},
    publisher = {Elsevier}
}

@article{ruan2024local,
    title     = {Local reduced subspaces of subchannel-inspired subdomains},
    author    = {Ruan, Shenhui and Yanez, Jorge and Class, Andreas G},
    journal   = {International Journal for Numerical Methods in Engineering},
    pages     = {e7552},
    year      = {2024},
    publisher = {Wiley Online Library}
}

@article{prusak2023optimisation,
    title     = {An optimisation--based domain--decomposition reduced order model for the incompressible Navier-Stokes equations},
    author    = {Prusak, Ivan and Nonino, Monica and Torlo, Davide and Ballarin, Francesco and Rozza, Gianluigi},
    journal   = {Computers \& Mathematics with Applications},
    volume    = {151},
    pages     = {172--189},
    year      = {2023},
    publisher = {Elsevier}
}

@phdthesis{sambataro2022component,
    title   = {Component-based model order reduction procedures for large scale THM systems},
    author  = {Sambataro, Giulia},
    year    = {2022},
    school  = {Bordeaux},
    address = {}
}

@article{iollo2023one,
    title     = {A one-shot overlapping Schwarz method for component-based model reduction: application to nonlinear elasticity},
    author    = {Iollo, Angelo and Sambataro, Giulia and Taddei, Tommaso},
    journal   = {Computer Methods in Applied Mechanics and Engineering},
    volume    = {404},
    pages     = {115786},
    year      = {2023},
    publisher = {Elsevier}
}

@incollection{buhr2020localized,
    title     = {Localized model reduction for parameterized problems},
    author    = {Buhr, Andreas and Iapichino, Laura and Ohlberger, Mario and Rave, Stephan and Schindler, Felix and Smetana, Kathrin},
    booktitle = {Handbook on Model Order Reduction},
    year      = {2020},
    publisher = {Walter De Gruyter}
}

@article{bernardi1990new,
    title   = {A new nonconforming approach to domain decomposition: the mortar element method},
    author  = {Bernardi, Christine},
    journal = {Technical Report, Universite Pierre at Marie Curie},
    year    = {1990}
}

@article{belgacem1999mortar,
    title     = {The mortar finite element method with {Lagrange} multipliers},
    author    = {Belgacem, Faker Ben},
    journal   = {Numerische Mathematik},
    volume    = {84},
    pages     = {173--197},
    year      = {1999},
    publisher = {Springer}
}

@article{maday2002reduced,
    title     = {A reduced-basis element method},
    author    = {Maday, Yvon and R{\o}nquist, Einar M},
    journal   = {Journal of scientific computing},
    volume    = {17},
    pages     = {447--459},
    year      = {2002},
    publisher = {Springer}
}

@article{maday2004reduced,
    title     = {The reduced basis element method: application to a thermal fin problem},
    author    = {Maday, Yvon and Ronquist, Einar M},
    journal   = {SIAM Journal on Scientific Computing},
    volume    = {26},
    number    = {1},
    pages     = {240--258},
    year      = {2004},
    publisher = {SIAM}
}

@article{lovgren2006stokes,
    title     = {A reduced basis element method for the steady Stokes problem},
    author    = {L{\o}vgren, Alf Emil and Maday, Yvon and R{\o}nquist, Einar M},
    journal   = {ESAIM: Mathematical Modelling and Numerical Analysis},
    volume    = {40},
    number    = {3},
    pages     = {529--552},
    year      = {2006},
    publisher = {EDP Sciences}
}

@incollection{lovgren2006reduced,
    title     = {The reduced basis element method for fluid flows},
    author    = {L{\o}vgren, Alf Emil and Maday, Yvon and R{\o}nquist, Einar M},
    booktitle = {Analysis and simulation of fluid dynamics},
    pages     = {129--154},
    year      = {2006},
    publisher = {Springer}
}

@book{quarteroni2009numerical,
    title     = {Numerical models for differential problems},
    author    = {Quarteroni, Alfio and Quarteroni, Silvia},
    volume    = {2},
    year      = {2009},
    publisher = {Springer}
}

@book{canuto2006spectral,
    title     = {Spectral methods},
    author    = {Canuto, Claudio and Hussaini, M Youssuff and Quarteroni, Alfio and Zang, Thomas A},
    volume    = {285},
    year      = {2006},
    publisher = {Springer}
}

@book{hesthaven2016certified,
    title     = {Certified reduced basis methods for parametrized partial differential equations},
    author    = {Hesthaven, Jan S and Rozza, Gianluigi and Stamm, Benjamin and others},
    volume    = {590},
    year      = {2016},
    publisher = {Springer}
}

@inproceedings{chen2011seamless,
    title        = {A seamless reduced basis element method for {2D} {Maxwell}'s problem: an introduction},
    author       = {Chen, Yanlai and Hesthaven, Jan S and Maday, Yvon},
    booktitle    = {Spectral and High Order Methods for Partial Differential Equations: {Selected} papers from the ICOSAHOM'09 conference, June 22-26, Trondheim, Norway},
    pages        = {141--152},
    year         = {2011},
    organization = {Springer}
}

@article{iapichino2012reduced,
    title     = {A reduced basis hybrid method for the coupling of parametrized domains represented by fluidic networks},
    author    = {Iapichino, Laura and Quarteroni, Alfio and Rozza, Gianluigi},
    journal   = {Computer Methods in Applied Mechanics and Engineering},
    volume    = {221},
    pages     = {63--82},
    year      = {2012},
    publisher = {Elsevier}
}

@article{eftang2012adaptive,
    title     = {Adaptive port reduction in static condensation},
    author    = {Eftang, JL and Huynh, DBP and Knezevic, DJ and Ronquist, EM and Patera, AT},
    journal   = {IFAC Proceedings Volumes},
    volume    = {45},
    number    = {2},
    pages     = {695--699},
    year      = {2012},
    publisher = {Elsevier}
}

@article{huynh2013static,
    title     = {A static condensation reduced basis element method: approximation and a posteriori error estimation},
    author    = {Huynh, Dinh Bao Phuong and Knezevic, David J and Patera, Anthony T},
    journal   = {ESAIM: Mathematical Modelling and Numerical Analysis},
    volume    = {47},
    number    = {1},
    pages     = {213--251},
    year      = {2013},
    publisher = {EDP Sciences}
}

@article{huynh2013complex,
    title     = {A static condensation reduced basis element method: {Complex} problems},
    author    = {Huynh, Dinh Bao Phuong and Knezevic, David J and Patera, Anthony T},
    journal   = {Computer Methods in Applied Mechanics and Engineering},
    volume    = {259},
    pages     = {197--216},
    year      = {2013},
    publisher = {Elsevier}
}

@article{eftang2013port,
    title     = {Port reduction in parametrized component static condensation: approximation and a posteriori error estimation},
    author    = {Eftang, Jens L and Patera, Anthony T},
    journal   = {International Journal for Numerical Methods in Engineering},
    volume    = {96},
    number    = {5},
    pages     = {269--302},
    year      = {2013},
    publisher = {Wiley Online Library}
}

@article{eftang2014port,
    title     = {A port-reduced static condensation reduced basis element method for large component-synthesized structures: approximation and a posteriori error estimation},
    author    = {Eftang, Jens L and Patera, Anthony T},
    journal   = {Advanced Modeling and Simulation in Engineering Sciences},
    volume    = {1},
    pages     = {1--49},
    year      = {2014},
    publisher = {Springer}
}

@phdthesis{iapichino2012thesis,
    author  = {Iapichino, Laura},
    title   = {Reduced basis methods for the solution of parametrized {PDEs} in repetitive and complex networks with application to {CFD}},
    school  = {EPFL},
    year    = {2012},
    address = {}
}

@article{martini2015reduced,
    title     = {Reduced basis approximation and a-posteriori error estimation for the coupled {Stokes-Darcy} system},
    author    = {Martini, Immanuel and Rozza, Gianluigi and Haasdonk, Bernard},
    journal   = {Advances in Computational Mathematics},
    volume    = {41},
    pages     = {1131--1157},
    year      = {2015},
    publisher = {Springer}
}

@phdthesis{rozza2005shape,
    title   = {Shape design by optimal flow control and reduced basis techniques: {Applications} to bypass configurations in haemodynamics},
    author  = {Rozza, Gianluigi},
    year    = {2005},
    school  = {EPFL},
    address = {}
}

@article{rozza2007stability,
    title     = {On the stability of the reduced basis method for {Stokes} equations in parametrized domains},
    author    = {Rozza, Gianluigi and Veroy, Karen},
    journal   = {Computer methods in applied mechanics and engineering},
    volume    = {196},
    number    = {7},
    pages     = {1244--1260},
    year      = {2007},
    publisher = {Elsevier}
}

@inproceedings{stenger1974convergence,
    title        = {On the convergence and error of the {Bubnov-Galerkin} method},
    author       = {Stenger, Frank},
    booktitle    = {Proceedings of the Conference on the Numerical Solution of Ordinary Differential Equations: 19, 20 October 1972, The University of Texas at Austin},
    pages        = {434--450},
    year         = {1974},
    organization = {Springer}
}

@incollection{lasaint1974finite,
    title     = {On a Finite Element Method for Solving the Neutron Transport Equation},
    author    = {P. Lasaint and P.A. Raviart},
    booktitle = {Mathematical Aspects of Finite Elements in Partial Differential Equations},
    publisher = {Academic Press},
    pages     = {89-123},
    year      = {1974}
}

@book{riviere2008discontinuous,
    title     = {Discontinuous {Galerkin} methods for solving elliptic and parabolic equations: theory and implementation},
    author    = {Rivi{\`e}re, B{\'e}atrice},
    year      = {2008},
    publisher = {SIAM}
}

@book{cockburn2012discontinuous,
    title     = {Discontinuous {Galerkin} methods: theory, computation and applications},
    author    = {Cockburn, Bernardo and Karniadakis, George E and Shu, Chi-Wang},
    volume    = {11},
    year      = {2012},
    publisher = {Springer Science \& Business Media}
}

@article{farhat2009domain,
    title     = {A domain decomposition method for discontinuous {Galerkin} discretizations of {Helmholtz} problems with plane waves and {Lagrange} multipliers},
    author    = {Farhat, Charbel and Tezaur, Radek and Toivanen, Jari},
    journal   = {International journal for numerical methods in engineering},
    volume    = {78},
    number    = {13},
    pages     = {1513--1531},
    year      = {2009},
    publisher = {Wiley Online Library}
}

@article{antonietti2016discontinuous,
    title   = {A discontinuous {Galerkin} reduced basis element method for elliptic problems},
    author  = {Antonietti, Paola F and Pacciarini, Paolo and Quarteroni, Alfio},
    journal = {ESAIM: Mathematical Modelling and Numerical Analysis},
    volume  = {50},
    number  = {2},
    pages   = {337--360},
    year    = {2016}
}

@article{chung2024train,
    title     = {Train small, model big: {Scalable} physics simulators via reduced order modeling and domain decomposition},
    author    = {Chung, Seung Whan and Choi, Youngsoo and Roy, Pratanu and Moore, Thomas and Roy, Thomas and Lin, Tiras Y and Nguyen, Du T and Hahn, Christopher and Duoss, Eric B and Baker, Sarah E},
    journal   = {Computer Methods in Applied Mechanics and Engineering},
    volume    = {427},
    pages     = {117041},
    year      = {2024},
    publisher = {Elsevier}
}

@phdthesis{riffaud2020reduced,
    title   = {Reduced-order models: convergence between scientific computing and data for fluid mechanics},
    author  = {Riffaud, S{\'e}bastien},
    year    = {2020},
    school  = {Universit{\'e} de Bordeaux},
    address = {}
}

@article{ferrero2018global,
    title     = {Global and local {POD} models for the prediction of compressible flows with {DG} methods},
    author    = {Ferrero, Andrea and Iollo, Angelo and Larocca, Francesco},
    journal   = {International Journal for Numerical Methods in Engineering},
    volume    = {116},
    number    = {5},
    pages     = {332--357},
    year      = {2018},
    publisher = {Wiley Online Library}
}

@article{pacciarini2016spectral,
    title     = {Spectral based discontinuous {Galerkin} reduced basis element method for parametrized {Stokes} problems},
    author    = {Pacciarini, Paolo and Gervasio, Paola and Quarteroni, Alfio},
    journal   = {Computers \& Mathematics with Applications},
    volume    = {72},
    number    = {8},
    pages     = {1977--1987},
    year      = {2016},
    publisher = {Elsevier}
}

@article{schleuss2022optimal,
    title     = {Optimal local approximation spaces for parabolic problems},
    author    = {Schleu{\ss}, Julia and Smetana, Kathrin},
    journal   = {Multiscale Modeling \& Simulation},
    volume    = {20},
    number    = {1},
    pages     = {551--582},
    year      = {2022},
    publisher = {SIAM}
}

@article{smetana2023localized,
    title     = {Localized model reduction for nonlinear elliptic partial differential equations: localized training, partition of unity, and adaptive enrichment},
    author    = {Smetana, Kathrin and Taddei, Tommaso},
    journal   = {SIAM Journal on Scientific Computing},
    volume    = {45},
    number    = {3},
    pages     = {A1300--A1331},
    year      = {2023},
    publisher = {SIAM}
}

@article{melenk1996partition,
    title     = {The partition of unity finite element method: basic theory and applications},
    author    = {Melenk, Jens M and Babu{\v{s}}ka, Ivo},
    journal   = {Computer methods in applied mechanics and engineering},
    volume    = {139},
    number    = {1-4},
    pages     = {289--314},
    year      = {1996},
    publisher = {Elsevier}
}

@article{babuvska1997partition,
    title     = {The partition of unity method},
    author    = {Babu{\v{s}}ka, Ivo and Melenk, Jens M},
    journal   = {International journal for numerical methods in engineering},
    volume    = {40},
    number    = {4},
    pages     = {727--758},
    year      = {1997},
    publisher = {Wiley Online Library}
}

@inproceedings{shah2020discontinuous,
    title        = {Discontinuous {Galerkin} model order reduction of geometrically parametrized {Stokes} equation},
    author       = {Shah, Nirav Vasant and Hess, Martin Wilfried and Rozza, Gianluigi},
    booktitle    = {Numerical Mathematics and Advanced Applications {ENUMATH} 2019: {European} Conference, Egmond aan Zee, The Netherlands, September 30-October 4},
    pages        = {551--561},
    year         = {2020},
    organization = {Springer}
}

@article{li2023simulation,
    title     = {Simulation of the interaction of light with {3-D} metallic nanostructures using a proper orthogonal decomposition-{Galerkin} reduced-order discontinuous {Galerkin} time-domain method},
    author    = {Li, Kun and Huang, Ting-Zhu and Li, Liang and Lanteri, St{\'e}phane},
    journal   = {Numerical Methods for Partial Differential Equations},
    volume    = {39},
    number    = {2},
    pages     = {932--954},
    year      = {2023},
    publisher = {Wiley Online Library}
}

@article{babuska2011optimal,
    title     = {Optimal local approximation spaces for generalized finite element methods with application to multiscale problems},
    author    = {Babuska, Ivo and Lipton, Robert},
    journal   = {Multiscale Modeling \& Simulation},
    volume    = {9},
    number    = {1},
    pages     = {373--406},
    year      = {2011},
    publisher = {SIAM}
}

@article{babuvska2020multiscale,
    title     = {Multiscale-spectral {GFEM} and optimal oversampling},
    author    = {Babu{\v{s}}ka, Ivo and Lipton, Robert and Sinz, Paul and Stuebner, Michael},
    journal   = {Computer Methods in Applied Mechanics and Engineering},
    volume    = {364},
    pages     = {112960},
    year      = {2020},
    publisher = {Elsevier}
}

@article{baiges2013domain,
    title     = {A domain decomposition strategy for reduced order models. Application to the incompressible {Navier}--{Stokes} equations},
    author    = {Baiges, Joan and Codina, Ramon and Idelsohn, Sergio},
    journal   = {Computer Methods in Applied Mechanics and Engineering},
    volume    = {267},
    pages     = {23--42},
    year      = {2013},
    publisher = {Elsevier}
}

@article{corigliano2013domain,
    title     = {Domain decomposition and model order reduction methods applied to the simulation of multi-physics problems in {MEMS}},
    author    = {Corigliano, Alberto and Dossi, Martino and Mariani, Stefano},
    journal   = {Computers \& Structures},
    volume    = {122},
    pages     = {113--127},
    year      = {2013},
    publisher = {Elsevier}
}

@article{corigliano2015model,
    title     = {Model order reduction and domain decomposition strategies for the solution of the dynamic elastic--plastic structural problem},
    author    = {Corigliano, Alberto and Dossi, Martino and Mariani, Stefano},
    journal   = {Computer Methods in Applied Mechanics and Engineering},
    volume    = {290},
    pages     = {127--155},
    year      = {2015},
    publisher = {Elsevier}
}

@inproceedings{corigliano2014combined,
    title        = {Combined domain decomposition and model order reduction methods for the solution of coupled and non-linear problems},
    author       = {Corigliano, Alberto and Dossi, Martino and Mariani, Stefano and others},
    booktitle    = {Proceedings of 11th. World Congress on Computational Mechanics (WCCM XI); 5th. European Congress on Computational Mechanics (ECCM V); 6th. European Congress on Computational Fluid Dynamics (ECFD VI)},
    pages        = {4115--4123},
    year         = {2014},
    organization = {}
}

@article{benaceur2022port,
    title  = {Port-reduced reduced-basis component method for steady state {Navier--Stokes} and passive scalar equations},
    author = {Benaceur, Amina and Patera, Anthony},
    year   = {2022}
}

@article{vallaghe2014static,
    title     = {The static condensation reduced basis element method for a mixed-mean conjugate heat exchanger model},
    author    = {Vallagh{\'e}, Sylvain and Patera, Anthony T},
    journal   = {SIAM Journal on Scientific Computing},
    volume    = {36},
    number    = {3},
    pages     = {B294--B320},
    year      = {2014},
    publisher = {SIAM}
}

@article{kaulmann2011new,
    title     = {A new local reduced basis discontinuous {Galerkin} approach for heterogeneous multiscale problems},
    author    = {Kaulmann, Sven and Ohlberger, Mario and Haasdonk, Bernard},
    journal   = {Comptes Rendus Mathematique},
    volume    = {349},
    number    = {23-24},
    pages     = {1233--1238},
    year      = {2011},
    publisher = {Elsevier}
}

@inproceedings{Felix2012localized,
    author       = {Schindler, Felix and Haasdonk, B. and Ohlberger, Mario and Kaulmann, Sven},
    year         = {2012},
    month        = {09},
    pages        = {393--403},
    title        = {The Localized Reduced Basis Multiscale Method},
    volume       = {9},
    journal      = {Algoritmy 2012},
    organization = {}
}

@article{kaulmann2015localized,
    title     = {The localized reduced basis multiscale method for two-phase flows in porous media},
    author    = {Kaulmann, Sven and Flemisch, Bernd and Haasdonk, Bernard and Lie, K-A and Ohlberger, Mario},
    journal   = {International Journal for Numerical Methods in Engineering},
    volume    = {102},
    number    = {5},
    pages     = {1018--1040},
    year      = {2015},
    publisher = {Wiley Online Library}
}

@article{ohlberger2015error,
    title     = {Error control for the localized reduced basis multiscale method with adaptive on-line enrichment},
    author    = {Ohlberger, Mario and Schindler, Felix},
    journal   = {SIAM Journal on Scientific Computing},
    volume    = {37},
    number    = {6},
    pages     = {A2865--A2895},
    year      = {2015},
    publisher = {SIAM}
}

@article{ohlberger2017true,
    title     = {True error control for the localized reduced basis method for parabolic problems},
    author    = {Ohlberger, Mario and Rave, Stephan and Schindler, Felix},
    journal   = {Model reduction of parametrized systems},
    pages     = {169--182},
    year      = {2017},
    publisher = {Springer}
}

@inproceedings{keil2023adaptive,
    title        = {Adaptive Localized Reduced Basis Methods for Large Scale {PDE}-Constrained Optimization},
    author       = {Keil, Tim and Ohlberger, Mario and Schindler, Felix},
    booktitle    = {International Conference on Large-Scale Scientific Computing},
    pages        = {108--116},
    year         = {2023},
    organization = {Springer}
}

@article{prusak2024optimisation,
    title   = {Optimisation--Based Coupling of Finite Element Model and Reduced Order Model for Computational Fluid Dynamics},
    author  = {Prusak, Ivan and Torlo, Davide and Nonino, Monica and Rozza, Gianluigi},
    journal = {arXiv preprint arXiv:2402.10570},
    year    = {2024}
}

@article{taddei2024non,
    title     = {A non-overlapping optimization-based domain decomposition approach to component-based model reduction of incompressible flows},
    author    = {Taddei, Tommaso and Xu, Xuejun and Zhang, Lei},
    journal   = {Journal of Computational Physics},
    volume    = {509},
    pages     = {113038},
    year      = {2024},
    publisher = {Elsevier}
}

@article{hoang2021domain,
    title     = {Domain-decomposition least-squares {Petrov--Galerkin} ({DD-LSPG}) nonlinear model reduction},
    author    = {Hoang, Chi and Choi, Youngsoo and Carlberg, Kevin},
    journal   = {Computer methods in applied mechanics and engineering},
    volume    = {384},
    pages     = {113997},
    year      = {2021},
    publisher = {Elsevier}
}

@article{carlberg2017galerkin,
    title     = {{Galerkin} v. least-squares {Petrov--Galerkin} projection in nonlinear model reduction},
    author    = {Carlberg, Kevin and Barone, Matthew and Antil, Harbir},
    journal   = {Journal of Computational Physics},
    volume    = {330},
    pages     = {693--734},
    year      = {2017},
    publisher = {Elsevier}
}

@article{diaz2023nonlinear,
    title   = {Nonlinear-manifold reduced order models with domain decomposition},
    author  = {Diaz, Alejandro N and Choi, Youngsoo and Heinkenschloss, Matthias},
    journal = {arXiv preprint arXiv:2312.00713},
    year    = {2023}
}

@article{mu2019domain,
    title     = {A domain decomposition model reduction method for linear convection-diffusion equations with random coefficients},
    author    = {Mu, Lin and Zhang, Guannan},
    journal   = {SIAM Journal on Scientific Computing},
    volume    = {41},
    number    = {3},
    pages     = {A1984--A2011},
    year      = {2019},
    publisher = {SIAM}
}

@article{kerfriden2013partitioned,
    title     = {A partitioned model order reduction approach to rationalise computational expenses in nonlinear fracture mechanics},
    author    = {Kerfriden, Pierre and Goury, Olivier and Rabczuk, Timon and Bordas, Stephane Pierre-Alain},
    journal   = {Computer methods in applied mechanics and engineering},
    volume    = {256},
    pages     = {169--188},
    year      = {2013},
    publisher = {Elsevier}
}

@article{buffoni2009iterative,
    title     = {Iterative methods for model reduction by domain decomposition},
    author    = {Buffoni, Marcelo and Telib, Haysam and Iollo, Angelo},
    journal   = {Computers \& Fluids},
    volume    = {38},
    number    = {6},
    pages     = {1160--1167},
    year      = {2009},
    publisher = {Elsevier}
}

@article{cinquegrana2011hybrid,
    title   = {A hybrid method based on {POD} and domain decomposition to compute the {2-D} aerodynamics flow field},
    author  = {Cinquegrana, D and Donelli, RS and Viviani, A},
    journal = {AIMETA, Bologna, Italy},
    pages   = {1--10},
    year    = {2011}
}

@article{barnett2022schwarz,
    title   = {The {Schwarz} alternating method for the seamless coupling of nonlinear reduced order models and full order models},
    author  = {Barnett, Joshua and Tezaur, Irina and Mota, Alejandro},
    journal = {arXiv preprint arXiv:2210.12551},
    year    = {2022}
}

@article{song2024reduced,
    title     = {A reduced-order {Schwarz} domain decomposition method based on {POD} for the convection-diffusion equation},
    author    = {Song, Junpeng and Rui, Hongxing},
    journal   = {Computers \& Mathematics with Applications},
    volume    = {160},
    pages     = {60--69},
    year      = {2024},
    publisher = {Elsevier}
}

@article{wentland2024role,
    title   = {The role of interface boundary conditions and sampling strategies for {Schwarz}-based coupling of projection-based reduced order models},
    author  = {Wentland, Christopher R and Rizzi, Francesco and Barnett, Joshua and Tezaur, Irina},
    journal = {arXiv preprint arXiv:2410.04668},
    year    = {2024}
}

@article{maier2014dirichlet,
    title     = {A {Dirichlet}--{Neumann} reduced basis method for homogeneous domain decomposition problems},
    author    = {Maier, Immanuel and Haasdonk, Bernard},
    journal   = {Applied Numerical Mathematics},
    volume    = {78},
    pages     = {31--48},
    year      = {2014},
    publisher = {Elsevier}
}

@mastersthesis{maier2011iterative,
    title   = {An iterative domain decomposition procedure for the reduced-basis-method},
    author  = {Maier, I and Haasdonk, B},
    year    = {2011},
    school  = {University of Stuttgart},
    address = {}
}

@article{reyes2020element,
    title     = {Element boundary terms in reduced order models for flow problems: {Domain} decomposition and adaptive coarse mesh hyper-reduction},
    author    = {Reyes, Ricardo and Codina, Ramon},
    journal   = {Computer Methods in Applied Mechanics and Engineering},
    volume    = {368},
    pages     = {113159},
    year      = {2020},
    publisher = {Elsevier}
}

@article{zappon2024reduced,
    title     = {A reduced order model for domain decompositions with non-conforming interfaces},
    author    = {Zappon, Elena and Manzoni, Andrea and Gervasio, Paola and Quarteroni, Alfio},
    journal   = {Journal of Scientific Computing},
    volume    = {99},
    number    = {1},
    pages     = {22},
    year      = {2024},
    publisher = {Springer}
}

@article{gervasio2018analysis,
    title     = {Analysis of the {INTERNODES} method for non-conforming discretizations of elliptic equations},
    author    = {Gervasio, Paola and Quarteroni, Alfio},
    journal   = {Computer Methods in Applied Mechanics and Engineering},
    volume    = {334},
    pages     = {138--166},
    year      = {2018},
    publisher = {Elsevier}
}

@article{chaturantabut2010nonlinear,
    title     = {Nonlinear model reduction via discrete empirical interpolation},
    author    = {Chaturantabut, Saifon and Sorensen, Danny C},
    journal   = {SIAM Journal on Scientific Computing},
    volume    = {32},
    number    = {5},
    pages     = {2737--2764},
    year      = {2010},
    publisher = {SIAM}
}

@article{aletti2017reduced,
    title     = {A reduced-order representation of the {Poincar{\'e}--Steklov} operator: an application to coupled multi-physics problems},
    author    = {Aletti, Matteo and Lombardi, Damiano},
    journal   = {International Journal for Numerical Methods in Engineering},
    volume    = {111},
    number    = {6},
    pages     = {581--600},
    year      = {2017},
    publisher = {Wiley Online Library}
}

@article{padula2024brief,
    title   = {A brief review of Reduced Order Models using intrusive and non-intrusive techniques},
    author  = {Padula, Guglielmo and Girfoglio, Michele and Rozza, Gianluigi},
    journal = {arXiv preprint arXiv:2406.00559},
    year    = {2024}
}

@article{bai2021non,
    title     = {Non-intrusive nonlinear model reduction via machine learning approximations to low-dimensional operators},
    author    = {Bai, Zhe and Peng, Liqian},
    journal   = {Advanced Modeling and Simulation in Engineering Sciences},
    volume    = {8},
    number    = {1},
    pages     = {28},
    year      = {2021},
    publisher = {Springer}
}

@book{rozza2024real,
    author    = {Rozza, Gianluigi and Ballarin, Francesco and Scandurra, Leonardo and Pichi, Federico},
    title     = {Real Time Reduced Order Computational Mechanics},
    publisher = {Springer},
    year      = {2024}
}

@article{heaney2022ai,
    title     = {An {AI}-based non-intrusive reduced-order model for extended domains applied to multiphase flow in pipes},
    author    = {Heaney, Claire E and Wolffs, Zef and T{\'o}masson, J{\'o}n Atli and Kahouadji, Lyes and Salinas, Pablo and Nicolle, Andr{\'e} and Navon, Ionel M and Matar, Omar K and Srinil, Narakorn and Pain, Christopher C},
    journal   = {Physics of Fluids},
    volume    = {34},
    number    = {5},
    year      = {2022},
    publisher = {AIP Publishing}
}

@article{discacciati2024overlapping,
    title     = {An overlapping domain decomposition method for the solution of parametric elliptic problems via proper generalized decomposition},
    author    = {Discacciati, Marco and Evans, Ben J and Giacomini, Matteo},
    journal   = {Computer Methods in Applied Mechanics and Engineering},
    volume    = {418},
    pages     = {116484},
    year      = {2024},
    publisher = {Elsevier}
}

@article{fidkowski2020output,
    title     = {Output-based mesh optimization for hybridized and embedded discontinuous {Galerkin} methods},
    author    = {Fidkowski, Krzysztof J and Chen, Guodong},
    journal   = {International Journal for Numerical Methods in Engineering},
    volume    = {121},
    number    = {5},
    pages     = {867--887},
    year      = {2020},
    publisher = {Wiley Online Library}
}

@article{cockburn2004discontinuous,
    title     = {Discontinuous {Galerkin} Methods for Computational Fluid Dynamics},
    author    = {Cockburn, Bernardo},
    journal   = {Encyclopedia of computational mechanics},
    year      = {2004},
    publisher = {Wiley Online Library}
}

@book{hesthaven2007nodal,
    title     = {Nodal discontinuous {Galerkin} methods: algorithms, analysis, and applications},
    author    = {Hesthaven, Jan S and Warburton, Tim},
    year      = {2007},
    publisher = {Springer Science \& Business Media}
}

@article{gunzburger1999optimization,
    title     = {An optimization based domain decomposition method for partial differential equations},
    author    = {Gunzburger, MD and Peterson, JS and Kwon, H},
    journal   = {Computers \& Mathematics with Applications},
    volume    = {37},
    number    = {10},
    pages     = {77--93},
    year      = {1999},
    publisher = {Elsevier}
}

@techreport{tezaur2022schwarz,
    title       = {The {Schwarz} Alternating Method for {ROM-FOM} and {ROM-ROM} Coupling.},
    author      = {Tezaur, Irina and Mota, Alejandro and Shimizu, Yukiko and Barnett, Joshua},
    year        = {2022},
    institution = {Sandia National Lab.(SNL-CA), Livermore, CA (United States); Stanford University},
    address     = {}
}

@article{henning2013oversampling,
    title     = {Oversampling for the multiscale finite element method},
    author    = {Henning, Patrick and Peterseim, Daniel},
    journal   = {Multiscale Modeling \& Simulation},
    volume    = {11},
    number    = {4},
    pages     = {1149--1175},
    year      = {2013},
    publisher = {SIAM}
}

@book{bordas2023partition,
    title     = {Partition of unity methods},
    author    = {Bordas, St{\'e}phane PA and Menk, Alexander and Natarajan, Sundararajan},
    year      = {2023},
    publisher = {John Wiley \& Sons}
}

@article{babuvska1994special,
    title     = {Special finite element methods for a class of second order elliptic problems with rough coefficients},
    author    = {Babu{\v{s}}ka, Ivo and Caloz, Gabriel and Osborn, John E},
    journal   = {SIAM Journal on Numerical Analysis},
    volume    = {31},
    number    = {4},
    pages     = {945--981},
    year      = {1994},
    publisher = {SIAM}
}

@book{efendiev2009multiscale,
    title     = {Multiscale finite element methods: theory and applications},
    author    = {Efendiev, Yalchin and Hou, Thomas Y},
    volume    = {4},
    year      = {2009},
    publisher = {Springer Science \& Business Media}
}

@article{hou1997multiscale,
    title     = {A multiscale finite element method for elliptic problems in composite materials and porous media},
    author    = {Hou, Thomas Y and Wu, Xiao-Hui},
    journal   = {Journal of computational physics},
    volume    = {134},
    number    = {1},
    pages     = {169--189},
    year      = {1997},
    publisher = {Elsevier}
}

@article{efendiev2004Multiscale,
    author  = {Y. Efendiev and T. Y. Hou and V. Ginting},
    title   = {Multiscale Finite Element Methods for Nonlinear Problems and Their Applications},
    journal = {Communications in Mathematical Sciences},
    volume  = {2},
    number  = {4},
    pages   = {553--589},
    year    = {2004}
}

@article{efendiev2007multiscale,
    title     = {Multiscale finite element methods for porous media flows and their applications},
    author    = {Efendiev, Yalchin and Hou, T23224321112},
    journal   = {Applied Numerical Mathematics},
    volume    = {57},
    number    = {5-7},
    pages     = {577--596},
    year      = {2007},
    publisher = {Elsevier}
}

@article{henning2014localizedSIAM,
    title     = {Localized orthogonal decomposition techniques for boundary value problems},
    author    = {Henning, Patrick and M{\aa}lqvist, Axel},
    journal   = {SIAM Journal on Scientific Computing},
    volume    = {36},
    number    = {4},
    pages     = {A1609--A1634},
    year      = {2014},
    publisher = {SIAM}
}

@article{henning2014localizedESAIM,
    title     = {A localized orthogonal decomposition method for semi-linear elliptic problems},
    author    = {Henning, Patrick and M{\aa}lqvist, Axel and Peterseim, Daniel},
    journal   = {ESAIM: Mathematical Modelling and Numerical Analysis},
    volume    = {48},
    number    = {5},
    pages     = {1331--1349},
    year      = {2014},
    publisher = {EDP Sciences}
}

@article{efendiev2014generalized,
    title     = {Generalized multiscale finite element methods: Oversampling strategies},
    author    = {Efendiev, Yalchin and Li, Guanglian and Presho, Michael and others},
    journal   = {International Journal for Multiscale Computational Engineering},
    volume    = {12},
    number    = {6},
    year      = {2014},
    publisher = {Begel House Inc.}
}

@article{efendiev2013generalized,
    title     = {Generalized multiscale finite element methods ({GMsFEM})},
    author    = {Efendiev, Yalchin and Galvis, Juan and Hou, Thomas Y},
    journal   = {Journal of computational physics},
    volume    = {251},
    pages     = {116--135},
    year      = {2013},
    publisher = {Elsevier}
}

@article{chung2014generalized,
    title     = {Generalized multiscale finite element methods for wave propagation in heterogeneous media},
    author    = {Chung, Eric T and Efendiev, Yalchin and Leung, Wing Tat},
    journal   = {Multiscale Modeling \& Simulation},
    volume    = {12},
    number    = {4},
    pages     = {1691--1721},
    year      = {2014},
    publisher = {SIAM}
}

@article{chung2016adaptive,
    title     = {Adaptive multiscale model reduction with generalized multiscale finite element methods},
    author    = {Chung, Eric and Efendiev, Yalchin and Hou, Thomas Y},
    journal   = {Journal of Computational Physics},
    volume    = {320},
    pages     = {69--95},
    year      = {2016},
    publisher = {Elsevier}
}

@article{gander2012best,
    title     = {Best {Robin} parameters for optimized {Schwarz} methods at cross points},
    author    = {Gander, Martin J and Kwok, Felix},
    journal   = {SIAM Journal on Scientific Computing},
    volume    = {34},
    number    = {4},
    pages     = {A1849--A1879},
    year      = {2012},
    publisher = {SIAM}
}

@article{xiao2019domainfluid,
    title     = {A domain decomposition method for the non-intrusive reduced order modelling of fluid flow},
    author    = {Xiao, Dunhui and Fang, Fangxin and Heaney, Claire E and Navon, IM and Pain, CC},
    journal   = {Computer Methods in Applied Mechanics and Engineering},
    volume    = {354},
    pages     = {307--330},
    year      = {2019},
    publisher = {Elsevier}
}

@article{xiao2019domaintur,
    title     = {A domain decomposition non-intrusive reduced order model for turbulent flows},
    author    = {Xiao, D and Heaney, CE and Fang, F and Mottet, L and Hu, R and Bistrian, DA and Aristodemou, E and Navon, IM and Pain, CC},
    journal   = {Computers \& Fluids},
    volume    = {182},
    pages     = {15--27},
    year      = {2019},
    publisher = {Elsevier}
}

@article{xiao2019non,
    title     = {Non-intrusive subdomain {POD-TPWL} for reservoir history matching},
    author    = {Xiao, Cong and Leeuwenburgh, Olwijn and Lin, Hai Xiang and Heemink, Arnold},
    journal   = {Computational Geosciences},
    volume    = {23},
    pages     = {537--565},
    year      = {2019},
    publisher = {Springer}
}

@article{xiao2019subdomain,
    title   = {Subdomain {POD-TPWL} with local parameterization for large-scale reservoir history matching problems},
    author  = {Xiao, Cong and Leeuwenburgh, Olwijn and Lin, Hai Xiang and Heemink, Arnold},
    journal = {arXiv preprint arXiv:1901.08059},
    year    = {2019}
}

@article{xiao2021efficient,
    title     = {Efficient estimation of space varying parameters in numerical models using non-intrusive subdomain reduced order modeling},
    author    = {Xiao, Cong and Leeuwenburgh, Olwijn and Lin, Hai Xiang and Heemink, Arnold},
    journal   = {Journal of Computational Physics},
    volume    = {424},
    pages     = {109867},
    year      = {2021},
    publisher = {Elsevier}
}

@article{wang2022mosaic,
    title     = {Mosaic flows: {A} transferable deep learning framework for solving {PDEs} on unseen domains},
    author    = {Wang, Hengjie and Planas, Robert and Chandramowlishwaran, Aparna and Bostanabad, Ramin},
    journal   = {Computer Methods in Applied Mechanics and Engineering},
    volume    = {389},
    pages     = {114424},
    year      = {2022},
    publisher = {Elsevier}
}

@article{discacciati2024model,
    title     = {Model reduction of coupled systems based on non-intrusive approximations of the boundary response maps},
    author    = {Discacciati, Niccol{\`o} and Hesthaven, Jan S},
    journal   = {Computer Methods in Applied Mechanics and Engineering},
    volume    = {420},
    pages     = {116770},
    year      = {2024},
    publisher = {Elsevier}
}

@inproceedings{iyengar2022nonlinear,
    title     = {Nonlinear reduced order modeling using domain decomposition},
    author    = {Iyengar, Nikhil and Rajaram, Dushhyanth and Decker, Kenneth and Perron, Christian and Mavris, Dimitri N},
    booktitle = {AIAA SciTech 2022 Forum},
    pages     = {1250},
    year      = {2022}
}

@article{chen2021reduced,
    title   = {A reduced order {Schwarz} method for nonlinear multiscale elliptic equations based on two-layer neural networks},
    author  = {Chen, Shi and Ding, Zhiyan and Li, Qin and Wright, Stephen J},
    journal = {arXiv preprint arXiv:2111.02280},
    year    = {2021}
}

@techreport{de2022lagrange,
    title       = {A Lagrange Multiplier Partitioned Scheme for Coupling Reduced Order Models.},
    author      = {de Castro, Amy and Kuberry, Paul and Tezaur, Irina and Bochev, Pavel},
    year        = {2022},
    institution = {Sandia National Laboratories and Stanford University}
}

@techreport{de2023partitioned,
    title       = {A Partitioned Method for Reduced OrderModel-Finite Element Model ({ROM-FEM}) and {ROM-ROM} Couplings with Separate Reduced Bases for Interior and Interface Variables},
    author      = {de Castro, Amy Grace and Kuberry, Paul Allen and Tezaur, Irina Kalashnikova and Bochev, Pavel B},
    year        = {2023},
    institution = {Sandia National Laboratories and Stanford University}
}

@techreport{mota2021schwarz,
    title       = {The Schwarz alternating method for multi-scale coupling and contact in solid mechanics.},
    author      = {Mota, Alejandro and Tezaur, Irina and Hoy, Jonathan},
    year        = {2021},
    institution = {Sandia National Laboratories and University of Southern California}
}

@techreport{koliesnikova2023dirichlet,
    title       = {The {Dirichlet-Neumann} {Schwarz} alternating method for contact problems in elastodynamics: techniques for reducing artificial oscillations},
    author      = {Koliesnikova, Daria and Mota, Alejandro and Tezaur, Irina Kalashnikova},
    year        = {2023},
    institution = {Sandia National Lab.(SNL-CA), Livermore, CA (United States)}
}

@article{mota2023fundamentally,
    title   = {Fundamentally New Coupled Approach to Contact Mechanics via the {Dirichlet-Neumann} {Schwarz} Alternating Method},
    author  = {Mota, Alejandro and Koliesnikova, Daria and Tezaur, Irina and Hoy, Jonathan},
    journal = {arXiv preprint arXiv:2311.05643},
    year    = {2023}
}

@article{romor2025friedrichs,
    title     = {{Friedrichs}' systems discretized with the {DGM}: domain decomposable model order reduction and {Graph Neural Networks} approximating vanishing viscosity solutions},
    author    = {Romor, Francesco and Torlo, Davide and Rozza, Gianluigi},
    journal   = {Journal of Computational Physics},
    pages     = {113915},
    year      = {2025},
    publisher = {Elsevier}
}

@phdthesis{jensen2004discontinuous,
    title  = {Discontinuous Galerkin methods for Friedrichs systems with irregular solutions},
    author = {Jensen, Max},
    year   = {2004},
    school = {University of Oxford}
}

@article{pan2024domain,
    title     = {Domain decomposition for physics-data combined neural network based parametric reduced order modelling},
    author    = {Pan, Xinyu and Xiao, Dunhui},
    journal   = {Journal of Computational Physics},
    volume    = {519},
    pages     = {113452},
    year      = {2024},
    publisher = {Elsevier}
}

@inproceedings{feeney2023breaking,
    title     = {Breaking Boundaries: Distributed Domain Decomposition with Scalable Physics-Informed Neural PDE Solvers},
    author    = {Feeney, Arthur and Li, Zitong and Bostanabad, Ramin and Chandramowlishwaran, Aparna},
    booktitle = {Proceedings of the International Conference for High Performance Computing, Networking, Storage and Analysis},
    pages     = {1--15},
    year      = {2023}
}

@article{snyder2023domain,
    title   = {Domain decomposition-based coupling of physics-informed neural networks via the Schwarz alternating method},
    author  = {Snyder, Will and Tezaur, Irina and Wentland, Christopher},
    journal = {arXiv preprint arXiv:2311.00224},
    year    = {2023}
}

@inproceedings{arcucci2020adaptive,
    title        = {Adaptive domain decomposition for effective data assimilation},
    author       = {Arcucci, Rossella and Mottet, Laetitia and Casas, C{\'e}sar A Quilodr{\'a}n and Guitton, Florian and Pain, Christopher and Guo, Yi-Ke},
    booktitle    = {Euro-Par 2019: Parallel Processing Workshops: Euro-Par 2019 International Workshops, G{\"o}ttingen, Germany, August 26--30, 2019, Revised Selected Papers 25},
    pages        = {583--595},
    year         = {2020},
    organization = {Springer}
}

@article{dwivedi2021distributed,
    title     = {Distributed learning machines for solving forward and inverse problems in partial differential equations},
    author    = {Dwivedi, Vikas and Parashar, Nishant and Srinivasan, Balaji},
    journal   = {Neurocomputing},
    volume    = {420},
    pages     = {299--316},
    year      = {2021},
    publisher = {Elsevier}
}

@article{mienye2024recurrent,
    title     = {Recurrent neural networks: A comprehensive review of architectures, variants, and applications},
    author    = {Mienye, Ibomoiye Domor and Swart, Theo G and Obaido, George},
    journal   = {Information},
    volume    = {15},
    number    = {9},
    pages     = {517},
    year      = {2024},
    publisher = {MDPI}
}

@article{pichi2024graph,
    title     = {A graph convolutional autoencoder approach to model order reduction for parametrized PDEs},
    author    = {Pichi, Federico and Moya, Beatriz and Hesthaven, Jan S},
    journal   = {Journal of Computational Physics},
    volume    = {501},
    pages     = {112762},
    year      = {2024},
    publisher = {Elsevier}
}

@article{chen2025tensor,
    title     = {Tensor decomposition-based neural operator with dynamic mode decomposition for parameterized time-dependent problems},
    author    = {Chen, Yuanhong and Lin, Yifan and Sun, Xiang and Yuan, Chunxin and Gao, Zhen},
    journal   = {Journal of Computational Physics},
    volume    = {533},
    pages     = {113996},
    year      = {2025},
    publisher = {Elsevier}
}

@article{diab2025temporal,
    title   = {Temporal Neural Operator for Modeling Time-Dependent Physical Phenomena},
    author  = {Diab, Waleed and Al-Kobaisi, Mohammed},
    journal = {arXiv preprint arXiv:2504.20249},
    year    = {2025}
}

@article{li2022smart,
    title     = {Smart and rapid design of nanophotonic structures by an adaptive and regularized deep neural network},
    author    = {Li, Renjie and Gu, Xiaozhe and Shen, Yuanwen and Li, Ke and Li, Zhen and Zhang, Zhaoyu},
    journal   = {Nanomaterials},
    volume    = {12},
    number    = {8},
    pages     = {1372},
    year      = {2022},
    publisher = {MDPI}
}

@article{moseley2023finite,
    title     = {Finite basis physics-informed neural networks (FBPINNs): a scalable domain decomposition approach for solving differential equations},
    author    = {Moseley, Ben and Markham, Andrew and Nissen-Meyer, Tarje},
    journal   = {Advances in Computational Mathematics},
    volume    = {49},
    number    = {4},
    pages     = {62},
    year      = {2023},
    publisher = {Springer}
}

@article{jagtap2020extended,
    title     = {Extended physics-informed neural networks (XPINNs): A generalized space-time domain decomposition based deep learning framework for nonlinear partial differential equations},
    author    = {Jagtap, Ameya D and Karniadakis, George Em},
    journal   = {Communications in Computational Physics},
    volume    = {28},
    number    = {5},
    year      = {2020},
    publisher = {Brown Univ., Providence, RI (United States)}
}

@article{trahan2024Quantum,
    title        = {Quantum {{Physics-Informed Neural Networks}}},
    author       = {Trahan, Corey and Loveland, Mark and Dent, Samuel},
    date         = {2024-08},
    journaltitle = {Entropy},
    volume       = {26},
    number       = {8},
    pages        = {649},
    publisher    = {Multidisciplinary Digital Publishing Institute},
    issn         = {1099-4300},
    doi          = {10.3390/e26080649},
    url          = {https://www.mdpi.com/1099-4300/26/8/649},
    urldate      = {2026-01-12},
    langid       = {english}
}

@article{jagtap2020conservative,
    title     = {Conservative physics-informed neural networks on discrete domains for conservation laws: Applications to forward and inverse problems},
    author    = {Jagtap, Ameya D and Kharazmi, Ehsan and Karniadakis, George Em},
    journal   = {Computer Methods in Applied Mechanics and Engineering},
    volume    = {365},
    pages     = {113028},
    year      = {2020},
    publisher = {Elsevier}
}

@article{kharazmi2021hp,
    title     = {hp-VPINNs: Variational physics-informed neural networks with domain decomposition},
    author    = {Kharazmi, Ehsan and Zhang, Zhongqiang and Karniadakis, George Em},
    journal   = {Computer Methods in Applied Mechanics and Engineering},
    volume    = {374},
    pages     = {113547},
    year      = {2021},
    publisher = {Elsevier}
}

@inproceedings{dolean2022finite,
    title        = {Finite basis physics-informed neural networks as a Schwarz domain decomposition method},
    author       = {Dolean, Victorita and Heinlein, Alexander and Mishra, Siddhartha and Moseley, Ben},
    booktitle    = {International Conference on Domain Decomposition Methods},
    pages        = {165--172},
    year         = {2022},
    organization = {Springer}
}

@article{li2019d3m,
    title     = {D3M: A deep domain decomposition method for partial differential equations},
    author    = {Li, Ke and Tang, Kejun and Wu, Tianfan and Liao, Qifeng},
    journal   = {Ieee Access},
    volume    = {8},
    pages     = {5283--5294},
    year      = {2019},
    publisher = {IEEE}
}

@inproceedings{li2020deep,
    title        = {Deep domain decomposition method: Elliptic problems},
    author       = {Li, Wuyang and Xiang, Xueshuang and Xu, Yingxiang},
    booktitle    = {Mathematical and Scientific Machine Learning},
    pages        = {269--286},
    year         = {2020},
    organization = {PMLR}
}

@article{shukla2021parallel,
    title     = {Parallel physics-informed neural networks via domain decomposition},
    author    = {Shukla, Khemraj and Jagtap, Ameya D and Karniadakis, George Em},
    journal   = {Journal of Computational Physics},
    volume    = {447},
    pages     = {110683},
    year      = {2021},
    publisher = {Elsevier}
}

@book{nielsen2015neural,
    title     = {Neural networks and deep learning},
    author    = {Nielsen, Michael A},
    volume    = {25},
    year      = {2015},
    publisher = {Determination press San Francisco, CA, USA}
}

@book{aggarwal2018neural,
    title     = {Neural networks and deep learning},
    author    = {Aggarwal, Charu C and others},
    volume    = {10},
    number    = {978},
    year      = {2018},
    publisher = {Springer}
}

@phdthesis{croft2015proper,
    title  = {Proper generalised decompositions: theory and applications},
    author = {Croft, Thomas Lloyd David},
    year   = {2015},
    school = {Cardiff University}
}

@article{he2015constraint,
    title     = {Constraint reduction procedures for reduced-order subsurface flow models based on POD--TPWL},
    author    = {He, Jincong and Durlofsky, Louis J},
    journal   = {International Journal for Numerical Methods in Engineering},
    volume    = {103},
    number    = {1},
    pages     = {1--30},
    year      = {2015},
    publisher = {Wiley Online Library}
}

@article{iyengar2024domain,
    title     = {Domain decomposition strategy for combining nonlinear and linear reduced-order models},
    author    = {Iyengar, Nikhil and Rajaram, Dushhyanth and Mavris, Dimitri},
    journal   = {AIAA Journal},
    volume    = {62},
    number    = {4},
    pages     = {1375--1389},
    year      = {2024},
    publisher = {American Institute of Aeronautics and Astronautics}
}

@article{buhmann2000radial,
    title     = {Radial basis functions},
    author    = {Buhmann, Martin Dietrich},
    journal   = {Acta numerica},
    volume    = {9},
    pages     = {1--38},
    year      = {2000},
    publisher = {Cambridge university press}
}

@book{williams2006gaussian,
    title     = {Gaussian processes for machine learning},
    author    = {Williams, Christopher KI and Rasmussen, Carl Edward},
    volume    = {2},
    number    = {3},
    year      = {2006},
    publisher = {MIT press Cambridge, MA}
}

@article{cuomo2022scientific,
    title     = {Scientific machine learning through physics--informed neural networks: Where we are and what’s next},
    author    = {Cuomo, Salvatore and Di Cola, Vincenzo Schiano and Giampaolo, Fabio and Rozza, Gianluigi and Raissi, Maziar and Piccialli, Francesco},
    journal   = {Journal of Scientific Computing},
    volume    = {92},
    number    = {3},
    pages     = {88},
    year      = {2022},
    publisher = {Springer}
}

@article{zhao2024comprehensive,
    title     = {A comprehensive review of advances in physics-informed neural networks and their applications in complex fluid dynamics},
    author    = {Zhao, Chi and Zhang, Feifei and Lou, Wenqiang and Wang, Xi and Yang, Jianyong},
    journal   = {Physics of Fluids},
    volume    = {36},
    number    = {10},
    year      = {2024},
    publisher = {AIP Publishing}
}

@article{raissi2017physics,
    title   = {Physics informed deep learning (part i): Data-driven solutions of nonlinear partial differential equations},
    author  = {Raissi, Maziar and Perdikaris, Paris and Karniadakis, George Em},
    journal = {arXiv preprint arXiv:1711.10561},
    year    = {2017}
}

@article{raissi2019physics,
    title     = {Physics-informed neural networks: A deep learning framework for solving forward and inverse problems involving nonlinear partial differential equations},
    author    = {Raissi, Maziar and Perdikaris, Paris and Karniadakis, George E},
    journal   = {Journal of Computational physics},
    volume    = {378},
    pages     = {686--707},
    year      = {2019},
    publisher = {Elsevier}
}

@article{hu2022extended,
    title     = {When Do Extended Physics-Informed Neural Networks (XPINNs) Improve Generalization?},
    author    = {Hu, Zheyuan and Jagtap, Ameya D and Karniadakis, George Em and Kawaguchi, Kenji},
    journal   = {SIAM Journal on Scientific Computing},
    volume    = {44},
    number    = {5},
    pages     = {A3158--A3182},
    year      = {2022},
    publisher = {SIAM}
}

@book{bonaccorso2018machine,
    title     = {Machine Learning Algorithms: Popular algorithms for data science and machine learning},
    author    = {Bonaccorso, Giuseppe},
    year      = {2018},
    publisher = {Packt Publishing Ltd}
}

@book{brunton2024data,
    title     = {Data-driven science and engineering: Machine learning, dynamical systems, and control},
    author    = {Brunton, Steven L and Kutz, J Nathan},
    year      = {2024},
    publisher = {Cambridge University Press}
}

@article{dolean2024multilevel,
    title     = {Multilevel domain decomposition-based architectures for physics-informed neural networks},
    author    = {Dolean, Victorita and Heinlein, Alexander and Mishra, Siddhartha and Moseley, Ben},
    journal   = {Computer Methods in Applied Mechanics and Engineering},
    volume    = {429},
    pages     = {117116},
    year      = {2024},
    publisher = {Elsevier}
}

@article{hrebenshchykova2023multilevel,
    title  = {Multilevel and distributed Physics-Informed Neural Networks for the Helmholtz equation},
    author = {Hrebenshchykova, Daria},
    year   = {2023}
}

@article{heinlein2021combining,
    title     = {Combining machine learning and domain decomposition methods for the solution of partial differential equations—a review},
    author    = {Heinlein, Alexander and Klawonn, Axel and Lanser, Martin and Weber, Janine},
    journal   = {GAMM-Mitteilungen},
    volume    = {44},
    number    = {1},
    pages     = {e202100001},
    year      = {2021},
    publisher = {Wiley Online Library}
}

@article{everson1995karhunen,
    title     = {Karhunen--Loeve procedure for gappy data},
    author    = {Everson, Richard and Sirovich, Lawrence},
    journal   = {Journal of the Optical Society of America A},
    volume    = {12},
    number    = {8},
    pages     = {1657--1664},
    year      = {1995},
    publisher = {Optical Society of America}
}

@book{griewank2008evaluating,
    title     = {Evaluating derivatives: principles and techniques of algorithmic differentiation},
    author    = {Griewank, Andreas and Walther, Andrea},
    year      = {2008},
    publisher = {SIAM}
}

@article{morris2008cfdbased,
    title     = {{CFD}-based optimization of aerofoils using radial basis functions for domain element parameterization and mesh deformation},
    volume    = {58},
    copyright = {http://onlinelibrary.wiley.com/termsAndConditions\#vor},
    issn      = {0271-2091, 1097-0363},
    doi       = {10.1002/fld.1769},
    number    = {8},
    journal   = {International Journal for Numerical Methods in Fluids},
    author    = {Morris, A. M. and Allen, C. B. and Rendall, T. C. S.},
    month     = nov,
    year      = {2008},
    pages     = {827--860}
}

@book{quarteroni2014reduced,
    address   = {Cham},
    title     = {Reduced {Order} {Methods} for {Modeling} and {Computational} {Reduction}},
    copyright = {http://www.springer.com/tdm},
    isbn      = {978-3-319-02089-1 978-3-319-02090-7},
    publisher = {Springer International Publishing},
    editor    = {Quarteroni, Alfio and Rozza, Gianluigi},
    year      = {2014},
    doi       = {10.1007/978-3-319-02090-7}
}

@article{manzoni2012model,
    title     = {Model reduction techniques for fast blood flow simulation in parametrized geometries},
    volume    = {28},
    copyright = {http://onlinelibrary.wiley.com/termsAndConditions\#vor},
    issn      = {2040-7939, 2040-7947},
    doi       = {10.1002/cnm.1465},
    number    = {6-7},
    journal   = {International Journal for Numerical Methods in Biomedical Engineering},
    author    = {Manzoni, Andrea and Quarteroni, Alfio and Rozza, Gianluigi},
    month     = jun,
    year      = {2012},
    pages     = {604--625}
}

@article{manzoni2012shape,
    title     = {Shape optimization for viscous flows by reduced basis methods and free-form deformation},
    volume    = {70},
    copyright = {http://onlinelibrary.wiley.com/termsAndConditions\#vor},
    issn      = {0271-2091, 1097-0363},
    doi       = {10.1002/fld.2712},
    number    = {5},
    journal   = {International Journal for Numerical Methods in Fluids},
    author    = {Manzoni, Andrea and Quarteroni, Alfio and Rozza, Gianluigi},
    month     = oct,
    year      = {2012},
    pages     = {646--670}
}

@article{koshakji2013free,
    title   = {Free {Form} {Deformation} techniques applied to {3D} shape optimization problems},
    volume  = {4},
    issn    = {2038-0909},
    journal = {Communications in Applied and Industrial Mathematics},
    author  = {Koshakji, Anwar and Quarteroni, Alﬁo and Rozza, Gianluigi},
    month   = dec,
    year    = {2013}
}

@article{ballarin2014shape,
    title      = {Shape {Optimization} by {Free}-{Form} {Deformation}: {Existence} {Results} and {Numerical} {Solution} for {Stokes} {Flows}},
    volume     = {60},
    copyright  = {http://www.springer.com/tdm},
    issn       = {0885-7474, 1573-7691},
    shorttitle = {Shape {Optimization} by {Free}-{Form} {Deformation}},
    doi        = {10.1007/s10915-013-9807-8},
    number     = {3},
    journal    = {Journal of Scientific Computing},
    author     = {Ballarin, Francesco and Manzoni, Andrea and Rozza, Gianluigi and Salsa, Sandro},
    month      = sep,
    year       = {2014},
    pages      = {537--563}
}

@article{ballarin2016fast,
    title   = {Fast simulations of patient-specific haemodynamics of coronary artery bypass grafts based on a {POD}-{Galerkin} method and a vascular shape parametrization},
    volume  = {315},
    issn    = {00219991},
    doi     = {10.1016/j.jcp.2016.03.065},
    journal = {Journal of Computational Physics},
    author  = {Ballarin, Francesco and Faggiano, Elena and Ippolito, Sonia and Manzoni, Andrea and Quarteroni, Alfio and Rozza, Gianluigi and Scrofani, Roberto},
    month   = jun,
    year    = {2016},
    pages   = {609--628}
}

@article{ballarin2019podselective,
    title   = {A {POD}-selective inverse distance weighting method for fast parametrized shape morphing},
    volume  = {117},
    issn    = {0029-5981, 1097-0207},
    doi     = {10.1002/nme.5982},
    number  = {8},
    journal = {International Journal for Numerical Methods in Engineering},
    author  = {Ballarin, Francesco and D'Amario, Alessandro and Perotto, Simona and Rozza, Gianluigi},
    month   = feb,
    year    = {2019},
    pages   = {860--884}
}

@phdthesis{forti2012Comparison,
    title      = {Comparison of Shape Parametrization Techniques for Fluid Structure Interaction Problems},
    author     = {Forti, Davide},
    year       = 2012,
    month      = oct,
    urldate    = {2026-01-07},
    langid     = {english},
    school     = {Politecnico di Milano}
}

@phdthesis{damario2014Reducedorder,
    title   = {A Reduced-Order Inverse Distance Weighting Technique for the Efficient Mesh-Motion of Deformable Interfaces and Moving Shapes in Computational Problems},
    author  = {D'Amario, Alessandro},
    year    = 2014,
    urldate = {2026-01-07},
    school  = {Politecnico di Milano}
}

@inproceedings{rozza2012reduction,
    address = {Vienna, Austria},
    title   = {Reduction {Strategies} for {Pde}-{Constrained} {Optimization} {Problems} in {Haemodynamics}},
    author  = {Rozza, Gianluigi and Manzoni, Andrea and Negri, Federico},
    month   = sep,
    year    = {2012}
}

@article{tezzele2021pygem,
    title      = {{PyGeM}: {Python} {Geometrical} {Morphing}},
    volume     = {7},
    issn       = {26659638},
    shorttitle = {{PyGeM}},
    doi        = {10.1016/j.simpa.2020.100047},
    journal    = {Software Impacts},
    author     = {Tezzele, Marco and Demo, Nicola and Mola, Andrea and Rozza, Gianluigi},
    month      = feb,
    year       = {2021},
    pages      = {100047}
}

@article{gadalla19bladex,
    author  = {Gadalla, Mahmoud and Tezzele, Marco and Mola, Andrea and Rozza, Gianluigi},
    title   = {{BladeX: Python Blade Morphing}},
    journal = {The Journal of Open Source Software},
    volume  = {4},
    number  = {34},
    pages   = {1203},
    year    = {2019},
    doi     = {https://doi.org/10.21105/joss.01203}
}

@misc{sissa2021ARGOS,
    title  = {{ARGOS}: Advanced Reduced Order Modeling Online Computational Web Server for Parametric Systems},
    author = {{SISSA MathLab}},
    year   = {2021},
    url    = {http://argos.sissa.it/}
}

@book{rozza2024rbnics,
    author    = {Rozza, Gianluigi and Ballarin, Francesco and Scandurra, Leonardo and Pichi, Federico},
    title     = {Real Time Reduced Order Computational Mechanics},
    subtitle  = {Parametric PDEs Worked Out Problems},
    publisher = {Springer Cham},
    year      = {2024},
    series    = {SISSA Springer Series},
    isbn      = {978-3-031-49891-6}
}

@article{Stabile2017CAF,
    title   = {{Finite volume POD-Galerkin stabilised reduced order methods for the parametrised incompressible Navier-Stokes equations}},
    author  = {Stabile, Giovanni and Rozza, Gianluigi},
    journal = {Computers \& Fluids},
    year    = {2018},
    doi     = {10.1016/j.compfluid.2018.01.035}
}

@article{coscia2023physics,
    title   = {Physics-informed neural networks for advanced modeling},
    author  = {Coscia, Dario and Ivagnes, Anna and Demo, Nicola and Rozza, Gianluigi},
    journal = {Journal of Open Source Software},
    volume  = {8},
    number  = {87},
    pages   = {5352},
    year    = {2023}
}

@article{eiximeno2025pylom,
    title     = {{PyLOM}: a {HPC} open source reduced order model suite for fluid dynamics applications},
    author    = {Eiximeno, Benet and Mir{\'o}, Arnau and Begiashvili, Beka and Valero, Eusebio and Rodriguez, Ivette and Lehmkuhl, Oriol},
    journal   = {Computer Physics Communications},
    volume    = {308},
    pages     = {109459},
    year      = {2025},
    publisher = {Elsevier}
}

@article{milk2016pymor,
    title     = {{pyMOR}-generic algorithms and interfaces for model order reduction},
    author    = {Milk, Ren{\'e} and Rave, Stephan and Schindler, Felix},
    journal   = {SIAM Journal on Scientific Computing},
    volume    = {38},
    number    = {5},
    pages     = {S194--S216},
    year      = {2016},
    publisher = {SIAM}
}

@article{demo18ezyrb,
    author  = {Demo, Nicola and Tezzele, Marco and Rozza, Gianluigi},
    title   = {{EZyRB: Easy Reduced Basis method}},
    journal = {The Journal of Open Source Software},
    volume  = {3},
    number  = {24},
    pages   = {661},
    year    = {2018},
    doi     = {https://doi.org/10.21105/joss.00661}
}

@misc{rozza2018advances,
    title     = {Advances in {Reduced} {Order} {Methods} for {Parametric} {Industrial} {Problems} in {Computational} {Fluid} {Dynamics}},
    doi       = {10.48550/arXiv.1811.08319},
    publisher = {arXiv},
    author    = {Rozza, Gianluigi and Malik, Haris and Demo, Nicola and Tezzele, Marco and Girfoglio, Michele and Stabile, Giovanni and Mola, Andrea},
    month     = nov,
    year      = {2018},
    note      = {arXiv:1811.08319 [math]}
}

@misc{tezzele2019shape,
    title     = {Shape optimization through proper orthogonal decomposition with interpolation and dynamic mode decomposition enhanced by active subspaces},
    doi       = {10.48550/arXiv.1905.05483},
    publisher = {arXiv},
    author    = {Tezzele, Marco and Demo, Nicola and Rozza, Gianluigi},
    month     = may,
    year      = {2019},
    note      = {arXiv:1905.05483 [math]}
}

@misc{mola2019efficient,
    title     = {Efficient {Reduction} in {Shape} {Parameter} {Space} {Dimension} for {Ship} {Propeller} {Blade} {Design}},
    doi       = {10.48550/arXiv.1905.09815},
    publisher = {arXiv},
    author    = {Mola, Andrea and Tezzele, Marco and Gadalla, Mahmoud and Valdenazzi, Federica and Grassi, Davide and Padovan, Roberta and Rozza, Gianluigi},
    month     = may,
    year      = {2019},
    note      = {arXiv:1905.09815 [cs]}
}

@article{demo2019non-intrusive,
    title   = {A non-intrusive approach for the reconstruction of {POD} modal coefficients through active subspaces},
    volume  = {347},
    issn    = {1873-7234},
    doi     = {10.1016/j.crme.2019.11.012},
    number  = {11},
    journal = {Comptes Rendus. Mécanique},
    author  = {Demo, Nicola and Tezzele, Marco and Rozza, Gianluigi},
    month   = nov,
    year    = {2019},
    pages   = {873--881}
}

@inproceedings{salmoiraghi2016advances,
    address    = {Crete Island, Greece},
    title      = {Advances in {Geometrical} {Parametrization} and {Reduced} {Order} {Models} and {Methods} for {Computational} {Fluid} {Dynamics} {Problems} in {Applied} {Sciences} and {Engineering}: {Overview} and {Perspectives}},
    isbn       = {978-618-82844-0-1},
    shorttitle = {{ADVANCES} {IN} {GEOMETRICAL} {PARAMETRIZATION} {AND} {REDUCED} {ORDER} {MODELS} {AND} {METHODS} {FOR} {COMPUTATIONAL} {FLUID} {DYNAMICS} {PROBLEMS} {IN} {APPLIED} {SCIENCES} {AND} {ENGINEERING}},
    doi        = {10.7712/100016.1867.8680},
    booktitle  = {Proceedings of the {VII} {European} {Congress} on {Computational} {Methods} in {Applied} {Sciences} and {Engineering} ({ECCOMAS} {Congress} 2016)},
    publisher  = {Institute of Structural Analysis and Antiseismic Research School of Civil Engineering National Technical University of Athens (NTUA) Greece},
    author     = {Salmoiraghi, F. and Ballarin, F. and Corsi, G. and Mola, A. and Tezzele, M. and Rozza, G.},
    year       = {2016},
    pages      = {1013--1031}
}

@article{demo2021hull,
    title     = {Hull {Shape} {Design} {Optimization} with {Parameter} {Space} and {Model} {Reductions}, and {Self}-{Learning} {Mesh} {Morphing}},
    volume    = {9},
    copyright = {https://creativecommons.org/licenses/by/4.0/},
    issn      = {2077-1312},
    doi       = {10.3390/jmse9020185},
    number    = {2},
    journal   = {Journal of Marine Science and Engineering},
    author    = {Demo, Nicola and Tezzele, Marco and Mola, Andrea and Rozza, Gianluigi},
    month     = feb,
    year      = {2021},
    pages     = {185}
}

@article{romor2024local,
    title   = {A {Local} {Approach} to {Parameter} {Space} {Reduction} for {Regression} and {Classification} {Tasks}},
    volume  = {99},
    issn    = {0885-7474, 1573-7691},
    doi     = {10.1007/s10915-024-02542-0},
    number  = {3},
    journal = {Journal of Scientific Computing},
    author  = {Romor, Francesco and Tezzele, Marco and Rozza, Gianluigi},
    month   = jun,
    year    = {2024},
    pages   = {83}
}

@article{tezzele2018dimension,
    title   = {Dimension reduction in heterogeneous parametric spaces with application to naval engineering shape design problems},
    volume  = {5},
    issn    = {2213-7467},
    doi     = {10.1186/s40323-018-0118-3},
    number  = {1},
    journal = {Advanced Modeling and Simulation in Engineering Sciences},
    author  = {Tezzele, Marco and Salmoiraghi, Filippo and Mola, Andrea and Rozza, Gianluigi},
    month   = dec,
    year    = {2018},
    pages   = {25}
}

@article{tezzele2020enhancing,
    title   = {Enhancing {CFD} predictions in shape design problems by model and parameter space reduction},
    volume  = {7},
    issn    = {2213-7467},
    doi     = {10.1186/s40323-020-00177-y},
    number  = {1},
    journal = {Advanced Modeling and Simulation in Engineering Sciences},
    author  = {Tezzele, Marco and Demo, Nicola and Stabile, Giovanni and Mola, Andrea and Rozza, Gianluigi},
    month   = dec,
    year    = {2020},
    pages   = {40}
}

@misc{siena2024accuracy,
    title      = {On the accuracy and efficiency of reduced order models: towards real-world applications},
    shorttitle = {On the accuracy and efficiency of reduced order models},
    doi        = {10.48550/arXiv.2407.03325},
    publisher  = {arXiv},
    author     = {Siena, Pierfrancesco and Africa, Paquale Claudio and Girfoglio, Michele and Rozza, Gianluigi},
    month      = sep,
    year       = {2024},
    note       = {arXiv:2407.03325 [math]}
}

@article{ivagnes2024shape,
    title   = {A shape optimization pipeline for marine propellers by means of reduced order modeling techniques},
    volume  = {125},
    issn    = {0029-5981, 1097-0207},
    doi     = {10.1002/nme.7426},
    number  = {7},
    journal = {International Journal for Numerical Methods in Engineering},
    author  = {Ivagnes, Anna and Demo, Nicola and Rozza, Gianluigi},
    month   = apr,
    year    = {2024},
    pages   = {e7426}
}

@article{padula2024generative,
    title      = {Generative models for the deformation of industrial shapes with linear geometric constraints: {Model} order and parameter space reductions},
    volume     = {423},
    issn       = {00457825},
    shorttitle = {Generative models for the deformation of industrial shapes with linear geometric constraints},
    doi        = {10.1016/j.cma.2024.116823},
    journal    = {Computer Methods in Applied Mechanics and Engineering},
    author     = {Padula, Guglielmo and Romor, Francesco and Stabile, Giovanni and Rozza, Gianluigi},
    month      = apr,
    year       = {2024},
    pages      = {116823}
}

@misc{padula2025generative,
    title     = {Generative {Models} for {Parameter} {Space} {Reduction} applied to {Reduced} {Order} {Modelling}},
    doi       = {10.48550/arXiv.2506.09721},
    publisher = {arXiv},
    author    = {Padula, Guglielmo and Rozza, Gianluigi},
    month     = jun,
    year      = {2025},
    note      = {arXiv:2506.09721 [math]}
}

@misc{padula2023generative,
    title      = {Generative {Models} for the {Deformation} of {Industrial} {Shapes} with {Linear} {Geometric} {Constraints}: model order and parameter space reductions},
    shorttitle = {Generative {Models} for the {Deformation} of {Industrial} {Shapes} with {Linear} {Geometric} {Constraints}},
    doi        = {10.48550/arXiv.2308.03662},
    publisher  = {arXiv},
    author     = {Padula, Guglielmo and Romor, Francesco and Stabile, Giovanni and Rozza, Gianluigi},
    month      = aug,
    year       = {2023},
    note       = {arXiv:2308.03662 [math]}
}

@article{karcher2022adaptive,
    title   = {Adaptive sampling strategies for reduced-order modeling},
    volume  = {13},
    issn    = {1869-5582, 1869-5590},
    doi     = {10.1007/s13272-022-00574-6},
    number  = {2},
    journal = {CEAS Aeronautical Journal},
    author  = {Karcher, Niklas and Franz, Thomas},
    month   = apr,
    year    = {2022},
    pages   = {487--502}
}

@article{bui-thanh2008model,
    title   = {Model {Reduction} for {Large}-{Scale} {Systems} with {High}-{Dimensional} {Parametric} {Input} {Space}},
    volume  = {30},
    issn    = {1064-8275, 1095-7197},
    doi     = {10.1137/070694855},
    number  = {6},
    journal = {SIAM Journal on Scientific Computing},
    author  = {Bui-Thanh, T. and Willcox, K. and Ghattas, O.},
    month   = jan,
    year    = {2008},
    pages   = {3270--3288}
}

@article{liu2024application,
    title   = {Application and comparison of several adaptive sampling algorithms in reduced order modeling},
    volume  = {10},
    issn    = {24058440},
    doi     = {10.1016/j.heliyon.2024.e34928},
    number  = {15},
    journal = {Heliyon},
    author  = {Liu, Xirui and Wang, Zhiyong and Ji, Hongjun and Gong, Helin},
    month   = aug,
    year    = {2024},
    pages   = {e34928}
}

@book{thompson2012sampling,
    title     = {Sampling},
    author    = {Thompson, Steven K},
    volume    = {755},
    year      = {2012},
    publisher = {John Wiley \& Sons}
}

@article{murthy1967sampling,
    title  = {Sampling theory and methods.},
    author = {Murthy, Mankal Narasinha},
    year   = {1967}
}

@article{benner2015survey,
    title     = {A survey of projection-based model reduction methods for parametric dynamical systems},
    author    = {Benner, Peter and Gugercin, Serkan and Willcox, Karen},
    journal   = {SIAM review},
    volume    = {57},
    number    = {4},
    pages     = {483--531},
    year      = {2015},
    publisher = {SIAM}
}

@article{clenshaw1960method,
    title     = {A method for numerical integration on an automatic computer},
    author    = {Clenshaw, Charles W and Curtis, Alan R},
    journal   = {Numerische Mathematik},
    volume    = {2},
    pages     = {197--205},
    year      = {1960},
    publisher = {Springer}
}

@article{mckay2000comparison,
    title     = {A comparison of three methods for selecting values of input variables in the analysis of output from a computer code},
    author    = {McKay, Michael D and Beckman, Richard J and Conover, William J},
    journal   = {Technometrics},
    volume    = {42},
    number    = {1},
    pages     = {55--61},
    year      = {2000},
    publisher = {Taylor \& Francis}
}

@article{ju2002probabilistic,
    title     = {Probabilistic methods for centroidal Voronoi tessellations and their parallel implementations},
    author    = {Ju, Lili and Du, Qiang and Gunzburger, Max},
    journal   = {Parallel Computing},
    volume    = {28},
    number    = {10},
    pages     = {1477--1500},
    year      = {2002},
    publisher = {Elsevier}
}

@article{veroy2005certified,
    title     = {Certified real-time solution of the parametrized steady incompressible Navier--Stokes equations: rigorous reduced-basis a posteriori error bounds},
    author    = {Veroy, Karen and Patera, Anthony T},
    journal   = {International Journal for Numerical Methods in Fluids},
    volume    = {47},
    number    = {8-9},
    pages     = {773--788},
    year      = {2005},
    publisher = {Wiley Online Library}
}

@inproceedings{veroy2003posteriori,
    title     = {A posteriori error bounds for reduced-basis approximation of parametrized noncoercive and nonlinear elliptic partial differential equations},
    author    = {Veroy, Karen and Prud'Homme, Christophe and Rovas, Dimitrios and Patera, Anthony},
    booktitle = {16th AIAA computational fluid dynamics conference},
    pages     = {3847},
    year      = {2003}
}

@article{rozza2008reduced,
    title     = {Reduced basis approximation and a posteriori error estimation for affinely parametrized elliptic coercive partial differential equations: application to transport and continuum mechanics},
    author    = {Rozza, Gianluigi and Huynh, Dinh Bao Phuong and Patera, Anthony T},
    journal   = {Archives of Computational Methods in Engineering},
    volume    = {15},
    number    = {3},
    pages     = {229--275},
    year      = {2008},
    publisher = {Springer}
}

@article{grepl2005posteriori,
    title     = {A posteriori error bounds for reduced-basis approximations of parametrized parabolic partial differential equations},
    author    = {Grepl, Martin A and Patera, Anthony T},
    journal   = {ESAIM: Mathematical Modelling and Numerical Analysis},
    volume    = {39},
    number    = {1},
    pages     = {157--181},
    year      = {2005},
    publisher = {EDP Sciences}
}

@article{negri2015reduced,
    title     = {Reduced basis approximation of parametrized optimal flow control problems for the Stokes equations},
    author    = {Negri, Federico and Manzoni, Andrea and Rozza, Gianluigi},
    journal   = {Computers \& Mathematics with Applications},
    volume    = {69},
    number    = {4},
    pages     = {319--336},
    year      = {2015},
    publisher = {Elsevier}
}

@article{chen2018greedy,
    title     = {Greedy nonintrusive reduced order model for fluid dynamics},
    author    = {Chen, Wang and Hesthaven, Jan S and Junqiang, Bai and Qiu, Yasong and Yang, Zhang and Tihao, Yang},
    journal   = {AIAA Journal},
    volume    = {56},
    number    = {12},
    pages     = {4927--4943},
    year      = {2018},
    publisher = {American Institute of Aeronautics and Astronautics}
}

@article{sleeman2025greedy,
    title   = {Greedy recursion parameter selection for the One-Way Navier-Stokes (OWNS) equations},
    author  = {Sleeman, Michael K and Colonius, Tim},
    journal = {arXiv preprint arXiv:2506.02320},
    year    = {2025}
}

@article{bond2007piecewise,
    title     = {A piecewise-linear moment-matching approach to parameterized model-order reduction for highly nonlinear systems},
    author    = {Bond, Bradley N and Daniel, Luca},
    journal   = {IEEE Transactions on Computer-Aided Design of Integrated Circuits and Systems},
    volume    = {26},
    number    = {12},
    pages     = {2116--2129},
    year      = {2007},
    publisher = {IEEE}
}

@phdthesis{gong2018data,
    title  = {Data assimilation with reduced basis and noisy measurement: Applications to nuclear reactor cores},
    author = {Gong, Helin},
    year   = {2018},
    school = {Sorbonne Universit{\'e}}
}

@article{fahl2000trust,
    title  = {Trust-region methods for flow control based on reduced order modelling},
    author = {Fahl, Marco},
    year   = {2000}
}

@article{bressloff2007parametric,
    title     = {Parametric geometry exploration of the human carotid artery bifurcation},
    author    = {Bressloff, Neil W},
    journal   = {Journal of Biomechanics},
    volume    = {40},
    number    = {11},
    pages     = {2483--2491},
    year      = {2007},
    publisher = {Elsevier}
}

@article{lee2008geometry,
    title     = {Geometry of the carotid bifurcation predicts its exposure to disturbed flow},
    author    = {Lee, Sang-Wook and Antiga, Luca and Spence, J David and Steinman, David A},
    journal   = {Stroke},
    volume    = {39},
    number    = {8},
    pages     = {2341--2347},
    year      = {2008},
    publisher = {Lippincott Williams \& Wilkins}
}

@article{ding2001flow,
    title     = {Flow field and oscillatory shear stress in a tuning-fork-shaped model of the average human carotid bifurcation},
    author    = {Ding, Zurong and Wang, Keqiang and Li, Jie and Cong, Xushen},
    journal   = {Journal of Biomechanics},
    volume    = {34},
    number    = {12},
    pages     = {1555--1562},
    year      = {2001},
    publisher = {Elsevier}
}

@article{schaback2007practical,
    title   = {A practical guide to radial basis functions},
    author  = {Schaback, Robert},
    journal = {Electronic Resource},
    volume  = {11},
    pages   = {1--12},
    year    = {2007}
}

@inproceedings{sederberg1986free,
    title     = {Free-form deformation of solid geometric models},
    author    = {Sederberg, Thomas W and Parry, Scott R},
    booktitle = {Proceedings of the 13th annual conference on Computer graphics and interactive techniques},
    pages     = {151--160},
    year      = {1986}
}

@book{constantine2015active,
    title     = {Active subspaces: Emerging ideas for dimension reduction in parameter studies},
    author    = {Constantine, Paul G},
    year      = {2015},
    publisher = {SIAM}
}

@article{ballarin2015reduced,
    title     = {Reduced-order models for patient-specific haemodynamics of coronary artery bypass grafts},
    author    = {Ballarin, Francesco},
    year      = {2015},
    publisher = {Politecnico di Milano}
}

@incollection{witteveen2009explicit,
    title     = {Explicit mesh deformation using inverse distance weighting interpolation},
    author    = {Witteveen, Jeroen and Bijl, Hester},
    booktitle = {19th AIAA computational fluid dynamics},
    pages     = {3996},
    year      = {2009}
}

@article{witteveen2010Explicit,
    title   = {Explicit and Robust Inverse Distance Weighting Mesh Deformation for {{CFD}}},
    author  = {Witteveen, Jeroen},
    year    = 2010,
    month   = jan,
    journal = {48th AIAA Aerospace Sciences Meeting Including the New Horizons Forum and Aerospace Exposition},
    doi     = {10.2514/6.2010-165}
}

@incollection{garotta2020reduced,
    title     = {Reduced order isogeometric analysis approach for pdes in parametrized domains},
    author    = {Garotta, Fabrizio and Demo, Nicola and Tezzele, Marco and Carraturo, Massimo and Reali, Alessandro and Rozza, Gianluigi},
    booktitle = {Quantification of Uncertainty: Improving Efficiency and Technology: QUIET Selected Contributions},
    pages     = {153--170},
    year      = {2020},
    publisher = {Springer}
}

@book{keyes2007domain,
    title     = {Domain decomposition methods in science and engineering XVI},
    author    = {Keyes, David E and Widlund, Olof B},
    year      = {2007},
    publisher = {Springer}
}

@article{forti2014efficient,
    title      = {Efficient geometrical parametrisation techniques of interfaces for reduced-order modelling: application to fluid-structure interaction coupling problems},
    volume     = {28},
    issn       = {1061-8562, 1029-0257},
    shorttitle = {Efficient geometrical parametrisation techniques of interfaces for reduced-order modelling},
    doi        = {10.1080/10618562.2014.932352},
    language   = {en},
    number     = {3-4},
    journal    = {International Journal of Computational Fluid Dynamics},
    publisher  = {Taylor \& Francis},
    author     = {Forti, Davide and Rozza, Gianluigi},
    month      = mar,
    year       = {2014},
    pages      = {158--169}
}

@misc{zappon2023staggered,
    title     = {A staggered-in-time and non-conforming-in-space numerical framework for realistic cardiac electrophysiology outputs},
    doi       = {10.48550/arXiv.2308.03884},
    language  = {en},
    publisher = {arXiv},
    author    = {Zappon, Elena and Manzoni, Andrea and Quarteroni, Alfio},
    month     = aug,
    year      = {2023},
    note      = {arXiv:2308.03884 [math]}
}

@article{zappon2023efficient,
    title    = {Efficient and certified solution of parametrized one-way coupled problems through {DEIM}-based data projection across non-conforming interfaces},
    volume   = {49},
    issn     = {1019-7168, 1572-9044},
    doi      = {10.1007/s10444-022-10008-w},
    language = {en},
    number   = {2},
    journal  = {Advances in Computational Mathematics},
    author   = {Zappon, Elena and Manzoni, Andrea and Quarteroni, Alfio},
    month    = apr,
    year     = {2023},
    pages    = {21}
}

@article{deparis2016internodes,
    title      = {{INTERNODES}: an accurate interpolation-based method for coupling the {Galerkin} solutions of {PDEs} on subdomains featuring non-conforming interfaces},
    volume     = {141},
    issn       = {00457930},
    shorttitle = {{INTERNODES}},
    doi        = {10.1016/j.compfluid.2016.03.033},
    language   = {en},
    journal    = {Computers \& Fluids},
    publisher  = {Elsevier},
    author     = {Deparis, Simone and Forti, Davide and Gervasio, Paola and Quarteroni, Alfio},
    month      = dec,
    year       = {2016},
    pages      = {22--41}
}

@article{ballarin2024projection-based,
    title    = {Projection-based reduced order modeling of an iterative scheme for linear thermo-poroelasticity},
    volume   = {21},
    issn     = {25900374},
    doi      = {10.1016/j.rinam.2023.100430},
    language = {en},
    journal  = {Results in Applied Mathematics},
    author   = {Ballarin, Francesco and Lee, Sanghyun and Yi, Son-Young},
    month    = feb,
    year     = {2024},
    pages    = {100430}
}

@article{pichi2022driving,
    title      = {Driving bifurcating parametrized nonlinear {PDEs} by optimal control strategies: application to {Navier}-{Stokes} equations with model order reduction},
    volume     = {56},
    copyright  = {https://creativecommons.org/licenses/by/4.0},
    issn       = {2822-7840, 2804-7214},
    shorttitle = {Driving bifurcating parametrized nonlinear {PDEs} by optimal control strategies},
    doi        = {10.1051/m2an/2022044},
    language   = {en},
    number     = {4},
    journal    = {ESAIM: Mathematical Modelling and Numerical Analysis},
    author     = {Pichi, Federico and Strazzullo, Maria and Ballarin, Francesco and Rozza, Gianluigi},
    month      = jul,
    year       = {2022},
    pages      = {1361--1400}
}

@article{pichi2020reduced,
    title      = {A {Reduced} {Order} {Modeling} {Technique} to {Study} {Bifurcating} {Phenomena}: {Application} to the {Gross}--{Pitaevskii} {Equation}},
    volume     = {42},
    issn       = {1064-8275, 1095-7197},
    shorttitle = {A {Reduced} {Order} {Modeling} {Technique} to {Study} {Bifurcating} {Phenomena}},
    doi        = {10.1137/20M1313106},
    language   = {en},
    number     = {5},
    journal    = {SIAM Journal on Scientific Computing},
    author     = {Pichi, Federico and Quaini, Annalisa and Rozza, Gianluigi},
    month      = jan,
    year       = {2020},
    pages      = {B1115--B1135}
}

@article{khamlich2022model,
    title    = {Model order reduction for bifurcating phenomena in fluid-structure interaction problems},
    volume   = {94},
    issn     = {0271-2091, 1097-0363},
    doi      = {10.1002/fld.5118},
    language = {en},
    number   = {10},
    journal  = {International Journal for Numerical Methods in Fluids},
    author   = {Khamlich, Moaad and Pichi, Federico and Rozza, Gianluigi},
    month    = oct,
    year     = {2022},
    pages    = {1611--1640}
}

@phdthesis{pichi2019reduced,
    address = {Trieste, Italy},
    title   = {Reduced order models for parametric bifurcation problems in nonlinear {PDEs}},
    school  = {SISSA},
    author  = {Pichi, Federico},
    year    = {2019}
}

@article{riva2024multi-physics,
    title      = {Multi-physics model bias correction with data-driven reduced order techniques: {Application} to nuclear case studies},
    volume     = {135},
    issn       = {0307904X},
    shorttitle = {Multi-physics model bias correction with data-driven reduced order techniques},
    doi        = {10.1016/j.apm.2024.06.040},
    language   = {en},
    journal    = {Applied Mathematical Modelling},
    author     = {Riva, Stefano and Introini, Carolina and Cammi, Antonio},
    month      = nov,
    year       = {2024},
    pages      = {243--268}
}

@article{vierendeels2007implicit,
    title     = {Implicit coupling of partitioned fluid-structure interaction problems with reduced order models},
    volume    = {85},
    copyright = {https://www.elsevier.com/tdm/userlicense/1.0/},
    issn      = {0045-7949},
    doi       = {10.1016/j.compstruc.2006.11.006},
    language  = {en},
    number    = {11-14},
    journal   = {Computers \& Structures},
    author    = {Vierendeels, Jan and Lanoye, Lieve and Degroote, Joris and Verdonck, Pascal},
    month     = jun,
    year      = {2007},
    note      = {Publisher: Elsevier BV},
    pages     = {970--976}
}

@article{ballarin2016podgalerkin,
    title     = {{POD}-{Galerkin} monolithic reduced order models for parametrized fluid-structure interaction problems},
    volume    = {82},
    copyright = {http://onlinelibrary.wiley.com/termsAndConditions\#vor},
    issn      = {0271-2091, 1097-0363},
    doi       = {10.1002/fld.4252},
    language  = {en},
    number    = {12},
    journal   = {International Journal for Numerical Methods in Fluids},
    author    = {Ballarin, Francesco and Rozza, Gianluigi},
    month     = dec,
    year      = {2016},
    note      = {Publisher: Wiley},
    pages     = {1010--1034}
}

@article{liberge2010reduced,
    title     = {Reduced order modelling method via proper orthogonal decomposition ({POD}) for flow around an oscillating cylinder},
    volume    = {26},
    copyright = {https://www.elsevier.com/tdm/userlicense/1.0/},
    issn      = {0889-9746},
    doi       = {10.1016/j.jfluidstructs.2009.10.006},
    language  = {en},
    number    = {2},
    journal   = {Journal of Fluids and Structures},
    author    = {Liberge, E. and Hamdouni, A.},
    month     = feb,
    year      = {2010},
    note      = {Publisher: Elsevier BV},
    pages     = {292--311}
}

@article{barzegar2025predictive,
    title      = {A predictive surrogate model based on linear and nonlinear solution manifold reduction in cardiovascular {FSI}: {A} comparative study},
    volume     = {189},
    copyright  = {https://www.elsevier.com/tdm/userlicense/1.0/},
    issn       = {0010-4825},
    shorttitle = {A predictive surrogate model based on linear and nonlinear solution manifold reduction in cardiovascular {FSI}},
    doi        = {10.1016/j.compbiomed.2025.109959},
    language   = {en},
    journal    = {Computers in Biology and Medicine},
    author     = {Barzegar Gerdroodbary, M. and Salavatidezfouli, Sajad},
    month      = may,
    year       = {2025},
    note       = {Publisher: Elsevier BV},
    pages      = {109959}
}

@article{zhang2022data-driven,
    title    = {Data-driven nonlinear reduced-order modeling of unsteady fluid-structure interactions},
    volume   = {34},
    issn     = {1070-6631, 1089-7666},
    doi      = {10.1063/5.0090394},
    language = {en},
    number   = {5},
    journal  = {Physics of Fluids},
    author   = {Zhang, Xinshuai and Ji, Tingwei and Xie, Fangfang and Zheng, Changdong and Zheng, Yao},
    month    = may,
    year     = {2022},
    note     = {Publisher: AIP Publishing}
}

@article{han2022deep,
    title    = {Deep neural network based reduced-order model for fluid-structure interaction system},
    volume   = {34},
    issn     = {1070-6631, 1089-7666},
    doi      = {10.1063/5.0096432},
    language = {en},
    number   = {7},
    journal  = {Physics of Fluids},
    author   = {Han, Renkun and Wang, Yixing and Qian, Weiqi and Wang, Wenzheng and Zhang, Miao and Chen, Gang},
    month    = jul,
    year     = {2022},
    note     = {Publisher: AIP Publishing}
}

@article{xiao2016non-intrusive,
    title     = {Non-intrusive reduced order modelling of fluid-structure interactions},
    volume    = {303},
    copyright = {https://www.elsevier.com/tdm/userlicense/1.0/},
    issn      = {0045-7825},
    doi       = {10.1016/j.cma.2015.12.029},
    language  = {en},
    journal   = {Computer Methods in Applied Mechanics and Engineering},
    author    = {Xiao, D. and Yang, P. and Fang, F. and Xiang, J. and Pain, C.C. and Navon, I.M.},
    month     = may,
    year      = {2016},
    note      = {Publisher: Elsevier BV},
    pages     = {35--54}
}

@article{colciago2014comparisons,
    title     = {Comparisons between reduced order models and full {3D} models for fluid-structure interaction problems in haemodynamics},
    volume    = {265},
    copyright = {https://www.elsevier.com/tdm/userlicense/1.0/},
    issn      = {0377-0427},
    doi       = {10.1016/j.cam.2013.09.049},
    language  = {en},
    journal   = {Journal of Computational and Applied Mathematics},
    author    = {Colciago, C.M. and Deparis, S. and Quarteroni, A.},
    month     = aug,
    year      = {2014},
    note      = {Publisher: Elsevier BV},
    pages     = {120--138}
}

@article{gobat2023reduced,
    title     = {Reduced {Order} {Modeling} of {Nonlinear} {Vibrating} {Multiphysics} {Microstructures} with {Deep} {Learning}-{Based} {Approaches}},
    volume    = {23},
    copyright = {https://creativecommons.org/licenses/by/4.0/},
    issn      = {1424-8220},
    doi       = {10.3390/s23063001},
    language  = {en},
    number    = {6},
    journal   = {Sensors},
    author    = {Gobat, Giorgio and Fresca, Stefania and Manzoni, Andrea and Frangi, Attilio},
    month     = mar,
    year      = {2023},
    note      = {Publisher: MDPI AG},
    pages     = {3001}
}

@article{lieu2006reduced-order,
    title     = {Reduced-order fluid/structure modeling of a complete aircraft configuration},
    volume    = {195},
    copyright = {https://www.elsevier.com/tdm/userlicense/1.0/},
    issn      = {0045-7825},
    doi       = {10.1016/j.cma.2005.08.026},
    language  = {en},
    number    = {41-43},
    journal   = {Computer Methods in Applied Mechanics and Engineering},
    author    = {Lieu, T. and Farhat, C. and Lesoinne, M.},
    month     = aug,
    year      = {2006},
    note      = {Publisher: Elsevier BV},
    pages     = {5730--5742}
}

@article{manzoni2018reduced,
    title     = {Reduced order modeling for cardiac electrophysiology and mechanics: New methodologies, challenges and perspectives},
    author    = {Manzoni, Andrea and Bonomi, Diana and Quarteroni, Alfio},
    journal   = {Mathematical and numerical modeling of the cardiovascular system and applications},
    pages     = {115--166},
    year      = {2018},
    publisher = {Springer}
}

@article{cicci2023projection,
    title   = {Projection-based reduced order models for parameterized nonlinear time-dependent problems arising in cardiac mechanics},
    author  = {Cicci, Ludovica and Fresca, Stefania and Pagani, Stefano and Manzoni, Andrea and Quarteroni, Alfio and others},
    journal = {Mathematics in Engineering},
    volume  = {5},
    number  = {2},
    pages   = {1--38},
    year    = {2023}
}

@article{nonino2021monolithic,
    title     = {A {Monolithic} and a {Partitioned}, {Reduced} {Basis} {Method} for {Fluid}-{Structure} {Interaction} {Problems}},
    volume    = {6},
    copyright = {https://creativecommons.org/licenses/by/4.0/},
    issn      = {2311-5521},
    doi       = {10.3390/fluids6060229},
    language  = {en},
    number    = {6},
    journal   = {Fluids},
    author    = {Nonino, Monica and Ballarin, Francesco and Rozza, Gianluigi},
    month     = jun,
    year      = {2021},
    note      = {Publisher: MDPI AG},
    pages     = {229}
}

@phdthesis{nonino2020application,
    address  = {Trieste, Italy},
    type     = {{PhD} {Thesis}},
    title    = {On the application of the {Reduced} {Basis} {Method} to {Fluid}-{Structure} {Interaction} problems},
    language = {en},
    school   = {SISSA},
    author   = {Nonino, Monica},
    year     = {2020}
}

@article{nonino2023projection,
    title     = {Projection {Based} {Semi}-{Implicit} {Partitioned} {Reduced} {Basis} {Method} for {Fluid}-{Structure} {Interaction} {Problems}},
    volume    = {94},
    copyright = {https://creativecommons.org/licenses/by/4.0},
    issn      = {0885-7474, 1573-7691},
    doi       = {10.1007/s10915-022-02049-6},
    language  = {en},
    number    = {1},
    journal   = {Journal of Scientific Computing},
    author    = {Nonino, Monica and Ballarin, Francesco and Rozza, Gianluigi and Maday, Yvon},
    month     = jan,
    year      = {2023},
    note      = {Publisher: Springer Science and Business Media LLC}
}

@article{nonino2023reduced,
    title   = {A reduced basis method by means of transport maps for a fluid-structure interaction problem with slowly decaying {Kolmogorov} n-width},
    volume  = {1},
    number  = {1},
    journal = {Advances in Computational Science \& Engineering (ACSE)},
    author  = {Nonino, Monica and Ballarin, Francesco and Rozza, Gianluigi and Maday, Yvon},
    year    = {2023}
}

@article{pitton2017computational,
    title     = {Computational reduction strategies for the detection of steady bifurcations in incompressible fluid-dynamics: Applications to Coanda effect in cardiology},
    author    = {Pitton, Giuseppe and Quaini, Annalisa and Rozza, Gianluigi},
    journal   = {Journal of Computational Physics},
    volume    = {344},
    pages     = {534--557},
    year      = {2017},
    publisher = {Elsevier}
}

@book{hale2012dynamics,
    title     = {Dynamics and bifurcations},
    author    = {Hale, Jack K and Ko{\c{c}}ak, H{\"u}seyin},
    volume    = {3},
    year      = {2012},
    publisher = {Springer Science \& Business Media}
}

@phdthesis{dossi2015combined,
    address = {Milan},
    type    = {{PhD} {Thesis}},
    title   = {Combined {Model} {Order} {Reduction} and {Domain} {Decomposition} strategies for the solution of non-linear and multi-physics structural p},
    school  = {Politecnico di Milano},
    author  = {Dossi, Martino},
    month   = jan,
    year    = {2015}
}

@article{lee2024parametric,
    title     = {Parametric model order reduction by machine learning for fluid-structure interaction analysis},
    volume    = {40},
    copyright = {https://creativecommons.org/licenses/by/4.0},
    issn      = {0177-0667, 1435-5663},
    doi       = {10.1007/s00366-023-01782-2},
    language  = {en},
    number    = {1},
    journal   = {Engineering with Computers},
    author    = {Lee, SiHun and Jang, Kijoo and Lee, Sangmin and Cho, Haeseong and Shin, SangJoon},
    month     = feb,
    year      = {2024},
    note      = {Publisher: Springer Science and Business Media LLC},
    pages     = {45--60}
}

@misc{farea2025learning,
    title     = {Learning {Fluid}-{Structure} {Interaction} {Dynamics} with {Physics}-{Informed} {Neural} {Networks} and {Immersed} {Boundary} {Methods}},
    doi       = {10.48550/arXiv.2505.18565},
    language  = {en},
    publisher = {arXiv},
    author    = {Farea, Afrah and Khan, Saiful and Daryani, Reza and Ersan, Emre Cenk and Celebi, Mustafa Serdar},
    month     = jun,
    year      = {2025},
    note      = {arXiv:2505.18565 [cs]}
}

@article{ma2022preliminary,
    title     = {A {Preliminary} {Study} on the {Resolution} of {Electro}-{Thermal} {Multi}-{Physics} {Coupling} {Problem} {Using} {Physics}-{Informed} {Neural} {Network} ({PINN})},
    volume    = {15},
    copyright = {https://creativecommons.org/licenses/by/4.0/},
    issn      = {1999-4893},
    doi       = {10.3390/a15020053},
    language  = {en},
    number    = {2},
    journal   = {Algorithms},
    author    = {Ma, Yaoyao and Xu, Xiaoyu and Yan, Shuai and Ren, Zhuoxiang},
    month     = feb,
    year      = {2022},
    note      = {Publisher: MDPI AG},
    pages     = {53}
}

@article{wu2023application,
    title     = {The {Application} of {Physics}-{Informed} {Machine} {Learning} in {Multiphysics} {Modeling} in {Chemical} {Engineering}},
    volume    = {62},
    copyright = {https://doi.org/10.15223/policy-029},
    issn      = {0888-5885, 1520-5045},
    doi       = {10.1021/acs.iecr.3c02383},
    language  = {en},
    number    = {44},
    journal   = {Industrial \& Engineering Chemistry Research},
    author    = {Wu, Zhiyong and Wang, Huan and He, Chang and Zhang, Bingjian and Xu, Tao and Chen, Qinglin},
    month     = nov,
    year      = {2023},
    note      = {Publisher: American Chemical Society (ACS)},
    pages     = {18178--18204}
}

@article{ryu2024multiphysics,
    title     = {Multiphysics generalization in a polymerization reactor using physics-informed neural networks},
    volume    = {298},
    copyright = {https://www.elsevier.com/tdm/userlicense/1.0/},
    issn      = {0009-2509},
    doi       = {10.1016/j.ces.2024.120385},
    language  = {en},
    journal   = {Chemical Engineering Science},
    author    = {Ryu, Yubin and Shin, Sunkyu and Lee, Won Bo and Na, Jonggeol},
    month     = oct,
    year      = {2024},
    note      = {Publisher: Elsevier BV},
    pages     = {120385}
}

@article{balajewicz2014Reduction,
    title        = {Reduction of Nonlinear Embedded Boundary Models for Problems with Evolving Interfaces},
    author       = {Balajewicz, Maciej and Farhat, Charbel},
    date         = {2014-10},
    journaltitle = {Journal of Computational Physics},
    shortjournal = {Journal of Computational Physics},
    volume       = {274},
    pages        = {489--504},
    issn         = {00219991},
    doi          = {10.1016/j.jcp.2014.06.038},
    url          = {https://linkinghub.elsevier.com/retrieve/pii/S0021999114004458},
    urldate      = {2026-01-05},
    langid       = {english}
}

@online{karatzas2018Reduced,
    title       = {A {{Reduced Order Approach}} for the {{Embedded Shifted Boundary FEM}} and a {{Heat Exchange System}} on {{Parametrized Geometries}}},
    author      = {Karatzas, E. N. and Stabile, G. and Atallah, N. and Scovazzi, G. and Rozza, G.},
    date        = {2018-12-17},
    eprint      = {1807.07753},
    eprinttype  = {arXiv},
    eprintclass = {math},
    doi         = {10.48550/arXiv.1807.07753},
    url         = {http://arxiv.org/abs/1807.07753},
    urldate     = {2026-01-05},
    pubstate    = {prepublished}
}

@article{karatzas2019Reduceda,
    title        = {A Reduced Basis Approach for {{PDEs}} on Parametrized Geometries Based on the Shifted Boundary Finite Element Method and Application to a {{Stokes}} Flow},
    author       = {Karatzas, Efthymios N. and Stabile, Giovanni and Nouveau, Leo and Scovazzi, Guglielmo and Rozza, Gianluigi},
    date         = {2019-04},
    journaltitle = {Computer Methods in Applied Mechanics and Engineering},
    shortjournal = {Computer Methods in Applied Mechanics and Engineering},
    volume       = {347},
    pages        = {568--587},
    issn         = {00457825},
    doi          = {10.1016/j.cma.2018.12.040},
    url          = {https://linkinghub.elsevier.com/retrieve/pii/S0045782518306479},
    urldate      = {2026-01-05},
    langid       = {english}
}

@article{karatzas2020Projectionbased,
    title        = {Projection-Based Reduced Order Models for a Cut Finite Element Method in Parametrized Domains},
    author       = {Karatzas, Efthymios N. and Ballarin, Francesco and Rozza, Gianluigi},
    date         = {2020-02},
    journaltitle = {Computers \& Mathematics with Applications},
    shortjournal = {Computers \& Mathematics with Applications},
    volume       = {79},
    number       = {3},
    pages        = {833--851},
    issn         = {08981221},
    doi          = {10.1016/j.camwa.2019.08.003},
    url          = {https://linkinghub.elsevier.com/retrieve/pii/S0898122119303931},
    urldate      = {2026-01-05},
    langid       = {english}
}

@article{karatzas2020ReducedOrder,
    title        = {A {{Reduced-Order Shifted Boundary Method}} for {{Parametrized}} Incompressible {{Navier-Stokes}} Equations},
    author       = {Karatzas, Efthymios N. and Stabile, Giovanni and Nouveau, Leo and Scovazzi, Guglielmo and Rozza, Gianluigi},
    date         = {2020-10},
    journaltitle = {Computer Methods in Applied Mechanics and Engineering},
    shortjournal = {Computer Methods in Applied Mechanics and Engineering},
    volume       = {370},
    eprint       = {1907.10549},
    eprinttype   = {arXiv},
    eprintclass  = {math},
    pages        = {113273},
    issn         = {00457825},
    doi          = {10.1016/j.cma.2020.113273},
    url          = {http://arxiv.org/abs/1907.10549},
    urldate      = {2026-01-05}
}

@online{karatzas2021Reduced,
    title       = {A Reduced Order Model for a Stable Embedded Boundary Parametrized {{Cahn-Hilliard}} Phase-Field System Based on Cut Finite Elements},
    author      = {Karatzas, Efthymios N. and Rozza, Gianluigi},
    date        = {2021-08-09},
    eprint      = {2009.01596},
    eprinttype  = {arXiv},
    eprintclass = {math},
    doi         = {10.48550/arXiv.2009.01596},
    url         = {http://arxiv.org/abs/2009.01596},
    urldate     = {2026-01-05},
    pubstate    = {prepublished}
}

@article{zeng2022Embedded,
    title        = {Embedded Domain {{Reduced Basis Models}} for the Shallow Water Hyperbolic Equations with the {{Shifted Boundary Method}}},
    author       = {Zeng, Xianyi and Stabile, Giovanni and Karatzas, Efthymios N. and Scovazzi, Guglielmo and Rozza, Gianluigi},
    date         = {2022-08},
    journaltitle = {Computer Methods in Applied Mechanics and Engineering},
    shortjournal = {Computer Methods in Applied Mechanics and Engineering},
    volume       = {398},
    pages        = {115143},
    issn         = {00457825},
    doi          = {10.1016/j.cma.2022.115143},
    url          = {https://linkinghub.elsevier.com/retrieve/pii/S0045782522003164},
    urldate      = {2026-01-05},
    langid       = {english}
}

@article{main2018Shifteda,
    title        = {The Shifted Boundary Method for Embedded Domain Computations. {{Part I}}: {{Poisson}} and {{Stokes}} Problems},
    shorttitle   = {The Shifted Boundary Method for Embedded Domain Computations. {{Part I}}},
    author       = {Main, A. and Scovazzi, G.},
    date         = {2018-11},
    journaltitle = {Journal of Computational Physics},
    shortjournal = {Journal of Computational Physics},
    volume       = {372},
    pages        = {972--995},
    issn         = {00219991},
    doi          = {10.1016/j.jcp.2017.10.026},
    url          = {https://linkinghub.elsevier.com/retrieve/pii/S0021999117307799},
    urldate      = {2026-01-06},
    langid       = {english}
}

@article{main2018Shiftedb,
    title        = {The Shifted Boundary Method for Embedded Domain Computations. {{Part II}}: {{Linear}} Advection–Diffusion and Incompressible {{Navier}}–{{Stokes}} Equations},
    shorttitle   = {The Shifted Boundary Method for Embedded Domain Computations. {{Part II}}},
    author       = {Main, A. and Scovazzi, G.},
    date         = {2018-11},
    journaltitle = {Journal of Computational Physics},
    shortjournal = {Journal of Computational Physics},
    volume       = {372},
    pages        = {996--1026},
    issn         = {00219991},
    doi          = {10.1016/j.jcp.2018.01.023},
    url          = {https://linkinghub.elsevier.com/retrieve/pii/S0021999118300330},
    urldate      = {2026-01-06},
    langid       = {english}
}

@article{wang2011Algorithms,
    title        = {Algorithms for Interface Treatment and Load Computation in Embedded Boundary Methods for Fluid and Fluid–Structure Interaction Problems},
    author       = {Wang, K. and Rallu, A. and Gerbeau, J.-F. and Farhat, C.},
    date         = {2011},
    journaltitle = {International Journal for Numerical Methods in Fluids},
    volume       = {67},
    number       = {9},
    pages        = {1175--1206},
    issn         = {1097-0363},
    doi          = {10.1002/fld.2556},
    url          = {https://onlinelibrary.wiley.com/doi/abs/10.1002/fld.2556},
    urldate      = {2026-01-06},
    langid       = {english}
}

@article{pasquariello2016Cutcell,
    title   = {A Cut-Cell Finite Volume -- Finite Element Coupling Approach for Fluid--Structure Interaction in Compressible Flow},
    author  = {Pasquariello, Vito and Hammerl, Georg and {\"O}rley, Felix and Hickel, Stefan and Danowski, Caroline and Popp, Alexander and Wall, Wolfgang A. and Adams, Nikolaus A.},
    year    = 2016,
    month   = feb,
    journal = {Journal of Computational Physics},
    volume  = {307},
    pages   = {670--695},
    issn    = {00219991},
    doi     = {10.1016/j.jcp.2015.12.013},
    urldate = {2026-01-07},
    langid  = {english}
}

\end{document}